\documentclass[12pt]{amsart}
\usepackage{amsfonts,latexsym,amsmath, amssymb}
\usepackage{mathrsfs,MnSymbol}
\usepackage{url,color}
\usepackage{upgreek}
\usepackage{fancyhdr}
\usepackage{hyperref}
\usepackage{times}
\usepackage{scalefnt}
\usepackage{url}
\usepackage{cite}

\newcommand{\bea}{\begin{eqnarray}}
\newcommand{\eea}{\end{eqnarray}}
\def\beaa{\begin{eqnarray*}}
\def\eeaa{\end{eqnarray*}}
\def\ba{\begin{array}}
\def\ea{\end{array}}
\def\be#1{\begin{equation} \label{#1}}
\def \eeq{\end{equation}}

\def\be{{\beta}}

\def\De{\Delta}

\def\la{\lambda}

\def\HH{{\mathcal H}}

\def\HH{\mathcal H}

\def\R{{\mathbb{R}}}
\def\C{{\mathbb{C}}}

\def\T{{\mathbb{T}}}
\def\N{{\mathbb N}}

\newcommand{\EQ}[1]{\begin{align*}\begin{split} #1 \end{split}\end{align*}}
\newcommand{\EQn}[1]{\begin{align}\begin{split} #1 \end{split}\end{align}}

\def\normo#1{\left\|#1\right\|}

\def\normbb#1{\bigg\|#1\bigg\|}

\def\brk#1{\left(#1\right)}
\def\fbrk#1{\left\lbrace#1\right\rbrace}
\def\jb#1{\langle#1\rangle}

\def\loe{\leqslant}

\def\lsm{\lesssim}

\def\half#1{\frac{#1}{2}}

\newtheorem{theorem}{Theorem}[section]
\newtheorem{lemma}[theorem]{Lemma}
\newtheorem{proposition}[theorem]{Proposition}
\newtheorem{corollary}[theorem]{Corollary}
\newtheorem{definition}[theorem]{Definition}
\newtheorem{remark}[theorem]{Remark}

\numberwithin{equation}{section}

\newcommand{\norm}[1]{\left\lVert #1 \right\rVert}
\newcommand{\abs}[1]{\left\vert#1\right\vert}

\numberwithin{equation}{section}

\topmargin - .23 in

\begin{document}
\title{Scattering of the three-dimensional cubic nonlinear Schr\"odinger equation with partial harmonic potentials}
\author[X. Cheng, C.-Y. Guo, Z. Guo, X. Liao and J. Shen]{Xing Cheng, Chang-Yu Guo, Zihua Guo, Xian Liao, and Jia Shen}

\address[Xing Cheng]{College of Science, Hohai University, Nanjing 210098, Jiangsu, P. R. China}
\email{chengx@hhu.edu.cn}

\address[Chang-Yu Guo]{Research Center for Mathematics and Interdisciplinary Sciences, Shandong University 266237,  Qingdao, P. R. China}
\email{changyu.guo@sdu.edu.cn}

\address[Zihua Guo]{School of Mathematical Sciences, Monash University, VIC 3800,  Australia}
\email{Zihua.Guo@monash.edu}

\address[Xian Liao]{Institute of Analysis, Karlsruhe Institute of Technology, Englestr. 2, 76131 Karlsruhe, Germany}
\email{xian.liao@kit.edu}

\address[Jia Shen]{School of Mathematical Sciencs, Peking University, No 5. Yiheyuan Road, Beijing 100871, P. R. China}
\email{shenjia@pku.edu.cn}

\date{\today}
\thanks{X. Cheng is partially supported by the NSF grant of China (No.~11526072). The corresponding author C.-Y. Guo is supported by the Qilu funding of Shandong University (No.~62550089963197). Z. Guo was supported by the ARC project (NO.~DP200101065 and No.~DP170101060). X. Liao is partially supported by DFG via Collaborative Research Centre 1173. J. Shen was supported by China Scholarship Council.}		
\begin{abstract}
	In this paper, we consider the following three dimensional defocusing cubic nonlinear Schr\"odinger equation (NLS) with partial harmonic potential
	\begin{equation*}\tag{NLS}
		\begin{cases}
			i\partial_t u + \left(\Delta_{\mathbb{R}^3 }-x^2 \right)  u   = |u|^2 u,\\
			u|_{t=0} = u_0.
		\end{cases}
	\end{equation*}
	Our main result shows that the solution $u$ scatters for any given initial data $u_0$ with finite mass and energy.
	
	The main new ingredient in our approach is to approximate (NLS) in the large-scale case by a relevant dispersive continuous resonant (DCR) system. The proof of global well-posedness and scattering of the new (DCR) system is greatly inspired by the fundamental works of Dodson \cite{D3,D1,D2} in his study of scattering for the mass-critical nonlinear Schr\"odinger equation. The analysis of (DCR) system allows us to utilize the additional regularity of the smooth nonlinear profile so that the celebrated concentration-compactness/rigidity argument of Kenig and Merle \cite{KM,KM1} applies.
\bigskip

\bigskip

\end{abstract}

\maketitle
\setcounter{tocdepth}{2}
\section{Introduction}

\subsection{Background and motivation}
Consider the Cauchy problem for the following family of \emph{nonlinear Schr\"odinger equations}
in $\R^{d}$, $d\in \N$,  with \emph{harmonic oscillators}:
\begin{equation}\label{eq1.2n}
\begin{cases}
i\partial_t u + \Delta_{\mathbb{R}^d}  u - \left(\omega^2  |y|^2 + |x|^2 \right) u = \mu |u|^{p-1} u,\\
u|_{t=0} = u_0,
\end{cases}
\end{equation}
where $1 < p < \infty$, $(y,x) \in \mathbb{R}^{d_1} \times \mathbb{R}^{d_2}$, $d= d_1+ d_2$, and $d_1, d_2\in \mathbb{N}$, $d_1 \ge 1$. The complex-valued function $u=u(t,y,x)\colon \R\times \R^d\to \C$ is the unknown wave function. The parameter $\omega = 0 \text{ or } 1$, with $\omega=1$ corresponding to the quadratic potential case, while $\omega=0$ corresponding to the partial harmonic oscillator on the left hand side. The parameter $\mu =1$ or $\mu=-1$ corresponds to the \emph{defocusing} or \emph{focusing} case respectively. Equation \eqref{eq1.2n} arises as models for diverse physical phenomena, including
Bose-Einstein condensates in a laboratory trap \cite{JP,PS} and  the envelope dynamics of a general dispersive wave in a weakly nonlinear medium. It can also be derived in the NLS  with  constant magnetic potential, see for example \cite{FO}.
The associated conserved mass and energy of equation \eqref{eq1.2n} read as
$$\mathcal{M}(u)(t)=\int_{\R^{d_1}\times \R^{d_2}}|u(t,y,x)|^2\mathrm{d}y \mathrm{d}x$$
and
$$\mathcal{E}^{d_1,d_2}_{\omega,\mu,p}(u)(t)=\int_{\mathbb{R}^{d_1}\times\mathbb{R}^{d_2}}
\frac12 |\nabla_{y, x} u(t, y, x )|^2
+ \frac12  \left(\omega^2|y|^2+|x|^2 \right) |u(t, y, x)|^2+  \frac{\mu} {p+1} |u(t, y, x)|^{p+1} \, \mathrm{d}y \mathrm{d}x.$$
It is natural to take the initial data from the following weighted Sobolev space
\begin{align*}
 u_0\in \left\{ f=f(y,x)  \in L_{y,x}^2(\mathbb{R}^d)\, : \,   \left\|\nabla_{y,x}  f \right\|_{L^2_{y,x}(\R^d)} +   \left\| |x|  f  \right\|_{L^2_{y,x} (\R^d)} +  \omega  \left\| |y| f  \right\|_{L_{y,x}^2(\R^d)}  + \|f  \|_{L_{y,x}^2(\R^d)} < \infty\right\}.
\end{align*}
In view of the Sobolev embedding
\begin{equation*}
H^1(\R^d)\hookrightarrow L^q(\R^d), \qquad
\begin{cases}
2\leq q\leq 2+\frac{4}{d-2} & \text{ if } d\geq 3\\
2\leq q<\infty & \text{ if } d=2\\
2\leq q\leq\infty & \text{ if } d=1,
\end{cases}
\end{equation*}
\begin{equation*}
\text{the initial data is of finite energy in the \emph{energy-subcritical }case} \qquad
\begin{cases}
 1 < p < 1 + \frac4{d-2} & \text{ if } d\geq 3\\
1 < p < \infty & \text{ if } d=1,2,
\end{cases}\qquad\qquad\qquad\qquad\qquad\qquad\quad
\end{equation*}
and we call the critical case $p = 1 + \frac4{d-2}$, $d\ge 3$ the \emph{energy critical} case.

The global well-posedness of \eqref{eq1.2n} has been established  in the energy-subcritical case by R. Carles \cite{C,C3} in the defocusing case $\mu=1$, and by J. Zhang \cite{Z} in the focusing case $\mu=-1$ when the initial energy  is assumed to be less than the energy of the ground state of the related elliptic equation. The Cauchy problem for the equation \eqref{eq1.2n}  with quadratic potential (that is, $\omega = 1$) in the energy-critical case  was considered  by R. Killip, M. Visan, and X. Zhang \cite{KVZ} in the radial case, and in the general case later by C. Jao \cite{J1,J2}.
 They proved the global well-posedness for the defocusing case and also for the focusing case when the initial energy (resp. kinetic energy) is less than the energy (resp. kinetic energy) of the ground state. We would also like to mention the work of C. Hao, L. Hsiao and H. Li \cite{HHL1,HHL2}, where the authors proved the global well-posedness for the equation \eqref{eq1.2n} (when $\omega = 1$) with an additional angular momentum rotational term.

It is well-known that solutions of the equation \eqref{eq1.2n} with a quadratic potential (i.e. $\omega=1$) can not scatter.
However, intuitively, in the defocusing case, if we turn off the confinement in \emph{some instead of all} of the directions, it should suffice for the condensate to evolve asymptotically freely: Indeed, if $\omega=0$, then the operator $i\partial_t+\Delta_y$ should yield large time dispersion and one expects a scattering theory for the equation~\eqref{eq1.2n}.
When $\omega= 0$, the scattering phenomena for the equation \eqref{eq1.2n} in the defocusing case has already been showed by P. Antonelli, R. Carles and J. D. Silva \cite{ACS} (see also \cite{CG}) in the fully weighted space when $\omega  =0$, $\mu=1$, $ d_1 = 1, 2, 3  $,   $d_2 = 1$ and $ 1+\frac4{d_1}<  p < 1+\frac4{d_1-1}$. The focusing case of \eqref{eq1.2n} has been investigated by A. H. Ardila and R. Carles \cite{AC} recently when the energy is strictly less than the static energy of the ground state. In this aspect, one expects the global-in-time well-posedness result for the defocusing/focusing (when energy is strictly less than the static energy of the ground state) energy-critical and subcritical cases for \eqref{eq1.2n}.
On the other hand, {the potential influences strongly the asymptotic dynamics of the solution}. In \eqref{eq1.2n}, the $x-$direction is not expected to have a global in time dispersive estimate in view of the Mehler's formula
\begin{align*}
e^{it \left( \Delta_x - |x|^2 \right)} f(y,x) = \left( 2 \pi i \sin( 2 t)  \right)^{- \frac{d_2}2} \int_{\mathbb{R}^{d_2}} e^{ \frac{i}{ \sin(2 t)} \left( \frac{ |x|^2 + |\tilde{x}|^2 }2 \cos( 2t) - x \cdot \tilde{x}  \right) } f(y,  \tilde{x}) \mathrm{d} \tilde{x} ,\  \forall\, y \in \mathbb{R}^{d_1}, x \in \mathbb{R}^{d_2},
\end{align*}
from which we can only derive the following periodically in time dispersive estimate
\begin{align*}
\left\|e^{it \left( \Delta_x - |x|^2 \right)} f(y, x)  \right\|_{L_x^\infty(\mathbb{R}^{d_2})} \lesssim |\sin (2 t) |^{- \frac{d_2}2} \left\|f (y,x) \right\|_{L_x^1(\mathbb{R}^{d_2})}, \ \forall\, t \not\in  \frac\pi2 \mathbb{Z}, \forall\, y \in \mathbb{R}^{d_1} .
\end{align*}
Nevertheless, we have the following global in time dispersive estimate in the $y-$direction:
\begin{align*}
\left\|e^{it \left( \Delta_x + \Delta_y   - |x|^2 \right) } f(y,x) \right\|_{L_y^\infty L_x^2( \mathbb{R}^d)} \lesssim |t|^{- \frac{d_1}2} \left\|f(y,x) \right\|_{L_y^1 L_x^2},
\end{align*}
where we used the dispersive estimate for the semigroup $e^{it\Delta_y}$ together with the $L^2$-norm conservation for the unitary of the operator $e^{it \left( \Delta_x - |x|^2 \right) }$. Thus, according to the scattering theory for the nonlinear Schr\"odinger equations without potential, see for instance \cite{Sta,T2}, one \emph{expects
a scattering result in the weighted Sobolev space} when $\omega=0$ in the case $ 1 + \frac4{d_1} \le p \le 1 + \frac4{d_1+ d_2-2}$ with $d_1 + d_2 \ge 2$. Generally, to obtain the scattering in the inter-critical case, one relies on the Morawetz estimate, see for instance \cite{ACS}. It is difficult to deal with the scattering on the two endpoints $p=1+ \frac4{d_1}$ and $p=1 + \frac4{d_1 + d_2 - 2}$, which correspond to the usual $d_1$ dimensional mass-critical and $d_1+d_2$ dimensional energy-critical nonlinear Schr\"odinger equation without potentials respectively. For the endpoint $p=1 + \frac4{d_1 + d_2 - 2}$, the scattering is a byproduct of the proof of the global well-posedness, and we need to use the induction on energy method \cite{CKSTT0} or the concentration-compactness/rigidity argument \cite{KM,KM1} to prove the global well-posedness. It seems for us one of the main difficulty is to establish a more delicate global in time Strichartz estimate which should be a lot combination of the local Strichartz estimate of 3 dimensional Schr\"odinger equations
as in \cite{B,HP}. We refer to \cite{Killip-Visan1} for more illustration on the proof of the scattering of the nonlinear Schr\"odinger equations at critical regularity Sobolev space. For the endpoint $p=1 + \frac4{d_1}$, global well-posedness is quite easy to get, and the main obstacle is to show the scattering. We cannot prove the scattering by the Morawetz estimate even when the initial data lies in a better regular Sobolev space $H^1_{y,x} $ because the Morawetz estimate only provides a priori estimate of the non-endpoint Strichartz norm on the $\dot{H}_y^\frac14 L_x^2-$level but cannot give a priori estimate of the Strichartz norm on the $L^2-$level, which is not enough to yield the scattering in this case. Therefore, to show the scattering, we still need to use the concentration-compactness/rigidity argument \cite{KM,KM1} and its mass-critical counterpart \cite{D1,D2,D3,KTV09,Killip-Visan1,KVZ0,TVZ,TVZ0} to show the finiteness of the $L^2-$level Strichartz norm. In the $L^2-$level Strichartz norm, we need to consider not only the space and time translations of the equation \eqref{eq1.2n} as in the case $1 + \frac4{d_1} < p \le 1 + \frac4{d_1+d_2 - 2}$, but also the following partial Galilean invariance
\begin{align*}
u(t,y,x) \mapsto e^{-it |\xi_0|^2} e^{iy\cdot \xi_0} u(t, y - 2\xi_0 t, x),
\end{align*}
where $\xi_0 \in \mathbb{R}^{d_1}$, of the equation \eqref{eq1.2n}.
In addition, by a limitation operation, it is realized that a new mass-critical nonlinear Schr\"odinger system can be embedded into \eqref{eq1.2n}, this new mass-critical nonlinear Schr\"odinger system inherits the above invariance and also has the scaling invariance in space-time, and its global well-posedness and scattering should be proven by the argument from \cite{D1,D2,D3,KTV09,Killip-Visan1,KVZ0,TVZ,TVZ0}.

In this paper, we will consider the following Cauchy problem for the defocusing cubic NLS on $\mathbb{R}^3$
\begin{equation}\label{eq1.1}
\begin{cases}
i\partial_t u + \left(\Delta_{\mathbb{R}^3 }-x^2 \right)  u   = |u|^2 u,\\
u|_{t=0} = u_0,
\end{cases}
\end{equation}
where $u=u(t,y,x)\colon \mathbb{R}\times \mathbb{R}^2 \times \mathbb{R}  \to \mathbb{C}$ is an unknown wave function. The following mass and energy quantities are conserved by the evolution of the equation \eqref{eq1.1}
\begin{align*}\tag{ME}\label{ME}
 \mathcal{M}(u(t)) &   = \int_{\mathbb{R}^2 \times \mathbb{R} } |u(t,y, x )|^2\, \mathrm{d}y \mathrm{d}x,\\
 \mathcal{E}(u(t)) &   = \int_{\mathbb{R}^2\times\mathbb{R}} \frac12 |\nabla_{y, x} u(t, y, x )|^2 + \frac12 x^2 |u(t, y, x)|^2+ \frac14 |u(t, y, x)|^4 \, \mathrm{d}y \mathrm{d}x.
\end{align*}
Motivated by the mass and energy formulations, we take the initial data in the following weighted Sobolev space
\begin{align}\label{Sigma}
&u_0\in\Sigma(\mathbb{R}^3) := \left\{ f\in L^2_{y,x}(\mathbb{R}^3)\,|\,\|f\|_{\Sigma(\R^3)} := \|\nabla_{y} f\|_{L^2_{y,x}(\R^3)}
+\|f\|_{L^2_y\mathcal{H}^1_x(\R^2\times\R)}  < \infty \right\}, \\
&\hbox{ with }\|f\|_{\mathcal{H}^1_x(\R)}=\|  f\|_{H^1_x(\R)}+\| x  f\|_{L^2_x(\R)} . \notag
\end{align}
By the Sobolev embedding $H^1(\R^3)\hookrightarrow L^q(\R^3)$, $2\leq q\leq 6$, the initial data is of finite mass and energy.

Observe that equation \eqref{eq1.1} is a special case of equation \eqref{eq1.2n}, namely, corresponding to $d_1 =2, d_2 = 1, \omega = 0, \mu = 1, p = 1+\frac{4}{d_1}=3$ in \eqref{eq1.2n}. In this case, the scattering phenomena is not yet clear. As we are in the energy subcritical case $1<p=3<5$, the equation \eqref{eq1.1} is global well-posed and the scattering of the solutions follows in the small initial data case $\norm{u_0}_\Sigma\ll1$, which is a by-product of the small data well-posedness theorem. We will briefly explore these results in Section \ref{se2} and outline the ideas of the proofs as we did not find them in the literature.

\subsection{Main results}
Our main result of this article is the following scattering result for solutions of the defocusing cubic NLS \eqref{eq1.1}. Recall that $\Sigma(\mathbb{R}^3)$ is defined in \eqref{Sigma}.

\begin{theorem}\label{th1.3}
For any initial data $u_0\in \Sigma(\mathbb{R}^3)$, there is a unique global solution $u\in C_t^0 (\mathbb{R}, \Sigma(  \mathbb{R}^3) ) $ of equation \eqref{eq1.1}. Moreover, the solution scatters, namely there exist $u_\pm \in \Sigma(\mathbb{R}^3)$ such that
\begin{equation*}
\left\| u(t)- e^{it \left(\Delta_{{\mathbb{R}}^3}- x^2 \right)}u_\pm \right\|_{\Sigma(\mathbb{R}^3)} \to 0,  \text{  as }  t\to \pm \infty.
\end{equation*}
\end{theorem}
In order to treat the general initial data with finite (but not necessarily small) $\Sigma$-norm $\norm{u_0}_\Sigma<\infty$, we turn to the celebrated concentration-compactness/rigidity argument developed by C. E. Kenig and F. Merle \cite{KM1,KM}, where one key ingredient is the linear and nonlinear profile decompositions for solutions with bounded $\Sigma-$norm. The proof of Theorem \ref{th1.3} shall rely on (a corollary of) Theorem \ref{as3.14v12} given below.

As for the nonlinear profile decomposition, we will consider a sequence of solutions exhibiting an extreme behavior to study the concentration of the data.
More precisely, we need to study the behavior of the nonlinear profile $u_\lambda$ when $\lambda \to \infty$. The (simplified) nonlinear profile $u_\lambda$, $\lambda > 0$, is the solution of the equation \eqref{eq1.1}
\begin{align}\label{eq1.5v18}
\begin{cases}
i\partial_t u_\lambda + \Delta_{y} u_\lambda + \left(\Delta_x - x^2 \right) u_\lambda = |u_\lambda|^2 u_\lambda,\\
u_\lambda(0,y,x) = \frac1\lambda \phi\left(\frac{y}\lambda, x\right),
\end{cases}
\end{align}
taking the initial data by rescaling the function $\phi$ only in the $y$-variable. Set
$$w_\lambda(t,y,x) = e^{-it\left(\Delta_x - x^2\right)} u_\lambda(t,y,x),$$
and we obtain from \eqref{eq1.5v18} the following  evolutionary equation for $w_\lambda$
\begin{align*}
\begin{cases}
\left(i\partial_t + \Delta_y\right) w_\lambda = e^{-it \left(\Delta_x - x^2 \right)} \left(\left|e^{it \left(\Delta_x - x^2 \right)} w_\lambda\right|^2 e^{it(\Delta_x - x^2)} w_\lambda \right),\\
w_\lambda(0,y,x) = \frac1\lambda \phi\left(\frac{y}\lambda, x\right).
\end{cases}
\end{align*}
If we define $w_\lambda(t,y,x) = \frac1\lambda \tilde{v}\left(\frac{t}{\lambda^2},\frac{y}\lambda, x \right)$, then $\tilde v$ satisfies
\begin{align*}
\left\{\begin{array}{l}
\left(i\partial_t + \Delta_y \right) \tilde{v} = e^{-i\lambda^2 t\left(\Delta_x - x^2 \right)} \left(\left|e^{i\lambda^2 t (\Delta_x - x^2)} \tilde{v}\right|^2 e^{i\lambda^2 t \left(\Delta_x - x^2 \right)} \tilde{v}\right),\\
\tilde v(0,y,x)=\phi(y,x).
\end{array}\right.
\end{align*}
Denote by $\Pi_n$ the orthogonal projector on the $n^{th}$ eigenspace of $-\Delta_x + x^2$ (see Section \ref{sec:not and pre} below for more details). Applying $\Pi_n$ to the equation for $\tilde v$, we arrive at the following equation for $\tilde v_n=\Pi_n\tilde v$:
\begin{align*}
\left\{\begin{array}{l}\left(i\partial_t + \Delta_y \right) \tilde{v}_n = e^{i\lambda^2 t(2n+1)} \Pi_n\left(\sum\limits_{n_1,n_2,n_3 \in \mathbb{N} } e^{-i\lambda^2 \left(2n_1-2n_2 + 2n_3 + 1\right) t } \, \tilde{v}_{n_1} \bar{\tilde{v}}_{n_2} \tilde{v}_{n_3}\right),\\
\tilde v_n(0,y,x)=\phi_n(y,x):=\Pi_n\phi(y,x).\end{array}\right.
\end{align*}
Letting $\lambda\to \infty$, we can formally get a limiting equation
\begin{align}\label{eq1.5v31}
\begin{cases}
\left(i\partial_t + \Delta_y \right) v_n(t,y,x) =
 \sum\limits_{\substack{n_1,n_2,n_3 \in \mathbb{N} , \\n_1-n_2+ n_3= n} } \Pi_n\left(v_{n_1} \bar{v}_{n_2} v_{n_3} \right)(t,y,x), \\
v_n(0,y,x) = \phi_n(y,x),
\end{cases}
\end{align}
By reversing the above process, we get an approximation solution of $u_\lambda$:
\begin{align}\label{eq1.6v38}
\tilde{u}_\lambda(t,y,x) = e^{it(\Delta_x - x^2)} \sum_{n\in \mathbb{N} } \left( \frac1\lambda v_n\left(\frac{t}{\lambda^2}, \frac{y}\lambda, x\right) \right), \ (t,y,x)\in \mathbb{R}\times \mathbb{R}^2 \times \mathbb{R},
\end{align}
where $v_n$ is the solution of \eqref{eq1.5v31}.

In the above deduction, the following equivalent form dispersive continuous resonant (DCR) system enters naturally
\begin{align}\label{eq1.3v28}\tag{DCR}
\begin{cases}
i\partial_t v  + \Delta_{\mathbb{R}^2} v  = F(v) ,\\
v(0,y,x) =\phi(y,x),
\end{cases}
\end{align}
where the nonlinear term $F(v)$ is defined by
\begin{align*}
F(v) : = \sum\limits_{\substack{n_1,n_2,n_3,n \in \mathbb{N} ,\\  n_1-n_2 + n_3 = n}} \Pi_n \left(v_{n_1} \bar{v}_{n_2} v_{n_3} \right).
\end{align*}
This \eqref{eq1.3v28} system can be viewed as a dispersive version of the (CR) system derived by E. Faou, P. Germain, and Z. Hani \cite{FGH} in their study of the weak turbulence of the nonlinear Schr\"odinger equations on compact domains; see also \cite{BGHS2,CKSTT,GHT,GHT1,DELRV,Fe}.
This new \eqref{eq1.3v28} system is very similar to the resonant nonlinear Schr\"odinger system arising in \cite{BBCE,CGYZ,CGZ,HP,HPTV}. It has nice local well-posedness theory, and also scatters for small data in $L_y^2 \mathcal{H}_x^1$.

In our second main result, we prove the following large data global well-posedness and scattering theorem for (DCR), which might be of independent interest.
\begin{theorem}\label{as3.14v12}
	For any  $\phi\in L_y^2\mathcal{H}_x^1(\mathbb{R}^2\times \mathbb{R})$, there exists a unique global solution $v$ of the equation \eqref{eq1.3v28} in $C_t^0 L_y^2\mathcal{H}_x^1(\R\times\mathbb{R}^2\times \mathbb{R})$ satisfying
	\begin{align*}
	\|v\|_{L_t^\infty L_y^2 \mathcal{H}_x^1 \cap L_{t,y}^4 \mathcal{H}_x^1 (\mathbb{R}\times \mathbb{R}^2 \times \mathbb{R})} \le C,
	\end{align*}
    where $C=C\left(\|\phi \|_{L_y^2 \mathcal{H}_x^1} \right)$ is a constant. Moreover, the solution scatters, namely there exist $v_\pm \in L_y^2\mathcal{H}_x^1$ such that
	\begin{align*}
	\left\|v(t) - e^{it\Delta_y} v_\pm \right\|_{L_y^2 \mathcal{H}_x^1(\mathbb{R}^2 \times \mathbb{R})} \to 0, \text{ as } t \to \pm \infty.
	\end{align*}
\end{theorem}
Theorem \ref{as3.14v12} shall be proved in the final two sections and it takes a vast bulk of the paper. We prove it again by the concentration-compactness/rigidity argument from \cite{KM,KM1}. The system \eqref{eq1.3v28} is essentially a defocusing mass-critical nonlinear Schr\"odinger system. In the proof, we follow the framework for scattering of mass-critical nonlinear Schr\"odinger equation  \cite{D1,D2,D3,TVZ} and our argument is also partly inspired by the scattering of the resonant Schr\"odinger system derived from the NLS on cylinders \cite{CGZ,CGYZ,HP,HPTV,YZ0,Z1}.

We would like to comment briefly on the relation between \eqref{eq1.3v28} and weak turbulence.
\begin{remark}[The \eqref{eq1.3v28} system and weak turbulence]
We can rewrite \eqref{eq1.5v31} in the
Hermite coordinate (see \eqref{hn} below for the definition of the Hermite functions): Taking the solution $v_n(t,y,x) = c_n(t,y) h_n(x)$ in the equation \eqref{eq1.5v31}, we get an equivalent but simplified equation
\begin{align}\label{eq1.7v41}
\left(i\partial_t + \Delta_{\mathbb{R}^2} \right) c_n(t,y) = \sum\limits_{\substack{n_1,n_2,n_3 \in \mathbb{N} ,\\ n_1-n_2+ n_3 = n} } D_{n_1,n_2,n_3,n } c_{n_1} \bar{c}_{n_2} c_{n_3},
\end{align}
where $D_{n_1,n_2,n_3,n}$ is the number such that $\Pi_n \left(h_{n_1} \bar{h}_{n_2} h_{n_3} \right)(x) = D_{n_1,n_2,n_3,n} h_n(x)$, $x\in \R$. It would be very interesting to understand the constant $D_{n_1,n_2,n_3,n}$ in \eqref{eq1.7v41}. Comparing with the success of the proof of the weak turbulence on cylinders given by Z. Hani, B. Pausader, N. Tzvetkov, and N. Visciglia \cite{HPTV}, the unclear expression of the nonlinear term of the (DCR) system seems to be one of the main obstacles to study the weak turbulence of the nonlinear Schr\"odinger equations with (partial) harmonic potentials; for more information we refer to \cite{HT}. However, there are some interesting recent attempts toward this direction in \cite{GGT,GT}.
\end{remark}

\begin{remark}[Focusing NLS equations with harmonic potentials]
In this paper, we only consider the scattering of the defocusing NLS with partial harmonic potentials.
It is an interesting problem to study the scattering of the focusing version of \eqref{eq1.1}. It seems difficult to find the threshold of the scattering of the focusing NLS. On the other hand, if we were able to find the threshold of the scattering, then most likely the scattering can be proven by following the argument in \cite{D1,D2,D3,D4,KTV09,TVZ,TVZ0}. We refer to \cite{AC,BBJV,FL,Z3,SS} and the references therein for the study of the instability/stability of soliton which may give some clues on the threshold of the scattering of the focusing NLS.
\end{remark}

 \subsection{Brief outline of the proofs}
The model with partial harmonic potential studied in this paper can be compared to the NLS on wave-guide $\R^2\times\T$, which was considered previously in \cite{YZ,CGYZ}. One key difference is that in our case, the linear operator has more complicated spectral theory, for example the eigenfunctions cannot be written explicitly.

The proof of this paper contains two main ingredients. In the first part, we prove that Theorem \ref{as3.14v12} implies Theorem \ref{th1.3}. The proof of Theorem \ref{th1.3} has a very standard skeleton based on the concentration-compactness/rigidity argument introduced by C. Kenig and F. Merle \cite{KM}, and it consists of three main steps: linear profile decomposition, the existence of an almost periodic solution to the defocusing cubic NLS \eqref{eq1.1}, and a rigidity theorem.

First of all, we establish the linear profile decomposition of Schr\"odinger operator with partial harmonic potentials, namely the linear solutions can be divided into several orthogonal bubbles modulo some transforms. This can be viewed as a vector-valued version of linear profile decomposition of the Schr\"odinger equation in $L^2$, which was first established by F. Merle and L. Vega \cite{MV} in 2D, and then extended to general dimensions; see for instance \cite{Killip-Visan1} for more details. The proof of this part is very similar to the wave-guide case in \cite{CGYZ}, and it is essentially related to the description of the lack of compactness of the embedding
$e^{it \left( \Delta_{\mathbb{R}^3} - x^2 \right)}: \Sigma(\mathbb{R}^3) \hookrightarrow L_{t,y}^4 \mathcal{H}_x^{1-\epsilon_0},$
for some fixed $0 < \epsilon_0 < \frac12$.

In the second step, we prove the existence of a critical element by the construction of approximation solutions. Since the non-linear flow is not commutable with the transform groups derived in the first step, in order to construct the approximation solutions, we need to assume that the limiting equations, which is exactly the \eqref{eq1.3v28} system, is globally well-posed and scatters, as stated in Theorem \ref{as3.14v12}. The idea of using limiting equations was first considered in \cite{IMN}, and was widely used in \cite{CGYZ,CGZ,HP,IP,J1}. Then, similarly as in \cite{CGYZ}, we use the normal form method to exploit additional decay to approximate the non-linear profile. In the wave-guide case \cite{CGYZ}, the eigenfunctions, which are the plain waves $e^{i y\cdot j}$, can be easily computed, thus the Fourier coefficients are summed naturally. The  difficulty in this step is that we need to sum up the spectral projections of the solution properly. To some extent, the main innovation of this paper is that we utilize the additional regularity of the smooth non-linear profile to update the $l^1$ summation of projections to $l^2$.

In the third step, we borrow the idea used in \cite{CGYZ} to prove the non-existence of non-trivial critical element. The key point is the use of the interaction Morawetz estimate developed by J. Colliander, M. Keel, G. Staffilani, H. Takaoka and T. Tao \cite{CKSTT1}, which is very important in the remarkable work \cite{CKSTT0} on scattering for energy-critical NLS in 3-D, and was further developed in \cite{PV,CGT}. Then, we can arrive at the contradiction similar to \cite{KM,KM1} using the compactness property of the critical element.

The second part of this paper is devoted to the proof of Theorem \ref{as3.14v12}. The proof is greatly inspired by the fundamental work of B. Dodson \cite{D1,D2,D3} in his study of mass critical NLS. We also refer to \cite{YZ0}, and the principal difference between \cite{YZ0} and this paper is that our system \eqref{eq1.3v28} involves the spectral projection of Sch\"odinger operator with harmonic potential. Here, one key observation is that the \eqref{eq1.3v28} system is scaling invariant, which indicates that the classical method as developed in \cite{CGZ,CGYZ,TVZ} could be potentially applied to our situation. Indeed, the linear profile decomposition developed for the Schr\"odinger propagator in $L_y^2 \mathcal{H}_x^1(\mathbb{R}^2 \times \mathbb{R})$ (see Theorem \ref{pro3.9v23}) can be directly applied here. The essential difficulty occuring in the proof of Theorem \ref{as3.14v12} lies in precluding the almost periodic solution to the \eqref{eq1.3v28} system.

There are two cases of the critical element: high-to-low frequency cascade and the quasi-soliton scenarios. We exclude these scenarios based on the rigidity argument of B. Dodson \cite{D3,D1,D2}.  The key tool is to establish a vector-valued version of 2D long time Strichartz estimate in \cite{D1}.
The long time Strichartz estimate is developed by B. Dodson to show the scattering of the mass-critical nonlinear Schr\"odinger equations and has been proved as an important technique in the scattering theory of nonlinear dispersive and wave equation. We refer to \cite{D6,DL,DLMM,Killip-Visan2,Visan,Mu,BMMZ,R} for more application of this powerful tool. The proof of the long time Strichartz estimate in our situation here is rather technical due to the spectral projection and the failure of 2D end-point Strichartz estimate. For the high-to-low frequency cascade scenario, it is more delicate and we have to exploit some additional regularity of the critical element through the long time Strichartz, and then preclude it using energy conservation law. For the quasi-soliton scenario, we  mainly use the long time Strichartz to control the error terms of low frequency cut-off of interaction Morawetz identity. With all these ingredients at hand, the contradiction argument of C.E. Kenig and F. Merle \cite{KM,KM1} allows us to conclude the proof.

The rest of the paper is organized as follows. Section \ref{sec:not and pre} contains some basic notations and preliminaries.
In Section \ref{se2}, we record the local well-posedness, the small data scattering result and the stability theory for system \eqref{eq1.1}. For convenience of the readers, we present the proofs in the Appendix. In Section \ref{se4}, we will give the linear profile decomposition for data in $\Sigma(\mathbb{R}^3)$ and also analyze the nonlinear profiles, therefore we reduce the non-scattering in $\Sigma(\mathbb{R}^3)$ to the existence of an almost-periodic solution. In Section \ref{se4v16}, we will show the extinction of such an almost-periodic solution. The scattering of the \eqref{eq1.3v28} system shall be proved in Section \ref{se6v55}, where the proofs of two auxiliary theorems are left to the final Section \ref{sec:two theorems}.

\section{Basic notations and preliminaries}\label{sec:not and pre}
In this section, we introduce some basic notations used in this paper. We will use the notation $X\lesssim Y$ whenever there exists some constant $C>0$ so that $X \le C Y$. Similarly, we will write $X \sim Y$ if
$X\lesssim Y \lesssim X$. We use $\mathbb{N} $ to denote the set of all non-negative integers.

Throughout the paper, we will take $\epsilon_0$ to be some small fixed number in $\left(0,\frac{1}{2} \right)$.

\subsection{Fourier transform and Sobolev spaces}
For any $a\in \mathbb{R}^d,  d \in \mathbb{N}$, the Japanese bracket $\langle a \rangle$ is defined to be
$\langle a \rangle = \left( 1 + |a|^2 \right)^\frac12$.
We define the Fourier transform $\hat{f}\colon \mathbb{R}^d \to \mathbb{C}$ of a function $f\colon  \mathbb{R}^d \to \mathbb{C}$, as
$$\hat{f}(\xi) =\frac1{(2\pi)^{\frac{d}2 }} \int_{\mathbb{R}^d} e^{- i z\cdot \xi} f(z )\,\mathrm{d}z.$$
For each $s \in \mathbb{R}$,
the fractional differential operator $|\nabla|^s$ is defined by $\widehat{|\nabla|^s f}(\xi)  = |\xi|^s \hat{f}(\xi).$
We also define ${ \langle \nabla \rangle}^s$ as an operator between function spaces by $\widehat{{\langle \nabla \rangle}^s f}(\xi)=  \left( 1 + |\xi|^2 \right)^{\frac{s}2}  \hat{f}(\xi)$. In the following we will use $\langle\nabla_x\rangle^s$ to emphasize the application of the operator on the $x$-variable.

We will frequently use the partial Fourier transform $\mathcal{F}_y f$ of a complex-valued function $f\colon \mathbb{R}^2  \times \mathbb{R} \to \mathbb{C}$ defined as
\begin{align*}
\mathcal{F}_y f(\xi,x ) = \frac1{2\pi} \int_{\mathbb{R}^2 } e^{-iy\cdot \xi} f(y,x ) \,\mathrm{d}y, \ \xi\in\R^2,
\end{align*}
where $x\in \R$ is viewed as a parameter.

We shall also use the the Littlewood-Paley projections. Take a cut-off function $\chi\in C^{\infty}\left((0,\infty)\right)$ such that $\chi(r)=1$ if $r\loe1$ and $\chi(r)=0$ if $r>2$. For $N\in2^\mathbb{Z}$, let $\chi_N(r) = \chi(N^{-1}r)$ and $\phi_N(r) =\chi_N(r)-\chi_{N/2}(r)$. We define the Littlewood-Paley dyadic operator $P_{\loe N} f := \mathcal{F}^{-1}\brk{ \chi_N(|\xi|) \hat{f}(\xi)}$ and $P_N f := \mathcal{F}^{-1}\brk{ \phi_N(|\xi|) \hat{f}(\xi)}$.
We also define the partial Littlewood-Paley projections to be $P_{\le N}^y f(y,x) : = \mathcal{F}_y^{-1} \left( \chi_N(\xi) \left(\mathcal{F}_y f \right)(\xi,x) \right)$ and $P_N^y f(y,x): = \mathcal{F}_y^{-1} \left( \phi_N(|\xi|) \left( \mathcal{F}_y f \right)(\xi, x) \right)$.

Next, we denote the usual Lebesgue space as $L^p(\mathbb{R}^d)$, and some time we write $\norm{f}_{p} = \norm{f}_{L^p(\mathbb{R}^d)}$ for abbreviation. For any $s\in\mathbb{R}$, we define the Sobolev space as
\EQ{
W^{s,p}(\mathbb{R}^d) := \fbrk{f \in L^p(\mathbb{R}^d):\norm{f}_{W^{s,p}(\mathbb{R}^d)} := \norm{\jb{\nabla}^s f}_{L^p(\mathbb{R}^d)}<+\infty}.
}
We also define $H^s(\mathbb{R}^d)=W^{s,2}(\mathbb{R}^d)$.

\subsection{Harmonic oscillator and Hermite-Sobolev spaces}
The harmonic oscillator $-\Delta_x + x^2$, $x\in \R$, has been studied by many authors, and we refer to the lecture notes of B. Helffer \cite{He} and also the seminal work of H. Koch and D. Tataru \cite{KT} and the references therein for a few basic facts that we shall record below. The harmonic oscillator admits a Hilbertian basis of eigenvectors for $L^2(\mathbb{R})$, and  for each $n\in \mathbb{N}$, we will denote the $n^{th}$
eigenspace by $E_n$ and the corresponding eigenvalue by $\lambda_n = 2n+ 1$. %Denote by $K_n$ the dimension of $E_n$, where $K_n\sim 1$.
Each eigenspace $E_n$ is spanned by the
Hermite functions $h_{n}$, where
\begin{align}\label{hn}
h_n(x) = \frac1{\sqrt{n!} 2^\frac{n}2 \pi^\frac14} (-1)^n e^{\frac{x^2}2} \frac{d^n}{dx^n} (e^{-x^2}),
\end{align}
for $n\in \mathbb{N}$.
We also let $\Pi_n$ be the orthogonal projector on the $n^{th}$ eigenspace $E_n$ of $-\Delta_x + x^2$.

For $s\in\R$ and $p\geq 1$, the Hermite-Sobolev space $\mathcal{W}^{s,p}(\mathbb{R}) $ is defined as follows:
\begin{align*}
\mathcal{W}^{s,p}(\mathbb{R})
 = \left\{ u\in L^p_x(\mathbb{R}):
 \|u\|_{\mathcal{W}^{s,p}}:= \left\|\langle \nabla\rangle^s u \right\|_{L_x^p} + \left\| |\cdot|^s u\right\|_{L_x^p}<\infty \right\}.
\end{align*}
In particular, if $p=2$, we denote $\mathcal{W}_x^{s,2}(\mathbb{R})$ by $\mathcal{H}_x^s(\mathbb{R})$, and the $\mathcal{H}^1_x(\R)$-norm was given in \eqref{Sigma}.
By \cite{YZ}, we have
\EQ{
	\|u\|_{\mathcal{W}^{s,p}_x} \sim \norm{\brk{-\De + x^2}^{s/2}u}_p + \norm{u}_p.
}
The Hermite-Sobolev space $L_y^p \mathcal{H}_x^s$ with $1 \le p <  \infty$ and $s \in \mathbb{R}$  is defined by
\begin{align*}
L_y^p \mathcal{H}^s_x   =  \Bigg\{ f\in L^p_{y}L^2_x(\mathbb{R}^2 \times \mathbb{R}) &:
  \|f\|_{L_y^p \mathcal{H}_x^s (\mathbb{R}^2\times\mathbb{R})  } := \brk{ \int_{\mathbb{R}^2} \|f(y,\cdot)\|_{\mathcal{H}^s_x(\mathbb{R})}^p \mathrm{d}y}^{1/p}\\
  &= \brk{ \int_{\mathbb{R}^2} \bigg\| \Big(\sum\limits_{n \in \mathbb{N}}  (2n+1)^s |f_n(y,x )|^2 \Big)^\frac12 \bigg\|_{L^2_x(\mathbb{R})}^p \mathrm{d}y}^{1/p}
< \infty   \Bigg\},
\end{align*}
where $f_n=\Pi_nf$. Similarly, for any time interval $I \subseteq \mathbb{R}$ and $u: I \times \mathbb{R}^2 \times \mathbb{R} \to \mathbb{C}$, we define the space-time norms $L_t^p W_y^{s, q } L_x^r$ and $L_t^p L_y^q \mathcal{H}_x^s$ of $u$ as
\begin{align*}
\|u\|_{L_t^p W_{y}^{s,q}  L_x^r(I \times \mathbb{R}^2 \times \mathbb{R})} & : = \brk{\int_{I} \brk{\int_{\mathbb{R}^2} \left\|\langle \nabla_y \rangle^s  u(t,y,\cdot) \right\|_{L_x^r(\mathbb{R})}^q\mathrm{d}y}^{p/q}\mathrm{d}t}^{1/p}, \\
\|u\|_{L_t^p L_y^q \mathcal{H}_x^s(I\times \mathbb{R}^2\times \mathbb{R})}  & : =
\brk{\int_{I} \brk{\int_{\mathbb{R}^2} \left\| u(t,y,\cdot) \right\|_{\mathcal{H}_x^s(\mathbb{R})}^q\mathrm{d}y}^{p/q}\mathrm{d}t}^{1/p},
\end{align*}
where $1 \le p, q , r \leq  \infty$, and $s\in \mathbb{R}$. When $s=0$ and $p = q = r$, we  shall write $L_{t,y,x}^p$ for $L_t^p W_y^{s,p} L_x^r$. Similarly, when $p = q $, we shall write $L_{t,y}^p \mathcal{H}_x^s$ for $L_t^p L_y^q \mathcal{H}_x^s$.
We also use the following space-time norm. For any $\left\{ u_n (t,y,x) \right\}_{n \in \mathbb{N}}$, with $(t,y,x ) \in I \times \mathbb{R}^2 \times \mathbb{R}$, we set
\begin{align*}
\left\|u_n \right\|_{L_t^p L_y^q L_x^r l^2_n (I \times \mathbb{R}^2 \times \mathbb{R}\times \mathbb{N})}
= \norm{\norm{u_n}_{l_n^2}}_{L_t^p L_y^q L_x^r\brk{I \times \mathbb{R}^2 \times \mathbb{R}}},
\end{align*}
where $1\le p, q, r\le \infty$.
\begin{lemma}\label{le1.27}
The Dirac function $\delta_0(x)$ belongs to $\mathcal{H}_x^{-1}(\mathbb{R})$.
\end{lemma}
\begin{proof}
By definition, we have
\begin{align}\label{eq1.6v31}
\left\|\delta_0(x)\right\|_{\mathcal{H}^{-1}_x}^2
= \sum_{n=0}^\infty {(2n+1)^{-1 }} |c_n|^2,
\end{align}
where $c_n=\left\langle\delta_0(x),  h_n(x) \right\rangle
= h_n(0)$. Since $e^{-x^2} = \sum\limits_{m=0}^\infty \frac{(-x^2)^m}{m!} = \sum\limits_{n=0}^\infty \frac{d^n}{dx^n}\Big|_{x=0} e^{-x^2} \cdot \frac{x^n}{n!}$,
we have
\begin{align*}
\frac{d^n}{dx^n}\bigg|_{x=0} e^{-x^2} =
\begin{cases}
0,   & \text{ $n$ is odd},\\
\frac{(-1)^\frac{n}2}{(\frac{n}2)!} n!, & \text{ $n$ is even}.
\end{cases}
\end{align*}
Thus
\begin{align*}
h_n(0)
& = \begin{cases}
0,  &  \text{ $n$ is odd},\\
\frac{(-1)^n}{\sqrt{n!} 2^\frac{n}2 \pi^\frac14}  \frac{(-1)^\frac{n}2}{(\frac{n}2)!} n!,  & \text{ $n$ is even}.
\end{cases}
\end{align*}
Together with \eqref{eq1.6v31}, this implies
\begin{align*}
 \left\|\delta_0(x) \right\|_{\mathcal{H}_x^{-1}}^2 \le  \pi^{-\frac14} \sum_{\substack{ n=0, \\ \text{ n even} } }^\infty \frac{n!} {2^n  \left( \left( \frac{n}2 \right)! \right)^2 (2n+1) } \lesssim \sum\limits_{m = 0}^\infty \frac1{2^m(4m + 1)} \lesssim 1.
\end{align*}
\end{proof}

\section{Local well-posedness and small data scattering}\label{se2}
In this section, we will review the local well-posedness {theorem and} the stability theorem {for solutions of \eqref{eq1.1},} which {shall be} crucial in proving the existence of the critical element, and then record another important theorem on the scattering norm in Theorem \ref{th2.3}, which says that a weak space-time norm $L_{t,y} ^4  \mathcal{H}_x^{1- \epsilon_0 }$ is sufficient to prove the scattering result. We shall only state these results in this section and leave the proofs to Appendix A. In fact, the results in this section can be proved by following the exact arguments as in  \cite[Section 2]{CGYZ} or \cite{CGZ}, upon noticing the embedding $\mathcal{H}^{\frac12+}(\mathbb{R})\hookrightarrow L^\infty(\mathbb{R})$.

Different from the Strichartz estimate for the harmonic oscillator, which is a local estimate, we have a global Strichartz estimate for the partial harmonic oscillator similar to the Schr\"odinger equation on waveguides \cite{CGZ,CGYZ,Ta,TV}. Before giving the Strichartz estimate, we first introduce the following definition.
\begin{definition}[Strichartz admissible pair]
We call a pair $(p,q)$ is Strichartz admissible if $2 < p \le \infty$, $2\le q  < \infty$, and $\frac1p + \frac1q= \frac12$.
\end{definition}
We can now state the Strichartz estimate. The proof is almost identical to \cite[Proof of Proposition 2.1]{TV}, we also refer to Proposition
3.1 in \cite{ACS}, and we omit the proof here.
\begin{proposition}[Strichartz estimate for the partial harmonic oscillator]\label{pr2.2v31}
For any Strichartz admissible pair $(p,q)$, we have
\begin{align*}
\left\|e^{it\left( \Delta_{\mathbb{R}^3}  -x^2\right)  } f(y,x)\right\|_{L_t^p L_y^q L_x^2(\mathbb{R}\times \mathbb{R}^2 \times \mathbb{R})} & \lesssim \|f\|_{L_{y,x}^2}.
\end{align*}
Meanwhile, for $\alpha = 0,1$, it holds
\begin{align*}
\left\|e^{it \Delta_y } f(y,x ) \right\|_{L_t^p L_y^q \mathcal{H}^\alpha_x ( \mathbb{R} \times \mathbb{R}^2 \times \mathbb{R}) } \lesssim \|f \|_{L_y^2 \mathcal{H}^\alpha_x(\mathbb{R}^2 \times \mathbb{R})}.
\end{align*}
\end{proposition}

The following nonlinear estimate, which follows from the H\"older and Sobolev inequalities, is useful in showing the local well-posedness result.
\begin{proposition}[Nonlinear estimate] \label{pr3.3}
For any $0 < \epsilon_0 < \frac12$, we have
\begin{align*}
\left\| u_1 u_2 u_3  \right\|_{L_{t,y}^\frac43   \mathcal{H}_x^{1-\epsilon_0}}   \lesssim  \left\| u_1\right\|_{L_{t,y}^4  \mathcal{H}_x^{1-\epsilon_0}}  \left\|  u_2\right\|_{L_{t,y}^4  \mathcal{H}_x^{1-\epsilon_0}} \left\|  u_3\right\|_{L_{t,y}^4   \mathcal{H}_x^{1-\epsilon_0}}.
\end{align*}
\end{proposition}

Using Proposition \ref{pr2.2v31} and Proposition \ref{pr3.3}, one can easily prove the following local well-posedness and small data scattering in $L_y^2 \mathcal{ H}_x^1(\mathbb{R}^2 \times \mathbb{R}  )$ and $\Sigma (\mathbb{R}^3)$. The local solution can be extended to be global by the conservation of mass and energy, we refer to \cite{C,T2}. The proof of the local well-posedness is given in the Appendix; see also \cite{ACS,AC,C,C2,CG,C4,C5} for a comparison.

\begin{theorem}[LWP and scattering in $L_y^2 \mathcal{H}_x^1 $ and $\Sigma$]\label{th2.3}
\
\begin{enumerate}
\item (Well-posedness)  Let $u_0 \in L_y^2 \mathcal{H}_x^1 $, there exits a unique solution $u\in C_t^0  L_y^2 \mathcal{H}_{ x}^{ 1 }(I\times \mathbb{R}^2 \times \mathbb{R})  $ of \eqref{eq1.1}, where $I\subseteq \mathbb{R}$ is the maximal lifespan.
Furthermore, if $u_0\in \Sigma (\mathbb{R}^3  )$, the solution $u$ can be extended to be global in $ C_t^0 \Sigma_{y,x} (\mathbb{R} \times \mathbb{R}^2  \times \mathbb{R})$.
\item  (Scattering norm) If the solution $u\in C_t^0 \Sigma_{y,x} (\mathbb{R} \times \mathbb{R}^3)$ of \eqref{eq1.1} satisfies
$\|u\|_{L_{t,y} ^4  \mathcal{H}_x^{1- \epsilon_0 } (\mathbb{R}  \times \mathbb{R}^2\times \mathbb{R})}  \le M$ for some positive constant $M$. Then $u$ scatters in $\Sigma (\mathbb{R}^3)$, that is there exist $u_\pm \in \Sigma_{y,x} (\mathbb{R}^2  \times \mathbb{R})$ such that
\begin{equation}\label{eq3.1v131}
\left\|u(t,y,x)- e^{it\left(\Delta_{\mathbb{R}^3}-x^2\right)} u_\pm(y,x)\right\|_{\Sigma_{y,x}  } \to 0, \ \text{ as } t\to \pm \infty.
\end{equation}
\end{enumerate}
\end{theorem}
We now state the stability theory in $ L_y^2 \mathcal{H}_x^{1  }(\mathbb{R}^2  \times \mathbb{R})$. The proof is again given in the Appendix. For a comparison, see \cite{CKSTT0,Killip-Visan1,KTV}, in particular, \cite[Theorem 3.7]{Killip-Visan1}. We also contend that the result in the following theorem can extended to $\Sigma(\mathbb{R}^3)$.
\begin{theorem}[Stability theorem]\label{th2.6}
Let $I$ be a compact interval and let $\tilde{u}$ be an approximate solution to \eqref{eq1.1} in the sense that $\tilde{u}$ satisfies
$i\partial_t \tilde{u} + \Delta_{\mathbb{R}^3 } \tilde{u}  - x^2  \tilde{u} = |\tilde{u}|^2 \tilde{u} + e$
for some function $e$.

Suppose
\begin{align*}
\left\|\tilde{u}\right\|_{ L_t^\infty L_y^2 \mathcal{H}_x^{1  }\cap L_{t,y}^4   \mathcal{H}_x^{1 }   } \le M
\end{align*}
for some positive constant $M$.

Let $t_0\in I$ and let $u(t_0)$ obey
\begin{equation}\label{eq3.32}
\left\|u(t_0)- \tilde{u}(t_0)\right\|_{ L_y^2 \mathcal{H}_x^{1 } } \le M'
\end{equation}
for some $M'> 0$. Assume in addition that the smallness condition holds
\begin{align}
\left\|e^{i(t-t_0)\left(\Delta_{\mathbb{R}^3}-x^2\right) } \left(u(t_0)-\tilde{u}(t_0) \right)\right\|_{L_{t,y}^4   \mathcal{H}_x^{1 }   }  + \|e\|_{L_{t,y}^\frac43   \mathcal{ H}_x^{1}   }  \le \epsilon, \label{eq3.33}
\end{align}
for some $0 < \epsilon \le \epsilon_1$, where $\epsilon_1 = \epsilon_1(M,M') > 0$ is a small constant.
Then, there exists a solution $u$ to \eqref{eq1.1} on $I\times \mathbb{R}^2  \times \mathbb{R}$ with an
initial data $u(t_0)$ at time $t=t_0$ satisfying
\begin{align*}
 &\  \left\|u-\tilde{u}\right\|_{L_{t,y}^4  \mathcal{H}_x^{1}}   \le C(M,M')\epsilon,\quad
\left\|u-\tilde{u}\right\|_{L_t^\infty L_y^2 \mathcal{H}_{x}^{ 1 }  }    \le C(M,M')M',   \\
& \text{ and }  \|u\|_{L_t^\infty L_y^2 \mathcal{H}_{x}^{ 1  }  \cap L_{t,y}^4   \mathcal{H}_x^{1}
}   \le C(M,M').
\end{align*}
\end{theorem}

\section{Existence of an almost-periodic solution} \label{se4}
In this section, we will show the existence of an almost-periodic solution by the profile decomposition and the nonlinear approximation.

\subsection{Linear profile decomposition}
In this subsection, we will establish the linear profile decomposition in $\Sigma(\R^3)$, which depends on the corresponding decomposition in $L^2(\mathbb{R}^2)$.
The linear profile decomposition in $L^2$ for the mass-critical nonlinear Schr\"odinger equation has been established by F. Merle and L. Vega \cite{MV}, R. Carles and S. Keraani \cite{CK}, and P. B\'egout and A. Vargas \cite{BV}. We also refer readers to \cite{Killip-Visan1,KTV} for other versions of the linear profile decomposition.
\begin{theorem}[Linear profile decomposition in $L_y^2 \mathcal{H}_x^1(\mathbb{R}^2 \times \mathbb{R})$ and $\Sigma $]\label{pro3.9v23}
Let $ \left\{u_k \right\}_{k \ge 1}$ be a bounded sequence in $L_y^2 \mathcal{H}_x^1(\mathbb{R}^2 \times \mathbb{R})$. Then after passing to a subsequence if necessary,
there exists $J^\ast \in \{0, 1, \cdots \} \cup \{ \infty \}$, so that for any $J \le J^\ast$, we have functions $\phi^j \in L_y^2 \mathcal{H}_x^1(\mathbb{R}^2 \times \mathbb{R})$, $1\leq j\leq J$, $r_k^J \in L_y^2 \mathcal{H}_x^1(\mathbb{R}^2 \times \mathbb{R})$, and mutually orthogonal frames $\left\{ \left( \lambda_k^j, t_k^j, y_k^j, \xi_k^j \right) \right\}_{k  \ge 1} \subseteq \mathbb{R}_+  \times \mathbb{R} \times \mathbb{R}^2 \times \mathbb{R}^2$ in the sense that for any $j \ne j'$,
\begin{align}\label{eq4.1v103}
\frac{\lambda^j_k}{\lambda_k^{j'}} + \frac{\lambda_k^{j'}}{\lambda_k^j} + \lambda_k^j \lambda^{j'}_k \left|\xi_k^j - \xi_k^{j'}\right|^2 + \frac{ \left|y_k^j - y_k^{j'} \right|^2}{\lambda_k^j \lambda_k^{j'}}
+ \frac{ \left| \left(\lambda_k^j \right)^2 t_k^j - \left(\lambda_k^{j'} \right)^2 t_k^{j'} \right|}{\lambda_k^j \lambda_k^{j'}} \to \infty, \text{ as } k \to \infty,
\end{align}
such that, for every $1 \le j \le J$, we have a decomposition
\begin{align*}
u_k (y,x) = \sum_{j=1}^J \frac1{\lambda_k^j} e^{iy\cdot \xi_k^j} \left(e^{it_k^j \Delta_{\mathbb{R}^2}} \phi^j \right)\left(\frac{y-y_k^j}{\lambda_k^j}, x\right) + r_k^J(y,x).
\end{align*}
In addition,
\begin{align}
& \lim_{k \to \infty} \left(\left\|u_k \right\|_{L_y^2 \mathcal{H}_x^1 }^2 - \sum_{j=1}^J \left\|\phi^j\right\|_{L_y^2 \mathcal{H}_x^1  }^2 - \left\|r_k^J\right\|_{L_y^2 \mathcal{H}_x^1  }^2 \right) = 0, \label{eq4.2v103}
\\
& \lambda_k^j e^{-it_k^j \Delta_{y}} \Big(e^{-i \left(\lambda_k^j y + y_k^j \right)\xi_n^j} r_k^J\left(\lambda_k^j y + y_k^j, x \right)\Big)\rightharpoonup
 0, \text{ in } L_y^2 \mathcal{H}_x^1, \text{ as } k \to \infty, \text{ for } j\le J, \label{eq4.3v103}
 \\
& \limsup_{k \to \infty} \left\|e^{it\left(\Delta_{\mathbb{R}^3}-x^2 \right)} r_k^J \right\|_{L_{t,y}^4 \mathcal{H}^{1-\epsilon_0}_x } \to 0,
\text{ as } J \to J^*. \label{eq4.4v103}
\end{align}
Furthermore, if $ \left\{u_k \right\}_{k \ge 1}$ is a bounded sequence in $\Sigma (\mathbb{R}^3 )$, then in the above conclusion, we can further take $\lambda_k^j\to 1$ or $\infty$, as $k \to \infty$, $ \left|\xi_k^j \right| \le C_j$, for every $1\le j\le J$. And we have a slight different decomposition
\begin{align*}
u_k (y,x) = \sum_{j=1}^J   \phi_k^j  (y,x) + r_k^J(y,x)
:
= \sum_{j=1}^J \frac1{\lambda_k^j} e^{iy\cdot \xi_k^j} \left(e^{it_k^j \Delta_{\mathbb{R}^2}} P_k^j \phi^j \right)\left(\frac{y-y_k^j}{\lambda_k^j}, x\right) + r_k^J(y,x),
\end{align*}
where
\begin{align*}
P_k^j \phi^j (y,x) =
\begin{cases}
\phi^j(y,x),  & \text{ if $\lim\limits_{k \to \infty}  \lambda_k^j  = 1 $,   }
\\
P_{\le ( \lambda_k^j)^\theta }^y \phi^j (y,x),  & \text{ if  $  \lim\limits_{k \to \infty} \lambda_k^j = \infty$,  }
\end{cases}
\end{align*}
and
$\theta$ is some fixed positive sufficiently small number.
In addition, we also have a slight different decoupling
\begin{align}
  \lim_{k\to \infty} \left(  \mathcal{E}\left( u_k \right)    - \sum_{j=1}^J \mathcal{E} \left( \phi_k^j \right)   - \mathcal{E} \left( r_k^J \right) \right)  = 0, \label{eq4.5v115}
  \intertext{ and }
  \lim_{k \to \infty} \left( \mathcal{M} \left( u_k \right)     - \sum_{j=1}^J \mathcal{M} \left(  \phi_k^j \right)   - \mathcal{M} \left(   r_k^J \right)  \right) = 0, \label{eq4.6v115}
\end{align}
where $\mathcal{E}$ and $\mathcal{M}$ are given in \eqref{ME}. Other conclusions \eqref{eq4.1v103}-\eqref{eq4.4v103} hold as before.
\end{theorem}
To prove the above theorem, we need to establish the inverse Strichartz estimate in Proposition \ref{pr4.2536} below. We first recall the following refined Strichartz estimate which is essentially established in \cite{CGYZ,CGZ}.

\begin{proposition}[Refined Strichartz estimate, \cite{CGYZ,CGZ}] \label{pr4.2436}
For any $f \in L_y^2 \mathcal{H}_x^{1- \frac{\epsilon_0}2}$, we have
\begin{align*}
\left\|e^{it\Delta_{\mathbb{R}^2 }} f \right\|_{L_{t,y,x}^4(\mathbb{R}  \times \mathbb{R}^2 \times \mathbb{R})} \lesssim
\|f\|_{L_{y}^2 \mathcal{H}_x^{1-\frac{\epsilon_0}2}}^\frac34 \left(\sup_{Q\in \mathcal{D}} |Q|^{-\frac3{22}} \left\|e^{it\Delta_{\mathbb{R}^2}} f_Q \right\|_{L_{t,y,x}^\frac{11}2}\right)^\frac14,
\end{align*}
where
\begin{align*}
\mathcal{D} = \bigcup\limits_{j\in \mathbb{Z}} \left\{ \left[2^j k_1, 2^j(k_1+1) \right) \times \left[2^j k_2, 2^j(k_2+1) \right): (k_1,k_2) \in \mathbb{Z}^2 \right\}
\end{align*}
is the collection of all dyadic cubes, and $f_Q$ is defined by $\mathcal{F}_y  ({f}_Q) = \chi_Q \,\mathcal{F}_y {f} $.
\end{proposition}

To prove the inverse Strichartz estimate, we shall need the following two facts:
\begin{proposition}[Local smoothing estimate, \cite{CS,Ve}] \label{pr4.1436}
For any given $\epsilon > 0$, we have
\begin{align*}
\int_{\mathbb{R}} \int_{\mathbb{R}^2\times \mathbb{R}} \left|\left(|\nabla_y |^\frac12 e^{it\Delta_{\mathbb{R}^2}} f\right)(y,x)\right|^2 \langle y\rangle^{-1-\epsilon} \,\mathrm{d}y\mathrm{d}x \mathrm{d}t \lesssim_\epsilon \|f\|_{L_{y,x}^2(\mathbb{R}^2 \times \mathbb{R})}^2.
\end{align*}
Furthermore, if $\epsilon \ge 1$, then we have
\begin{align*}
\int_{\mathbb{R}} \int_{\mathbb{R}^2\times \mathbb{R}} \left|\left(\langle \nabla_y  \rangle ^\frac12 e^{it\Delta_{\mathbb{R}^2}} f \right)(y,x)\right|^2 \langle y\rangle^{-1-\epsilon} \,\mathrm{d}y\mathrm{d}x \mathrm{d}t \lesssim_\epsilon \|f\|_{L_{y,x}^2(\mathbb{R}^2 \times \mathbb{R})}^2.
\end{align*}
\end{proposition}

\begin{lemma}\label{le3.4v47}
For each $f\in \mathcal{H}_x^{1}(\mathbb{R})$ and any $R> 0$, we have
\begin{align*}
\|f\|_{L_x^\infty(|x|\ge R)} \lesssim R^{-\frac12} \left(\|f(x)\|_{L_x^2} + \|x f(x)\|_{L_x^2}^\frac12 \|  f'(x)\|_{L_x^2}^\frac12 \right).
\end{align*}
\end{lemma}
\begin{proof}
For any $f \in \mathcal{H}_x^{1}(\mathbb{R})$, we have
\begin{align*}
x f^2(x) = \int_0^x \left(z f^2(z) \right)' \,\mathrm{d}z = \int_0^x f^2(z) + 2 z f(z) f'(z)\,\mathrm{d}z,
\end{align*}
then by H\"older's inequality, we get
\begin{align*}
\left\|x f^2(x) \right\|_{L_x^\infty} \lesssim \|f\|_{L_x^2}^2 + \|x f(x) \|_{L_x^2} \| f'(x) \|_{L_x^2}.
\end{align*}
Therefore, for any $R>0$,
\begin{align*}
\|f(x)\|_{L_x^\infty(|x|\ge R)} \lesssim R^{-\frac12} \left(\|f\|_{L_x^2} + \|x f(x)\|_{L_x^2}^\frac12 \|f'(x)\|_{L_x^2}^\frac12\right).
\end{align*}
\end{proof}
We also have the following estimate.
\begin{lemma}\label{re3.441}
By interpolation, the H\"older inequality, the embedding $L_x^4(\mathbb{R})\hookrightarrow \mathcal{H}^{-1}(\mathbb{R})$,
and Proposition \ref{pr2.2v31}, we have
\begin{align}\label{eq4.7v139}
 %\limsup_{k \to \infty}
 \left\|e^{it \Delta_{y} } f \right\|_{L_{t,y}^4 \mathcal{H}_x^{1-\epsilon_0}}
\lesssim   \left\|e^{it\Delta_y} f \right\|_{L_{t,y}^4 \mathcal{H}_x^{-  1 }}^\frac{\epsilon_0}2
\left\| e^{it\Delta_y} f  \right\|_{L_{t,y}^4 \mathcal{H}_x^1}^{1-\frac{\epsilon_0}2}
\lesssim   \left\|e^{it\Delta_y}  f \right\|_{L_{t,y,x}^4}^\frac{\epsilon_0}2 \left\|
f  \right\|_{L_{y}^2 \mathcal{H}_x^1 }^{1-\frac{\epsilon_0}2} .
\end{align}
\end{lemma}
We can now prove the inverse Strichartz estimate.
\begin{proposition}[Inverse Strichartz estimate] \label{pr4.2536}
For $\left\{f_k \right\}_{k \ge 1} \subseteq L_y^2 \mathcal{H}_x^{1}(\mathbb{R}^2 \times \mathbb{R} )  $ satisfying
\begin{align}\label{eq3.948}
\lim_{k \to \infty} \left\|f_k \right\|_{L_y^2 \mathcal{H}_x^{1} } = A \quad\text{ and }\quad \lim_{k \to \infty} \left\|e^{it\Delta_{\mathbb{R}^2}} f_k \right\|_{L_{t,y}^4 \mathcal{H}_x^{1- \epsilon_0} } = \epsilon,
\end{align}
there exist $\phi \in L_y^2 \mathcal{H}_x^{1} $ and $ \left(\lambda_k,t_k,  \xi_k,  y_k \right) \in \mathbb{R}_+ \times \mathbb{R}  \times \mathbb{R}^2 \times \mathbb{R}^2$, so that passing to a further subsequence of if necessary, we have
\begin{align}
& \lambda_k  e^{-i\xi_k \cdot (\lambda_k   y + y_k ) } \left(e^{it_k \Delta_{\mathbb{R}^2}} f_k \right)\left(\lambda_k  y + y_k , x \right) \rightharpoonup  \phi(y,x)\  \text{ in } L_y^2 \mathcal{H}_x^{1} , \text{ as } k \to \infty, \notag\\
 & \lim_{k \to \infty} \left(\|f_k  \|_{L_y^2 \mathcal{H}_x^{1}}^2 - \|f_k -\phi_k \|_{L_y^2 \mathcal{H}_x^{1}}^2 \right) = \|\phi\|
_{L_y^2 \mathcal{H}_x^{1}}^2
\gtrsim A^2 \left(\frac\epsilon A\right)^{\frac{48}{\epsilon_0} },\label{eq4.4936} \\
&  \limsup_{k \to \infty} \left\|e^{it\Delta_{\mathbb{R}^2}} (f_k -\phi_k ) \right\|_{L_{t,y}^4 \mathcal{H}_x^{1- \epsilon_0} }^4 \le \epsilon^{\frac8{\epsilon_0} } A^{4 - \frac8{\epsilon_0}}  \left( 1- c A^{2\beta} \left( \frac{ \epsilon}{A} \right)^{\frac{2\beta}{\epsilon_0} } \right), \label{eq4.5036}
\end{align}
where $c$ and $\beta$ are small positive constants, and
\begin{align*}
\phi_k (y,x) = \frac1{\lambda_k } e^{iy\cdot \xi_k } \left(e^{-i\frac{t_k }{\lambda_k^2} \Delta_{\mathbb{R}^2}} \phi\right)\left(\frac{y-y_k }{\lambda_k },x\right).
\end{align*}
Moreover, if $\{f_k \}_{k \ge 1} $ is bounded in $\Sigma (\mathbb{R}^3)$, and also
\begin{align}\label{eq4.8v115}
\lim_{k \to \infty} \left\|f_k \right\|_{\Sigma } = A \quad\text{ and }\quad \lim_{k \to \infty} \left\|e^{it\Delta_{\mathbb{R}^2}} f_k \right\|_{L_{t,y}^4 \mathcal{H}_x^{1- \epsilon_0} } = \epsilon,
\end{align}
then we can take $\lambda_k  \ge 1$, $|\xi_k | \lesssim 1$ and
$\phi \in L_y^2 \mathcal{H}_x^{1 }(\mathbb{R}^2 \times \mathbb{R}) $ such that
\begin{align}
\lambda_k  e^{-i\xi_k\cdot (\lambda_k  y + y_k ) } \left(e^{it_k \Delta_{\mathbb{R}^2}} f_k  \right)\left(\lambda_k  y+ y_k , x \right) \rightharpoonup  \phi(y,x)\  \text{ in } L_y^2 \mathcal{H}_x^{1 } , \text{ as } k \to \infty, \label{eq4.9v115}
 \\
\intertext{ and }
\lim_{k \to \infty} \left( \left\|f_k \right\|_{\Sigma}^2 - \left\|f_k -\phi_k \right\|_{\Sigma}^2 \right) = \lim\limits_{k \to \infty} \|\phi_k \|
_{\Sigma }^2
\gtrsim A^2 \left(\frac\epsilon A\right)^{\frac{48}{\epsilon_0} }.\label{eq3.1348}
\end{align}
\end{proposition}
\begin{proof}
{\it Case 1. $ \left\{f_k \right\}_{k \ge 1} $ is bounded in $L_y^2 \mathcal{H}_x^1$.}
By Proposition \ref{pr4.2436}, \eqref{eq4.7v139} and \eqref{eq3.948}, there exists $ \left\{Q_k \right\}_{k \ge 1}  \subseteq \mathcal{D}$ so that
\begin{align}\label{eq4.5236}
  \epsilon^\frac8{\epsilon_0} A^{1 - \frac8{\epsilon_0} }   \lesssim \liminf_{k \to \infty}  \left|Q_k \right|^{-\frac3{22}} \left\|e^{it\Delta_{\mathbb{R}^2}} \left(f_k \right)_{Q_k }\right\|_{L_{t,y,x}^\frac{11}2}.
\end{align}
Let $\lambda_k  $ be the inverse of the side-length and $\xi_k $ be the center of the cube $Q_k $. By H\"older's inequality and \eqref{eq3.948}, we have
\begin{align*}
\liminf_{k \to \infty} \left|Q_k \right|^{-\frac3{22}} \left\|e^{it\Delta_{\mathbb{R}^2}} (f_k )_{Q_k } \right\|_{L_{t,y,x}^\frac{11}2}
 \lesssim \liminf_{k \to \infty} \lambda_k^\frac3{11} \left( \epsilon^\frac2{\epsilon_0} A^{1 - \frac2{\epsilon_0} } \right)^\frac8{11} \left\|e^{it\Delta_{\mathbb{R}^2}} \left(f_k \right)_{Q_k } \right\|_{L_{t,y,x}^\infty}^\frac3{11}.
\end{align*}
Together with \eqref{eq4.5236}, this implies
\begin{align*}
\liminf_{k \to \infty} \lambda_k  \left\|e^{it\Delta_{\mathbb{R}^2}} \left(f_k \right)_{Q_k } \right\|_{L_{t,y,x}^\infty(\mathbb{R}\times \mathbb{R}^2\times \mathbb{R})} \gtrsim   \epsilon^\frac{24}{\epsilon_0} A^{1 - \frac{24}{\epsilon_0} }  .
\end{align*}
Then by Lemma \ref{le3.4v47} and Bernstein's inequality, we have
\begin{align*}
& \liminf_{k \to \infty} \lambda_k  \left\|e^{it\Delta_{\mathbb{R}^2}} \left(f_k \right)_{Q_k } \right\|_{L_{t,y,x}^\infty(\mathbb{R}\times \mathbb{R}^2\times \{  |x|\ge R \} )}\\
\lesssim & R^{-\frac12} \liminf_{k  \to \infty} \lambda_k  \left(\left\| \left(f_k \right)_{Q_k } \right\|_{L_{t,y}^\infty L_x^2} + \left\||x| \left(f_k \right)_{Q_k }(x)  \right\|_{L_{t,y}^\infty L_x^2}^\frac12
\left\|\partial_x \left( \left(f_k \right)_{Q_k }  \right) \right\|_{L_{t,y}^\infty L_x^2}^\frac12\right)\\
\lesssim & R^{-\frac12} \liminf_{k \to \infty} \lambda_k  \left( \left|Q_k \right|^\frac12 \|f_k \|_{L_t^\infty L_{y,x}^2  }
+ |Q_k |^\frac12 \|x f_k \|_{L_t^\infty L_{y,x}^2  }^\frac12 \left\|\partial_x f_k  \right\|_{L_t^\infty L_{y,x}^2}^\frac12\right)\\
\sim & R^{-\frac12} \liminf_{k \to \infty} \left(\|f_k \|_{L_t^\infty L_{y,x}^2  } + \|x f_k \|_{L_t^\infty L_{y,x}^2 }^\frac12 \left\|\partial_x f_k \right\|_{L_t^\infty L_{y,x}^2}^\frac12\right) \to 0, \ \text{ as } R\to \infty.
\end{align*}
Therefore, we can take $R$ large enough such that
\begin{align*}
\liminf_{k \to \infty} \lambda_k  \left\|e^{it\Delta_{\mathbb{R}^2}} \left(f_k \right)_{Q_k } \right\|_{L_{t,y,x}^\infty(|x|\ge R)} \lesssim \frac12
  \epsilon^\frac{24}{\epsilon_0} A^{1- \frac{24}{\epsilon_0} }  .
\end{align*}
As a consequence, there exists $\left(t_k,y_k ,x_k \right) \in \mathbb{R} \times \mathbb{R}^2 \times \mathbb{R}$ with $ \left|x_k \right|\le R$, so that
\begin{align}\label{eq4.5336}
\liminf_{k \to \infty} \lambda_k  \left|\left(e^{it_k  \Delta_{\mathbb{R}^2}} \left(f_k \right)_{Q_k } \right)(y_k ,x_k )\right|\gtrsim 
 \epsilon^\frac{24}{\epsilon_0} A^{1- \frac{24}{\epsilon_0} }  .
\end{align}
Since $ \left|x_k \right| \le R$, we may assume, up to a subsequence, $x_k \to x^*$, as $k \to \infty$, with $|x^*|\lesssim 1$.

By the weak compactness of $L_y^2 \mathcal{H}_x^{1}$, we have
\begin{align*}
\lambda_k  e^{-i\xi_k (\lambda_k  y+y_k )} \left(e^{it_k  \Delta_{\mathbb{R}^2}} f_k \right)\left(\lambda_k  y + y_k , x \right) \rightharpoonup  \phi(y,x) \text{ in } L_y^2 \mathcal{H}_x^{1}, \text{ as } k \to \infty.
\end{align*}
By the very basic fact in Hilbert space $H$ that
\begin{align*}
g_k \rightharpoonup g \text{ in } H \Rightarrow \|g_k \|_H^2 - \|g_k - g\|_H^2 \to \|g\|_H^2,
\end{align*}
we have
\begin{align*}
\lim_{k \to \infty} \left(\|f_k  \|_{L_y^2 \mathcal{H}_x^{1}}^2 - \|f_k -\phi_k \|_{L_y^2 \mathcal{H}_x^{1}}^2 \right) = \|\phi\|
_{L_y^2 \mathcal{H}_x^{1}}^2.
\end{align*}
We now turn to the remaining part \eqref{eq4.4936}. Define $h$ so that $\mathcal{F}_y {h}$ is the characteristic function of the cube $ \left[-\frac12,\frac12 \right]^2$. By Lemma \ref{le1.27},  the function $(x,y)\mapsto h(y) \delta_0(x)\in L_y^2 \mathcal{H}_x^{-1}(\mathbb{R}^2 \times \mathbb{R})$.
From \eqref{eq4.5336}, we obtain
\begin{align}\label{eq4.5436}
 \left|\left\langle h(y) \delta_0(x), \phi \left(y,x+ x^* \right)\right\rangle_{y,x} \right|
= &\lim_{k \to \infty} \left|\left\langle \delta_0(x), \int_{\R^2}\bar{h}(y) \lambda_k  e^{-i\xi_k\cdot (\lambda_k  y + y_k )} \left(e^{it_k  \Delta_{\mathbb{R}^2}} f_k  \right) \left(\lambda_k  y +y_k , x+ x_k \right) \,\mathrm{d}y \right\rangle_x \right| \nonumber\\
=  & \lim_{k \to \infty} \lambda_k  \left|\left(e^{it_k \Delta_{\mathbb{R}^2}} \left(f_k \right)_{Q_k }\right)(y_k ,x_k )\right|
\gtrsim   \epsilon^\frac{24}{\epsilon_0} A^{1- \frac{24}{\epsilon_0} }  ,
\end{align}
from which it follows
\begin{align*}
\left\|\phi \left(y,x+ x^* \right) \right\|_{L_y^2 \mathcal{H}_x^1} \gtrsim \epsilon^\frac{24}{\epsilon_0} A^{1- \frac{24}{\epsilon_0} }  .
\end{align*}
At the same time, since
\begin{align*}
\|\phi(y,x)\|_{L_y^2 \mathcal{H}_x^1}
& \ge \left\|\phi \left(y,x+x^* \right) \right\|_{L_y^2 H_x^1} + \left\||x| \phi \left(y,x + x^* \right) \right\|_{L_{y,x}^2}-  \left\| \left|x^* \right| \phi \left(y,x + x^* \right) \right\|_{L_{y,x}^2}\\
& = \left\|\phi \left(y,x+x^* \right) \right\|_{L_y^2 \mathcal{H}_x^1} -  \left|x^* \right| \| \phi(y,x)\|_{L_{y,x}^2},
\end{align*}
we get
 \begin{align*}
 \left\|\phi \left(y,x+x^* \right)  \right\|_{L_y^2 \mathcal{H}_x^1}\le  \|\phi\|_{L_y^2 \mathcal{H}_x^1} +  \left|x^* \right| \|\phi\|_{L_{y,x}^2} \lesssim \|\phi\|_{L_y^2 \mathcal{H}_x^1}.
\end{align*}
Therefore $\|\phi\|_{L_y^2 \mathcal{H}_x^1} \gtrsim    \epsilon^{\frac{24}{\epsilon_0} } A^{1- \frac{24}{\epsilon_0} }  $ and \eqref{eq4.4936} follows. We turn to \eqref{eq4.5036}, by Proposition \ref{pr4.1436} and the Rellich-Kondrashov theorem, we have
\begin{align*}
e^{it\Delta_{\mathbb{R}^2}} \left( \lambda_k  e^{-i \xi_k\cdot (\lambda_k  y + y_k ) } (e^{it_k \Delta_{\mathbb{R}^2}} f_k )(\lambda_k  y + y_k , x + x_k )\right) \to
e^{it\Delta_{\mathbb{R}^2}} \phi(y,x), \text{ as $k \to \infty$},
\end{align*}
for almost every $(t,y,x) \in \mathbb{R} \times \mathbb{R}^2\times \mathbb{R}$.
By the refined Fatou's lemma \cite{Le}, we obtain
\begin{align*}
\left\|e^{it \Delta_{\mathbb{R}^2}} f_k \right\|_{L_{t,y}^4 \mathcal{H}_x^{1- \epsilon_0} }^4 - \left\|e^{it\Delta_{\mathbb{R}^2}} (f_k -\phi_k )\right\|_{L_{t,y}^4 \mathcal{H}_x^{1- \epsilon_0} }^4 - \left\|e^{it\Delta_{\mathbb{R}^2}} \phi_k  \right\|_{L_{t,y}^4 \mathcal{H}_x^{1- \epsilon_0} }^4 \to 0, \text{ as } k \to \infty.
\end{align*}
Thus, by the invariance of Galilean transform, we have
\begin{align}\label{eq4.16v133}
\limsup\limits_{k \to \infty} \left\|e^{it\Delta_{\mathbb{R}^2}} (f_k -\phi_k )\right\|_{L_{t,y}^4 \mathcal{H}_x^{1- \epsilon_0} }^4
= &  \limsup\limits_{k \to \infty} \left(  \left\|e^{it \Delta_{\mathbb{R}^2}} f_k \right\|_{L_{t,y}^4 \mathcal{H}_x^{1- \epsilon_0} }^4  - \left\|e^{it\Delta_{\mathbb{R}^2}} \phi_k
 \right\|_{L_{t,y}^4 \mathcal{H}_x^{1- \epsilon_0} }^4 \right)\\
= &  \left( \epsilon^\frac2{\epsilon_0} A^{1- \frac2{\epsilon_0} } \right)^4 -   \left\|e^{it\Delta_{\mathbb{R}^2}} \phi   \right\|_{L_{t,y}^4 \mathcal{H}_x^{1- \epsilon_0} }^4. \notag
\end{align}
We now take $c(t) \in C^\infty$ which has compact support on $[-1,1]$, such that
\begin{align*}
\|c(t) e^{it \Delta} h \|_{L_{t,y}^\frac43} = 1.
\end{align*}
Then by \eqref{eq4.5436}, we have
\begin{align*}
& \left|\int_{\mathbb{R}} \left\langle c(t) h(y) \delta_0(x), \phi \left(y,x+ x^* \right)\right\rangle_{y,x} \,\mathrm{d}t  \right|
\gtrsim  \epsilon^\frac{24}{\epsilon_0} A^{1- \frac{24}{\epsilon_0}}  .
\end{align*}
On the other hand, by H\"older's inequality, Sobolev's inequality and Lemma \ref{le1.27},
\begin{align*}
 \left| \int_{\mathbb{R}} \left\langle c(t) h(y) \delta_0(x), \phi \left(y, x+ x^* \right) \right\rangle_{y,x} \,\mathrm{d}t \right|
= &  \left|\int_{\mathbb{R}}  \left\langle e^{it \Delta_y} \left( c(t) h(y) \delta_0(x) \right) , e^{it \Delta_y} \phi \left(y, x+x^* \right) \right\rangle_{y,x} \,\mathrm{d}t \right|
\\
\lesssim &  \left\|e^{it \Delta_y} \left( c(t) h(y)  \right) \right\|_{L_{t,y}^\frac43   }  \left\|e^{it \Delta_y} \phi \left(y, x \right) \right\|_{L_{t,y}^4 \mathcal{H}_x^{1- \epsilon_0} } \lesssim  \left\|e^{it \Delta_y} \phi(y,x) \right\|_{L_{t,y}^4 \mathcal{H}_x^{1- \epsilon_0} }.
\end{align*}
Therefore, by the above two estimates and \eqref{eq4.16v133}, we get \eqref{eq4.5036}.

{\it Case 2. $ \left\{f_k \right\}_{k \ge 1} $ is bounded in $\Sigma (\mathbb{R}^3)$.} In this case, we have
\begin{align*}
\limsup\limits_{k \to \infty} \left\|P_{\ge R}^y f_k \right\|_{L_y^2 \mathcal{H}_x^{1- {\epsilon_0} }}
 & \lesssim \langle R\rangle^{- {\epsilon_0} }  \limsup\limits_{k \to \infty} \left\|f_k \right\|_{\Sigma(\mathbb{R}^3)  } \to 0, \text{ as  }  R\to \infty.
\end{align*}
For $R\in 2^{\mathbb{Z}}$ large enough depending on $A$ and $\epsilon$, by \eqref{eq4.8v115}, Sobolev embedding, and Strichartz estimate, $P_{\le R}^y f_k $ satisfies
\begin{align*}
\lim\limits_{k \to \infty} \left\|e^{it\Delta_{\mathbb{R}^2}} P_{\le R}^y f_k \right\|_{L_{t,y}^4 \mathcal{H}_x^{1- \epsilon_0} } & \ge \lim\limits_{k \to \infty}
\left\|e^{it\Delta_{\mathbb{R}^2}} f_k \right\|_{L_{t,y}^4 \mathcal{H}_x^{1- \epsilon_0} } - \lim\limits_{k \to \infty} \left\|e^{it\Delta_{\mathbb{R}^2}} P_{\ge R}^y  f_k \right\|_{L_{t,y}^4 \mathcal{H}_x^{1- \epsilon_0} }\\
& \ge  \lim\limits_{k \to \infty}\left\|e^{it\Delta_{\mathbb{R}^2}} f_k \right\|_{L_{t,y}^4 \mathcal{H}_x^{1- \epsilon_0} } - C   \lim\limits_{k \to \infty} \left\| P_{\ge R}^y f_k \right\|_{L_y^2 \mathcal{H}_x^{1- {\epsilon_0} }} \ge \frac12  \epsilon^{\frac2{\epsilon_0} } A^{1- \frac2{\epsilon_0} }  .
\end{align*}
So we can replace $f_k $ by $P_{\le R}^y f_k $ in the above case, and for $R = R(A,\epsilon) > 0$ large enough, we may take $ \left\{Q_k \right\}_{k  \ge 1}  \subseteq \mathcal{D}$ and $ \left|Q_k\right |\lesssim R^2$ such that $\lambda_k  \gtrsim R^{-1}$, and $|\xi_k |\lesssim R$.
{As in the proof of {\it Case 1}, we still have \eqref{eq4.9v115} and also \eqref{eq4.4936}, \eqref{eq4.5036}. Furthermore, if $\limsup\limits_{k \to \infty} \lambda_k < \infty$, then
\begin{align*}
\lambda_k  e^{-i\xi_k \cdot(\lambda_k  y + y_k )} \left(e^{it_k  \Delta_{\mathbb{R}^2}} f_k\right ) \left(\lambda_k  y + y_k , x \right) \rightharpoonup  \phi(y,x)
\end{align*}
holds for some $\phi \in \Sigma (\mathbb{R}^2 \times \mathbb{R})$}. To show \eqref{eq3.1348}, we just need to consider the case when $\lambda_k \to \infty$ because the situation when $\limsup \limits_{k \to \infty } \lambda_k < \infty$ is as in {\it Case 1}. We note
\begin{align*}
{ \lim\limits_{k \to \infty}\left \|\phi_k \right\|_{\Sigma}^2 \ge \lim\limits_{k \to \infty} \left\|P^y_{\le \lambda_k^\theta} \phi \right\|_{L_y^2 \mathcal{H}_x^1}^2 }\gtrsim A^2 \left(\frac\epsilon A\right)^{\frac{48}{\epsilon_0} }.
\end{align*}
Then the decoupling of the $\Sigma-$norm comes from { $P_{\leq \lambda_k^\theta} \to Id$} in $L_y^2 \mathcal{H}_x^1$ and  \eqref{eq4.9v115}.
\end{proof}
\begin{proof}[Proof of Theorem \ref{pro3.9v23}]
The conclusion follows by applying Proposition \ref{pr4.2536} repeatedly until the asymptotically linear evolution of the remainder is trivial in $L_{t,y}^4 \mathcal{H}_x^{1- \epsilon_0} $. The decoupling \eqref{eq4.5v115} and \eqref{eq4.6v115} follow from \eqref{eq3.1348} and the orthogonality \eqref{eq4.1v103}.
\end{proof}
\begin{remark}
For a linear profile decomposition for the Schr\"odinger propagator of the Schr\"odinger operator $-\Delta + |x|^2$ in $L^2$, we refer to the work of C. Jao, R. Killip, and M. Visan \cite{JKV} and C. Jao \cite{J4}, we believe that some part of their argument can be applied in our equation. We also refer to the linear profile decomposition proved by A. Ardila and R. Carles \cite{AC}.
\end{remark}

\subsection{Approximation of the nonlinear profile - The case of concentrated initial data}
In this section, we will show that
the nonlinear profile $u_\lambda$ given in \eqref{eq1.5v18}
\begin{align*}
\begin{cases}
i\partial_t u_\lambda + \Delta_{\mathbb{R}^3} u_\lambda - x^2 u_\lambda = |u_\lambda|^2 u_\lambda,\\
u_\lambda(0,y,x) = \frac1\lambda \phi(\frac{y}\lambda, x),
\end{cases}
\end{align*}
can be approximated by $\tilde u_\lambda$ given in  \eqref{eq1.6v38}
\begin{align*}
\tilde u_\lambda(t,y,x)   = e^{it\left(\Delta_{\mathbb{R}} -x^2\right)} \sum_{n\in \mathbb{N} } \left( \frac1\lambda v_n\left(\frac{t}{\lambda^2}, \frac{y}\lambda, x\right) \right) , \ (t,y,x)\in \mathbb{R}\times \mathbb{R}^2 \times \mathbb{R},
\end{align*}
when $\lambda$ is sufficiently large.
Here $v_n$ is the solution of the (DCR) system \eqref{eq1.5v31}
\begin{align*}
\begin{cases}
\left(i\partial_t + \Delta_y \right) v_n(t,y,x) =
 \sum\limits_{\substack{n_1,n_2,n_3 \in \mathbb{N} , \\n_1-n_2+ n_3= n} } \Pi_n\left(v_{n_1} \bar{v}_{n_2} v_{n_3} \right)(t,y,x), \\
v_n(0,y,x) = \phi_n(y,x)= \Pi_n\phi(y,x).
\end{cases}
\end{align*}
The following corollary can be proven from Theorem \ref{as3.14v12} by following the argument in \cite{CKSTT0,Killip-Visan1}. In particular, we refer to \cite[Lemma 3.12]{CKSTT1}.

\begin{corollary}[Corollary of Theorem \ref{as3.14v12}: Preservation of higher regularity]\label{co4.7v65}
Suppose that $\phi \in L_y^2 \mathcal{H}_x^1 (\mathbb{R}^2 \times \mathbb{R})$ and $v$ is the global solution of \eqref{eq1.3v28} given as in Theorem \ref{as3.14v12}. For any $s_1\ge 0 $ and $ s_2 \ge 1 $, if we assume further $v|_{t=0}\in H_y^{s_1 } \mathcal{H}_x^{s_2} (\mathbb{R}^2\times \mathbb{R})$, then the solution $v\in C_t^0 H_y^{s_1 } \mathcal{H}_x^{s_2} (\mathbb{R}\times \mathbb{R}^2 \times \mathbb{R})$ and satisfies
\begin{align*}
\|v\|_{L_t^\infty H_y^{s_1 } \mathcal{H}_x^{s_2}  \cap L_{t}^4 W_y^{s_1 ,4} \mathcal{H}_x^{s_2}  (\mathbb{R}\times \mathbb{R}^2 \times \mathbb{R})} \le C\left(\|\phi \|_{H_y^{s_1 } \mathcal{H}_x^{s_2} (\mathbb{R}^2\times \mathbb{R})} \right).
\end{align*}
\end{corollary}

Relying on Corollary \ref{co4.7v65}, we can now prove the  following general result on approximation of the non-linear profile in the large-scale case. We will prove it with the help of Theorem \ref{th2.6}.
\begin{theorem}\label{pr3.12}
For any $\phi \in L_y^2 \mathcal{H}_x^1$, $0 < \theta << 1$, $( \lambda_k, t_k, y_k, \xi_k) \in \mathbb{R}_+ \times \mathbb{R} \times \mathbb{R}^2 \times \mathbb{R}^2$, $ | \xi_k  | \lesssim 1$ and $\lambda_k  \to \infty$ when $k  \to \infty$. There exists a global solution $u_k \in C_t^0 L_y^2 \mathcal{H}_x^1$ of
\begin{align*}
\begin{cases}
i \partial_t u_k + \Delta_y u_k + \Delta_x u_k - x^2 u_k  = |u_k|^2 u_k , \\
u_k (0,y,x) = \lambda_k^{-1} e^{iy\cdot \xi_k } \left( e^{it_k  \Delta_y} { P_{\le \lambda_k^\theta} \phi} \right) \left( \frac{ y - y_k}{\lambda_k }, x \right),
\end{cases}
\end{align*}
for $k $ large enough, satisfying
\begin{align*}
\|u_k  \|_{L_t^\infty L_y^2 \mathcal{H}_x^1 \cap L_{t,y}^4 \mathcal{H}_x^1( \mathbb{R}\times \mathbb{R}^2 \times \mathbb{R})} \lesssim_{\|\phi\|_{L_y^2 \mathcal{H}_x^1} }
1.
\end{align*}
Furthermore, assume that $\epsilon_4 = \epsilon_4 \left( \|\phi\|_{L_y^2 \mathcal{H}_x^1} \right) $ is a sufficiently small positive constant and $\psi  \in H_y^{10 } \HH_x^{10} $ such that
\begin{align*}
\| \phi - \psi  \|_{L_y^2 \mathcal{H}_x^1} \le \epsilon_4.
\end{align*}
Then there exists a solution $v \in C_t^0 H_y^2 \mathcal{H}_x^1( \mathbb{R} \times \mathbb{R}^2 \times \mathbb{R})$ of (DCR), with
\begin{align*}
v(0, y ,x )  = \psi(y,x) ,   & \text{ if } t_k = 0, \\
\lim\limits_{t\to \pm \infty } \| v(t,y,x) - e^{it \Delta_y} \psi \|_{L_y^2 \mathcal{H}_x^1}  =  0,  & \text{ if  } t_k \to \pm \infty,
\end{align*}
such that for $k $ large enough, we have
$ \left\|u_k  \right\|_{L_t^\infty L_y^2 \mathcal{H}_x^1 \cap L_{t,y}^4 \mathcal{H}_x^1( \mathbb{R} \times \mathbb{R}^2 \times \mathbb{R})}  \lesssim 1,$
with
\begin{align*}
\left\|u_k (t) - w_{\lambda_k } (t) \right\|_{L_t^\infty L_y^2 \mathcal{H}_x^1 \cap L_{t,y}^4 \mathcal{H}_x^1( \mathbb{R} \times \mathbb{R}^2 \times \mathbb{R})}    \to 0, \text{ as } k \to \infty,
\end{align*}
where
\begin{align*}
w_{\lambda_k}  (t,y,x) = e^{- i(t - t_k ) |\xi_k |^2} e^{iy\cdot \xi_k } \lambda_k^{-1} e^{it \left (  \Delta_x - x^2 \right) } v\left( \frac{t}{\lambda_k^2} + t_k , \frac{ y - y_k - 2 \xi_k ( t-t_k)}{ \lambda_k } , x \right).
\end{align*}
\end{theorem}
\begin{proof}[Proof of Theorem \ref{pr3.12}]
By translation invariance, we may take $y_k = 0$. By Galilean transformation and $ \left|\xi_k \right|$ is bounded, we may take $\xi_k  = 0$. Then
\begin{align*}
w_{\lambda_k}  (t,y,x) = \lambda_k^{-1} e^{it \left(  \Delta_x - x^2 \right) } v\left( \frac{t}{\lambda_k^2} + t_k, \frac{ y  }{ \lambda_k } , x \right)
\end{align*}
When $t_k= 0$, we will show $w_{\lambda_k}$ is an approximate solution to \eqref{eq1.1}.
After a simple computation, we see
\begin{align}\label{eq3.19v11}
e_{\lambda_k} & := \left(i\partial_t + \Delta_y + \Delta_x - x^2 \right) w_{\lambda_k}  - |w_{\lambda_k} |^2 w_{\lambda_k} \notag\\
 & = -  {\lambda_k^{- 3} } \sum\limits_{ n \in \mathbb{N} }  e^{-it(2n+1)} \sum\limits_{\substack{n_1,n_2,n_3\in \mathbb{N} , \\ n_1-n_2+ n_3 \ne n}} e^{- 2 it \left(n_1-n_2 + n_3 -n \right)} \left( \Pi_n \left(v_{n_1} \bar{v}_{n_2} v_{n_3} \right) \right)\left(\frac{t}{\lambda_k^2}, \frac{y}{\lambda_k} , x\right).
\end{align}
We will show this error term is small in the dual Strichartz space. Dividing the right hand side of \eqref{eq3.19v11} into three terms:
\begin{align*}
e_{\lambda_k} (t,y,x) = & -{\lambda_k^{- 3 } } \sum\limits_{n \in \mathbb{N} }  e^{-it(2n+1)} \sum\limits_{n_1,n_2,n_3 \in \mathbb{N}  } e^{- 2 i t  \left(n_1-n_2 + n_3 -n \right)} P_{\ge 2^{-10}}^y \left( \Pi_n \left(v_{n_1} \bar{v}_{n_2} v_{n_3} \right)\left(\frac{t}{\lambda_k^2}, \frac{y}{\lambda_k} , x\right)\right)\\
& + {\lambda_k^{- 3}}
 \sum\limits_{n \in \mathbb{N} }  e^{-it(2n+1)} \sum\limits_{\substack{n_1,n_2,n_3 \in \mathbb{N} , \\ n_1-n_2+ n_3 = n}} e^{- 2 i t \left(n_1-n_2 + n_3 -n \right)} P_{\ge 2^{-10}}^y \left( \Pi_n \left(v_{n_1} \bar{v}_{n_2} v_{n_3} \right)\left(\frac{t}{\lambda_k^2}, \frac{y}{\lambda_k} , x\right)\right)\\
& - {\lambda_k^{- 3}}  \sum\limits_{n \in \mathbb{N} }  e^{-it(2n+1)} \sum\limits_{\substack{n_1,n_2,n_3 \in \mathbb{N} , \\ n_1-n_2+ n_3 \ne n}} e^{- 2 i t \left(n_1-n_2 + n_3 -n \right)} P_{\le 2^{-10}}^y \left( \Pi_n\left(v_{n_1} \bar{v}_{n_2} v_{n_3} \right)\left(\frac{t}{\lambda_k^2}, \frac{y}{\lambda_k} , x\right)\right)\\
&=:   e_{\lambda_k}^1 + e_{\lambda_k}^2 + e_{\lambda_k}^3.
\end{align*}
We first consider $e_{\lambda_k}^1$ and shall use Bernstein's inequality, Leibnitz's rule, Plancherel's identity and H\"older's inequality to estimate as follows:
\begin{equation}\label{eq3.11v47}
\begin{aligned}
 \left\|e_{\lambda_k}^1(t,y,x)\right\|_{L_{t,y}^\frac43 \mathcal{H}_x^1}
\lesssim & \lambda_k^{-1 }
 \left\|   \sum\limits_{ n \in \mathbb{N} } e^{-i\lambda_k^2 t(2n+1)} \sum\limits_{n_1,n_2,n_3 \in \mathbb{N}  } e^{- 2 i\lambda_k^2 t  \left(n_1-n_2 + n_3 -n \right)}    \Pi_n\left( \nabla_y v_{n_1} \cdot \bar{v}_{n_2} \cdot v_{n_3} \right)\left(t, y, x\right)\right \|_{L_{t,y}^\frac43 \mathcal{H}_x^1}
 + \cdots \\ \sim &  \lambda_k^{-1 }
 \left \|e^{i\lambda_k^2 t \left(\Delta_x - x^2 \right)}\left( e^{i\lambda_k^2 t \left(\Delta_x - x^2 \right)} \nabla_y v \cdot \overline{e^{ i\lambda_k^2 t \left(\Delta_x - x^2 \right)} v} e^{i\lambda_k^2 t \left(\Delta_x -x^2 \right)} v \right)(t,y,x) \right\|_{L_{t,y}^\frac43 \mathcal{H}_x^1}\\
\lesssim &  \lambda_k^{-1 }
\left \|   \nabla_y v \right\|_{L_{t,y}^4 \mathcal{H}_x^1} \left\|v \right\|^2_{L_{t,y}^4 \mathcal{H}_x^1}
 \to 0, \text{ as }  k \to \infty,
\end{aligned}
\end{equation}
where $\cdots$ are the missing two terms with $\nabla_y$ acting on $\bar{v}_{n_2}$ and $ {v}_{n_3}$.

We now turn to the estimate of $e_{\lambda_k}^2$. Using Bernstein's inequality and Leibniz's rule as above, we have
\EQn{\label{eq3.12v42}
\left\|e_{\lambda_k}^2 \right\|_{L_{t,y}^\frac43 \mathcal{H}_x^1}
\lesssim & \lambda_k^{-1 }
 \bigg\|   \sum\limits_{ n \in \mathbb{N} }  e^{-i\lambda_k^2 t(2n+1)} \sum\limits_{\substack{n_1,n_2,n_3 \in \mathbb{N}  , \\ n_1-n_2+ n_3 = n}}   \Pi_n( \nabla_y v_{n_1} \cdot \bar{v}_{n_2} v_{n_3} )\left(t, y, x\right)  \bigg\|_{L_{t,y}^\frac43 \mathcal{H}_x^1}+ \cdots \\
\sim &\lambda_k^{-1 }
 \left\|   \bigg(\sum\limits_{ n \in \mathbb{N} }   \Big|  \langle n \rangle^{\half1} \sum\limits_{\substack{n_1,n_2,n_3 \in \mathbb{N}  , \\ n_1-n_2+ n_3 = n}}   \Pi_n( \nabla_y v_{n_1} \cdot \bar{v}_{n_2} v_{n_3} )\left(t, y, x\right) \Big|^2\bigg)^\frac12 \right\|_{L_{t,y}^\frac43 L_x^2}+ \cdots , \\}
where $\cdots$ are the missing two terms with $\nabla_y$ acting on $\bar{v}_{n_2}$ and $ {v}_{n_3}$.

We observe the following elementary inequality:  for $n=n_1-n_2+n_3$
\EQ{
	\jb{n}^{\frac12} \loe \jb{n}^{-1}\jb{n}^2 \loe \jb{n}^{-1}\jb{n_1}^2\jb{n_2}^2\jb{n_3}^2.
}
Using the fact $\left\{ \jb{n}^{-1} \right\}_{n \in \mathbb{N}   }  \in l_n^2$, the Minkowski inequality and boundedness of $\Pi_n$, we have
\begin{align*}
&\bigg\| \jb{n}^{\half1} \sum_{\substack{n_1,n_2,n_3 \in \mathbb{N}  , \\ n_1-n_2+ n_3 = n}} \Pi_n\brk{\nabla_y v_{n_1} \cdot \bar{v}_{n_2} v_{n_3}}(t,y,x)
\bigg\|_{L_{t,y}^\frac43 L_x^2 l_n^2} \\
	\lsm & \bigg\| \jb{n}^{-1} \sum_{\substack{n_1,n_2,n_3 \in \mathbb{N}  , \\ n_1-n_2+ n_3 = n}} \jb{n_1}^2\jb{n_2}^2\jb{n_3}^2
 \left|\Pi_n\brk{\nabla_y v_{n_1} \cdot \bar{v}_{n_2} v_{n_3}}(t,y,x) \right|
\bigg\|_{L_{t,y}^\frac43 L_x^2 l_n^2} \\
\lsm &  \bigg\| \sum_{n_1,n_2,n_3 \in \mathbb{N}  } \jb{n_1}^2\jb{n_2}^2\jb{n_3}^2 \normo{\brk{\nabla_y v_{n_1} \cdot \bar{v}_{n_2} v_{n_3}}(t,y,x)}_{L_x^2}
 \bigg\|_{L_{t,y}^\frac43}.
\end{align*}
By H\"older's inequality and the embedding $\HH^1(\mathbb{R}) \subseteq L^\infty(\mathbb{R})$, we find
\EQ{
	\norm{\brk{\nabla_y v_{n_1} \cdot \bar{v}_{n_2} v_{n_3}}(t,y,x )}_{L_x^2(\mathbb{R}) } \lsm \norm{\nabla_y v_{n_1}(t,y,x)}_{L^2_x(\mathbb{R}) } \norm{ v_{n_2}(t,y,x )}_{\HH^1_x(\mathbb{R}) } \norm{ v_{n_3}(t,y,x)}_{\HH^1_x(\mathbb{R}) }.
}
Similar arguments can be applied to the other two terms on the right hand side of \eqref{eq3.12v42}. All together lead to the estimate:
\EQn{\label{eq3.14v47}
\left\|e_{\lambda_k}^2 \right\|_{L_{t,y}^\frac43 \mathcal{H}_x^1}
\lsm & \lambda_k^{-1 }
 \norm{\jb{n_1}^2\jb{n_2}^2\jb{n_3}^2 \norm{\brk{\nabla_y v_{n_1} \cdot \bar{v}_{n_2} v_{n_3}}(t,y,x)}_{L_x^2}}_{L_{t,y}^\frac43 l_{n_1}^1l_{n_2}^1l_{n_3}^1} \\
& + \lambda_k^{-1}
 \norm{\jb{n_1}^2\jb{n_2}^2\jb{n_3}^2 \norm{\brk{v_{n_1} \cdot  \overline{ \nabla_y  {v}_{n_2}}  v_{n_3}}(t,y,x)}_{L_x^2}}_{L_{t,y}^\frac43 l_{n_1}^1l_{n_2}^1l_{n_3}^1} \\
& + \lambda_k^{-1 }
 \norm{\jb{n_1}^2\jb{n_2}^2\jb{n_3}^2 \norm{\brk{ v_{n_1} \cdot \bar{v}_{n_2} \nabla_y  v_{n_3}}(t,y,x)}_{L_x^2}}_{L_{t,y}^\frac43 l_{n_1}^1l_{n_2}^1l_{n_3}^1} \\
\lsm & \lambda_k^{-1 }
 \norm{\jb{n_1}^3\jb{n_2}^3\jb{n_3}^3 \norm{\nabla_y v_{n_1}(t,y,\cdot)}_{L^2} \norm{ v_{n_2}(t,y,\cdot)}_{\HH^1} \norm{ v_{n_3}(t,y,\cdot)}_{\HH^1}}_{L_{t,y}^\frac43 l_{n_1}^2l_{n_2}^2l_{n_3}^2} \\
\lsm & \lambda_k^{-1 }
 \norm{\nabla_y v}_{L_{t,y}^4\HH_x^6} \norm{ v}_{L_{t,y}^4\HH_x^7}^2
\lsm   \lambda_k^{-1 }  C\brk{\norm{\psi
%v_0
}_{H_{y}^1\HH_x^6}} C\brk{\norm{ \psi
%v_0
}_{L_{y}^2\HH_x^7}} \to 0, \text{ as }  k \to \infty.
}

Now, we only need to deal with $e_{\lambda_k}^3$. We will use the normal form transform to exploit additional decay of $\lambda_k $, since it possesses time non-resonance property. Integrating by parts and direct computation imply
\begin{align*}
& \int_0^t e^{i(t-\tau) \left(\Delta_y + \Delta_x - x^2 \right)} e^3_{\lambda_k} (\tau) \,\mathrm{d}\tau\\
= & - {\lambda_k^{- 3 } } \sum\limits_{\substack{n_1,n_2,n_3,n \in \mathbb{N}  \\ n_1-n_2 + n_3 \ne n}} \int_0^t e^{it \left(\Delta_y - 2n-1 \right)}
e^{ - i \tau \tilde{\Delta}_y} P_{\le 2^{-10}}^y \left( \Pi_n\left(v_{n_1 } \bar{v}_{n_2 } v_{n_3 }\right)\left(\frac{\tau}{\lambda_k^2}, \frac{y}{\lambda_k } ,x \right) \right)  \,\mathrm{d}\tau\\
= & - {\lambda_k^{- 3 } } \sum\limits_{\substack{n_1,n_2,n_3, n\in \mathbb{N}  ,\\ n_1 -n_2 + n_3 \ne n}} e^{it \left(\Delta_y - 2n-1 \right)}
{e^{ -it \tilde\Delta_y}}{ \left( - i \tilde\Delta_y \right)^{- 1}  }
 P_{\le 2^{-10}}^y\Pi_n \left(v_{n_1 } \bar{v}_{n_2 } v_{n_3 } \right) \left(\frac{t}{\lambda_k^2}, \frac{y}{\lambda_k},x \right)\\
& + {\lambda_k^{- 3 } } \sum\limits_{\substack{n_1,n_2,n_3, n \in \mathbb{N}  ,\\ n_1 -n_2 + n_3 \ne n}} e^{it \left(\Delta_y - 2n-1 \right)}
{ \left( - i \tilde\Delta_y \right)^{- 1}}
 P_{\le 2^{-10}}^y\Pi_n\left( v_{n_1 } \bar{v}_{n_2 } v_{n_3 } \right) \left(0, \frac{y}{\lambda_k} ,x \right)\\
& +{\lambda_k^{- 3 } }  \sum\limits_{\substack{n_1,n_2,n_3, n \in \mathbb{N}  ,\\ n_1 -n_2 + n_3 \ne n}} \int_0^t
 e^{it \left(\Delta_y - 2n-1 \right)} {e^{- i \tau \tilde\Delta_y }}{\left( -i \tilde\Delta_y  \right)^{- 1} }
 \partial_\tau\left( P_{\le 2^{-10}}^y\Pi_n \left(v_{n_1 } \bar{v}_{n_2 } v_{n_3 } \right) \left( \frac{\tau}{\lambda_k^2} , \frac{y}{\lambda_k} ,x \right) \right) \,\mathrm{d}\tau\\
= &{\lambda_k^{- 3 } }  \sum\limits_{\substack{n_1,n_2,n_3, n \in \mathbb{N}  ,\\ n_1 -n_2 + n_3 \ne n}}
 e^{it \left(\Delta_y + \Delta_x - x^2 \right)}  { \left( -i   \tilde\Delta_y \right)^{- 1}}
P_{\le 2^{-10}}^y\Pi_n \left( v_{n_1 } \bar{v}_{n_2 } v_{n_3 } \right) \left(0,\frac{y}{\lambda_k} , x\right)  \\
& -{\lambda_k^{- 3 } } \sum\limits_{\substack{n_1, n_2,n_3, n \in \mathbb{N}  ,\\ n_1 -n_2 + n_3 \ne n}}
 {e^{-2it \left(n_1-n_2 + n_3 \right) -it}} { \left( -i \tilde\Delta_y  \right)^{- 1} } P_{\le 2^{-10}}^y\Pi_n \left( v_{n_1 } \bar{v}_{n_2 } v_{n_3 } \right) \left( \frac{t}{\lambda_k^2} , \frac{y}{\lambda_k},x \right)\\
& +{\lambda_k^{- 3 } } \sum\limits_{\substack{n_1,n_2,n_3, n \in \mathbb{N}  ,\\ n_1 -n_2 + n_3 \ne n}} e^{it \left(\Delta_y - 2n-1 \right)} \int_0^t
{e^{-i \tau \tilde\Delta_y}}{ \left( -i \tilde\Delta_y   \right)^{- 1}}
\partial_\tau P_{\le 2^{-10}}^y\Pi_n \left( v_{n_1 } \bar{v}_{n_2 } v_{n_3 } \right) \left( \frac{\tau}{\lambda_k^2} , \frac{y}{\lambda_k} ,x \right) \,\mathrm{d}\tau,
\end{align*}
where the operator $ \tilde\Delta_y   $ is defined to be
\begin{align*}
\tilde\Delta_y   : =  2 (n_1 -n_2 +n_3 -n) +  \Delta_y.
\end{align*}
This is a perturbation of the Laplacian operator and we suppress the parameters $n_1, n_2,n_3,n$. The inverse operator $\left( - \tilde\Delta_y \right)^{- 1} $ is defined by the Fourier transform
\begin{align*}
\mathcal{F}_y \left( \left( - i \tilde\Delta_y \right)^{- 1} f \right) (\xi, x) =  \frac{i  ( \mathcal{F}_y f) ( \xi, x)  }{  2(n_1 -n_2 + n_3 - n) -  |\xi|^2} .
\end{align*}
This operator is invertible when $n_1 -n_2 +n_3 -n \ne 0$ and $|\xi|\leq 2^{-10}$. We will use this expression in the remaining of the proof.

Denote
\begin{align*}
e_{\la_k}^{3,1}:= & \bigg\| {\lambda_k^{- 3 } }
\sum\limits_{\substack{n_1,n_2,n_3, n \in \mathbb{N}  ,\\ n_1 -n_2 + n_3 \ne n}} e^{it  \left(\Delta_y + \Delta_x - x^2 \right)} P_{\le 2^{-10}}^y \left(
{ \left(  -i \tilde\Delta_y   \right)^{- 1}}
\left( \Pi_n \left( v_{n_1 } \bar{v}_{n_2 } v_{n_3 } \right)\left(0, \frac{y}{\lambda_k} ,x \right) \right) \right)  \bigg\|_{L_{t,y}^4 \mathcal{H}_x^1}, \\
e_{\la_k}^{3,2}:= & \bigg\|
{\lambda_k^{- 3} }
 \sum\limits_{\substack{n_1,n_2,n_3, n \in \mathbb{N}  ,\\ n_1 -n_2 + n_3 \ne n}}
  {e^{-2it \left(n_1-n_2 + n_3 \right) -it}}{ \left(  -i \tilde\Delta_y  \right)^{- 1}  }
   \left( P_{\le 2^{-10}}^y \Pi_n \left( v_{n_1 } \bar{v}_{n_2} v_{n_3 } \right)\left( \frac{t}{\lambda_k^2} , \frac{y}{\lambda_k} ,x \right) \right) \bigg\|_{L_{t,y}^4 \mathcal{H}_x^1},
\end{align*}
and
\begin{align*}
e_{\la_k}^{3,3} & :=\bigg\| {\lambda_k^{- 3 } }
\sum\limits_{\substack{n_1,n_2,n_3, n \in \mathbb{N}  ,\\ n_1 -n_2 + n_3 \ne n}} e^{it \left(\Delta_y - 2n-1 \right)}
  \int_0^t {e^{ -i \tau \tilde\Delta_y }}{ \left( -i \tilde\Delta_y \right)^{- 1}}  \left( \partial_\tau P_{\le 2^{-10}}^y \Pi_n \left( v_{n_1 } \bar{v}_{n_2 } v_{n_3 } \right)\left( \frac{\tau}{\lambda_k^2} , \frac{y}{\lambda_k} ,x \right) \right)  \,\mathrm{d}\tau \bigg\|_{L_{t,y}^4 \mathcal{H}_x^1}.
\end{align*}
Then, we have
\begin{align} \label{eq3.18v48}
 \left\| \int_0^t e^{i(t-\tau) \left(\Delta_y + \Delta_x - x^2 \right)} e_{\lambda_k}^3(\tau,y,x) \mathrm{d}\tau \right\|_{L_{t,y}^4 \mathcal{H}_x^1}
\sim &\ e_{\la_k}^{3,1} + e_{\la_k}^{3,2} + e_{\la_k}^{3,3}.
\end{align}
First, we consider the term $e_{\la_k}^{3,1}$. By the boundedness of the operator ${ P_{\le 2^{-10}}^y }{ \left( -i \tilde\Delta_y   \right)^{- 1}} $ when $n_1 - n_2 + n_3 \ne n$ and Minkowski's inequality, we may estimate as follows:
\EQn{ \label{eq3.19v48}
e_{\la_k}^{3,1} \lsm & \bigg\|\normbb{ \jb{n}^{\half 1} \sum\limits_{\substack{n_1,n_2,n_3 \in \mathbb{N}  \\ n_1-n_2 + n_3 \ne n}}
{\lambda_k^{- 3 } } P_{\le 2^{-10}}^y { \left( -i \tilde\Delta_y  \right)^{- 1}}   \Pi_n \left( \left(v_{n_1} \bar{v}_{n_2} v_{n_3} \right) \left(0,\frac{y}{\lambda_k} , x \right) \right) }_{l_n^2} \bigg\|_{L_{y,x}^2 }\\
\lsm &  {\lambda_k^{- 3 } } \bigg\| \jb{n}^{\half 1} \sum\limits_{\substack{n_1,n_2,n_3 \in \mathbb{N} \\ n_1-n_2 + n_3 \ne n}} \norm{ \Pi_n \left( (v_{n_1} \bar{v}_{n_2} v_{n_3})\left(0,\frac{y}{\lambda_k} , x \right) \right) }_{L_y^2}  \bigg\|_{L_x^2 l_n^2}\\
\lsm & {\lambda_k^{- 2 } } \bigg\| \jb{n}^{\half 1} \sum\limits_{\substack{n_1,n_2,n_3 \in \mathbb{N}  }} \norm{ \Pi_n\left(  \left(v_{n_1} \bar{v}_{n_2} v_{n_3} \right)(0,y, x) \right) }_{L_y^2}  \bigg\|_{L_x^2 l_n^2} \\
\lsm & {\lambda_k^{- 2 } }   \sum\limits_{\substack{n_1,n_2,n_3 \in \mathbb{N} }} \norm{\jb{n}^{\half 1} \Pi_n \left(  \left(v_{n_1} \bar{v}_{n_2} v_{n_3} \right)(0,y, x) \right) }_{L_{y,x}^2 l_n^2}
\lsm {\lambda_k^{- 2}}   \sum\limits_{\substack{n_1,n_2,n_3 \in \mathbb{N}  }} \norm{ \left( v_{n_1}  \bar{v}_{n_2} v_{n_3} \right) (0,y, x) }_{L_y^2 \HH_x^1}  \\
\lsm & {\lambda_k^{- 2 } }   \sum\limits_{\substack{n_1,n_2,n_3 \in \mathbb{N}  }} \norm{v_{n_1}(0,y, x) }_{L_y^6 \HH_x^1} \norm{v_{n_2}(0,y, x) }_{L_y^6 \HH_x^1} \norm{v_{n_3}(0,y, x) }_{L_y^6 \HH_x^1}
\\
 \lsm &   {\lambda_k^{- 2 } } C \left( \norm{v(0,y, x) }_{H_y^1 \HH_x^3} \right)^3 \to 0, \text{ as }  k \to \infty. }
Next, we consider the term $e_{\la_k}^{3,2}$. As in the estimate of $e_{\la_k}^{3,1}$. By the boundedness of the operator ${ P_{\le 2^{-10}}^y }{ \left( -i \tilde\Delta_y   \right)^{- 1}} $ when $n_1 - n_2 + n_3 \ne n$, Minkowski's inequality, the fractional Leibniz rule, Sobolev's inequality and H\"older's inequality, we have
\begin{align} \label{eq3.25v48}
e_{\la_k}^{3,2} \lsm & \la_k^{-3} \bigg\|\jb{n}^{\half 1}\sum_{\substack{n_1,n_2,n_3 \in \mathbb{N}  \\ n_1-n_2 + n_3 \ne n}}
{e^{-2it \left(n_1 -n_2+ n_3 \right) -it}}{ \left( -i \tilde\Delta_y  \right)^{- 1}}
P_{\le 2^{-10}}^y \Pi_n\Big( \left( v_{n_1}\bar{v}_{n_2} v_{n_3} \right) \Big(\frac{t}{\lambda_k^2},\frac{y}{\lambda_k} , x\Big)\Big) \bigg\|_{L_{t,y}^4 L_x^2l_n^2}   \\
\lsm & \la_k^{-3} \bigg\|\jb{n}^{\half 1}\sum_{\substack{n_1,n_2,n_3 \in \mathbb{N}  \\ n_1-n_2 + n_3 \ne n}} \norm{ \Pi_n\Big( \left( v_{n_1}\bar{v}_{n_2} v_{n_3} \right) \Big(\frac{t}{\lambda_k^2},\frac{y}{\lambda_k} , x\Big)\Big)}_{L_{y}^4} \bigg\|_{L_t^4L_x^2l_n^2} \notag \\
\lsm &     \la_k^{-\frac32} \bigg\|\sum_{\substack{n_1,n_2,n_3 \in \mathbb{N}  }} \norm{ v_{n_1}\bar{v}_{n_2} v_{n_3} }_{H_y^{\frac12}\HH_x^1} \bigg\|_{L_t^4}
\lsm   \la_k^{-\frac32} \bigg\| \norm{ v(t,y,x)}_{W_y^{\frac34,4}\HH_x^5}^3 \bigg\|_{L_t^4} \notag \\
	\lsm  & \la_k^{- \frac32} \norm{ v(t,y,x)}_{L_t^{12}W_y^{\frac34,4}\HH_x^5}^3
 \lsm   \la_k^{-\frac32} C\brk{\norm{ v(0,y,x)}_{H_y^{ \frac{13}{12} }\HH_x^5}}^3\to 0, \text{ as }  k \to \infty. \notag
\end{align}
Finally, we are left to consider the term $e_{\la_k}^{3,3}$. Applying the Strichartz estimate, we obtain
\EQn{ \label{eq3.26v42}
& \bigg\| {\lambda_k^{- 3 } }
\sum\limits_{\substack{n_1,n_2,n_3, n \in \mathbb{N}  ,\\ n_1 -n_2 + n_3 \ne n}} e^{it \left(\Delta_y - 2n-1 \right)} \int_0^t {e^{-i\tau  \tilde\Delta_y}}{ \left( -i \tilde\Delta_y  \right)^{-1}}  \left( \partial_\tau P_{\le 2^{-10}}^y \Pi_n \left( v_{n_1 } \bar{v}_{n_2 } v_{n_3 }\right)\left( \frac{\tau}{\lambda_k^2} , \frac{y}{\lambda_k} ,x \right) \right)  \,\mathrm{d}\tau  \bigg\|_{L_{t,y}^4 \mathcal{H}_x^1}\\
\lesssim & \bigg\|\left(i\partial_t + \Delta_y + \Delta_x - x^2 \right) \sum\limits_{\substack{n_1,n_2,n_3,n \in \mathbb{N} \\ n_1-n_2 + n_3 \ne n }} e^{it(\Delta_y -2n -1)} \\
& \ \int_0^t P_{\le 2^{-10}}^y \bigg({e^{-i \tau \tilde\Delta_y}}{ \left( -i \tilde\Delta_y   \right)^{- 1}}
 \partial_\tau \Pi_n\left(  {\lambda_k^{- 3 } }
v_{n_1}\left(\frac{\tau}{\lambda_k^2}, \frac{y}{\lambda_k} , x\right) \overline{v_{n_2}\left(\frac{\tau}{\lambda_k^2}, \frac{y}{\lambda_k} , x \right) }  v_{n_3} \left(\frac{\tau}{\lambda_k^2}, \frac{y}{\lambda_k} ,x \right) \right) \bigg) \mathrm{d}\tau \bigg\|_{L_t^1 L_y^2 \mathcal{H}_x^1}.
}
We observe, after some computation, that
\begin{align*}
& \left(i\partial_t + \Delta_y + \Delta_x - x^2 \right) \sum\limits_{\substack{n_1,n_2,n_3,n \in \mathbb{N}  ,\\ n_1-n_2 + n_3  \ne n}}
 e^{it \left(\Delta_y - 2n-1 \right)} \int_0^t
 {e^{-i \tau  \tilde\Delta_y }}{ \left( -i \tilde\Delta_y \right)^{- 1} P_{\le 2^{-10}}^y } \\
&\qquad  \cdot \partial_\tau \Pi_n\left({\lambda_k^{- 3 } } v_{n_1}\left(\frac{\tau}{\lambda_k^2}, \frac{y}{\lambda_k} , x\right) \overline{v_{n_2}\left(\frac{\tau}{\lambda_k^2},\frac{y}{\lambda_k} ,x\right)  } v_{n_3}\left(\frac{\tau}{\lambda_k^2}, \frac{y}{\lambda_k} , x\right)\right) \mathrm{d}\tau\\
= & \sum\limits_{\substack{n_1,n_2,n_3,n \in \mathbb{N}  \\ n_1-n_2 + n_3 \ne n}}
{e^{-2it \left (n_1-n_2 + n_3 \right) - it}}{ \left( -i \tilde\Delta_y \right)^{- 1}P_{\le 2^{-10}}^y} \partial_t \Pi_n\left({\lambda_k^{- 3 } }
 v_{n_1}\left(\frac{t}{\lambda_k^2}, \frac{y}{\lambda_k}, x\right) \overline{v_{n_2}\left(\frac{t}{\lambda_k^2}, \frac{y}{\lambda_k} , x \right) } v_{n_3}\left(\frac{t}{\lambda_k^2}, \frac{y}{\lambda_k}, x\right) \right).
\end{align*}
Therefore, by the above observation, Plancherel's theorem and Leibniz's rule, we have
\EQn{ \label{eq3.27v45}
\eqref{eq3.26v42} \lesssim
& \bigg\| \bigg(\sum\limits_{ n \in \mathbb{N} }  \bigg(\sum\limits_{\substack{n_1,n_2,n_3 \in \mathbb{N}  ,\\ n_1-n_2 + n_3 \ne n}} \bigg|
{e^{-2it \left(n_1-n_2 + n_3-n \right) -it}}{ \left( -i \tilde\Delta_y \right)^{- 1}} P_{\le 2^{-10}}^y \partial_t \Pi_n\bigg({\lambda_k^{- 3 } }
 v_{n_1}\left(\frac{t}{\lambda_k^2}, \frac{y}{\lambda_k} , x \right)  \\
& \qquad \quad \cdot   \overline{v_{n_2}\left(\frac{t}{\lambda_k^2},\frac{y}{\lambda_k} , x \right) }
v_{n_3}\left(\frac{t}{\lambda_k^2}, \frac{y}{\lambda_k} ,x \right) \bigg)\bigg| \bigg)^2
 \langle n\rangle\bigg)^\frac12 \bigg\|_{L_t^1 L_y^2 L_x^2}\\
\lesssim & \bigg\|\bigg(\sum\limits_{n \in \mathbb{N} }  \bigg(\sum\limits_{\substack{n_1,n_2,n_3 \in \mathbb{N}  ,\\ n_1-n_2 + n_3 \ne n}} \Big\|\partial_t \Pi_n\Big({\lambda_k^{- 3 } } v_{n_1}\Big(\frac{t}{\lambda_k^2}, \frac{y}{\lambda_k} , x\Big) \overline{v_{n_2}\Big(\frac{t}{\lambda_k^2}, \frac{y}{\lambda_k} , x\Big)}  v_{n_3}\Big(\frac{t}{\lambda_k^2}, \frac{y}{\lambda_k} , x\Big) \Big) \bigg\|_{L_{y,x}^2}\bigg)^2\langle n\rangle\bigg)^\frac12 \bigg\|_{L_t^1}\\
\lesssim & {\lambda_k^{- 2 } } \sum\limits_{n_1,n_2,n_3 \in \mathbb{N} } \norm{ \partial_t v_{n_1} \cdot \bar{v}_{n_2} v_{n_3} }_{L_t^1 L_y^2 \mathcal{H}_x^1}
+ {\lambda_k^{- 2 } } \sum\limits_{n_1,n_2,n_3 \in \mathbb{N}  } \norm{   v_{n_1} \cdot  \overline{ \partial_t  {v}_{n_2} }  v_{n_3} }_{L_t^1 L_y^2 \mathcal{H}_x^1}\\
& \ +{\lambda_k^{- 2 } } \sum\limits_{n_1,n_2,n_3 \in \mathbb{N}  } \norm{  v_{n_1} \cdot \bar{v}_{n_2} \partial_t  v_{n_3} }_{L_t^1 L_y^2 \mathcal{H}_x^1}.
}
We shall only show how to estimate the first term on the right hand side of \eqref{eq3.27v45} as the other two terms can be estimated similarly. By H\"older's inequality, and the fact that $v$ satisfies \eqref{eq1.3v28}, we have
\begin{align*}
\sum\limits_{\substack{n_1,n_2,n_3 \in \mathbb{N}   }} \norm{ \partial_t v_{n_1} \bar{v}_{n_2} v_{n_3} }_{L_t^1 L_y^2 \mathcal{H}_x^1}
\lesssim & \|v\|_{L_t^3 L_y^6 \mathcal{H}_x^3}^2 \left\|\partial_t v \right\|_{L_t^3 L_y^6 \mathcal{H}_x^3}\\
\lesssim & \|v\|_{L_t^3 L_y^6 \mathcal{H}_x^3}^2 \|\Delta_y v \|_{L_t^3 L_y^6 \mathcal{H}_x^3} +
\|v\|_{L_t^3 L_y^6 \mathcal{H}_x^3}^2 \bigg\|\sum\limits_{\substack{n_4,n_5,n_6,n_1 \in \mathbb{N}  ,\\ n_4-n_5 + n_6 = n_1}} \Pi_{n_1}\left(v_{n_4} \bar{v}_{n_5} v_{n_6} \right) \bigg\|_{L_t^3 L_y^6 \mathcal{H}_x^3}.\\
\end{align*}
Applying H\"older's inequality and the Sobolev embedding, we have
\begin{align*}
 \  \normbb{\sum\limits_{\substack{n_4,n_5,n_6,n_1 \in \mathbb{N}  ,\\ n_4-n_5 + n_6 = n_1}} &\Pi_{n_1} \left(v_{n_4} \bar{v}_{n_5} v_{n_6}  \right) }_{L_t^3 L_y^6 \mathcal{H}_x^3}
 \lsm  \normbb{\jb{n_1}^{\frac32}  \sum\limits_{\substack{n_4,n_5,n_6 \in \mathbb{N}  ,\\ n_4-n_5 + n_6 = n_1}} \Pi_{n_1}(v_{n_4} \bar{v}_{n_5} v_{n_6} )}_{L_t^3 L_y^6 L_x^2 l_{n_1}^2} \\
\lsm &  \normbb{\sum\limits_{\substack{n_4,n_5,n_6 \in \mathbb{N}  ,\\ n_4-n_5 + n_6 = n_1}} \jb{n_1}^{-1} \jb{n_4}^{3} \jb{n_5}^{3} \jb{n_6}^{3}  \abs{\Pi_{n_1}(v_{n_4} \bar{v}_{n_5} v_{n_6} ) }}_{L_t^3 L_y^6 L_x^2 l_{n_1}^2} \\
\lsm & \norm{\sum\limits_{\substack{n_4,n_5,n_6 \in \mathbb{N}  }} \jb{n_4}^{3} \jb{n_5}^{3} \jb{n_6}^{3}  \norm{v_{n_4} \bar{v}_{n_5} v_{n_6} }_{L_x^2} }_{L_t^3 L_y^6} \\
\lsm & \norm{\sum\limits_{\substack{n_4,n_5,n_6 \in \mathbb{N}  }} \jb{n_4}^{3} \jb{n_5}^{3} \jb{n_6}^{3}  \norm{v_{n_4}}_{\HH_x^1} \norm{v_{n_5}}_{\HH_x^1} \norm{v_{n_6}}_{\HH_x^1}}_{L_t^3 L_y^6}\lsm   \norm{v}_{L_t^9 L_y^{18} \HH_x^9}^3 \le
  C\brk{\norm{ \psi }_{ H_y^{\frac23} \HH_x^9}}^3,
\end{align*}
where in the last inequality we use the fact that by Corollary \ref{co4.7v65}, we have
\begin{align*}
\|v\|_{L_t^3 L_y^6 \mathcal{H}_x^3} \le C\brk{\left\|\psi \right\|_{ L_y^2 \mathcal{H}_x^3}} ,
\end{align*}
and
\begin{align*}
\left\|\Delta_y v \right\|_{L_t^3 L_y^6 \mathcal{H}_x^3} \le C\brk{ \left\|\Delta_y \psi \right\|_{L_y^2 \mathcal{H}_x^3}}.
\end{align*}
Combining all these estimates together, we finally obtain
\EQn{ \label{eq3.32v48}
e_{\la_k}^{3,3} \lsm {\la_k^{- 2 } } \brk{C\brk{ \left\| \psi \right\|_{ L_y^2 \mathcal{H}_x^3}}^2C\brk{ \left\|\Delta_y \psi \right\|_{L_y^2 \mathcal{H}_x^3}} + C\brk{ \left\| \psi  \right\|_{ L_y^2 \mathcal{H}_x^3}}^2C\brk{\norm{\psi  }_{ H_y^{\frac23} \HH_x^9}}^3}\to 0, \text{ as $ k \to\infty$. }}
To apply Theorem \ref{th2.6}, we see
\begin{align*}
\lim_{k\to\infty}\| w_{\lambda_k} (0,y,x)  - u_{\lambda_k} (0,y,x) \|_{L_y^2 \mathcal{H}_x^1} = \|\phi - \psi \|_{L_y^2 \mathcal{H}_x^1} \le \epsilon_4,
\end{align*}
\begin{align*}
\|w_{\lambda_k} \|_{L_t^\infty L_y^2 \mathcal{H}_x^1(\mathbb{R}\times \mathbb{R}^2 \times \mathbb{R})}
& = \left\| \frac1{\lambda_k}  v\left(\frac{t}{\lambda_k^2}, \frac{y}{\lambda_k} , x \right) \right\|_{L_t^\infty L_y^2 \mathcal{H}_x^1}
= \|v\|_{L_t^\infty L_y^2 \mathcal{H}_x^1},
\end{align*}
and
\begin{align*}
\|w_{\lambda_k} \|_{L_{t,y}^4 \mathcal{H}_x^1 (\mathbb{R}\times \mathbb{R}^2 \times \mathbb{R})}
= & \left\|\frac1{\lambda_k}  v\left(\frac{t}{\lambda_k^2},\frac{y}{\lambda_k} , x\right) \right\|_{L_{t,y}^4 \mathcal{H}_x^1}
= \|v\|_{L_{t,y}^4 \mathcal{H}_x^1}.
\end{align*}
These together with the estimates \eqref{eq3.11v47}, \eqref{eq3.14v47}, \eqref{eq3.18v48}, \eqref{eq3.19v48}, \eqref{eq3.25v48}, \eqref{eq3.32v48} and Theorem \ref{as3.14v12}  yields Theorem \ref{pr3.12} when $t_k = 0$.

If $t_k \to \pm \infty$, as $k \to \infty$, $v$ is the solution of (DCR) with
\begin{align*}
\lim\limits_{t\to \pm \infty } \| v(t,y,x) - e^{it \Delta_y} \psi \|_{L_y^2 \mathcal{H}_x^1}  =  0.
\end{align*}
By the argument in the case $t_k =0$, we can also obtain Theorem \ref{pr3.12} in this case.
\end{proof}

\subsection{Existence of an almost-periodic soltuion}
Define
\begin{align*}
\Lambda(L) = \sup \|u\|_{L_{t,y}^4 \mathcal{H}_x^{1 - \epsilon_0}  ( \mathbb{R} \times \mathbb{R}^2 \times \mathbb{R})},
\end{align*}
where the supremum is taken over all global solutions $u\in C_t^0( \mathbb{R},  \Sigma(\mathbb{R}^3)) $ of \eqref{eq1.1} with
\begin{align*}
\mathcal{E}(u(t)) + \mathcal{M}(u(t)) \le L.
\end{align*}
The proof of Theorem \ref{th2.3} implies $\Lambda(L) < \infty$ for sufficiently small $L$. Let
\begin{align}\label{eq:def L-max}
L_{\max} = \sup \left\{ L \ge 0 : \Lambda(L) < \infty \right\}.
\end{align}
If $L_{\max}<\infty$, then following the arguments in \cite{CGYZ,CGZ}, one can show the existence of an almost periodic solution with the help of Theorems \ref{pro3.9v23} and \ref{pr3.12}. The proof is rather standard, we refer to \cite{CGYZ,CGZ,Killip-Visan1,KM1,KM,TVZ} and omit the proof here.

\begin{theorem}[Existence of an almost-periodic solution]\label{th4.3}
Assume that $L_{max} < \infty$. Then there exists a solution $u_c \in C_t^0( \mathbb{R} ,  \Sigma ( \mathbb{R}^3) )$ of the defocusing cubic NLS with partial harmonic potential \eqref{eq1.1} satisfying
\begin{equation}\label{eq3.31v40}
\mathcal{E}(u_c) + \mathcal{M}(u_c)  = L_{max}\quad\text{ and }\quad
\|u_c\|_{L_{t,y}^4 \mathcal{H}_x^{ 1 - \epsilon_0}     (\mathbb{R} \times \mathbb{R}^2 \times \mathbb{R})} = \infty.
\end{equation}
Furthermore, $u_c$ is almost periodic in the sense that for any $\eta > 0$, there is a Lipschitz function $t\mapsto y(t)$ and a sufficiently large positive number
 $C(\eta) $ such that
\begin{align}\label{eq3.259}
\int_{|y+ y(t)|\ge C(\eta)} \left\|u_c(t,y,x) \right\|_{\mathcal{H}_x^1}^2 \,\mathrm{d}y < \eta, \ \forall\, t \in \mathbb{R}.
\end{align}
\end{theorem}

\section{Rigidity theorem}\label{se4v16}
In this section, we will exclude the almost-periodic solution in Theorem \ref{th4.3} by the interaction Morawetz estimate with an appropriately chosen weight function. Once the almost-periodic solution is excluded, we can finish the proof of Theorem \ref{th1.3}.

\begin{proposition}[Non-existence of the almost-periodic solution]\label{pr7.3}
The almost-periodic solution $u_c$ as in Theorem \ref{th4.3} does not exist.
\end{proposition}
\begin{proof}
For each $r_0>0$, we define the interaction Morawetz action
\begin{align*}
M_{r_0}(t) = \int_{\mathbb{R}^2 \times \mathbb{R}} \int_{\mathbb{R}^2\times \mathbb{R}} \Im\left(\overline{ {u}_c(t,y,x) }  \nabla_y u_c (t,y,x) \right) \cdot \nabla_y \psi_{r_0} \left(|y-\tilde{y}| \right) \left|u_c \left(t,\tilde{y},\tilde{x} \right)\right|^2 \,\mathrm{d}y\mathrm{d}x \mathrm{d}\tilde{y} \mathrm{d}\tilde{x},
\end{align*}
where $\Im=\text{Im}$ denotes the imaginary part of a complex number and $\psi_{r_0}\colon \mathbb{R} \to \mathbb{R}$ is a radial function defined as in \cite{CGT,PV} with
\begin{align*}
\Delta \psi_{r_0}(r) = \int_r^\infty s \log\left(\frac{s}r\right) w_{r_0}(s) \,\mathrm{d}s,
\end{align*}
where
\begin{align*}
w_{r_0}(s) =
\begin{cases}
\frac1{s^3}, \ \text{if } s\ge r_0,\\
0, \  \text{ if } s < r_0.
\end{cases}
\end{align*}
It is straightforward to verify that $\psi_{r_0}$ is convex and $ \left| \nabla \psi_{r_0} \right| $ is uniformly bounded (independent of $r_0$), with
\begin{align*}
-\Delta^2 \psi_{r_0} (r) = \frac{2\pi}{r_0} \delta_0(r) - w_{r_0}(r).
\end{align*}
Using the above properties of the weight function $\psi_{r_0}$, one can show (see \cite[Section 3.3]{CGT}) for all $T_0 > 0$,
\begin{align}\label{eq5.032}
\int_{-T_0}^{T_0} \int_{\mathbb{R}^2} \left| \left|\nabla_y \right|^\frac12 \left( \left\|u_c(t,y,x) \right\|_{L_x^2(\mathbb{R}) }^2\right) \right|^2 \,\mathrm{d}y\mathrm{d}t
 \lesssim \left\|u_c \right\|_{L_t^\infty L_{y,x}^2}^3 \left\|\nabla_y u_c \right\|_{L_t^\infty L_{y,x}^2} \lesssim 1.
\end{align}
By \eqref{eq3.259} and the conservation of mass, we have
\begin{align}\label{eq5.232}
\frac{\left\|u_c \right\|_{L_{y,x}^2}^2  }2 \le \int_{|y + y(t)|\le C\left(\frac{m_c }{100}\right)}  \left\|  u_c(t,y,x) \right\|_{L_x^2}^2 \,\mathrm{d}y.
\end{align}
where $m_c  := \mathcal{M}(u_c ) > 0$ by \eqref{eq3.31v40}.

Therefore, for each $  T_0 > 0$, by \eqref{eq5.232}, Sobolev's inequality, and \eqref{eq5.032}, we deduce
\begin{align*}
 \frac{ m_c^2 T_0} 2 &  \le \int_{-T_0}^{T_0} \left(\int_{|y + y(t)|\le C\left(\frac{m_c }{100}\right)} \left\|u_c(t,y,x) \right\|_{L_x^2}^2 \,\mathrm{d}y \right)^2 \,\mathrm{d}t\\
& \lesssim C\left(\frac{m_c }{100}\right)  \int_{-T_0}^{T_0}  \left( \left(\int_{\mathbb{R}^2} \left(\left\|u_c(t,y,x)\right\|_{L_x^2}^2\right)^4 \,\mathrm{d}y \right)^\frac14 \right)^2 \,\mathrm{d}t\\
& \lesssim   C\left(\frac{m_c }{100}\right) \int_{-T_0}^{T_0} \int_{\mathbb{R}^2} \left| \left|\nabla_y \right|^\frac12 \left(\left\|u_c(t,y,x)\right\|_{L_x^2}^2\right) \right|^2 \,\mathrm{d}y \mathrm{d}t
  \lesssim C\left(\frac{m_c  }{100}\right).
\end{align*}
 Letting $T_0 \to \infty$, we obtain a contradiction, and this concludes the proof.
\end{proof}

Finally, we can now reach to Theorem \ref{th1.3}.
\begin{proof}[Proof of Theorem \ref{th1.3} ]
By Theorems \ref{as3.14v12} and \ref{th2.3}, to prove the scattering of solutions to \eqref{eq1.1}, it suffices to show the finiteness of the $L_{t,y}^4 \mathcal{H}_x^{1- \epsilon_0}-$norm of the solution of \eqref{eq1.1}.

To this end, let $L_{\max}$ be given as in \eqref{eq:def L-max}. Then, equivalently, we need show that $L_{\max} = \infty$. Suppose for a contradiction that $L_{\max} < \infty$. Then Theorem \ref{th4.3} would yield an almost-periodic solution of \eqref{eq1.1}, which is impossible in view of Proposition \ref{pr7.3}. This completes the proof.
\end{proof}

\section{Scattering of the \eqref{eq1.3v28}} \label{se6v55}
In this section, we will prove Theorem \ref{as3.14v12}, that is the global well-posedness and scattering of  the \eqref{eq1.3v28} system:
\begin{align*}
	\begin{cases}
		i\partial_t v  + \Delta_{\mathbb{R}^2} v  = F(v), \\
		v(0,y,x) =\phi(y,x),
	\end{cases}
\end{align*}
where
\begin{align*}
	F(v) : = \sum\limits_{ \substack{ n_1,n_2,n_3,n\in \mathbb{N},\\  n_1 - n_2 + n_3 = n } } \Pi_n \left(\Pi_{n_1} v \overline{\Pi_{n_2} v } \Pi_{n_3} v \right)
\end{align*}
and $\Pi_n$ is the orthogonal projector on the $n^{th}$ eigenspace $E_n$ of $-\Delta_x + x^2$.

We will mainly follow the approach to the global well-posedness and scattering of the two-dimensional mass-critical nonlinear Schr\"odinger equation as in \cite{D1}. The main ingredient is to establish an infinite dimensional vector-valued version of 2D long-time Strichartz estimate, which helps us to preclude certain almost periodic solutions.

The \eqref{eq1.3v28} system is Hamiltonian with an energy functional
\begin{align*}
\mathcal{E}(v) = \frac12 \sum\limits_{n\in \mathbb{N} } \int_{\mathbb{R}^2\times \mathbb{R}} |\nabla_y v_n|^2\,\mathrm{d}y\mathrm{d}x
+ \frac14 \sum\limits_{\substack{n, n_1,n_2,n_3,n_4 \in \mathbb{N} ,\\ n_1-n_2 + n_3 - n_4 = n}} \int_{\mathbb{R}^2\times\mathbb{R}} v_{n_1} \bar{v}_{n_2} v_{n_3} \bar{v}_{n_4} \,\mathrm{d}y\mathrm{d}x,
\end{align*}
under the symplectic structure on $L_{y,x}^2(\mathbb{R}^2\times \mathbb{R})$ given by $\omega(f,g):=\Im \int_{\mathbb{R}^2 \times \mathbb{R}} f(y,x) \overline{g(y,x)} \,\mathrm{d}y \mathrm{d}x$.
It also conserves the following mass $\mathcal{M}$ and kinetic energy $\mathcal{E}_0$:
\begin{align*}
\mathcal{M}(v) & = \int_{\mathbb{R}^2\times \mathbb{R}} |v(t,y,x)|^2\,\mathrm{d}y\mathrm{d}x,\\
\mathcal{E}_0(v) & = \int_{\mathbb{R}^2\times \mathbb{R}}  \left|\left(-\Delta_x + x^2\right)^\frac12 v(t,y,x)  \right|^2 \,\mathrm{d}y\mathrm{d}x
 = \sum\limits_{n \in \mathbb{N} }  (2n+1) \|v_n\|_{L_{y,x}^2 (\mathbb{R}^2 \times \mathbb{R})}^2  = \|v\|_{L_y^2 \mathcal{H}_x^1(\mathbb{R}^2 \times \mathbb{R})}^2.
\end{align*}

We shall divide this section into three subsections. In Section \ref{subsec:6.1}, we establish the local well-posedness theory for \eqref{eq1.3v28} and reduce the scattering to the exclusion of almost periodic solutions. In Section \ref{subsec:6.2}, we derive the long time Strichartz estimate and in Section \ref{subsec:6.3}, we exclude the almost periodic solution.

\subsection{Local well-posedness and reduction to the almost periodic solution}\label{subsec:6.1}
In this subsection, we will present the well-posedness theory of the \eqref{eq1.3v28}. Then following similar ideas as in \cite{TVZ,CGYZ,YZ0,CGZ}, we shall prove that there is an almost periodic solution of \eqref{eq1.3v28} if the system is not global well-posed and if the solution does not scatter in $L_y^2 \mathcal{H}_x^1$. That is we reduce the global well-posedness and scattering of \eqref{eq1.3v28} to the exclusion of this almost periodic solution.
\subsubsection{Local well-posedness theory and the existence of almost periodic solution}
The local well-posedness theory of the \eqref{eq1.3v28} system follows from a more or less standard argument: the Strichartz estimate in Proposition \ref{pr2.2v31} and the nonlinear estimate in Lemma \ref{le4.8v56}. The proof of the nonlinear estimate relies on the following Strichartz estimate for the harmonic oscillator.

\begin{lemma}[Strichartz estimate for the harmonic oscillator,\cite{C,KT0}]\label{le2.1v27}
For $2 \le q, r\le \infty $ with $\frac2q + \frac1r = \frac12$, we have the following estimate
\begin{align*}
\left\|e^{it\left(\Delta_x - x^2 \right)} f\right\|_{L_t^q \mathcal{W}_x^{s,r}([-T_1,T_1] \times \mathbb{R})} \lesssim \|f\|_{\mathcal{H}^s_x(\mathbb{R})}.
\end{align*}
holds for any $T_1  > 0$ and $s\ge 0$.
\end{lemma}
We can now give the nonlinear estimate.
\begin{lemma}\label{le4.8v56}
For functions $F_1,F_2,F_3$ defined on $\mathbb{R}^2 \times \mathbb{R}$, we have
\begin{align}\label{eq3.8v41}
\bigg\|\sum\limits_{\substack{n_1,n_2,n_3,n\in \mathbb{N} ,\\ n_1-n_2 + n_3 = n}} \Pi_n\left( \Pi_{n_1} F_1 \overline{\Pi_{n_2} F_2} \Pi_{n_3} F_3 \right)\bigg\|_{L_y^\frac43 L_x^2}
\lesssim  \left\|F_1  \right\|_{L_y^4 L_x^2} \left\|F_2 \right\|_{L_y^4 L_x^2 } \left\|F_3 \right\|_{L_y^4 L_x^2 },
\end{align}
and consequently, for any $\beta \ge 0$,
\begin{align}\label{eq3.9v41}
\bigg\|\sum\limits_{\substack{n_1,n_2,n_3,n\in \mathbb{N} ,\\ n_1-n_2 + n_3 = n}} \Pi_n\left( \Pi_{n_1} F_1 \overline{\Pi_{n_2} F_2} \Pi_{n_3} F_3 \right)\bigg\|_{L_y^\frac43 \mathcal{H}_x^\beta}
\lesssim \min\limits_{\tau \in \mathcal{\sigma}_3} \left\|F_{\tau(1)}\right\|_{L_y^4 \mathcal{H}_x^\beta} \left\|F_{\tau(2)} \right\|_{L_y^4 L_x^2} \left\|F_{\tau(3)}\right\|_{L_y^4 L_x^2},
\end{align}
where $\sigma_3$ is a permutation of the set $\{1,2,3\}$.
\end{lemma}
\begin{proof}
Let $F_0\in L_y^4 L_x^2$, by H\"older's inequality and Lemma \ref{le2.1v27}, we have
\begin{align*}
& \left\langle \sum\limits_{\substack{n_1,n_2,n_3,n \in \mathbb{N} ,\\ n_1-n_2 + n_3 =n}} \Pi_n\left(\Pi_{n_1} F_1 \overline{\Pi_{n_2} F_2} \Pi_{n_3} F_3 \right), F_0\right\rangle \\
= & \frac1\pi \sum\limits_{n_1,n_2,n_3,n \in \mathbb{N} } \int_0^\pi e^{2it(n_1-n_2 +n_3 -n)} \int_{\mathbb{R}^2 \times \mathbb{R}} \Pi_{n_1} F_1 \overline{\Pi_{n_2} F_2} \Pi_{n_3} F_3 \overline{\Pi_n F_0} \,\mathrm{d}y\mathrm{d}x \mathrm{d}t\\
= & \frac1\pi \int_0^\pi \int_{\mathbb{R}^2 \times \mathbb{R}} e^{it \left(-\Delta_x + x^2 \right)} F_1(y,x) \overline{e^{it(-\Delta_x + x^2)} F_2(y,x)} e^{it\left(-\Delta_x + x^2 \right)} F_3(y,x) \overline{e^{it \left(-\Delta_x + x^2 \right)} F_0(y,x) } \,\mathrm{d}y\mathrm{d}x \mathrm{d}t\\
\lesssim & \int_{\mathbb{R}^2} \left\|e^{it \left(-\Delta_x + x^2 \right)} F_0(y,x) \right\|_{L_t^\infty L_x^2([0,\pi]\times \mathbb{R})} \left\|e^{it \left(-\Delta_x + x^2 \right)} F_1(y,x) \right\|_{L_t^2 L_x^4 ([0,\pi] \times \mathbb{R})} \left\|e^{it\left(-\Delta_x + x^2 \right)} F_2(y,x) \right\|_{L_t^4 L_x^8([0,\pi] \times \mathbb{R})} \\
& \quad \left\|e^{it\left(-\Delta_x +  x^2 \right)} F_3(y,x) \right\|_{L_t^4 L_x^8([0,\pi]\times \mathbb{R})} \,\mathrm{d}y\\
\lesssim & \int_{\mathbb{R}^2} \left\|F_0(y,x) \right\|_{L_x^2(   \mathbb{R})} \left\|e^{it(-\Delta_x + x^2)} F_1(y,x) \right\|_{L_t^8 L_x^4 ([0,\pi]\times \mathbb{R})} \left\|e^{it \left(-\Delta_x + x^2 \right)} F_2(y,x) \right\|_{L_t^\frac{16}3 L_x^8([0,\pi]\times \mathbb{R})} \\
& \qquad \left\|e^{it \left(-\Delta_x +  x^2 \right)} F_3(y,x)  \right\|_{L_t^\frac{16}3 L_x^8([0,\pi]\times \mathbb{R})} \,\mathrm{d}y\\
\lesssim & \int_{\mathbb{R}^2} \left\|F_0(y,x)\right\|_{L_x^2} \left\|F_1(y,x)\right\|_{L_x^2} \left\|F_2(y,x) \right\|_{L_x^2} \left\|F_3(y,x) \right\|_{L_x^2}\,\mathrm{d}y.
\end{align*}
Therefore,
\begin{align*}
\bigg\|\sum\limits_{\substack{n_1,n_2,n_3,n \in \mathbb{N} ,\\ n_1-n_2 + n_3 = n}} \Pi_n \left(\Pi_{n_1} F_1 \overline{\Pi_{n_2} F_2} \Pi_{n_3} F_3 \right) \bigg\|_{L_y^\frac43 L_x^2}
\lesssim \left\|F_1 \right\|_{L_y^4 L_x^2} \left\|F_2\right\|_{L_y^4 L_x^2} \left\|F_3\right\|_{L_y^4 L_x^2},
\end{align*}
which is \eqref{eq3.8v41}. One can similarly prove \eqref{eq3.9v41} using the fractional Leibniz rule.
\end{proof}
Lemma \ref{le4.8v56} provides the following estimate for the nonlinearity $F(v)$.
\begin{lemma}\label{pr2.3v51}
For each solution $v$ of \eqref{eq1.3v28}, we have
\begin{align*}
\| F(v) \|_{L_{t,y}^\frac43 \mathcal{H}_x^\alpha }
\lesssim \|v\|_{L_{t,y}^4 \mathcal{H}_x^\alpha }^3, \text{  where } \alpha = 0, 1.
\end{align*}
Thus by Proposition \ref{pr2.2v31}, the solution $v$ of \eqref{eq1.3v28} satisfies the Strichartz estimate
\begin{align}\label{eq2.6v51}
\|v\|_{L_t^p L_y^q \mathcal{H}_x^\alpha (I \times \mathbb{R}^2 \times \mathbb{R})}
\lesssim \|v_0 \|_{L_y^2 \mathcal{H}_x^\alpha (\mathbb{R}^2\times \mathbb{R})}
+ \|v\|_{L_{t,y}^4 \mathcal{H}_x^\alpha(I\times \mathbb{R}^2 \times \mathbb{R})}^3, \ \text{  for }\, \alpha = 0,1,
\end{align}
where $I \subseteq \mathbb{R}$, and $(p,q)$ is Strichartz admissible pair.
\end{lemma}
As a consequence of Lemma \ref{pr2.3v51} and  \eqref{eq3.9v41}, we obtain the following well-posedness theory. Since the proof is well-known (see for instance \cite{CGYZ,CGZ,T2,Killip-Visan1}), we omit it.
\begin{theorem}[Well-posedness and scattering of the equation \eqref{eq1.3v28}] \label{th2.4v51}
\
\begin{enumerate}
\item (Local well-posedness)
Assume $ \left\|v_0 \right\|_{L_y^2 \mathcal{H}^1_x} <\infty$. The \eqref{eq1.3v28} admits a unique solution
\begin{align*}
v \in \left(C_t^0 L_y^2 \mathcal{H}_x^1 \cap L_{t,y}^4 \mathcal{H}_x^1\right) \left((-T,T) \times \mathbb{R}^2 \times \mathbb{R}\right)
\end{align*}
for some $T > 0$.
\item (Small data scattering)
There is a sufficient small constant $\delta >0$, such that when $\|v_0 \|_{L_y^2 \mathcal{H}_x^1} \le \delta$, \eqref{eq1.3v28} admits a unique global solution $v
$ with $v(0) = v_0$, which scatters in $L_y^2 \mathcal{H}_x^1$ in the sense that there exist $v^\pm \in L_y^2 \mathcal{H}_x^1(\mathbb{R}^2 \times \mathbb{R})$, such that
\begin{align*}
\left\|v(t) - e^{it \Delta_y} v^\pm \right\|_{L_y^2 \mathcal{H}_x^1 } \to 0, \text{ as } t \to \pm \infty.
\end{align*}
\item (Scattering norm) Suppose $v$ is a maximal lifespan solution on $I$ with
$\|v\|_{L_{t,y}^4 L_x^2(I \times \mathbb{R}^2 \times \mathbb{R})} < \infty,$ then $v$ globally exists and scatters in $L_y^2 \mathcal{H}_x^1$.
\end{enumerate}
\end{theorem}

We also have the stability theorem by Lemmas \ref{le4.8v56} and \ref{pr2.3v51}. The argument is similar to the proof of Theorem \ref{th2.6}, and we also refer to \cite{CKSTT0,Killip-Visan1}.
\begin{theorem}[Stability] \label{th2.6v51}
Let $ l   \in \{0,1\}$, $I$ be a compact interval and $\tilde{v}\in \left(C_t^0 L_y^2 \mathcal{H}_x^1 \cap L_{t,y}^4 \mathcal{H}_x^1 \right) \left(I \times \mathbb{R}^2 \times \mathbb{R}\right)$ be an approximate solution of \eqref{eq1.3v28} with the error term
$e = i \partial_t \tilde{v}  + \Delta_y \tilde{v} - F( \tilde{v})$.
Then for any $\epsilon > 0$, there is $\delta > 0$ such that if
\begin{align*}
\|e\|_{L_{t,y}^\frac43 \mathcal{H}_x^l (I \times \mathbb{R}^2 \times \mathbb{R})} + \left\|\tilde{v} (t_0) - v_0  \right\|_{L_y^2 \mathcal{H}_x^1} \le \delta,
\end{align*}
then \eqref{eq1.3v28} admits a solution $v \in \left(L_t^\infty L_y^2 \mathcal{H}_x^1 \cap L_{t,y}^4 \mathcal{H}_x^l  \right) (I \times \mathbb{R}^2 \times \mathbb{R})$ with $v(t_0) = v_0$ and %satisfying
\begin{align*}
\left\| \tilde{v} - v \right\|_{L_{t,y}^4 \mathcal{H}_x^l  \cap L_t^\infty L_y^2 \mathcal{H}_x^1(I \times \mathbb{R}^2 \times \mathbb{R})} <  \epsilon.
\end{align*}
\end{theorem}
To prove \eqref{eq1.3v28} is globally well-posed and scatters for large data, by Theorem \ref{th2.4v51}, we need to prove
\begin{align*}
\|v \|_{L_{t,y}^4 L_x^2(\mathbb{R} \times \mathbb{R}^2 \times \mathbb{R})} < \infty,
\end{align*}
where $v$ is a solution to \eqref{eq1.3v28} with initial data $v_0 \in L_y^2 \mathcal{H}_x^1(\mathbb{R}^2 \times \mathbb{R})$.
For the solution $v$ of \eqref{eq1.3v28} with maximal lifespan interval $I$, let
\begin{align*}
S(m) = \sup \left\{ \|v \|_{L_{t,y}^4 L_x^2(I \times \mathbb{R}^2 \times \mathbb{R})}: \|v (0)\|_{L_y^2 \mathcal{H}_x^1(\mathbb{R}^2 \times \mathbb{R})} \le m \right\},
\end{align*}
and
\begin{align*}
m_0 = \sup \left\{ m: S ( \tilde{m} ) < \infty, \forall\,  \tilde{m} < m \right\} > 0.
\end{align*}
If have $m_0 = \infty$, then the global well-posedness and scattering in $L_y^2 \mathcal{H}_x^1$ of \eqref{eq1.3v28}  hold. Following the argument in \cite{TVZ,Killip-Visan1}, and using Theorems \ref{th2.6v51} and \ref{pro3.9v23} during the proof,
we have
\begin{theorem}[Existence of an almost periodic solution to \eqref{eq1.3v28}]\label{th3.3v51}
Assume $m_0 < \infty$. Then there exists an non-zero almost periodic solution $v \in C_t^0 L_y^2 \mathcal{H}_x^1 \cap L_{t,y}^4 L_x^2(I \times \mathbb{R}^2 \times \mathbb{R})$ to \eqref{eq1.3v28} with $I$ the maximal lifespan interval such that $\mathcal{M}( v ) = m_0$. In addition,  for any $\eta > 0$, there exists   $C(\eta) > 0$ and $\left( y (t), \xi(t), N(t) \right) \in \mathbb{R}^2 \times \mathbb{R}^2 \times \mathbb{R}_+ $ such that
\begin{align}\label{eq6.3v55}
\int_{| y - y (t)| \ge \frac{C(\eta)}{N(t)}} \left\| v (t,y,x) \right\|_{\mathcal{H}_x^1}^2 \, \mathrm{d}y  + \int_{ |\xi - \xi(t) | \ge C(\eta) N(t)} \left\|(\mathcal{F}_y  v )(t, \xi, x) \right\|_{\mathcal{H}_x^1}^2 \,\mathrm{d}\xi <  \eta, \  \forall\, t \in I .
\end{align}
Furthermore, we can take $[0, \infty) \subseteq I$, and $N(0) = 1, \xi(0) = y(0) = 0$, with
\begin{align*}
 N(t) \le 1 ,  \ \left|N'(t) \right| + \left|\xi'(t) \right| \lesssim N(t)^3, \ \forall\, t \in [0, \infty).
\end{align*}
\end{theorem}

As in \cite{D1,D2,D3,Killip-Visan1,R}, we see the almost periodic solution in Theorem \ref{th3.3v51} has the following property:
\begin{theorem}\label{th3.9v51}
(1)  If $J \subseteq I$ is an interval which is partitioned into small intervals $J_k$ in the sense that $\| v \|_{L_{t,y}^4 L_x^2(J_k \times \mathbb{R}^2 \times \mathbb{R})} = 1$, then we have
\begin{align}\label{eq6.4v58}
N(J_k) \sim \int_{J_k} N(t)^3 \,\mathrm{d}t \sim \inf\limits_{t \in J_k} N(t), \text{ and } \sum\limits_{J_k \subseteq J } N(J_k) \sim \int_{J} N(t)^3 \,\mathrm{d}t,
\end{align}
 where $N(J_k) = \sup\limits_{t \in J_k} N(t)$.

(2) For any interval $J \subseteq [0, \infty)$, we have
\begin{align}\label{eq6.6v122}
\int_J N(t)^2 \,\mathrm{d}t \lesssim \| v \|_{L_{t,y}^4 L_x^2(J\times \mathbb{R}^2 \times \mathbb{R})}^4 \lesssim 1 + \int_J N(t)^2 \,\mathrm{d}t.
\end{align}
\end{theorem}
{
\begin{proof}[Proof of Theorems \ref{th3.3v51} and \ref{th3.9v51}]
With Theorems \ref{pro3.9v23} and \ref{th2.6v51} at hand, one can follow the arguments in \cite{CGZ,CGYZ,D1,D2,D3,R,TVZ,Killip-Visan1}.
\end{proof}}

\subsubsection{Some functional spaces and bilinear Strichartz estimates}
{
As in \cite{D1}, due to the failure of the endpoint Strichartz estimate in 2-D, we need to utilize the function spaces $U^p_\Delta$ and $V^p_\Delta$ introduced originally in the seminal work of Koch and Tataru \cite{KT3}; see also \cite{HHK,KT1,KTV} for more detailed study on these spaces.
The structure of our  \eqref{eq1.3v28} system motivates us to introduce the Banach spaces $U_\Delta^p(L_x^2)$ and $V_\Delta^p(L_x^2)$ as follows.}
\begin{definition}[$U_\Delta^p(L_x^2)$ space]\label{de4.1v51}
For  $1 \le p < \infty$, let $U_\Delta^p(L_x^2)$ be an atomic space, where an atom $ v^\gamma$ is defined to be
\begin{align*}
 v^\gamma  (t,y,x) = \sum\limits_{k = 0}^N  \chi_{[t_k, t_{k+1})}(t) e^{it \Delta_y} v_k^\gamma (y,x), \text{ with }
\sum\limits_{k = 0}^N   \left\| v_k^\gamma(y,x) \right\|_{L_{y,x}^2 }^p = 1.
\end{align*}
In the expansion of $v^\gamma$, $N$ may be finite or infinite, $t_0 = - \infty$, and $t_{N+1} = \infty$ if $N$ is finite.
We impose a norm on $\| \cdot \|_{U_\Delta^p(L^2_x)}$ as
\begin{align*}
\|v  \|_{U_\Delta^p (L^2_x)} = \inf \left\{ \sum\limits_\gamma  |c_\gamma |:  v = \sum\limits_\gamma c_\gamma v^\gamma , \text{ where }  v^\gamma  \text{ are } U_\Delta^p(L_x^2) \text{ atoms} \right\}.
\end{align*}
For a time interval $I \subseteq \mathbb{R}$, we define
\begin{align*}
\|v\|_{U_\Delta^p(L^2_x, I) } = \| v 1_I \|_{U_\Delta^p(L^2_x)}.
\end{align*}
Let $DU_\Delta^p(L_x^2)$ be the space
\begin{align*}
DU_\Delta^p(L_x^2) = \left\{ \left( i \partial_t+ \Delta_y \right) v  : \  v  \in U_\Delta^p(L_x^2) \right\},
\end{align*}
endowed with the following norm
\begin{align*}
\left\| \left(i\partial_t + \Delta_y \right)  v (t,y,x) \right\|_{DU_\Delta^p(L_x^2)}  = \left\| \int_0^t e^{i(t-s)\Delta_y } (i \partial_s + \Delta_y )  v (s,y,x) \,\mathrm{d}s \right\|_{U_\Delta^p(L_x^2)}.
\end{align*}
For each time interval $I \subseteq \mathbb{R}$, we can similarly define the restriction space $DU^p_\Delta(L^2_x,I)$.
\end{definition}

\begin{definition}[$V_\Delta^p(L_x^2)$ space]\label{de4.2v51}
For $1 \le p < \infty$, $V_\Delta^p(L_x^2)$ is defined to be the space of right continuous functions $v\in L_t^\infty L_{y,x}^2  $ such that
\begin{align*}
\|v\|_{V_\Delta^p(L^2_x)}^p = \|v\|_{L_t^\infty L_{y,x}^2 }^p + \sup\limits_{\{ t_k\}_k \nearrow } \sum\limits_k \left\|e^{-it_{k+1} \Delta_y } v(t_{k+1} ) - e^{-it_k \Delta_y } v(t_k ) \right\|_{L_{y,x}^2 }^p<\infty.
\end{align*}
When the time is restricted to $I \subseteq \mathbb{R}$, We can similarly define the function space  $V_\Delta^p(L_x^2, I)$. Then we have
\begin{align}\label{eq6.5v57}
	\left( DU_\Delta^p(L_x^2) \right)^\ast  =   V_\Delta^{p'}(L_x^2) .
\end{align}
\end{definition}

The following basic properties are strightfoward to verify. For the proofs, see \cite{HHK,KTV}.
\begin{remark}[Basic properties of $U_\Delta^p(L_x^2)$ and $V_\Delta^p(L_x^2)$]\label{pr4.3v51}
For any $  1 < p < q < \infty$ and $a\le b\le c$, we have
\begin{align}
 & U_\Delta^p(L_x^2) \subseteq V_\Delta^p(L_x^2) \subseteq U_\Delta^q(L_x^2),\label{eq6.4v57}\\
%\end{align}
%\begin{align}
& \| v  \|_{U_\Delta^p \left(L_x^2,  [a,b] \right)} \le \| v \|_{U_\Delta^p \left(L_x^2,  [a, c] \right)},
\text{ and }
   \| v  \|_{U_\Delta^p \left(L_x^2,  [ a, c]  \right)}^p \le \| v  \|_{U_\Delta^p \left(L_x^2,  [ a,b]  \right)}^p + \| v  \|_{U_\Delta^p \left(L_x^2,  [b, c]  \right)}^p, \label{eq6.6v57}\\
\intertext{ and }
& \| v \|_{U_\Delta^p(L_x^2)} \lesssim  \| v|_{t=0} \|_{L_{y,x}^2  } + \left\|\left(i\partial_t + \Delta_y \right) v  \right\|_{DU_\Delta^p(L_x^2)}. \label{eq6.8v57}
\end{align}
Moreover,
\begin{align}
 & L_t^{p'}  L_y^{r'}  L_x^2 \subseteq DU_\Delta^2(L^2_x), \text{ and }
 U^p_\Delta  (L_x^2) \subseteq  L_t^p L_y^r L_x^2, \text{ where  $(p,r )$ is Strichartz admissible.} \label{eq6.7v57}
\end{align}
\end{remark}

Following the argument in \cite{D1}, we also have
\begin{lemma}\label{le4.4v51}
Suppose $I = \bigcup_{j =1 }^m  I^j $, where $I^j = [a_j,b_j ]$, $a_{j +1} = b_j$. If $ f \in L_t^1 L_{y,x}^2 (I  \times \mathbb{R}^2 \times \mathbb{R})$, then $\forall\,  t_0 \in I  $,
\begin{align*}
\left\|\int_{t_0}^t e^{i(t-\tau) \Delta_y } f (\tau, y,x) \,\mathrm{d}\tau  \right\|_{U_\Delta^2(L_x^2, I )}
\lesssim \sum\limits_{j = 1}^m   \left\|\int_{I^j } e^{-i\tau \Delta_y}  f (\tau) \,\mathrm{d}\tau \right\|_{L_{y,x}^2 }
+ \left( \sum\limits_{j = 1}^m  \| f \|_{DU_\Delta^2(L_x^2, I^j  )}^2 \right)^\frac12,
\end{align*}
where
\begin{align*}
\|f \|_{DU_\Delta^2(L_x^2, I^j  )}
= \sup\limits_{ \|w  \|_{V_\Delta^2 \left(L_x^2, I^j  \right)} = 1} \int_{ I^j } \int_{\mathbb{R}^2 \times \mathbb{R}}  f(\tau, y, x) \overline{  w (\tau, y,x) }  \,\mathrm{d}\tau \mathrm{d} y \mathrm{d}x.
\end{align*}
\end{lemma}

By the bilinear Strichartz estimate in \cite{Bo}, Minkowski's inequality, H\"older's inequality, and interpolation, we have the following two propositions. The proofs are similar to the bilinear Strichartz estimates in \cite{D1,D2}.
\begin{proposition}[Bilinear Strichartz estimate I] \label{pr4.5v51}
Let $(p,q)$ satisfy $1< p , q < \infty$, $\frac1p + \frac1q = 1$.
For $M\ll N$, assume $supp \, \mathcal{F}_y u_0 \subseteq \left\{ \xi: |\xi| \sim N \right\}$ and $supp \,\mathcal{F}_y v_0 \subseteq \left\{ \xi: |\xi| \sim M \right\}$. Then we have
\begin{align}\label{eq6.11v122}
\left\| \left\|e^{it \Delta_y} u_0 \right\|_{L_x^2} \left\|e^{it\Delta_y} v_0 \right\|_{L_x^2} \right\|_{L_t^p L_y^q (\mathbb{R}\times \mathbb{R}^2)}
\lesssim \left(\frac{M}{N} \right)^\frac1p \|u_0\|_{L_{y,x}^2 } \|v_0 \|_{L_{y,x}^2}.
\end{align}
Furthermore, suppose that $g\left(t, y - \tilde{y} \right)$ and $h \left(t, y - \tilde{\tilde{y}} \right)$ are convolution kernels {with respect to $y$-variable} and
\begin{align*}
\left\| \sup\limits_{t \in \mathbb{R}} |g(t,y)| \right\|_{L_y^1(\mathbb{R}^2)} +
\left\|\sup\limits_{t \in \mathbb{R}} |h(t,y)| \right\|_{L_y^1(\mathbb{R}^2)} \lesssim 1,
\end{align*}
we also have
\begin{align*}
\left\| \left\|g\ast_y e^{it\Delta_y} u_0 \right\|_{L_x^2}\left \|h \ast_y e^{it \Delta_y} v_0 \right\|_{L_x^2} \right\|_{L_t^p L_y^q(\mathbb{R} \times \mathbb{R}^2)}
\lesssim \left(\frac{M}N \right)^\frac1p \left\|u_0 \right\|_{L_{y,x}^2 } \|v_0 \|_{L_{y,x}^2 }.
\end{align*}
\end{proposition}
Similar to the argument in the proof of Lemma 3.5 in \cite{D1}, we can transfer the estimate \eqref{eq6.11v122} to the $U_\Delta^p$ space. Therefore, we have
\begin{proposition}[Bilinear Strichartz estimate II]\label{pr4.7v51}
Let $(p,q)$ satisfy $1< p , q < \infty$, $\frac1p + \frac1q = 1$.
For $M \ll N$, assume $supp\, \mathcal{F}_y u \subseteq \left\{ \xi: |\xi| \sim N \right\}$ and $supp \, \mathcal{F}_y v \subseteq  \left\{ \xi: |\xi| \sim M \right\}$. Then we have
\begin{align*}
\left\| \|u\|_{L_x^2} \|v\|_{L_x^2} \right\|_{L_t^p L_y^q} \lesssim \left( \frac{M}N \right)^\frac1p \|u\|_{U_\Delta^p(L_x^2)} \|v\|_{U_\Delta^p(L_x^2)}.
\end{align*}
\end{proposition}

\subsection{Long time Strichartz estimate}\label{subsec:6.2}
From now on, we shall take our following setting as standarding assumpmtions. Fix
\begin{align}\label{eq6.13v109}
	0< \epsilon_3 \ll \epsilon_2 \ll \epsilon_1 < 1  \text{ and } \epsilon_3< \epsilon_2^{10}.
\end{align}
By Theorem \ref{th3.3v51}, we can take
\begin{align}
	& |N'(t)| + |\xi'(t)| \le 2^{-20} \epsilon_1^{- \frac12}  {N(t)^3}, \label{eq5.1v51}\\
	\intertext{ and }
	& \int_{|y- y(t)| \ge \frac{2^{-20} \epsilon_3^{- \frac14}}{N(t)} } \|  v (t,y,x)\|_{\mathcal{H}_x^1}^2 \,\mathrm{d}y + \int_{ |\xi - \xi(t)| \ge 2^{-20} \epsilon_3^{- \frac14 } N(t)} \left\| \left( \mathcal{F}_y  v \right) (t, \xi,x) \right\|_{\mathcal{H}_x^1}^2 \,\mathrm{d}\xi \le \epsilon_2^2. \label{eq5.2v51}
\end{align}

If $[0,T]$ is an interval with
\begin{align}\label{eq6.15v122}
	\| v \|_{L_{t,y}^4 L_x^2([0,T] \times \mathbb{R}^2 \times \mathbb{R})}^4 = 2^{k_0}
	\text{ and }   \int_0^T N(t)^3 \,\mathrm{d}t = \epsilon_3 2^{k_0} , \text{ for some }  k_0 \ge 0,
\end{align}
then we can partition $[0,T] = \bigcup\limits_{ \alpha =0}^{M-1} J^\alpha $, where $J^\alpha $ are intervals that satisfy
\begin{align}\label{eq5.3v51}
		\int_{J^\alpha} \left(  N(t)^3 + \epsilon_3 \| v (t)\|_{L_y^4 L_x^2 \left(\mathbb{R}^2 \times \mathbb{R} \right)}^4  \right) \,\mathrm{d}t = 2 \epsilon_3.
\end{align}
We can define the interval $G_k^j$ now.
\begin{definition}\label{de5.1v51}
	For any nonnegative integer $j < k_0$, and nonnegative integer $ k < 2^{k_0 - j}$, we can define
	\begin{align}\label{eq6.13v58}
		G_k^j = \bigcup_{\alpha= k2^j}^{(k+1)2^j - 1} J^\alpha.
	\end{align}
For $j \ge k_0$, we simply define  $G_k^j = [0,T]$. We let $\xi \left(G_k^j \right) = \xi \left(t_k^j \right) $, where $t_k^j$ is the left endpoint of $G_k^j$.
\end{definition}
On the time interval $G_k^j$ defined above, we have
\begin{lemma} \label{re5.2v51}
	(1) Let $J_l $ be the small intervals contained in $G_k^j$. By \eqref{eq6.4v58} and \eqref{eq5.3v51}, the following estimate holds:
	\begin{align}\label{eq5.4v51}
		\sum\limits_{J_l  \subseteq G_k^j} N \left(J_l  \right) \lesssim \int_{G_k^j} N(t)^3 \,\mathrm{d}t \lesssim \sum\limits_{\alpha = k 2^j}^{(k+1) 2^j- 1} \int_{J^\alpha} N(t)^3 \,\mathrm{d}t \lesssim 2^j \epsilon_3.
	\end{align}
	(2) By \eqref{eq5.1v51} and Definition \ref{de5.1v51}, we have { for each $t\in G_k^j$},
	\begin{align}\label{eq6.20v124}
		\left| \xi(t) - \xi \left(G_k^j \right) \right| \le 2^{j- 19} \epsilon_3 \epsilon_1^{- \frac12}.
	\end{align}
	Thus, for any $t \in G_k^j$, and $i \ge j$,
	\begin{align}\label{eq6.17v58}
		\left\{ \xi : 2^{i-1} \le \left|\xi - \xi(t) \right| \le 2^{i+ 1} \right\}\subseteq  \left\{ \xi: 2^{i-2} \le \left|\xi - \xi \left(G_k^j \right) \right| \le 2^{i+2} \right\}\subseteq  \left\{ \xi: 2^{i-3} \le \left|\xi - \xi(t) \right| \le 2^{i+3}  \right\},
	\end{align}
	and also
	\begin{align}\label{eq6.21v122}
		\left\{ \xi: |\xi - \xi(t)| \le 2^{i+1} \right\} \subseteq  \left\{ \xi: \left|\xi- \xi \left(G_k^j \right) \right| \le 2^{i+2} \right\} \subseteq \left\{ \xi: |\xi - \xi(t) | \le 2^{i+3} \right\}.
	\end{align}
\end{lemma}

\begin{lemma}\label{le5.3v51}
	For the almost periodic solution $v(t)$  to \eqref{eq1.3v28}, and assume $\| v \|_{L_{t,y}^4 L_x^2 \left(J \times \mathbb{R}^2 \times \mathbb{R} \right)} \le 1$ on $J \subseteq \mathbb{R}$, then we have
	\begin{equation*}
		\| v  \|_{U_\Delta^2  \left(L_x^2, J \right)} \lesssim 1,
		\text{ and }
		\bigg\| P^y_{> {2^{- 4 }  \epsilon_3^{- \frac14 } }    {N(J)} }  v  \bigg\|_{U_\Delta^2 \left(L_x^2, J \right)} \lesssim \epsilon_2,
	\end{equation*}
	where $N(J) = \sup\limits_{t \in J} N(t)$.
\end{lemma}
\begin{proof}
	Let $J = [t_0 , t_1 ]$, by \eqref{eq6.8v57}, \eqref{eq6.4v57}, \eqref{eq6.5v57}, and \eqref{eq6.7v57}, we have
	\begin{align*}
		\|v \|_{U_\Delta^2(L_x^2, J)} \lesssim  \|v(t_0)  \|_{L_{y,x}^2 } + \| v \|_{L_{t,y}^4 L_x^2(J \times \mathbb{R}^2 \times \mathbb{R})}^3 \lesssim 1.
	\end{align*}
	By \eqref{eq5.2v51}, we have
	\begin{align*}
		\bigg\|P^y_{>   { 2^{- 20} \epsilon_3^{- \frac14} }       {N(J)}   } v \bigg\|_{L_t^\infty L_{y,x}^2
			 \left(J \times \mathbb{R}^2 \times \mathbb{R} \right)}
		\le \bigg \|P^y_{>     { 2^{- 20} \epsilon_3^{- \frac14} }    {N(t)}  } v  \bigg\|_{L_t^\infty L_{y,x}^2
			\left(J \times \mathbb{R}^2 \times \mathbb{R} \right)} \le \epsilon_2.
	\end{align*}
	Therefore, by Strichartz estimate, we have
	\begin{align*}
		\bigg\|P^y_{> { 2^{- 4 }  \epsilon_3^{- \frac14 } }   {N(J)}    } v \bigg\|_{U_\Delta^2 \left(L_x^2, J \right)}
		& \lesssim \bigg \|P^y_{> {2^{- 20} \epsilon_3^{- \frac14} }  {N(J)}  } v(t_0)  \bigg\|_{ L_{y,x}^2
			 \left(J \times \mathbb{R}^2 \times \mathbb{R} \right)} + \bigg\|P^y_{>   {2^{- 4}  \epsilon_3^{- \frac14} }    {N(J)}  } F(v ) \bigg\|_{L_t^\frac32 L_y^\frac65 L_x^2 \left(J \times \mathbb{R}^2 \times \mathbb{R} \right)} \\
		& \lesssim \epsilon_2 +  \bigg\|P^y_{>  {2^{- 20} \epsilon_3^{- \frac14} }   {N(J)}  } v\bigg\|_{L_t^\infty L_{y,x}^2
		} \|v  \|_{L_t^3 L_y^6 L_x^2 \left(J \times \mathbb{R}^2 \times \mathbb{R} \right)}^2\\
		&\lesssim \epsilon_2 + \epsilon_2 \left( \left\|v ( t_0) \right\|_{L_t^\infty L_{y,x}^2
			\left(J \times \mathbb{R}^2 \times \mathbb{R} \right)} + \|v \|_{L_{t,y}^4 L_x^2(J \times \mathbb{R}^2 \times \mathbb{R} )}^3 \right)^2 \lesssim \epsilon_2.
	\end{align*}
\end{proof}

We also have the following fact as a consequence of the above lemma.
\begin{remark}\label{re5.4v51}
	If $N(J) < 2^{i-5} \epsilon_3^\frac12$, we have
	\begin{align*}
		\left\|P^y_{\xi \left(G_\alpha^i \right), 2^{i- 2}  \le \cdot \le 2^{ i+2}  } F(v ) \right\|_{L_t^\frac32 L_y^\frac65 L_x^2 \left(J\times \mathbb{R}^2 \times \mathbb{R} \right)}
		\lesssim \bigg\|P^y_{> { 2^{- 20} \epsilon_3^{- \frac14} }   {N(J)}   } v  \bigg\|_{L_t^\infty L_{y,x}^2
			 \left(J \times \mathbb{R}^2 \times \mathbb{R}  \right)}  \| v \|_{L_t^3 L_y^6 L_x^2 \left(J \times \mathbb{R}^2 \times \mathbb{R} \right)}^2 \lesssim \epsilon_2,
	\end{align*}
   where the operator $P^y_{\xi \left(G_\alpha^i \right), 2^{i- 2}  \le \cdot \le 2^{ i+2}  }$ is  given in Definition \ref{de5.5v51} below.
	Thus, for $0 \le i \le 11$, and $N \left(G_\alpha^i \right) < 2^{i-5} \epsilon_3^\frac12$, by the fact that $G_\alpha^i$ is a union of at most $ 2^{11}$ small intervals, we have
	\begin{align*}
		\left\|P^y_{\xi \left(G_\alpha^i \right), 2^{ i- 2 } \le \cdot \le 2^{ i+2} }   F(v ) \right\|_{L_t^\frac32 L_y^\frac65 L_x^2 \left(G_\alpha^i \times \mathbb{R}^2 \times \mathbb{R} \right)} \lesssim \epsilon_2.
	\end{align*}
\end{remark}

We can now define the long-time Strichartz estimate norm as in \cite{D1,D2,D3}; see also \cite{CGZ,CGYZ}.
\begin{definition}[Long-time Strichartz estimate norm] \label{de5.5v51}
	For any $G_k^j \subseteq [0,T]$, let
	\begin{align*}
		\|v  \|_{X \left(G_k^j \right)}^2 = \sum\limits_{0 \le i < j} 2^{i-j} \sum\limits_{G_\alpha^i \subseteq G_k^j} \left\|P^y_{\xi \left(G_\alpha^i \right), 2^{ i - 2 }  \le \cdot \le 2^{  i+2} }  v \right\|_{U_\Delta^2 \left(L_x^2, G_\alpha^i  \right)}^2 +
		\sum\limits_{i \ge j} \left\|P^y_{\xi \left(G_k^j \right), 2^{ i - 2 } \le \cdot \le  2^{ i+2  } }  v \right \|_{U_\Delta^2 \left(L_x^2, G_k^j \right)}^2 ,
	\end{align*}
	where
	\begin{align*}
		P^y_{\xi(t), 2^{ i - 2 }  \le \cdot\leq  2^{ i+2 } } v  = e^{iy\cdot \xi(t)} P^y_{ 2^{ i - 2 } \le \cdot \le 2^{ i+2} }  \left(e^{-iy\cdot \xi(t)} v  \right).
	\end{align*}
	We define the $\tilde{X}_{k_0}$ norm to be
	\begin{align*}
		\| v \|_{\tilde{X}_{k_0}([0,T])}^2 = \sup\limits_{0 \le j \le k_0} \sup\limits_{G_k^j \subseteq [0,T]} \| v \|_{X \left(G_k^j \right)}^2.
	\end{align*}
	For any nonnegative integer $ k_\ast \le k_0$, we take
	\begin{align}\label{eq6.19v122}
		\|v \|_{\tilde{X}_{k_\ast}([0,T])}^2 = \sup\limits_{ 0 \le j \le k_\ast} \sup\limits_{G_k^j \subseteq [0,T]} \| v  \|_{X \left(G_k^j \right)}^2 .
	\end{align}
\end{definition}

To close our bootstrap argument in the proof of the long time Strichartz estimate, we also need to introduce the following norm to measure $\tilde{X}_{k_0}$ norm of $v $ at scales much higher than $N(t)$.

\begin{definition}\label{de5.6v51}
	Let
	\begin{align*}
		\| v \|_{Y \left(G_k^j \right)}^2
		= \sum\limits_{0 \le i < j} 2^{i-j} \sum\limits_{ \substack{ G_\alpha^i \subseteq G_k^j, \\
				N \left(G_\alpha^i \right) \le 2^{i-5} \epsilon_3^\frac12 } } \left\|P^y_{\xi \left(G_\alpha^i \right), 2^{ i - 2 }  \le \cdot \le  2^{ i+2} }   v  \right\|_{U_\Delta^2 \left(L_x^2, G_\alpha^i \right)}^2
		+ \sum\limits_{\substack{ i \ge j ,\\
				i > 0, \\
				N \left(G_k^j \right) \le 2^{i-5} \epsilon^\frac12_3}}  \left\|P^y_{\xi \left(G_k^j \right), 2^{ i - 2 }  \le \cdot \le  2^{ i+2 }  }  v  \right\|_{U_\Delta^2 \left(L_x^2, G_k^j \right)}^2.
	\end{align*}
	We can define the norm  $\| v \|_{\tilde{Y}_{k_\ast}([0,T])}$ similar to \eqref{eq6.19v122} in definition \ref{de5.5v51}.
\end{definition}
For $i < j$, and the solution $v$ on the time interval $G_k^j$, we can define the Littlewood-Paley projector around $\xi(t)$ of $v$ as
\begin{align*}
P^y_{\xi(t), 2^{ i}   } v  = e^{iy\cdot \xi(t)} P^y_{ 2^{ i  } }  \left(e^{-iy\cdot \xi(t)} v  \right), \
P^y_{\xi(t), > 2^{ j }   } v  = e^{iy\cdot \xi(t)} P^y_{  > 2^{ j   } }  \left(e^{-iy\cdot \xi(t)} v  \right),
\end{align*}
Then, as a consequence of \eqref{eq6.5v57}, \eqref{eq6.4v57}, \eqref{eq6.7v57}, the Littlewood-Paley theorem and Proposition \ref{pr2.2v31}, we have the following estimates which reveal the relationship between the Strichartz norm $L_t^p L_y^q L_x^2$ of the Littlewood-Paley projector around $\xi(t)$ of $v$  and the long time Strichartz norm of $v$. We still refer to \cite{D1,D2} for the argument and without presenting the proof here.

\begin{lemma}\label{le5.7v51}
	For $i < j$, we have
	\begin{align}
		& \left\|P^y_{\xi(t), 2^i}  v  \right\|_{L_t^p L_y^q L_x^2\left(G_k^j \times \mathbb{R}^2 \times \mathbb{R} \right)} \lesssim 2^{\frac{j-i}p} \| v \|_{\tilde{X}_j \left(G_k^j \right)}, \label{eq5.10v51} \\
		& \left\|P^y_{\xi(t), \ge 2^j}  v \right \|_{L_t^p L_y^q L_x^2 \left(G_k^j \times \mathbb{R}^2 \times \mathbb{R} \right)} \lesssim \| v \|_{X(G_k^j)},  \label{eq5.11v51}
	\end{align}
	where $(p,q)$ is Strichartz admissible pair.
\end{lemma}

Our aim is to prove the long time Strichartz estimate.
\begin{theorem}[Long time Strichartz estimate] \label{th5.8v51}
For the almost periodic solution $v$ in Theorem \ref{th3.3v51}, which satisfies \eqref{eq6.13v109}, \eqref{eq5.1v51} and \eqref{eq5.2v51}, there exists a positive constant $C = C(v) $, such that for any nonnegative integer $ k_0$,
$v$ and $N(t)$ satisfy \eqref{eq6.15v122}, we have
\begin{align*}
\| v  \|_{\tilde{X}_{k_0}([0,T])} \le C .
\end{align*}
\end{theorem}

To prove Theorem \ref{th5.8v51}, it suffices to show for any $0 \le j \le k_0$ and $G_k^j \subseteq [0,T]$,
\begin{align*}
\sum\limits_{0 \le i < j} 2^{i-j} \sum\limits_{G_\alpha^i \subseteq G_k^j} \left\|P^y_{\xi \left(G_\alpha^i \right), 2^{ i- 2 } \le \cdot \le 2^{ i+2  } } v  \right\|_{U_\Delta^2 \left(L_x^2, G_\alpha^i \right)}^2
+ \sum\limits_{i \ge j}  \left\|P^y_{\xi \left(G_k^j \right), 2^{ i - 2  } \le \cdot \le  2^{ i+2 }  }  v  \right\|_{U_\Delta^2 \left(L_x^2, G_k^j \right )}^2 \le C.
\end{align*}
To reach the above estimate, we will perform an induction argument on $0 \le k_\ast \le k_0$ and then a bootstrap argument in Subsections \ref{subsec:basic inductive estimates} and \ref{subsec:bootstrap estimate}, respectively.

\subsubsection{Basic inductive estimates.}\label{subsec:basic inductive estimates}
First we show the basic estimates to start up our induction.
\begin{lemma}[Basic inductive estimate] \label{le5.10v51}
\begin{align}\label{eq6.21v59}
\| v  \|_{\tilde{X}_{0} ([0, T])} \le C, \text{ and }
\| v  \|_{\tilde{Y}_{0} ([0,T])} \le C \epsilon_2^\frac34.
\end{align}
For $0 \le k_\ast \le k_0$, we have
\begin{align} \label{eq6.22v59}
\| v \|_{\tilde{X}_{k_\ast + 1}([0,T])}^2 \le 2 \|  v  \|_{\tilde{X}_{k_\ast}([0,T])}^2 , \text{ and }
\| v \|_{ \tilde{Y}_{k_\ast + 1}([0,T])}^2 \le 2 \| v  \|_{\tilde{Y}_{k_\ast}([0,T])}^2.
\end{align}
\end{lemma}
\begin{proof}
By Lemma \ref{le5.3v51}, we have
\begin{align} \label{eq6.26v122}
\| v  \|_{U_\Delta^2(L_x^2, J^\alpha)} \lesssim 1, \text{ for any  } J^\alpha \text{ in the decomposition of }  G_k^j \text{ in } \eqref{eq6.13v58} .
\end{align}
Therefore, by Strichartz estimate, \eqref{eq6.8v57}, \eqref{eq6.7v57}, \eqref{eq6.4v57}, we have for $t_\alpha \in J^\alpha$,
\begin{align*}
 \left( \sum\limits_{i \ge 0} \left\| P^y_{\xi(J^\alpha), 2^{  i- 2 } \le \cdot \le  2^{ i+2 } }   v  \right\|_{U_\Delta^2(L_x^2, J^\alpha)}^2 \right)^\frac12
 & \lesssim \| v (t_\alpha) \|_{L_{y,x}^2 } + \| v  \|_{L_t^3 L_y^6 L_x^2(J^\alpha \times \mathbb{R}^2 \times \mathbb{R})}^3\\
&\lesssim \| v (t_\alpha ) \|_{L_{y,x}^2 } + \|  v \|_{U_\Delta^2(L^2_x, J^\alpha)}^3 \lesssim 1.
\end{align*}
Thus, $\| v \|_{\tilde{X}_{0}([0,T])} \le C$.

At the same time, by \eqref{eq5.2v51}, the conservation of mass, and \eqref{eq6.26v122}, we infer that
\begin{align*}
& \bigg( \sum\limits_{ \substack{ i \ge 0, \\
N(J^\alpha) \le \epsilon_3^\frac12 2^{i- 5}} }
\left\|P^y_{\xi(J^\alpha), 2^{ i - 2 } \le \cdot \le  2^{ i+2} }  v  \right\|_{U_\Delta^2 (L_x^2 ,J^\alpha)}^2 \bigg)^\frac12\\
 \lesssim & \left\|P^y_{\xi(J^\alpha), \ge 8 \epsilon_3^{- \frac12} N(J^\alpha)}   v (t_\alpha ) \right\|_{L_{y,x}^2  }
+ \left\|P^y_{\xi(J^\alpha), \ge 8 \epsilon_3^{- \frac12} N(J^\alpha) } F( v ) \right\|_{L_t^1 L_{y,x}^2
(J^\alpha \times \mathbb{R}^2 \times \mathbb{R})} \\
\lesssim  & \left\|P^y_{\xi(t), \ge 4 \epsilon_3^{- \frac12} N(t)}  v  \right\|_{L_t^\infty L_{y,x}^2
(J^\alpha \times \mathbb{R}^2 \times \mathbb{R})}
+  \left\|P^y_{\xi(t), \ge 4 \epsilon_3^{- \frac12} N(t)} F( v ) \right\|_{L_t^1 L_{y,x}^2
 (J^\alpha \times \mathbb{R}^2 \times \mathbb{R})} \\
\lesssim & \left\|P^y_{\xi(t), \ge \epsilon_3^{- \frac12} N(t)}  v  \right\|_{L_t^\infty L_{y,x}^2
(J^\alpha \times \mathbb{R}^2 \times \mathbb{R} )}^\frac34 \left( \|  v \|_{L_t^\infty L_{y,x}^2
 (J^\alpha \times \mathbb{R}^2 \times \mathbb{R} )}^\frac14 + \| v  \|_{U_\Delta^2 \left(L^2_x, J^\alpha \right)}^\frac94  \right) \lesssim \epsilon_2^\frac34 .
\end{align*}
Thus, by Definition \ref{de5.6v51}, we have
$\| v \|_{\tilde{Y}_{ 0}([0,T])} \le C\epsilon_2^\frac34$.

By Definition \ref{de5.1v51}, we see $G_k^{j+1} = G_{2k}^j \cup G_{2k+1}^j$ with $G_{2k}^j \cap G_{2k+1}^j = \emptyset$, then for $0 \le i \le j$, if $G_\alpha^i \subseteq G_k^{j+1}$, we have $G_\alpha^i \subseteq G_{2k}^j$ or $G_\alpha^i \subseteq G_{2k+1}^j$. Thus
\begin{align}\label{eq5.13v51}
& \ \sum\limits_{0 \le i < j+1} 2^{i- (j+1)} \sum\limits_{G_\alpha^i \subseteq G_k^{j+1}} \left\|P^y_{\xi \left(G_\alpha^i \right), 2^{ i- 2 }  \le \cdot \le  2^{ i+2 }  }  v \right\|_{U^2_\Delta \left(L_x^2, G_\alpha^i \right)}^2 \\
\le &  2^{-1} \sum\limits_{0 \le i < j} 2^{i- j} \bigg( \sum\limits_{G_\alpha^i \subseteq G_{2k}^j} \left\|P^y_{\xi \left(G_\alpha^i \right), 2^{  i - 2  } \le \cdot \le 2^{ i+2} }   v  \right\|_{U_\Delta^2 \left(L_x^2, G_\alpha^i \right)}^2
+ \sum\limits_{G_\alpha^i \subseteq G_{2k+1}^j} \left\|P^y_{\xi \left(G_\alpha^i \right) , 2^{  i - 2 }  \le \cdot \le  2^{ i+2} }   v  \right\|_{U_\Delta^2 \left(L_x^2, G_\alpha^i \right)}^2 \bigg)
\notag \\
 & \ + 2^{-1} \left( \left\|P^y_{\xi \left(G_{2k}^j \right),  2^{ j- 2 }  \le \cdot \le 2^{ j+2  } }  v  \right\|_{U_\Delta^2 \left(L_x^2, G_{2k}^j \right)}^2 + \left\|P^y_{\xi \left(G_{2k+1}^j \right), 2^{ j- 2 }  \le \cdot \le  2^{ j+2} }  v \right\|_{U_\Delta^2 \left(L_x^2, G_{2k+1}^j \right)}^2 \right)\notag\\
&\le \frac12 \left( \| v  \|_{X \left(G_{2k}^j \right)}^2 + \| v  \|_{X \left(G_{2k+1}^j \right)}^2 \right). \notag
\end{align}
At the same time, by \eqref{eq6.17v58} and \eqref{eq6.6v57}, we see
\begin{align}\label{eq5.15v51}
 & \ \sum\limits_{i \ge j + 1} \left\|P^y_{\xi(G_k^{j+1}), 2^{ i - 2 } \le \cdot \le 2^{  i + 2} }  v  \right\|_{U_\Delta^2(L_x^2, G_k^{j+1} )}^2\\
\le &  \sum\limits_{i \ge j+ 1} \left( \left\|P^y_{\xi(G_k^{j+1}), 2^{  i - 2 }  \le \cdot \le  2^{ i +2 } }  v  \right\|_{U_\Delta^2(L_x^2, G_{2k}^j)}^2
+ \left\|P^y_{\xi(G_k^{j+1}),  2^{ i - 2 }  \le \cdot \le  2^{ i+2 }  }  v  \right\|_{U_\Delta^2 (L_x^2, G_{2k+1}^j )}^2 \right)
\notag \\
\le &  \sum\limits_{i \ge j+ 1} \left( \left\|P^y_{ \xi(G_{2k}^j), 2^{ i - 3 }  \le \cdot \le  2^{ i +3} }   v  \right\|_{U_\Delta^2(L_x^2 , G_{2k}^j)}^2
+ \left\|P^y_{\xi(G_{2k+1}^j),  2^{ i - 3 }  \le \cdot \le  2^{ i +3 } }  v  \right\|_{U_\Delta^2(L_x^2, G_{2k+1}^j)}^2 \right).   \notag
\end{align}
Therefore, by \eqref{eq5.13v51} and \eqref{eq5.15v51}, and Definition \ref{de5.5v51}, we get
\begin{align*}
\| v  \|_{\tilde{X}_{k_\ast + 1}([0,T])}^2 \le 2 \|  v  \|_{ \tilde{X}_{k_\ast}([ 0, T])}^2 .
\end{align*}
By a similar argument, we can deduce
\begin{align*}
\| v  \|_{\tilde{Y}_{k_\ast  + 1} ([0,T])}^2 \le 2 \| v \|_{\tilde{Y}_{k_\ast([0,T])}}^2 .
\end{align*}
\end{proof}

\subsubsection{The bootstrap estimate}\label{subsec:bootstrap estimate}
In the following, we will establish the bootstrap estimate, which is necessary for the proof of Theorem \ref{th5.8v51}. For $0 \le j \le k_0$ and $G_k^j  \subseteq [0,T]$. By Duhamel's formula, we have for $0 \le i < j$,
\begin{align}\label{eq5.17v51}
\left\|P^y_{\xi \left(G_\alpha^i \right), 2^{ i - 2  }  \le \cdot \le  2^{ i+2 } } v  \right\|_{U_\Delta^2 \left(L_x^2, G_\alpha^i \right)}
&\le \left\|P^y_{\xi \left(G_\alpha^i \right),  2^{ i - 2 }  \le \cdot \le  2^{ i+2} }   v \left(t_\alpha^i \right) \right\|_{L_{y,x}^2 }\\
& + \left\| \int_{t_\alpha^i}^t e^{i(t-\tau)\Delta_y } P^y_{\xi \left(G_\alpha^i \right) ,  2^{ i - 2 } \le \cdot \le  2^{ i+2 } } F( v (\tau)) \,\mathrm{d}\tau \right\|_{U_\Delta^2 \left(L_x^2, G_\alpha^i \right)}\notag.
\end{align}
Here we take $t_\alpha^i$ to satisfy
\begin{align*}
\left\|P^y_{\xi \left(G_\alpha^i \right), 2^{  i - 2 }  \le\cdot \le  2^{ i+2} }   v \left(t_\alpha^i \right) \right\|_{L_{y,x}^2 }
= \inf\limits_{t \in G_\alpha^i} \left\|P^y_{\xi \left (G_\alpha^i \right),  2^{ i - 2  } \le \cdot \le 2^{  i+2} }   v (t) \right\|_{L_{y,x}^2}.
\end{align*}

We now consider the first term on the right hand side of \eqref{eq5.17v51}. By \eqref{eq5.3v51} and Lemma \ref{re5.2v51}, we have
\begin{align*}
& \sum\limits_{0 \le i < j} 2^{i- j} \sum\limits_{G_\alpha^i \subseteq G_k^j} \left\|P^y_{\xi \left(G_\alpha^i \right),  2^{ i - 2 } \le \cdot \le  2^{ i+2} }   v \left(t_\alpha^i \right) \right\|_{L_{y,x}^2 }^2\\
\lesssim &  2^{-j} \epsilon_3^{-1} \int_{G_k^j} \left(  N(t)^3 + \epsilon_3 \| v (t) \|_{L_y^4 L_x^2(\mathbb{R}^2 \times \mathbb{R})}^4 \right) \sum\limits_{ 0 \le i < j} \left\|P^y_{\xi(t),  2^{ i - 3 } \le \cdot \le  2^{ i+3} }   v (t) \right\|_{L_{y,x}^2 }^2 \,\mathrm{d}t \lesssim 1.
\end{align*}
For $i \ge j$, we can just take $t_k^j$ to be the left endpoint of $G_k^j$, then we have
\begin{align*}
\sum\limits_{i \ge j} \left\|P^y_{\xi \left(G_k^j \right),  2^{ i - 2 } \le \cdot \le 2^{  i + 2 } }  v \left(t_k^j  \right) \right\|_{L_{y,x}^2 }^2 \lesssim \left\| v \left(t_k^j  \right) \right\|_{L_{y,x}^2
}^2 \lesssim 1.
\end{align*}
Thus
\begin{align}\label{eq5.20v51}
\sum\limits_{0 \le i < j} 2^{i-j} \sum\limits_{G_\alpha^i \subseteq G_k^j}  \left\|P^y_{\xi \left(G_\alpha^i \right),  2^{ i - 2 }  \le \cdot \le  2^{ i + 2} }   v \left(t_\alpha^i \right) \right\|_{L_{y,x}^2 }^2
+ \sum\limits_{ i \ge j} \left\|P^y_{\xi \left(G_k^j \right),  2^{ i - 2 } \le \cdot \le  2^{ i + 2} }   v \left(t_k^j  \right) \right\|_{L_{y,x}^2 }^2 \lesssim 1.
\end{align}

We next consider the second term on the right hand side of \eqref{eq5.17v51}. Observe that there are at most two small intervals, called for instance $J_1$ and $J_2$, which intersect $G_k^j$ but are not contained in $G_k^j$. Then by Lemma \ref{le5.3v51} and \eqref{eq6.7v57}, we have
\begin{align}\label{eq5.23v51}
\sum\limits_{0 \le i < j} 2^{i-j} \sum\limits_{G_\alpha^i \subseteq G_k^j} \|F( v )\|_{L_t^1 L_{y,x}^2
( ( G_\alpha^i \cap (J_1 \cup J_2 ) ) \times \mathbb{R}^2 \times \mathbb{R} )}^2
& \lesssim \sum\limits_{0 \le i < j} 2^{i-j} \|F( v ) \|_{L_t^1 L_{y,x}^2
 ( (J_1 \cup J_2 )  \times \mathbb{R}^2 \times \mathbb{R} )}^2\\
& \lesssim \| v \|_{L_t^3 L_y^6 L_x^2(J_1 \times \mathbb{R}^2 \times \mathbb{R} )}^6 + \| v \|_{L_t^3 L_y^6 L_x^2 (J_2 \times \mathbb{R}^2 \times \mathbb{R} )}^6 \lesssim 1. \notag
\end{align}
Then by \eqref{eq6.6v57}, \eqref{eq6.7v57}, \eqref{eq5.23v51}, \eqref{eq5.1v51}, \eqref{eq5.4v51} and Definition \ref{de5.1v51}, we obtain
\begin{align}\label{eq6.28v58}
& \ \sum\limits_{0 \le i <  j} 2^{i - j} \sum\limits_{ \substack{ G_\alpha^i \subseteq G_k^j ,\\  N(G_\alpha^i) \ge 2^{i - 5} \epsilon_3^\frac12}}
\left\| \int_{t_\alpha^i}^t e^{i(t- \tau) \Delta_y } P^y_{\xi \left(G_\alpha^i \right), 2^{ i - 2 } \le \cdot \le  2^{ i+2} } F( v (\tau)) \,\mathrm{d}\tau  \right\|_{U_\Delta^2 \left(L_x^2, G_\alpha^i \right)}^2  \\
\lesssim &  \sum\limits_{0 \le i <  j} 2^{i-j} \sum\limits_{ \substack{ G_\alpha^i \subseteq G_k^j, \\N(G_\alpha^i) \ge 2^{i - 5} \epsilon_3^\frac12} }  \sum\limits_{J_l \cap G_k^j \ne \emptyset} \left\|P^y_{\xi \left(G_\alpha^i \right), 2^{ i- 2 } \le \cdot \le  2^{ i + 2  } } F( v ) \right\|_{DU_\Delta^2 \left(L_x^2, G_\alpha^i \cap J_l \right)}^2 \notag \\
\lesssim & \sum\limits_{0 \le i < j} 2^{i - j} \sum\limits_{ \substack{  G_\alpha^i \subseteq G_k^j, \\N \left(G_\alpha^i \right) \ge 2^{ i - 5} \epsilon_3^\frac12} }  \sum\limits_{J_l \cap G_k^j \ne \emptyset} \|F( v )\|_{L_t^1 L_{y,x}^2
 \left( \left(G_\alpha^i \cap J_l \right) \times \mathbb{R}^2 \times \mathbb{R} \right) }^2  \notag \\
\lesssim & 1 + \sum\limits_{0 \le i < j} 2^{i - j} \bigg( \sum\limits_{ \substack{ J_l \subseteq G_k^j, \\ N(J_l) \ge 2^{i-6} \epsilon_3^\frac12} }  \|F( v )\|_{L_t^1 L_{y,x}^2(J_l \times \mathbb{R}^2 \times \mathbb{R} )}^2 \bigg)
\lesssim   1 + \sum\limits_{J_l \subseteq G_k^j} \sum\limits_{ \substack{ 0 \le i <j,\\  2^i \le 2^6 \epsilon_3^{- \frac12} N(J_l)} } 2^{i-j}
\lesssim 1. \notag
\end{align}
On the interval $G_k^j$ with $N\left(G_k^j \right) \ge 2^{i -5} \epsilon_3^\frac12$, by \eqref{eq5.1v51} and \eqref{eq5.3v51}, we have
\begin{align}\label{eq6.35v122}
\int_{G_k^j} N(t)^2 \,\mathrm{d}t \lesssim 1.
\end{align}
Thus, by Minkowski's inequality, \eqref{eq2.6v51},  \eqref{eq5.4v51}, \eqref{eq6.6v122}, and \eqref{eq6.35v122}, we have
\begin{align}\label{eq6.29v58}
\sum\limits_{ \substack{ i \ge j , \\ N \left(G_k^j \right) \ge 2^{i  - 5} \epsilon_3^\frac12} } \left\|P^y_{\xi \left(G_k^j \right), 2^{ i -2 }  \le \cdot \le  2^{ i+2  } } F( v ) \right\|_{L_t^1 L_{y,x}^2 \left(G_k^j  \times \mathbb{R}^2 \times \mathbb{R} \right)}^2
\lesssim \|F( v ) \|_{L_t^1 L_{y,x}^2
\left(G_k^j  \times \mathbb{R}^2 \times \mathbb{R} \right)}^2 \lesssim \| v  \|_{L_t^3 L_y^6 L_x^2 \left(G_k^j \times \mathbb{R}^2 \times \mathbb{R} \right)}^6 \lesssim 1.
\end{align}
Thus, by \eqref{eq5.17v51}, \eqref{eq5.20v51}, \eqref{eq6.28v58}, and \eqref{eq6.29v58}, we infer
\begin{align}\label{eq5.26v51}
\|v  \|_{X \left(G_k^j \right)}^2  \lesssim  &
1 + \sum\limits_{ \substack{ i \ge j,\\  N \left(G_k^j \right) \le2^{i -5} \epsilon_3^\frac12} }
\left\|\int_{t_k^j}^t e^{i(t- \tau)\Delta_y } P^y_{\xi \left(G_k^j \right), 2^{  i - 2 }  \le \cdot \le  2^{ i+2 } } F( v (\tau)) \,\mathrm{d}\tau  \right\|_{U_\Delta^2 \left(L_x^2, G_k^j \right)}^2\\
& + \sum\limits_{0 \le i <  j} 2^{i-j} \sum\limits_{ \substack{ G_\alpha^i \subseteq G_k^j, \\ N \left(G_\alpha^i \right) \le 2^{i - 5} \epsilon_3^\frac12} }
\left\| \int_{t_\alpha^i}^t e^{i (t- \tau) \Delta_y } P^y_{\xi \left(G_\alpha^i \right),  2^{ i - 2  } \le \cdot \le  2^{ i + 2 } }F( v (\tau)) \,\mathrm{d}\tau \right\|_{U_\Delta^2 \left(L_x^2,
G_\alpha^i \right)}^2 .\notag
\end{align}
We can further get
\begin{align}\label{eq5.28v51}
\| v  \|_{X(G_k^j)}^2
\lesssim &  1 + \sum\limits_{ \substack{ i \ge j,\\  N \left(G_k^j \right) \le 2^{i - 10} \epsilon_3^\frac12}  } \left\|\int_{t_k^j}^t e^{i(t - \tau) \Delta_y } P^y_{\xi \left(G_k^j \right),  2^{ i - 2 } \le\cdot \le 2^{  i+2} }  F( v (\tau)) \,\mathrm{d}\tau \right\|_{U_\Delta^2 \left(L_x^2, G_k^j \right)}^2 \\
& + \sum\limits_{ 0 \le i <  j} 2^{i - j} \sum\limits_{ \substack{ G_\alpha^i \subseteq G_k^j,\\  N \left(G_\alpha^i \right) \le 2^{ i- 10} \epsilon_3^\frac12} }
\left\| \int_{t^i_\alpha}^t e^{i (t- \tau)\Delta_y } P^y_{\xi \left(G_\alpha^i \right),  2^{ i - 2 }  \le \cdot \le  2^{ i+2} }  F( v (\tau)) \,\mathrm{d}\tau \right\|_{U_\Delta^2 \left(L_x^2, G_\alpha^i \right)}^2, \notag
\end{align}
because the contribution of those terms for $i$ satisfying $ 2^{i  - 10} \epsilon_3^\frac12 \le N \left(t \right) \le 2^{i - 5} \epsilon_3^\frac12$ in the right hand side of \eqref{eq5.26v51} is small by similar argument as in the proof of \eqref{eq5.26v51}.

By a similar argument as above for  \eqref{eq5.28v51}, we also refer to \cite{D1} for more expatiation. Then, we have
\begin{align}\label{eq5.29v51}
\| v  \|_{Y \left(G_k^j \right)}^2
\lesssim & \epsilon_2^\frac32
+ \sum\limits_{{  \substack{ i \ge j,\\  N \left(G_k^j \right) \le2^{i - 10} \epsilon_3^\frac12}} }
\left\| \int_{t_k^j}^t e^{i(t- \tau)\Delta_y } P^y_{\xi \left(G_k^j \right), 2^{ i - 2 }  \le  \cdot \le  2^{ i + 2}  } F( v (\tau)) \,\mathrm{d}\tau \right\|_{U_\Delta^2 \left(L_x^2, G_k^j \right)}^2\\
& + \sum\limits_{0 \le i <  j} 2^{i-j} \sum\limits_{ \substack{ G_\alpha^i \subseteq G_k^j, \\ N \left(G_\alpha^i \right) \le 2^{i- 10} \epsilon_3^\frac12} }
\left\| \int_{t_\alpha^i}^t e^{i(t- \tau) \Delta_y } P^y_{\xi \left(G_\alpha^i \right),  2^{ i - 2 }  \le \cdot \le  2^{ i+2} }  F( v (\tau)) \,\mathrm{d}\tau \right\|_{U_\Delta^2 \left(L_x^2, G_\alpha^i \right)}^2. \notag
\end{align}

\begin{remark}\label{re5.11v51}
By Lemma \ref{le5.10v51}, Remark \ref{re5.4v51}, and \eqref{eq6.7v57}, we have
\begin{align*}
& \sum\limits_{ \substack{ i \ge j, \\
0 \le i \le 11,\\  N \left(G_k^j \right) \le 2^{i  - 5} \epsilon_3^\frac12} }   \left\| \int_{t_k^j}^t e^{i(t- \tau) \Delta_y } P^y_{\xi \left(G_k^j \right), 2^{ i - 2 } \le \cdot \le 2^{  i+2} }  F( v (\tau)) \,\mathrm{d}\tau \right\|_{U_\Delta^2 \left(L_x^2, G_k^j \right)}^2
\\
& + \sum\limits_{ 0 \le i \le 11} 2^{i - j} \sum\limits_{ \substack{ G_\alpha^i \subseteq G_k^j, \\ N \left(G_\alpha^i \right) \le 2^{i - 5} \epsilon_3^\frac12} }
 \left\| \int_{t_\alpha^i}^t e^{i(t - \tau) \Delta_y } P^y_{\xi \left(G_\alpha^i \right), 2^{ i- 2  } \le \cdot\le  2^{ i +2} }  F( v (\tau)) \,\mathrm{d}\tau \right\|_{U_\Delta^2 \left(L_x^2, G_\alpha^i \right)}^2 \lesssim 1.
\end{align*}
So we can further reduce the summation over $i$ on the right hand side of \eqref{eq5.28v51} and \eqref{eq5.29v51} to $i > 11$.
\end{remark}

Therefore, we have reduced to the proof of the following estimate.
\begin{theorem}[Reduced estimate]\label{th5.12v51}
\begin{align}\label{eq5.30v51}
& \sum\limits_{{  \substack{ i \ge j,\\ N\left(G_k^j \right) \le 2^{i - 10} \epsilon_3^\frac12} }}  \left\| \int_{t_k^j}^t e^{i(t - \tau) \Delta_y } P^y_{\xi \left(G_k^j \right), 2^{ i - 2 }  \le \cdot \le  2^{ i +2 } } F( v (\tau)) \,\mathrm{d}\tau \right\|_{U_\Delta^2 \left(L_x^2, G_k^j \right)}^2 \\
& + \sum\limits_{0 \le i <  j} 2^{i - j} \sum\limits_{ \substack{ G_\alpha^i \subseteq G_k^j,\\ N \left(G_\alpha^i \right) \le 2^{i - 10} \epsilon_3^\frac12}}  \left\| \int_{t_\alpha^i}^t e^{i (t- \tau) \Delta_y } P^y_{\xi \left(G_\alpha^i \right), 2^{  i - 2 }  \le \cdot \le  2^{ i + 2 } } F( v (\tau)) \,\mathrm{d}\tau \right\|_{U_\Delta^2 \left(L_x^2, G_\alpha^i \right)}^2\notag \\
\lesssim & \epsilon_2^\frac13 \| v \|_{\tilde{X}_j \left([0,T] \right)}^\frac53 \| v  \|_{\tilde{Y}_j \left([0,T] \right)}^2 + \epsilon_2^2 \| v  \|_{\tilde{Y}_j \left([0,T] \right)}^2 + \| v \|_{\tilde{Y}_j \left([0,T] \right)}^4 \left( 1 + \| v  \|_{\tilde{X}_j \left([0,T] \right)}^8  \right). \notag
\end{align}
\end{theorem}
Once this theorem is proved, we can close the proof of Theorem \ref{th5.8v51} by a bootstrap argument. In the proof given below, we shall assume Theorem \ref{th5.12v51} holds, while leaving its proof to Section \ref{subsec:theorem 1}.
\begin{proof}[Proof of Theorem \ref{th5.8v51}]
	Suppose
	\begin{align*}
		\| v  \|^2_{\tilde{X}_{k_\ast}([0,T])}\le C_0, \text{ and  } \| v \|_{\tilde{Y}_{k_\ast}([0,T])}^2 \le C \epsilon_2^\frac32 \le \epsilon_2,
	\end{align*}
	and from \eqref{eq6.22v59}, we have
	\begin{align*}
		\| v \|_{\tilde{X}_{k_\ast + 1}([0,T])}^2 \le2C_0, \text{ and }
		\| v  \|_{\tilde{Y}_{k_\ast + 1}([0,T])}^2 \le 2 \epsilon_2.
	\end{align*}
Then, by \eqref{eq5.28v51}, \eqref{eq5.29v51}, and \eqref{eq5.30v51}, we can further get
	\begin{align*}
		\|  v \|_{\tilde{X}_{k_\ast + 1} ([0,T]) } \le C \left( 1 + \epsilon_2^\frac23 \left(2 C_0 \right)^\frac56 + \epsilon_2^\frac32 +\epsilon_2 \left( 1 + 2C_0 \right)^8  \right), \intertext{ and }
		\| v  \|_{\tilde{Y}_{k_\ast + 1}([0,T])}
		\le C \left( \epsilon_2^\frac34 + \epsilon_2^\frac23 \left( 2 C_0 \right)^\frac56 + \epsilon_2^\frac32 + \epsilon_2 \left( 1 + 2 C_0 \right)^8  \right).
	\end{align*}
	If we choose $C_0 = 2^6 C$, and $\epsilon_2$ small enough, then we may deduce
	\begin{align*}
		\|v  \|_{\tilde{X}_{k_\ast + 1}([0,T])} \le C_0^{\frac{1}{2}},
		\text{ and  } \| v \|_{\tilde{Y}_{k_\ast + 1}([0,T]) } \le \epsilon_2^\frac12.
	\end{align*}
	Theorem \ref{th5.8v51} now follows from this and \eqref{eq6.21v59} by performing an induction on $k_\ast$.
\end{proof}

\subsubsection{The low frequency localized interaction Morawetz estimate}

As an application of the long time Strichartz estimate, we can obtain the low frequency localized interaction Morawetz estimate of the \eqref{eq1.3v28}. The Morawetz estimate is a very important tool to prove the scattering of the nonlinear dispersive equations for the radial case \cite{Lin-S, Morawetz}. In the non-radial case, J. Colliander, M. Keel, G. Staffilani, H. Takaoka, and T. Tao \cite{CKSTT1} develop the interaction Morawetz estimate, which is used to prove the scattering of the nonlinear Schr\"odinger equation \cite{CKSTT0,TVZ0,TVZ1,D1,D2,D3} in the non-radial case. The low frequency localized interaction Morawetz estimate will be used to preclude the soliton-like solution in Theorem \ref{th7.2v51}.
\begin{theorem}[Low frequency localized interaction Morawetz estimate] \label{th7.1v51}
Let  $v (t,y,x)$ be the almost periodic solution in Theorem \ref{th3.3v51} on $[0,T]$ with $\int_0^T N(t)^3 \,\mathrm{d}t = K$. Then we have
\begin{align}\label{eq6.97v77}
\left\| \int_{\mathbb{R}} |\nabla_y |^\frac12 \left( \Big| P^y_{\le  {10 \epsilon^{-1}_1  K}} v (t,y,x) \Big|^2 \right) \,\mathrm{d}x \right\|_{L_{t,y}^2([0,T] \times \mathbb{R}^2)}^2 \lesssim o(K).
\end{align}
\end{theorem}

The proof of this theorem follows from similar arguments in \cite{D1,D2,D3} and relies on Theorem \ref{th5.12v51} (and also some part of the proof).
In our \eqref{eq1.3v28} system, the interaction Morawetz quantity is
\begin{align*}
		M_0(t) = \int_{\mathbb{R}} \int_{\mathbb{R}} \int_{\mathbb{R}^2} \int_{\mathbb{R}^2} | { v } (t, \tilde{y}, \tilde{x})|^2 \frac{ y - \tilde{y}}{ | y - \tilde{y}|} \Im \left( \overline{  { v  } } \nabla_y    { v } \right)(t,y,x) \,\mathrm{d} y \mathrm{d} \tilde{y}  \mathrm{d} x \mathrm{d} \tilde{x},
	\end{align*}
which is invariant under the Galilean transform in the $\mathbb{R}^2$ component.
Following the argument in \cite{CGT,PV}, we can get
\begin{align*}
\left\| \int_{\mathbb{R}} |\nabla_y|^\frac12 \left( |v(t,y,x)|^2 \right) \,\mathrm{d}x \right\|_{L_{t,y}^2([0, T] \times \mathbb{R}^2)}^2 \lesssim \left|  M_0(T) - M_0(0) \right| .
\end{align*}
Replacing $v$ by its low frequency cut-off $P^y_{\le  {10 \epsilon^{-1}_1  K}} v$, we then get the low frequency localized interaction Morawetz quantity
\begin{align*}
		M (t) = \int_{\mathbb{R}} \int_{\mathbb{R}} \int_{\mathbb{R}^2} \int_{\mathbb{R}^2} \left| {P^y_{\le  {10 \epsilon^{-1}_1  K}} v } (t, \tilde{y}, \tilde{x}) \right|^2 \frac{ y - \tilde{y}}{ \left| y - \tilde{y} \right|} \Im \left( \overline{  { P^y_{\le  {10 \epsilon^{-1}_1  K}} v  } } \nabla_y    { P^y_{\le  {10 \epsilon^{-1}_1  K}} v } \right)(t,y,x) \,\mathrm{d} y \mathrm{d} \tilde{y}  \mathrm{d} x \mathrm{d} \tilde{x}.
	\end{align*}
Because for any $\eta > 0$ independent of $\epsilon_1$,
by Theorem \ref{th3.3v51} and Bernstein's inequality, we have
\begin{align*}
|M(T)|+ |M(0)| \lesssim \eta K,
\end{align*}
we then obtain
\begin{align*}
\left\| \int_{\mathbb{R}} |\nabla_y|^\frac12 \left( |v(t,y,x)|^2 \right) \,\mathrm{d}x \right\|_{L_{t,y}^2([0, T] \times \mathbb{R}^2)}^2 \lesssim \eta K + \mathscr{E},
\end{align*}
where $\mathscr{E}$ are the error terms coming from the low frequency cut-off of the solution of the \eqref{eq1.3v28}. These error terms can be proven to be $o(K)$, using Theorem \ref{th5.12v51} and also some estimates from the proof of it. We shall leave the detailed proof of this theorem to Section \ref{subsec:theorem 2}.

\subsection{Exclusion of the almost periodic solution}\label{subsec:6.3}
\begin{theorem}\label{th7.2v51}
The almost periodic solution to \eqref{eq1.3v28} in Theorem \ref{th3.3v51} does not exist.
\end{theorem}
\begin{proof}
We will preclude two scenarios in the following.
\medskip

{ \textit{Case I. } $\int_0^\infty N(t)^3 \,\mathrm{d}t < \infty.$ }
\medskip

By the proof of Theorem \ref{th5.8v51}, as in \cite{D1,D2}, we have
\begin{align}\label{eq6.42v122}
\| v(t,y,x) \|_{L_t^\infty \dot{H}_y^3 L_x^2([ 0, \infty) \times \mathbb{R}^2 \times \mathbb{R} )}
\lesssim_{m_0} \left( \int_0^\infty N(t)^3 \,\mathrm{d}t \right)^3.
\end{align}
By \eqref{eq6.42v122} and \eqref{eq6.3v55}, we have
\begin{align*}
\left\|e^{i y \xi(t)} v \right\|_{\dot{H}_y^1 L_x^2} \lesssim N(t) C(\eta(t)) + \eta(t)^\frac12 \to 0, \text{ as } t \to \infty.
\end{align*}
 Thus, for any $ \epsilon > 0$, we can take sufficiently large positive constant $t_0$ such that $\left\| e^{iy \xi(t_0)} v(t_0) \right\|_{\dot{H}_y^1 L_x^2} \le \epsilon$. In the following, We can assume $t_0 = 0$ because of the Galilean invariance. By Minkowski's inequality, the Gargliardo-Nirenberg inequality, and H\"older's inequality, we have
\begin{align*}
\mathcal{E}( v (t)) = \mathcal{E}( v (0)) \lesssim \| v (0) \|_{\dot{H}_y^1 L_x^2}^2 \lesssim \epsilon^2.
\end{align*}
Because we can take $\epsilon$ as small as we wish, this scenario does not exist.
\medskip

{ \textit{Case II. }
$\int_0^\infty N(t)^3 \,\mathrm{d} t = \infty$.  }
\medskip

By H\"older's inequality and Sobolev's inequality, we have
\begin{align*}
& \int_{|y - y(t)|\le \frac{C\Big(\frac{ \| v(0)  \|_{L_{y,x}^2}^2}{100} \Big)}{N(t)}} \int_{\mathbb{R}} \left|P^y_{\le 10 \epsilon_1^{-1} K}  v (t,y,x) \right|^2 \,\mathrm{d}y \mathrm{d}x
\lesssim \bigg( \frac{C \Big(\frac{ \| v(0)  \|_{L_{y,x}^2}^2}{100} \Big)}{ N(t)}\bigg)^\frac32 \left\| \int_{\mathbb{R}} |\nabla_y|^\frac12 \left( \left|P^y_{\le 10 \epsilon_1^{-1} K}  v (t,y,x) \right|^2 \right) \,\mathrm{d} x \right\|_{L_y^2}.
\end{align*}
By Theorem \ref{th3.3v51}, we have for $K \ge C \Big( \frac{ \| v \|_{L_{y,x}^2}^2}{100} \Big)$,
\begin{align*}
\frac{\| v \|_{L_{y,x}^2}^2}2 \le \int_{\mathbb{R}} \int_{|y - y(t)| \le \frac{ C\Big(\frac{ \| v \|_{L_{y,x}^2}^2}{100} \Big)}{ N(t)} } \left|P^y_{\le 10 \epsilon_1^{-1} K} v (t,y,x) \right|^2 \,\mathrm{d}y \mathrm{d}x.
\end{align*}
By the above two estimates, together with Theorem \ref{th7.1v51} and the conservation of mass, we have the following contradiction when $K$ is sufficiently large,
\begin{align*}
\|  v \|_{L_{y,x}^2}^4 K  \lesssim &   \int_0^T N(t)^3 \bigg( \int_{| y - y(t)| \le \frac{ C \Big( \frac{\| v \|_{L_{y,x}^2}^2} {100 } \Big) } {N(t)} } \int_{\mathbb{R}} \left|P^y_{\le 10 \epsilon_1^{-1} K} v (t,y,x) \right|^2 \,\mathrm{d}x \,\mathrm{d} y  \bigg)^2 \,\mathrm{d}t\\
\lesssim &  \left\| \int_{\mathbb{R}} |\nabla_y|^\frac12 \left( \left|P^y_{\le 10 \epsilon_1^{-1} K}  v (t,y,x) \right|^2 \right) \,\mathrm{d} x \right\|_{L_{t,y}^2([0,T] \times \mathbb{R}^2 )}^2 \lesssim o(K).
\end{align*}
This completes the proof of Theorem \ref{th7.2v51}.
\end{proof}
\begin{proof}[Proof of Theorem \ref{as3.14v12}]
	This is an immediate consequence of Theorem \ref{th3.3v51} and Theorem \ref{th7.2v51}.
\end{proof}

\section{Proof of Theorems \ref{th5.12v51} and \ref{th7.1v51}}\label{sec:two theorems}

\subsection{Proof of Theorem \ref{th5.12v51}}\label{subsec:theorem 1}
In this section, we complete the proof of  Theorem \ref{th5.12v51}. To prove this theorem, we decompose the nonlinear term $P^y_{\xi(G_\alpha^i), 2^{ i - 2 }  \le \cdot \le  2^{ i+2} }  F(v) $ and also use the fact that on the time interval $G_\alpha^i$, $\xi(t)$ can replace $\xi(G_\alpha^i)$ up to $2^{i - 20} $ by \eqref{eq6.20v124}. Then, we can see it is enough to prove the estimate the left hand side of \eqref{eq5.30v51} with $P^y_{\xi(G_\alpha^i), 2^{ i - 2 }  \le \cdot \le  2^{ i+2} }  F(v) $ being replaced by
\begin{align}\label{eq6.41v65}
  & P^y_{\xi(G_\alpha^i),  2^{ i - 2  } \le \cdot \le  2^{ i+2} }  \mathcal{O} \bigg( \sum\limits_{ \substack{ n_1,n_2,n_3,n\in \mathbb{N}, \\n_1- n_2 + n_3 = n} }
	\Pi_n \left(  P^y_{\xi(G_\alpha^i),  \ge 2^{ i - 5} }  v_{n_1}  \overline{ P^y_{\xi(t ), \ge 2^{ i -10 } }  v_{n_2} }  v_{n_3}  \right) \bigg)  \\
+ &
P^y_{\xi( t ), 2^{  i - 2  } \le \cdot \le  2^{ i+2} }  \mathcal{O} \bigg(\sum\limits_{ \substack{ n_1,n_2,n_3,n \in \mathbb{N}, \\ n_1 - n_2 + n_3 = n}  }
	\Pi_n \left( P^y_{\xi( t ), \le 2^{  i - 10 } }  v_{n_1}
	\overline{ P^y_{\xi(G_\alpha^i), \le 2^{ i - 10}  }   v_{n_2} }  P^y_{\xi(G_\alpha^i),  2^{ i - 5 }  \le \cdot \le 2^{  i +5} }  v_{n_3}  \right) \bigg),  \label{eq6.40v65},
\end{align}
we also has similar fact of the nonlinear term $P^y_{\xi(G_k^j), 2^{ i - 2 }  \le \cdot \le  2^{ i+2} }  F(v) $, where the symbol $``\mathcal{O}"$ represents the different frequencies will be located in different $v_{n_l}$s', $l= 1,2,3$. Since their estimates are almost identical, we denote them as a single  $``\mathcal{O}"$. The estimate of the Duhamel propagator of the term \eqref{eq6.41v65} is very short and easy, and mainly relies on the bilinear Strichartz estimate in Proposition \ref{pr4.7v51}.
The estimate of the Duhamel propagator of the term \eqref{eq6.40v65} is lengthy. This is because to prove the estimate of the Duhamel propagator of the term \eqref{eq6.40v65}, we need to prove the bilinear Strichartz estimates on the union of the small intervals. It turns out the proof of these bilinear Strichartz estimates cannot be proven just by the harmonic analysis but also rely heavily on the structure of the \eqref{eq1.3v28} system or more precisely the corresponding interaction Morawetz estimate of \eqref{eq1.3v28}. During the proof of this part, some terms can be estimated by the following  bilinear Strichartz estimate established recently \cite{Candy} instead of the interaction Morawetz estimate as in \cite{D1}. This new bilinear Strichartz estimate is very useful in \cite{SW}.
\begin{lemma}[Bilinear Strichartz estimate, \cite{Candy}] \label{le2.2v118}
Let $1 \le q, r \le 2$, $\frac1q + \frac3{2r} < \frac32$, and suppose $M, N \in 2^{\mathbb{Z}}$ satisfy $M<<N$, then for any $\phi, \psi \in L^2(\mathbb{R}^2)$, \begin{align*}
\left\|e^{it \Delta}P_N \phi e^{i t \Delta} P_M \psi \right\|_{L_t^q L_x^r(\mathbb{R} \times \mathbb{R}^2) } \lesssim \frac{M^{3- \frac2q - \frac3r}}{N^{1- \frac1r}} \left\|P_N \phi \right\|_{L^2( \mathbb{R}^2) } \left\|P_M \psi \right\|_{L^2 ( \mathbb{R}^2) }.
\end{align*}
\end{lemma}

\subsubsection{Estimate of \eqref{eq6.41v65}}
We first deal with \eqref{eq6.41v65}.
\begin{theorem}\label{th5.13v51}
	For any fixed $G_k^j \subseteq [0,T]$, $j> 0$, we have
	\begin{align*}
 &   \sum\limits_{0 \le i <  j} 2^{i - j} \sum\limits_{ \substack{ G_\alpha^i \subseteq G_k^j, \\ N \left(G_\alpha^i \right) \le 2^{i-10} \epsilon_3^\frac12} }
		\Bigg\| \int_{t_\alpha^i}^t e^{i(t- \tau) \Delta_y } P^y_{\xi \left(G_\alpha^i \right),  2^{ i - 2 }  \le \cdot \le 2^{  i+2} } F^{high}_{\alpha, i}( v( \tau ))
 \,\mathrm{d}\tau \Bigg\|_{U_\Delta^2 \left(L_x^2, G_\alpha^i \right)}^2 \\
		& +   \sum\limits_{ \substack{ i \ge j, \\ N \left(G_k^j \right) \le 2^{i- 10} \epsilon_3^\frac12} }
		\Bigg\|\int_{t_k^j}^t e^{i(t-\tau) \Delta_y } P^y_{\xi \left(G_k^j \right), 2^{ i - 2  } \le \cdot \le  2^{ i+2 } }
F^{high}_{k, j }( v( \tau )  ) \,\mathrm{d}\tau  \Bigg\|_{U_\Delta^2 \left(L_x^2, G_k^j \right)}^2	\lesssim   \epsilon_2^\frac13 \|v \|_{\tilde{X}_j([0,T])}^\frac53 \|v \|_{\tilde{Y}_j([0,T])}^2,
	\end{align*}	
where
\begin{align*}
F^{high}_{k , j }(v(t) )  &  =   \sum\limits_{ \substack{ n_1,n_2,n_3,n \in \mathbb{N}, \\ n_1 - n_2 + n_3 = n} }
		\Pi_n \left( P^y_{\xi \left(G_k^j  \right),  \ge 2^{ i - 5} }  v_{n_1} \overline{  P^y_{\xi( t ), \ge  2^{ i - 10} }  v_{n_2} }  v_{n_3} \right)(t),
\intertext{ and }
F^{high}_{\alpha, i}(v(t) )  &  =   \sum\limits_{ \substack{ n_1,n_2,n_3,n \in \mathbb{N}, \\ n_1 - n_2 + n_3 = n} }
		\Pi_n \left( P^y_{\xi \left(G_\alpha^i \right),  \ge 2^{ i - 5} }  v_{n_1} \overline{  P^y_{\xi( t ), \ge  2^{ i - 10} }  v_{n_2} }  v_{n_3} \right)(t).
\end{align*}

\end{theorem}
\begin{proof}
On the time interval $G_\alpha^i$ with $N\left(G_\alpha^i \right) \le 2^{i - 10} \epsilon_3^\frac12$, we take $w \in V_\Delta^2 \left(L_x^2, G_\alpha^i \right) $ be normalized {so that} $ \left( \mathcal{F}_y w  \right)(t, \xi, x) $ is supported on
	$$\left\{ \xi: 2^{i - 2} \le \left|\xi - \xi \left(G_\alpha^i \right) \right| \le 2^{i+2} \right\}$$
	for any $(t, x) \in \mathbb{R} \times \mathbb{R}$. By the Cauchy-Schwarz inequality, \eqref{eq3.8v41}, Proposition \ref{pr4.7v51}, the conservation of mass, \eqref{eq6.4v57}, \eqref{eq6.7v57}, Lemma \ref{le5.7v51}, and \eqref{eq5.2v51}, we infer
	\begin{align}\label{eq5.35v51}
		& \int_{G_\alpha^i} \left\langle w (  \tau ) , F_{\alpha,i}^{high} (v(\tau) )  \right\rangle   \,\mathrm{d}\tau\\
		\le & \sum\limits_{l \ge i - 5} \left\| \| w \|_{L_x^2}
		\left\|P^y_{\xi \left(G_\alpha^i \right), 2^l } v \right\|_{L_x^2} \right \|_{L_t^\frac52 L_y^\frac53 \left(G_\alpha^i \times \mathbb{R}^2 \right)}^\frac12
		\left\|P^y_{\xi \left(G_\alpha^i \right), 2^l } v \right\|_{L_t^\frac52 L_y^{10} L_x^2 \left(G_\alpha^i \times \mathbb{R}^2 \times \mathbb{R} \right)}^\frac12
		\| w  \|_{L_t^\frac52 L_y^{10} L_x^2(G_\alpha^i \times \mathbb{R}^2 \times \mathbb{R})}^\frac12  \notag\\
		& \quad \left\|P^y_{\xi(t), \ge 2^{ i - 10} }  v \right\|_{L_t^\frac52 L_y^{10} L_x^2 \left(G_\alpha^i \times \mathbb{R}^2 \times \mathbb{R} \right)}
		\|v \|_{L_t^\infty L_{y,x}^2
			( G_\alpha^i \times \mathbb{R}^2 \times \mathbb{R})} \notag \\
		\lesssim & \sum\limits_{l \ge i - 5} 2^{\frac{i - l}5} \left\|P^y_{\xi(G_\alpha^i), 2^l } v \right\|_{U_\Delta^2 \left(L_x^2, G_\alpha^i \right)}
		\left\| P^y_{\xi(t), \ge  2^{ i - 10} }  v \right\|_{L_t^\infty L_{y,x}^2
			\left(G_\alpha^i \times \mathbb{R}^2 \times \mathbb{R} \right)}^\frac16 \left\|P^y_{\xi(t), \ge  2^{ i - 10}  } v \right\|_{L_t^\frac{25}{12} L_y^{50} L_x^2 \left(G_\alpha^i \times \mathbb{R}^2 \times \mathbb{R} \right)}^\frac56  \notag \\
		\lesssim & \epsilon_2^\frac16 \|v \|_{\tilde{X}_i([0, T])}^\frac56 \sum\limits_{l \ge i - 5} 2^{ \frac{i - l}5} \left\|P^y_{\xi \left(G_\alpha^i \right), 2^l} v \right\|_{U_\Delta^2 \left(L^2_x, G_\alpha^i \right)}
		\lesssim  \epsilon_2^\frac16 \| v \|_{\tilde{X}_j([0,T])}^\frac56 \left( \sum\limits_{l \ge i - 5} 2^{\frac{i - l}5} \left\|P^y_{\xi \left(G_\alpha^i \right), 2^l} v \right\|_{U_\Delta^2 \left(L^2_x, G_\alpha^i \right)}^2 \right)^\frac12. \notag
	\end{align}
As in \cite{D1}, we see for any $0 \le l \le j$, $G_k^j$ overlaps $2^{j-l}$ intervals $G_\beta^l$ and for $0 \le i \le l$, $G_\beta^l$ overlaps $2^{l-i}$ intervals $G_\alpha^i$. In addition, every $G_\alpha^i$ is contained in one $G_\beta^l$. Thus, we can divide the summation in the left hand side of the following \eqref{eq5.36v51} and \eqref{eq7.5v122} into different groups according to $l \ge j$ and $0 \le l < j$. Then by some easy calculation and reordering the summation of $i$ and $l$, we have
\begin{align}\label{eq5.36v51}
		& \sum\limits_{0 \le i <  j} 2^{i - j} \sum\limits_{ \substack{ G_\alpha^i \subseteq G_k^j,\\ N \left(G_\alpha^i \right) \le 2^{i - 10} \epsilon_3^\frac12} }
		\left( \sum\limits_{l \ge i - 5} 2^{\frac{i-l}5} \left\|P^y_{\xi \left(G_\alpha^i \right),2^l} v  \right\|_{U_\Delta^2 \left(L_x^2, G_\alpha^i \right)}^2  \right)
\intertext{ and }
		& \sum\limits_{ \substack{ i \ge j, \\
				N \left(G_k^j \right) \le 2^{i-10} \epsilon_3^\frac12}} \left( \sum\limits_{l \ge i- 5} 2^{\frac{i-l}5} \left\|P^y_{\xi \left(G_k^j \right), 2^l}  v  \right\|_{U_\Delta^2 \left(L_x^2, G_k^j \right)}^2 \right) \lesssim \| v \|_{\tilde{Y}_j([0,T])}^2. \label{eq7.5v122}
	\end{align}
	Theorem \ref{th5.13v51} follows from \eqref{eq6.5v57}, \eqref{eq5.35v51} and \eqref{eq5.36v51}.
\end{proof}

\subsubsection{Estimate of \eqref{eq6.40v65}}
Now we turn to the estimate of \eqref{eq6.40v65}. Denote
\begin{align*}
F^{low}_{k, j }(v(t))  & :=   \sum\limits_{ \substack{ n_1,n_2,n_3,n \in \mathbb{N}, \\n_1-n_2 + n_3 = n } } \Pi_n \left( P^y_{\xi(t ), \le 2^{ i - 10 } } v_{n_1} \overline{ P^y_{\xi( t ), \le  2^{ i -10}  } v_{n_2} } P^y_{\xi \left(G_k^j \right),  2^{ i - 5}  \le \cdot \le 2^{ i+5} } v_{n_3}  \right) ,
\intertext{ and }
F^{low}_{\alpha, i}(v(t))  & :=   \sum\limits_{ \substack{ n_1,n_2,n_3,n \in \mathbb{N}, \\n_1-n_2 + n_3 = n } } \Pi_n \left( P^y_{\xi(t ), \le 2^{ i - 10 } } v_{n_1} \overline{ P^y_{\xi( t ), \le  2^{ i -10}  } v_{n_2} } P^y_{\xi \left(G_\alpha^i \right),  2^{ i - 5}  \le \cdot \le 2^{ i+5} } v_{n_3}  \right) .
\end{align*}
Then, we have
\begin{theorem}\label{th5.14v51}
	For any $0 \le i \le j$, on the time interval $G_\alpha^i \subseteq G_k^j$ with $N \left(G_\alpha^i \right) \le 2^{i - 10} \epsilon_3^\frac12$, we have
	\begin{align}\label{eq5.38v51}
		& \left\| \int_{t^i_\alpha}^t e^{i(t- \tau)\Delta_y } P^y_{\xi \left(G_\alpha^i \right),  2^{ i - 2 }  \le \cdot \le  2^{ i+2} }  F^{low}_{\alpha, i}(v(\tau)) \,\mathrm{d}\tau \right\|_{U_\Delta^2 \left(L_x^2, G_\alpha^i \right)} \\
		\lesssim  &
		\left\|P^y_{\xi \left(G_\alpha^i \right),  2^{ i- 5 }  \le \cdot \le  2^{ i+5} }   v  \right\|_{U_\Delta^2 \left(L_x^2, G_\alpha^i \right)} \left( \epsilon_2 + \| v \|_{\tilde{Y}_i([0,T]) } \left( 1 + \| v \|_{\tilde{X}_i([0,T])} \right)^4  \right).\notag
	\end{align}
In addition, for $i \ge j$, $N(G_k^j) \le 2^{i-10} \epsilon_3^\frac12$, we have
	\begin{align}\label{eq5.39v51}
& \left\| \int_{t_k^j}^t e^{i(t- \tau) \Delta_y } P^y_{\xi \left(G_k^j \right), 2^{ i - 2 } \le \cdot \le  2^{ i+2} } F_{k,j }^{low}(v(\tau)) \,\mathrm{d}\tau \right\|_{U_\Delta^2 \left(L_x^2, G_k^j \right)}	\\
		\lesssim & 2^{\frac{3(j-i)}4} \left\|P^y_{\xi \left(G_k^j \right), 2^{  i- 5 } \le \cdot \le 2^{ i+5} } v \right\|_{U_\Delta^2 \left(L_x^2, G_k^j \right)} \left( \epsilon_2+ \| v \|_{\tilde{Y}_j([0,T])} \left( 1 + \| v \|_{\tilde{X}_j([0,T])} \right)^4 \right).  \notag
	\end{align}	
\end{theorem}
\begin{proof}[Proof of Theorem \ref{th5.14v51}.]
We will only prove \eqref{eq5.38v51}, as \eqref{eq5.39v51} follows by a similar argument.
Fix $G_\alpha^i$ with $N \left(G_\alpha^i \right) \le 2^{i - 10} \epsilon_3^\frac12$. We can see there are no more than two small intervals $J_1$ and $J_2$ which overlap $G_\alpha^i$ but are not contained in $G_\alpha^i$.
Let $\tilde{G}_\alpha^i= G_\alpha^i \backslash (J_1 \cup J_2)$, by \eqref{eq6.4v57}, \eqref{eq6.5v57}, \eqref{eq6.6v57}, and \eqref{eq6.7v57}, we have
\begin{align}
	& \left\| \int_{t_\alpha^i}^t e^{i(t- \tau )\Delta_y } P^y_{\xi \left(G_\alpha^i \right),  2^{ i - 2 } \le \cdot \le  2^{ i+2} }
F^{low}_{\alpha,i}(v(\tau)) \,\mathrm{d}\tau \right\|_{U_\Delta^2 \left(L_x^2, G_\alpha^i \right)}  \label{eq7.8v123} \\
	\lesssim & \left\|\int_{t_\alpha^i}^t e^{i(t- \tau) \Delta_y } P^y_{\xi \left(G_\alpha^i \right), 2^{  i - 2  } \le \cdot \le 2^{  i+2}  }
F^{low}_{\alpha,i}(v(\tau)) \,\mathrm{d}\tau \right\|_{U_\Delta^2 \left(L_x^2, \tilde{G}_\alpha^i \right)} \notag \\
& + \left\| F^{low}_{\alpha,i}(v(t)) \right\|_{L_{t,y}^\frac43 L_x^2 \left( \left( J_1 \cap G_\alpha^i \right) \times \mathbb{R}^2 \times \mathbb{R} \right) }
 + \left\|F^{low}_{\alpha,i}(v(t))\right\|_{L_{t,y}^\frac43 L_x^2 \left(  \left( J_2 \cap G_\alpha^i \right) \times \mathbb{R}^2 \times \mathbb{R} \right)}. \notag
\end{align}
Here, we may assume $t_\alpha^i \in \tilde{G}_\alpha^i$, because if $t_\alpha^i \not\in \tilde{G}_\alpha^i$, we may move $t_\alpha^i$ into $\tilde{G}_\alpha^i$ with the errors being absorbed by the last two terms on the right hand side of the above inequality. We can show the last two terms on the right hand side of \eqref{eq7.8v123} is small in the following.

On the intervals $J_l$ for $l = 1,2$, by Propositions \ref{pr4.5v51} and \ref{pr4.7v51}, \eqref{eq6.7v57}, the fact $N(t) \le 2^{i-5} \epsilon_3^\frac12$ on $G_\alpha^i$, \eqref{eq5.1v51}, \eqref{eq5.2v51}, Lemma \ref{le5.3v51}, and \eqref{eq6.7v57}, we can get
\begin{align}\label{eq5.43v51}
	& \left\| F^{low}_{\alpha,i}(v(t))\right\|_{L_{t,y}^\frac43 L_x^2 \left(  \left( J_l \cap G_\alpha^i \right) \times \mathbb{R}^2 \times \mathbb{R} \right)} \\
	\lesssim &  \bigg\|  \Big\| P^y_{\xi(J_l), \le \epsilon_3^{- \frac14} N(J_l)} v \Big \|_{L_x^2} \left\|P^y_{\xi \left(G_\alpha^i \right),  2^{ i - 5 }  \le \cdot \le 2^{ i+5 } } v \right\|_{L_x^2} \bigg\|_{L_{t,y}^2  \left(  \left( J_l \cap G_\alpha^i \right) \times \mathbb{R}^2 \times \mathbb{R} \right)} \left\|P^y_{\xi( t ), \le 2^{ i - 10} }  v \right\|_{L_{t,y}^4 L_x^2 \left(  \left(  J_l \cap G_\alpha^i \right)  \times \mathbb{R}^2 \times \mathbb{R} \right)  }\notag \\
	+ & \left\|P^y_{\xi(J_l), \ge \epsilon_3^{- \frac14} N(J_l)} v \right\|_{L_t^\infty L_{y,x}^2
		\left(  \left( J_l \cap G_\alpha^i \right)  \times \mathbb{R}^2 \times \mathbb{R} \right) } \left\|P^y_{\xi \left(G_\alpha^i \right), 2^{ i - 5  } \le \cdot \le 2^{  i+5  } } v \right\|_{L_t^\frac83 L_y^8 L_x^2 \left(  \left(J_l \cap G_\alpha^i \right) \times \mathbb{R}^2 \times \mathbb{R} \right) } \left\|P^y_{\xi(t ), \le 2^{ i - 10} }  v \right\|_{L_t^\frac83 L_y^8 L_x^2 \left(  \left(J_l \cap G_\alpha^i  \right) \times \mathbb{R}^2 \times \mathbb{R} \right) }\notag  \\
	\lesssim & \sum\limits_{2^k \le \epsilon_3^{- \frac14} N(J_l)}2^{\frac{k-i}2} \left\|P^y_{\xi \left(G_\alpha^i \right),  2^{ i - 5 }  \le \cdot \le 2^{  i+5} }  v\right\|_{U_\Delta^2 \left(L_x^2, G_\alpha^i \right)}
	\left\|P^y_{\xi(J_l),2^k} v \right\|_{U_\Delta^2 \left(L_x^2, J_l \right)}
	+ \epsilon_2 \left\|P^y_{\xi \left(G_\alpha^i \right),  2^{ i - 5 } \le \cdot \le  2^{ i+5}  } v \right\|_{U_\Delta^2 \left(L_x^2, G_\alpha^i \right)} \notag \\
	\lesssim  & \epsilon_2 \left\|P^y_{\xi \left(G_\alpha^i \right),  2^{ i - 5 } \le \cdot \le  2^{ i+5} }  v \right\|_{U_\Delta^2 \left(L_x^2, G_\alpha^i \right)}. \notag
\end{align}
Thus, we can simplify the estimate of \eqref{eq7.8v123} to the case that $G_\alpha^i$ is the union of finite many small intervals $J_l$'s. (If not, we just need to add the right hand side of \eqref{eq5.43v51}). Let
\begin{align*}
F^{low, l_2}_{\alpha, i}(v(t))  =  \sum\limits_{ \substack{ n_1,n_2,n_3,n \in \mathbb{N},\\
			n_1 -n_2 + n_3 = n}  }
	\Pi_n \Big( P^y_{\xi(t), 2^{ l_2} }  v_{n_1}    \overline{ P^y_{\xi(t), \le 2^{ l_2 } }  v_{n_2}  }
	P^y_{\xi \left(G_\alpha^i \right),  2^{ i -5 }  \le \cdot \le 2^{ i+5 } } v_{n_3} \Big),
\end{align*}
then by Lemma \ref{le4.4v51}, we have
\begin{align*}
\text{ LHS of } \eqref{eq7.8v123}  \lesssim A_1 + A_2 + A_3 + A_4,
\end{align*}
where
\begin{align}
A_1 =  &  \sum\limits_{0 \le l_2 \le i- 10} \Bigg( \sum\limits_{ \substack{ J_l \subseteq G_\alpha^i,\\ N(J_l) \ge \epsilon_3^\frac12 2^{l_2- 5}} }
	\left\|P^y_{\xi \left(G_\alpha^i \right), 2^{ i- 2 } \le \cdot \le 2^{  i+2 } }
F^{low, l_2}_{\alpha, i}(v(t))\right\|_{DU_\Delta^2 \left(L_x^2,J_l \right)}^2 \Bigg)^\frac12 , \label{eq5.45v51}
\end{align}
\begin{align}
A_2 =   & \sum\limits_{ \substack{ 0 \le l_2 \le i - 10,\\ J_l \subseteq G_\alpha^i, \\
			N(J_l) \ge \epsilon_3^\frac12 2^{l_2- 5}  } }  \left\| \int_{J_l} e^{-it \Delta_y } P^y_{\xi \left(G_\alpha^i \right), 2^{ i - 2}  \le \cdot \le 2^{ i+2 } } F^{low, l_2}_{\alpha, i}(v(t))\,\mathrm{d} t \right\|_{L_{y,x}^2  }, \label{eq5.46v51}
\end{align}
\begin{align}%\\
 A_3 = &  \sum\limits_{0 \le l_2 \le i - 10} \bigg( \sum\limits_{ \substack{ G_\beta^{l_2} \subseteq G_\alpha^i,\\ N \left(G_\beta^{l_2} \right) \le \epsilon_3^\frac12 2^{l_2 - 5}} }
	\left\|P^y_{\xi \left(G_\alpha^i \right),   2^{ i - 2 } \le \cdot \le  2^{ i+2} }   F^{low,l_2}_{\alpha, i}(v(t)) \right\|_{DU_\Delta^2 \left(L_x^2, G_\beta^{l_2} \right)}^2 \bigg)^\frac12,   \label{eq5.47v51}
\end{align}
and
\begin{align}
A_4 =  &  \sum\limits_{ \substack{0 \le l_2 \le i - 10,\\ G_\beta^{l_2} \subseteq G_\alpha^i,\\ N \left(G_\beta^{l_2} \right) \le \epsilon_3^\frac12 2^{l_2 -5}} }\left\| \int_{G_\beta^{l_2}} e^{- it \Delta_y } P^y_{\xi \left(G_\alpha^i \right),  2^{ i - 2 }  \le \cdot \le  2^{ i+2 } }
F^{low,l_2}_{\alpha, i}(v(t))\,\mathrm{d} t \right\|_{L_{y,x}^2 }.   \label{eq5.48v51}
\end{align}
The proof of the first two terms are easy. We first prove the following auxiliary estimate.
\begin{lemma}\label{le6.28v118}
Let $(p_0, q_0)$ be Strichartz admissible with $  q_0 \ge 20  $. Suppose that $w(t,y,x) \in L_t^{p_0} L_y^{q_0}  L_x^2 \left(J_l \times \mathbb{R}^2 \times \mathbb{R} \right)$ satisfies that $\mathcal F w(t, \cdot , x)$ is supported on $ \left\{ \xi: 2^{i - 2}  \le \left|\xi - \xi \left(G_\alpha^i \right) \right| \le 2^{i + 2}  \right\}$. If $N(J_l) \ge \epsilon_3^\frac12 2^{l_2 - 5}$, then we have
\begin{align}\label{eq6.50v118}
 & \left|\int_{J_l} \int_{\mathbb{R}^2} \int_{\mathbb{R}} \overline{  w(t,y,x) }  F^{low, l_2}_{\alpha,i}(v( t,y,x) )  \,\mathrm{d}y \mathrm{d}x \mathrm{d}t \right|\\
\lesssim &  2^{ \frac{2l_2}{q_0}    } \epsilon_3^{\frac14- \frac1{q_0}   } 2^{- \frac{i}2 } \|w \|_{L_t^{p_0 }  L_y^{q_0}  L_x^2 \left(J_l \times \mathbb{R}^2 \times \mathbb{R} \right)} \left\|P_{\xi \left(G_\alpha^i \right), 2^{i - 5} \le \cdot \le 2^{i +5}}^y v(t,y,x) \right\|_{U_\Delta^2 \left(G_\alpha^i, L_x^2 \right)}. \notag
\end{align}
\end{lemma}
\begin{proof}
 By \eqref{eq5.1v51}, we see $|\xi - \xi(t)| \le 2^{l_2 + 2}$ implies $ |\xi - \xi(J_l)| \le \epsilon_3^{- \frac12} N(J_l)$ for $t \in J_l$. By the argument in the proof of Lemma \ref{le4.8v56} and  H\"older's inequality, we have
\begin{align}\label{eq6.53v118}
& \bigg| \int_{J_l} \int_{\mathbb{R}^2} \int_{\mathbb{R}} \overline{  w (t,y ,x ) }  F^{low,l_2}_{\alpha, i} v( (t,y ,x ) )\,\mathrm{d}y \mathrm{d}x \mathrm{d}t \bigg| \\
\le & \left\| e^{i \tau \left( - \Delta_x  + x^2 \right) } w P_{\xi(t), 2^{l_2}}^y e^{i \tau \left( - \Delta_x  + x^2 \right)}  v \right\|_{L_{\tau, x}^{p_0 } L_t^{p_0 } L_y^2 \left(  [0 ,\pi ] \times \mathbb{R} \times J_l \times \mathbb{R}^2 \right)} \notag\\
& \qquad
\cdot  \left\|P_{\xi(J_l), \lesssim \epsilon_3^{- \frac12} N(J_l)}^y e^{i \tau \left( - \Delta_x + x^2 \right)} v P_{\xi \left(G_\alpha^i \right), 2^{i - 5} \le \cdot 2^{i + 5}}^y e^{i \tau \left( - \Delta_x  + x^2 \right)} v \right\|_{L_{\tau, x}^{\frac{p_0}{p_0 - 1} } L_t^{\frac{p_0}{p_0 - 1 }} L_y^2 \left( [0 , \pi ] \times \mathbb{R} \times
J_l \times \mathbb{R}^2  \right)} . \notag
\end{align}
By Minkowski's inequality, H\"older's inequality, Lemma \ref{le2.1v27},
 Bernstein's inequality and the conservation of mass, we have
\begin{align}\label{eq6.54v118}
& \left\|e^{i \tau \left(- \Delta_x  + x^2 \right)} w P_{\xi(t), 2^{l_2}}^y e^{i \tau \left( - \Delta_x + x^2 \right) }v \right\|_{L_{\tau, x}^{p_0 }  L_t^{p_0 } L_y^2 \left( [ 0, \pi ] \times \mathbb{R} \times  J_l \times \mathbb{R}^2 \right)} \\
\lesssim  & \left\| \left\| e^{i \tau \left( - \Delta_x  + x^2 \right)} w \right\|_{L_{\tau, x}^{2 p_0  }  ([0, \pi] \times \mathbb{R})} \left\|P_{\xi(t), 2^{l_2}}^y e^{i \tau \left( - \Delta_x  + x^2 \right)} v \right\|_{L_{\tau, x}^{2 p_0  } ([0, \pi ] \times \mathbb{R})} \right\|_{L_t^{p_0 } L_y^2(J_l \times \mathbb{R}^2)} \notag \\
\lesssim & \|w \|_{L_t^{p_0 } L_y^{ q_0 }  L_x^2(J_l \times \mathbb{R}^2 \times \mathbb{R})}  \left\| P_{\xi(t), 2^{l_2}}^y v \right\|_{L_t^\infty L_y^{p_0} L_x^2(J_l \times \mathbb{R}^2 \times \mathbb{R})}
\lesssim 2^{ \frac{2   l_2} {q_0} } \|w\|_{L_t^{p_0 } L_y^{q_0}  L_x^2(J_l \times \mathbb{R}^2 \times \mathbb{R})} . \notag
\end{align}
Next, we use the vector-valued version of transference principle to estimate
\begin{align*}
\bigg\| P_{\xi(J_l), \lesssim \epsilon_3^{- \frac12} N(J_l)}^y e^{i \tau \left( - \Delta_x + x^2 \right)} v P_{\xi \left(G_\alpha^i \right), 2^{i - 5} \le \cdot \le 2^{i +5}}^y e^{i \tau \left( - \Delta_x  + x^2 \right)} v \bigg\|_{L_{\tau, x}^{p_0 }  L_t^{p_0 } L_y^2( [0, \pi ] \times \mathbb{R} \times
J_l \times \mathbb{R}^2)}.
\end{align*}
Then, using the similar argument of Corollary 1.6 in \cite{Candy}, we are reduced to consider $P_{\xi(J_l), \lesssim \epsilon_3^{- \frac12} N(J_l)}^y v$ and $P_{\xi \left(G_\alpha^i \right), 2^{i - 5} \le \cdot \le 2^{i+5}}^y v$ are $U_\Delta^2(L_x^2)-$atoms. Let
\begin{align*}
P_{\xi(J_l), \lesssim \epsilon_3^{- \frac12} N(J_l)}^yv = \sum\limits_{I \in \mathcal{I}} \chi_I e^{it \Delta_y } f_I(y,x), \
P_{\xi(G_\alpha^i), 2^{i - 5} \le \cdot \le 2^{i +5}}^y v = \sum\limits_{J \in \mathcal{J}} \chi_J e^{it \Delta_y  } g_J(y,x),
\end{align*}
where $\mathcal{I}$ and $\mathcal{J}$ are partitions as in the definition of $U_\Delta^2(L_x^2) $.
We see  by Lemma \ref{le2.2v118}, H\"older's inequality and Lemma \ref{le2.1v27},
\begin{align}
& \left\| e^{i t \Delta_y } e^{i\tau \left( - \Delta_x  + x^2 \right)} f_I  e^{it \Delta_y} e^{i \tau \left( - \Delta_x  + x^2 \right)} g_J \right\|_{L_{\tau, x}^{\frac{p_0}{p_0 - 1 } } L_t^{\frac{p_0}{p_0 - 1} } L_y^2([0, \pi ] \times \mathbb{R} \times
J_l \times \mathbb{R}^2)} \notag \\
\lesssim & 2^{- \frac{i}2} \left( \epsilon_3^{- \frac12} N(J_l) \right)^{\frac12 - \frac2{q_0}  }  \left\|e^{i \tau \left( - \Delta_x  + x^2 \right) } f_I  \right\|_{L_y^2 L_{\tau, x}^{\frac{2p_0}{p_0 - 1 } } ( \mathbb{R}^2 \times [0, \pi ] \times \mathbb{R} )} \left\|e^{i \tau \left( - \Delta_x  + x^2 \right)} g_J  \right\|_{ L_y^2 L_{\tau, x}^{\frac{2 p_0 } {p_0  - 1 } } ( \mathbb{R}^2 \times [0, \pi ] \times \mathbb{R})} \notag \\
\lesssim & 2^{- \frac{i}2} \left( \epsilon_3^{- \frac12} N(J_l) \right)^{\frac12 - \frac2{q_0}  } \|f_I  \|_{L_{y,x}^2  } \|g_J \|_{L_{y,x}^2  }. \notag
\end{align}
Then, we have
\begin{align}\label{eq6.60v118}
&  \left\| P_{\xi(J_l), \lesssim \epsilon_3^{- \frac12} N(J_l)}^y e^{i \tau \left( - \Delta_x + x^2 \right)} v P_{\xi \left(G_\alpha^i \right), 2^{i - 5} \le \cdot \le 2^{i +5}}^y e^{i \tau \left( - \Delta_x  +x^2 \right)} v \right\|_{L_{\tau, x}^{\frac{p_0 }{p_0 - 1  } } L_t^{ \frac{p_0 }{p_0  - 1 } } L_y^2 \left([0, \pi ] \times \mathbb{R} \times
J_l \times \mathbb{R}^2 \right)} \\
\lesssim & 2^{- \frac{i}2}  \left( \epsilon_3^{- \frac12} N(J_l) \right)^{\frac12 - \frac2{q_0}  }  \left\|P_{\xi(J_l), \lesssim \epsilon_3^{- \frac12} N(J_l)}^y v \right\|_{U_\Delta^2(J_l, L_x^2)} \left\|P_{\xi \left(G_\alpha^i \right), 2^{i-5} \le \cdot \le 2^{i+5}}^y v \right\|_{U_\Delta^2(J_l, L_x^2)}. \notag
\end{align}
Therefore, \eqref{eq6.50v118} follows from \eqref{eq6.53v118}, \eqref{eq6.54v118} and \eqref{eq6.60v118}.
\end{proof}
We first consider \eqref{eq5.45v51}. By duality, we have
\begin{align*}
\eqref{eq5.45v51} = \sum\limits_{0\le l_2 \le i - 10} \bigg( \sum\limits_{ \substack{ J_l \subseteq G_\alpha^i,\\
 N(J_l) \ge \epsilon_3^\frac12 2^{l_2 - 5}} } \sup\limits_{\|w\|_{V_\Delta^2 \left(J_l, L_x^2 \right) } = 1} \left|\int_{J_l} \int_{\mathbb{R}^2} \int_{\mathbb{R}} \overline{  w(t,y ,x ) }   P_{\xi \left(G_\alpha^i \right), 2^{i - 2} \le \cdot \le 2^{i+2}}^y  F^{low, l_2}_{\alpha, i}(v(t,y, x))  \,\mathrm{d}y \mathrm{d}x \mathrm{d}t \right|^2 \bigg)^\frac12.
\end{align*}
By $\|w\|_{L_t^{p_0 } L_y^{q_0}  L_x^2 \left(J_l \times \mathbb{R}^2 \times \mathbb{R} \right)} \lesssim \|w\|_{V_\Delta^2 \left(J_l, L_x^2 \right)} \lesssim 1$, and \eqref{eq6.50v118}, we get
\begin{align}
\eqref{eq5.45v51} \lesssim & \sum\limits_{0 \le l_2 \le i -10} \bigg( \sum\limits_{ \substack{ J_l \subseteq G_\alpha^i, \\
N(J_l) \ge \epsilon_3^\frac12 2^{l_2 - 5} } } 2^{ \frac{4  l_2}{q_0} } \epsilon_3^{\frac12 - \frac2{q_0}  } 2^{- i} \bigg)^\frac12 \left\|P_{\xi \left(G_\alpha^i \right), 2^{i - 5} \le \cdot \le 2^{i+ 5}}^y v(t,y,x) \right\|_{U_\Delta^2 \left(G_\alpha^i, L_x^2 \right)} \notag \\
\lesssim &  \left( \epsilon_3^{- \frac12} N(J_l) \right) \epsilon_3^{ \frac14- \frac1{q_0}   }  \left\|P_{\xi \left(G_\alpha^i \right), 2^{i - 5} \le \cdot \le 2^{i +5}}^y v \right\|_{U_\Delta^2 \left(G_\alpha^i, L_x^2 \right)}
\lesssim
 \epsilon_2^2  \left\|P_{\xi \left(G_\alpha^i \right), 2^{i - 5} \le \cdot \le 2^{i +5}}^y v \right\|_{U_\Delta^2 \left(L_x^2, G_\alpha^i \right)}. \notag
\end{align}
We now consider \eqref{eq5.46v51}. By duality and Lemma \ref{le6.28v118}, we have
\begin{align}\label{eq6.65v118}
 & \bigg\|\int_{J_l} e^{- it \Delta_y} P_{\xi \left(G_\alpha^i \right), 2^{i - 2} \le \cdot \le 2^{i+2}}^y F^{low, l_2}_{\alpha, i}(v(t,y,x))
  \,\mathrm{d}t \bigg\|_{L_{y,x}^2} \\
\lesssim & \sup\limits_{ \left \|w_0 \right\|_{L_{y,x}^2} = 1} \left\| e^{it \Delta_y} w_0 \right\|_{L_t^{p_0} L_y^{q_0}  L_x^2 \left(J_l \times \mathbb{R}^2 \times \mathbb{R} \right)} 2^{ \frac{2   l_2 }{q_0} } \epsilon_3^{\frac14- \frac1{q_0}  } 2^{- \frac{i}2} \left\|P_{\xi \left(G_\alpha^i \right), 2^{i - 5} \le \cdot \le 2^{i+5}}^y v \right\|_{U_\Delta^2 \left(G_\alpha^i, L_x^2 \right)}, \notag
\end{align}
where $\mathcal{F}_y w_0$ is supported on $\left\{ \xi: 2^{i - 5} \le  \left|\xi - \xi \left(G_\alpha^i \right) \right| \le 2^{i +5} \right\}$ in the above estimate. For fixed $i$, we take $q_0 = 20 + 2i $, then $2^{\frac{ i}{q_0} } \lesssim 1$. For the right hand side of \eqref{eq6.65v118}, by H\"older's inequality, Young's inequality, \eqref{eq5.1v51}, \eqref{eq6.4v58}, \eqref{eq5.4v51}, and the conservation of mass, we have
\begin{align*}
\eqref{eq5.46v51} & \lesssim \sum\limits_{ \substack{ 0 \le l_2 \le i - 10,\\
 J_l \subseteq G_\alpha^i,\\
  N(J_l) \ge \epsilon_3^\frac12 2^{l_2 - 5}}}  \sup\limits_{  \left\|w_0 \right\|_{L_{y,x}^2} = 1} 2^{ \frac{2   l_2}{q_0 } } \epsilon_3^{\frac14- \frac1{q_0}   } 2^{- \frac{i}2}  \left\|e^{i t\Delta_y} w_0 \right\|_{L_t^{p_0 } L_y^{q_0}  L_x^2(J_l \times \mathbb{R}^2 \times \mathbb{R})} \left\|P_{\xi \left(G_\alpha^i \right), 2^{i - 5} \le \cdot \le 2^{i +5}}^y v \right\|_{U_\Delta^2 \left(G_\alpha^i, L_x^2 \right)} \\
& \lesssim \sum\limits_{0 \le l_2 \le i - 10} 2^{ \frac{2  l_2}{q_0} }  \epsilon_3^{\frac14- \frac1{q_0}  } \bigg( \sum\limits_{
\substack{ J_l \subseteq G_\alpha^i,\\
 N(J_l) \ge \epsilon_3^\frac12 2^{l_2 - 5}}}
  \left\|e^{it \Delta_y} w_0 \right\|_{L_t^{p_0  } L_y^{q_0}  L_x^2(J_l \times \mathbb{R}^2 \times \mathbb{R})}^{p_0} \bigg)^\frac1{p_0 } \\
&\qquad  \qquad \cdot \bigg( \sum\limits_{J_l \subseteq G_\alpha^i} \left( 2^{- \frac{i}2} \right)^{\frac{p_0}{p_0 - 1  } }
\bigg)^{\frac{p_0 - 1 }{p_0}}
\left\|P_{\xi \left(G_\alpha^i \right), 2^{i - 5} \le \cdot \le 2^{i +5} }^y v \right\|_{U_\Delta^2 \left(G_\alpha^i, L_x^2 \right)} \\
& \lesssim \sum\limits_{ \substack{ 0 \le l_2\le i - 10, \\
2^{l_2 - 5} \le \epsilon_3^{- \frac12} N(J_l)}}
 2^{ \frac{2   l_2}{q_0}}  \epsilon_3^{\frac14- \frac1{q_0} } \left\|e^{it \Delta_y} w_0 \right\|_{L_t^{p_0 } L_y^{q_0}  L_x^2 \left(G_\alpha^i \times \mathbb{R}^2 \times \mathbb{R} \right)} \left( 2^{ \frac{2 i}{q_0 + 2} }  \right)^{ \frac12 + \frac1{q_0}   }
 \left\|P_{\xi \left(G_\alpha^i \right), 2^{i - 5} \le \cdot \le 2^{i+5}}^y v \right\|_{U_\Delta^2 \left(G_\alpha^i, L_x^2 \right)} \\
& \lesssim  \sum\limits_{ \substack{ 0 \le l_2 \le i - 10,\\
 2^{l_2 - 5} \le \epsilon_3^{- \frac12} N(J_l)} }
  2^{ \frac{2   l_2}{q_0} } \epsilon_3^{\frac14- \frac1{q_0}   } 2^{\frac{  i}{q_0} } \left\|P_{\xi \left(G_\alpha^i \right), 2^{i - 5} \le \cdot \le 2^{i +5}}^y v \right\|_{U_\Delta^2 \left(G_\alpha^i, L_x^2 \right)}  \lesssim \epsilon_3^\frac14 \left\|P_{\xi \left(G_\alpha^i \right), 2^{i - 5} \le \cdot \le 2^{i + 5}}^y v \right\|_{U_\Delta^2 \left( L_x^2, G_\alpha^i \right)}.
\end{align*}
For the estimates of \eqref{eq5.47v51} and \eqref{eq5.48v51}, we separate the proofs in the next section using two bilinear Strichartz estimates.
\end{proof}

\subsubsection{Two bilinear Strichartz estimates}

We have the following two bilinear Strichartz estimates.
\begin{theorem}[First bilinear Strichartz estimate] \label{th5.16v51}
	Let $w_0 \in L_{y,x}^2 (\mathbb{R}^2 \times \mathbb{R})$ with $supp \, \mathcal{F}_y w_0$ is supported on $\left\{ \xi: 2^{i - 5} \le  \left|\xi - \xi \left(G_\alpha^i \right) \right| \le 2^{i+5} \right\}$. Then for any $0 \le l_2 \le i - 10$, we have on $G_\beta^{l_2} \subseteq G_\alpha^i$,
	\begin{align}\label{eq5.58v51}
		\left\| \left\|e^{it \Delta_y } w_0 \right\|_{L_x^2}  \left\|P^y_{\xi(t), \le  2^{ l_2} }   v  \right\|_{L_x^2} \right\|_{L_{t,y}^2 \left(G_\beta^{l_2} \times \mathbb{R}^2 \right)}^2
		\lesssim \|w_0 \|_{L_{y,x}^2 }^2 \left( 1 + \| v  \|_{\tilde{X}_i \left(G_\alpha^i \right)}^4 \right).
	\end{align}
\end{theorem}
\begin{theorem}[Second bilinear Strichartz estimate]\label{th5.17v51}
Let $w_0 \in L_{y,x}^2 (\mathbb{R}^2 \times \mathbb{R})$ with $supp\,  \mathcal{F}_y w_0$ is supported on $\left\{ \xi: 2^{i - 5} \le \left|\xi - \xi \left(G_\alpha^i \right) \right| \le 2^{i+5} \right\}$, we have
	\begin{align}\label{eq5.63v51}
		\sum\limits_{0 \le l_2 \le i-10}\left\| \left\|e^{it \Delta_y } w_0 \right\|_{L_x^2}  \left\|P^y_{\xi(t), \le 2^{ l_2 }  } v  \right\|_{L_x^2} \right\|_{L_{t,y}^2 \left(G_\alpha^i \times \mathbb{R}^2 \right)}^2
		\lesssim \|w_0 \|_{L_{y,x}^2 }^2 \left( 1 + \| v  \|_{ \tilde{X}_i \left(G_\alpha^i \right)}^6 \right).
	\end{align}
\end{theorem}
With the above two bilinear Strichartz estimates, we can now estimate \eqref{eq5.47v51} and \eqref{eq5.48v51}.
\begin{proof}[Estimate of \eqref{eq5.47v51}]
For any $0 \le l_2 \le i - 10$, by the fact that $G_\alpha^i$ consists of $2^{10}$ subintervals $G_\beta^{i - 10}$, Proposition \ref{pr4.5v51} and Theorem \ref{th5.16v51} on the subintervals $G_\beta^{i-10}$,
we get
\begin{align}\label{eq5.60v51}
	\left\|\left \|P^y_{\xi\left(G_\alpha^i \right), 2^{  i - 5 }  \le \cdot \le 2^{  i+5} }  v  \right\|_{L_x^2} \left\|P^y_{\xi(t), \le  2^{ l_2}  }  v  \right\|_{L_x^2} \right\|_{L_{t,y}^2 \left(G_\alpha^i \times \mathbb{R}^2 \right)}
	\lesssim \left\|P^y_{\xi \left(G_\alpha^i \right),  2^{ i - 5 }  \le \cdot \le 2^{  i+5} }  v  \right\|_{U_\Delta^2 \left(L_x^2, G_\alpha^i \right)} \left( 1 + \| v \|_{\tilde{X}_i \left(G_\alpha^i \right)}^2 \right).
\end{align}
For any $G_\beta^{l_2} \subseteq G_\alpha^i$, choose $w_\beta^{l_2} \in V^2_\Delta \left(L_x^2, G_\beta^{l_2} \right)$ with $supp \,\mathcal{F}_y w_\beta^{l_2} \subseteq \left\{ \xi: 2^{i - 2} \le  \left|\xi - \xi \left(G_\alpha^i \right) \right| \le 2^{i+2} \right\}$ and $ \left\|w_\beta^{l_2} \right\|_{V_\Delta^2 \left(L_x^2, G_\beta^{l_2} \right) } = 1$. Then by H\"older's inequality, \eqref{eq3.8v41}, and \eqref{eq5.60v51}, we have
\begin{align}\label{eq6.59v67}
	& \bigg( \sum\limits_{ \substack{ G_\beta^{l_2} \subseteq G_\alpha^i, \\
			N \left(G_\beta^{l_2} \right) \le \epsilon_3^\frac12 2^{l_2-5} } }
	\bigg\| \overline{ w_\beta^{l_2} }  \sum\limits_{ \substack{ n_1,n_2,n_3,n \in \mathbb{N},\\
			n_1 -n_2 + n_3 = n}  }
	\Pi_n \left( P^y_{\xi(t), 2^{  l_2} }  v_{n_1} \overline{ P^y_{\xi(t), \le 2^{  l_2} }   {v}_{n_2} }
	P^y_{\xi \left(G_\alpha^i \right),  2^{ i - 5  } \le \cdot \le  2^{ i+5} }  v_{n_3}  \right) \bigg\|_{L_{t,y,x }^1 \left(G_\beta^{l_2} \times \mathbb{R}^2 \times \mathbb{R}
		\right)}^2 \bigg)^\frac12 \\
	\lesssim &  \bigg( \sup\limits_{ \substack{ G_\beta^{l_2} \subseteq G_\alpha^i, \\ N \left(G_\beta^{l_2} \right) \le \epsilon_3^\frac12 2^{l_2- 5}} }
	\left\| \left\| w_\beta^{l_2} \right\|_{L_x^2} \left \|P^y_{\xi(t),  2^{ l_2} }  v  \right\|_{L_x^2} \right\|_{L_{t,y}^2 \left(G_\beta^{l_2} \times \mathbb{R}^2
		\right)} \bigg)
	\left\|P^y_{\xi \left(G_\alpha^i \right),  2^{ i - 5 }  \le \cdot \le 2^{  i+5} }  v \right\|_{U_\Delta^2 \left(L_x^2, G_\alpha^i \right)} \left( 1 + \| v  \|_{\tilde{X}_i \left(G_\alpha^i \right)}^2  \right). \notag
\end{align}
By Proposition \ref{pr4.7v51}, \eqref{eq6.4v57}, \eqref{eq6.7v57}, and $N \left(G_\beta^{l_2} \right) \le \epsilon_3^\frac12 2^{l_2- 5}$, we can estimate the term in the first bracket on the right hand side of \eqref{eq6.59v67} as follows
\begin{align*}
	& \left\| \left\| w_\beta^{l_2} \right\|_{L_x^2} \left\|P^y_{\xi(t), 2^{ l_2} }  v  \right\|_{L_x^2} \right\|_{L_{t,y}^2 \left(G_\beta^{l_2} \times \mathbb{R}^2 \right)}\\
	\lesssim &  \left\| \left\|w_\beta^{l_2} \right\|_{L_x^2} \left\|P^y_{\xi(t),  2^{ l_2}  } v  \right\|_{L_x^2} \right\|_{L_t^3 L_y^\frac32 \left(G_\beta^{l_2} \times \mathbb{R}^2 \right)}^\frac12 \left\|w_\beta^{l_2} \right\|_{L_t^3 L_y^6 L_x^2 \left(G_\beta^{l_2}  \times \mathbb{R}^2 \times \mathbb{R} \right)}^\frac12 \left\|P^y_{\xi(t), 2^{ l_2} }   v  \right\|_{L_t^3 L_y^6 L_x^2 \left(G_\beta^{l_2} \times \mathbb{R}^2 \times \mathbb{R} \right)}^\frac12
	\lesssim 2^{\frac{l_2- i}6} \| v \|_{\tilde{Y}_i \left(G_\alpha^i \right)}.
\end{align*}
Thus, by the above inequalities, we obtain
\begin{align*}
	\eqref{eq5.47v51} \lesssim \left\|P^y_{\xi \left(G_\alpha^i \right),  2^{ i - 5 }  \le \cdot \le  2^{ i+5}  } v  \right\|_{U_\Delta^2 \left(L_x^2, G_\alpha^i \right)} \| v  \|_{\tilde{Y}_i \left(G_\alpha^i \right)} \left( 1 + \| v \|_{\tilde{X}_i \left(G_\alpha^i \right)}^2 \right).
\end{align*}
\end{proof}
\begin{proof}[Estimate of \eqref{eq5.48v51}]
Let $w_0 \in L_{y,x}^2$ have unit norm with $\mathcal{F}_y w_0$ is supported on $\left\{ \xi: 2^{i - 2} \le  \left|\xi - \xi\left(G_\alpha^i \right) \right| \le 2^{i+2} \right\}$. By the H\"older inequality and Proposition \ref{pr4.5v51}, we have
\begin{align*}
	& \left\| \int_{G_\beta^{l_2}} e^{-it \Delta_y } P^y_{\xi \left(G_\alpha^i \right), 2^{ i - 2 }  \le \cdot \le 2^{  i+2} } F^{low, l_2}_{\alpha , i}(v(t,y,x)) \,\mathrm{d}t \right\|_{L_{y,x}^2 } \\
	\lesssim & \sup\limits_{\|w_0 \|_{L_{y,x}^2 } = 1} \left\|e^{it \Delta_y } w_0 \cdot  F^{low, l_2}_{\alpha, i}(v(t,y,x))\right\|_{L_{t,y,x}^1 \left(G_\beta^{l_2} \times \mathbb{R}^2 \times \mathbb{R} \right)} \\
	\lesssim & \sup\limits_{ \|w_0 \|_{L_{y,x}^2} = 1} \left\| \left\|e^{it \Delta_y } w_0\right \|_{L_x^2} \left\|P^y_{\xi(t), 2^{  l_2} }  v  \right\|_{L_x^2} \left\|P^y_{\xi(t), \le 2^{  l_2 } } v  \right\|_{L_x^2} \left\|P^y_{\xi \left(G_\alpha^i \right),  2^{ i- 5 }  \le \cdot \le 2^{  i+5} }  v  \right\|_{L_x^2} \right\|_{L_{t,y}^1 \left(G_\beta^{l_2} \times \mathbb{R}^2  \right)} \\
	\lesssim & 2^{\frac{l_2-i}2} \left\|P^y_{\xi \left(G_\beta^{l_2} \right), 2^{  l_2 - 2 }  \le \cdot \le  2^{ l_2 + 2 }  } v  \right\|_{U_\Delta^2 \left(L_x^2, G_\beta^{l_2} \right)}
	\left\| \left \|P^y_{\xi(t), \le 2^{  l_2} }  v  \right\|_{L_x^2} \left\|P^y_{\xi \left(G_\alpha^i \right),  2^{ i- 5 } \le \cdot \le  2^{ i+5} }  v  \right\|_{L_x^2} \right\|_{L_{t,y}^2 \left(G_\beta^{l_2} \times \mathbb{R}^2 \right)}.
\end{align*}
Then by the Cauchy-Schwarz inequality and \eqref{eq5.63v51}, we have
\begin{align*}
	\eqref{eq5.48v51} & \lesssim
	\|  v \|_{\tilde{Y}_i(G_\alpha^i)} \bigg( \sum\limits_{0 \le l_2 \le i-10} \left\| \left\|P^y_{\xi(t), \le  2^{ l_2} }  v  \right\|_{L_x^2} \left\|P^y_{\xi(G_\alpha^i),  2^{ i-5 }  \le \cdot \le 2^{ i+5} } v \right\|_{L_x^2} \right\|_{L_{t,y}^2(G_\alpha^i \times \mathbb{R}^2 )}^2 \bigg)^\frac12\\
	& \lesssim \| v  \|_{\tilde{Y}_i(G_\alpha^i)} \left\|P^y_{\xi(G_\alpha^i), 2^{ i - 5  } \le \cdot \le  2^{ i+5} }  v  \right\|_{U_\Delta^2(L_x^2, G_\alpha^i)} \left( 1 + \| v  \|_{\tilde{X}_i(G_\alpha^i)}^3 \right).
\end{align*}
\end{proof}
Therefore, this completes the proof of Theorem \ref{th5.14v51}. Then we can prove Theorem \ref{th5.12v51} by summation with respect to $i$ in the same way as \eqref{eq5.38v51} and \eqref{eq5.39v51} in Theorem \ref{th5.14v51} and Theorem \ref{th5.13v51}.

\subsubsection{Proofs of the bilinear Strichartz estimates}
It remains to prove the two bilinear Strichartz estimates, that is Theorem \ref{th5.16v51} and Theorem \ref{th5.17v51}. The proofs of these results are basically the same and rely on the interaction Morawetz estimates of the \eqref{eq1.3v28} system, the argument here follow from that in \cite{D1}. We shall only present the proof of Theorem \ref{th5.16v51} here, because argument of the proof of Theorem \ref{th5.17v51} is similar to the proof of Theorem \ref{th5.16v51} and also rely on the result of Theorem \ref{th5.16v51} as the proof of the corresponding bilinear Strichartz estimate in \cite{D1}.
\begin{proof}[Proof of Theorem \ref{th5.16v51}]
Let $ w = e^{it \Delta_y } w_0$ and $ \tilde{w} = P^y_{\xi(t), \le  2^{ l_2} }  v $. Then $ w $ and $ \tilde{w} $ satisfy $i \partial_t w  + \Delta_y  w  = 0$, and
\begin{align*}
	i \partial_t \tilde{w}  + \Delta_y  \tilde{w}  = F( \tilde{w} ) + N_1 + N_2= F( \tilde{w} )+ N,
\end{align*}
where
\begin{align*}
	N_1 = P^y_{\xi(t), \le  2^{ l_2} }  F( v ) - F( \tilde{w} ) ,
\end{align*}
and $N_2= \left( \frac{d}{dt} P^y_{\xi(t), \le 2^{  l_2}  }  \right)  v $ with $\tfrac{d}{dt} P^y_{\xi(t), \le 2^{  l_2} } $ being given by the Fourier multiplier $- \nabla \phi\left( \frac{\xi- \xi(t)}{2^{l_2}} \right) \frac{\xi'(t)}{2^{l_2}}$.

We define the interaction Morawetz action
\begin{align*}
	M(t) & = \int_{\mathbb{R}^2\times \mathbb{R}} \int_{\mathbb{R}^2 \times \mathbb{R}} \left| \tilde{ w }
	(t, \tilde{y}, \tilde{x} ) \right|^2 \frac{ y - \tilde{y}}{ \left| y - \tilde{y} \right|}
	\Im \left( \bar{ w } \nabla_y   w \right)(t,y,x) \,\mathrm{d}y \mathrm{d}\tilde{y} \mathrm{d}x \mathrm{d} \tilde{x}\\
	& \ + \int_{\mathbb{R}^2\times \mathbb{R}} \int_{\mathbb{R}^2 \times \mathbb{R}} \left| w \left(t, \tilde{y}, \tilde{x} \right) \right|^2 \frac{  y - \tilde{y}} { \left|y - \tilde{y} \right|}
	\Im \left( \overline{ \tilde{w } } \nabla_y  \tilde{ w}  \right)(t,y,x) \,\mathrm{d} y \mathrm{d}\tilde{y} \mathrm{d}x \mathrm{d} \tilde{x}.
\end{align*}
After some tedious calculation, we get
\begin{align}
	& \int_{G_\beta^{l_2}}  \int_{\mathbb{R}^2} \int_{\mathbb{R}} \int_{\mathbb{R}} \left| \overline{ \tilde{ w } \left(t,y, \tilde{x} \right)  }
	w (t,y,x) \right|^2 \, \mathrm{d}x \mathrm{d} \tilde{x} \mathrm{d}y \mathrm{d}t\notag  \\
	\lesssim & 2^{l_2 - 2i} \sup\limits_{t \in G_\beta^{l_2} } |M(t)|  \notag\\
	& +2^{l_2 - 2i} \left| \int_{G_\beta^{l_2}} \int_{\mathbb{R}^2} \int_{\mathbb{R}^2} \int_{\mathbb{R}} \int_{\mathbb{R}} \left|  w \left(t, \tilde{y}, \tilde{x} \right) \right|^2 \frac{ y - \tilde{y}}{ \left|y - \tilde{y} \right|} \Im \left( \bar{N} \left( \nabla_y  - i \xi(t) \right)  \tilde{ w}  \right)(t,y,x) \,\mathrm{d}x \mathrm{d} \tilde{x} \mathrm{d}y \mathrm{d}\tilde{y} \mathrm{d}t \right| \label{eq6.3v51} \\
	& +  2^{l_2 - 2i} \left| \int_{G_\beta^{l_2}} \int_{\mathbb{R}^2} \int_{\mathbb{R}^2} \int_{\mathbb{R}} \int_{\mathbb{R}} \left| w (t, \tilde{y}, \tilde{x} )\right|^2 \frac{ y - \tilde{y}} { \left| y - \tilde{y} \right|} \Im \left( \overline{ \tilde{ w } } \left( \nabla_y   - i \xi(t) \right) N \right)(t,y,x) \,\mathrm{d} x \mathrm{d} \tilde{x} \mathrm{d}y \mathrm{d} \tilde{y} \mathrm{d}t \right| \label{eq6.4v51} \\
	& + 2^{l_2 - 2i} \left| \int_{G_\beta^{l_2}} \int_{\mathbb{R}^2} \int_{\mathbb{R}^2} \int_{\mathbb{R}} \int_{\mathbb{R}} \Im ( \bar{ w } \left( \nabla_y  - i \xi(t) \right) w )(t, \tilde{y}, \tilde{x}) \frac{ y - \tilde{y}}{ \left|y - \tilde{y} \right|} \Im \left(\overline{ \tilde{ w } } N \right)(t,y,x) \,\mathrm{d}x \mathrm{d} \tilde{x} \mathrm{d} y \mathrm{d} \tilde{y} \mathrm{d}t \right|. \label{eq6.5v51}
\end{align}
By the invariance of the Galilean transformation of $M(t)$, H\"older's inequality, and the conservation of mass, we infer that $2^{l_2 - 2i } \sup\limits_{ t\in G_\beta^{l_2}} |M(t)|$ can be bounded by the right hand side of \eqref{eq5.58v51}.
\medskip

\textbf{Estimate of \eqref{eq6.3v51}.}
\medskip

By  \eqref{eq5.1v51}, \eqref{eq5.3v51}, Bernstein's inequality, the conservation of mass and the Strichartz estimate, we have
\begin{align*}
	| \eqref{eq6.3v51} | \lesssim &  2^{l_2 - 2i} \left\|N_1 \right\|_{L_{t,y}^\frac43 L_x^2 \left(G_\beta^{l_2} \times \mathbb{R}^2 \times \mathbb{R} \right)}
	\left\| \left( \nabla_y  - i \xi(t) \right)  \tilde{ w }  \right\|_{L_{t,y}^4 L_x^2 \left(G_\beta^{l_2} \times \mathbb{R}^2 \times \mathbb{R} \right)} \| w  \|_{L_t^\infty L_{y,x}^2
		 \left(G_\beta^{l_2} \times \mathbb{R}^2 \times \mathbb{R} \right)}^2\\
	& + 2^{-2 i } \|  w  \|_{L_t^\infty L_{y,x}^2  \left(G_\beta^{l_2} \times \mathbb{R}^2 \times \mathbb{R} \right)}^2 \int_{\mathbb{R}} \int_{G_\beta^{l_2}} |\xi'(t)| \left\| P^y_{\xi(t),  2^{ l_2 - 3 }  \le \cdot \le 2^{  l_2 + 3} }   v (t,y,x) \right\|_{L_y^2} \\
	& \qquad \left\|  \left( \nabla_y  - i \xi(t) \right) P^y_{\xi(t), \le  2^{ l_2}  }  v  \left(t, \tilde{y}, x \right) \right\|_{L_{\tilde{y}}^2} \,\mathrm{d}x\mathrm{d}t\\
	\lesssim  & 2^{l_2 - 2i} \|N_1 \|_{L_{t,y}^\frac43 L_x^2(G_\beta^{l_2} \times \mathbb{R}^2 \times \mathbb{R} )} \left\| ( \nabla_y  - i \xi(t))  \tilde{ w }  \right\|_{L_{t,y}^4 L_x^2 \left(G_\beta^{l_2} \times \mathbb{R}^2 \times \mathbb{R} \right)} \left\|  w_0 \right\|_{L_{y,x}^2 }^2 + 2^{2l_2 - 2i} \left\| w_0 \right \|_{L_{y,x}^2}^2.
\end{align*}
Let $m(t, \xi) = \frac{ \xi - \xi(t)}{ 2^{l_1}} \phi \left( \frac{ \xi - \xi(t)}{ 2^{l_1}} \right)  $, by Minkowski's inequality, Young's inequality, $\sup\limits_t \left\| \left( \mathcal{F}_\xi^{-1} m \right) (t,y) \right\|_{L_y^1}  \lesssim 1$ and \eqref{eq5.10v51}, we get
\begin{align*}
	\left\| \left(\nabla_y  - i \xi(t) \right)  \tilde{ w }  \right\|_{ L_{t,y}^4 L_x^2 \left(G_\beta^{l_2} \times \mathbb{R}^2 \times \mathbb{R} \right)}
	&\lesssim  \sum\limits_{0 \le l_1 \le l_2} \left\| \left( \nabla_y  - i \xi(t) \right) P^y_{\xi(t), 2^{  l_1} }   v  \right\|_{L_{t,y}^4 L_x^2 \left(G_\beta^{l_2} \times \mathbb{R}^2 \times \mathbb{R} \right)}\\
	\lesssim   \sum\limits_{ 0 \le l_1 \le l_2 } &2^{l_1}  \left\| \int \left| \left( \mathcal{F}_\xi^{-1} m \right)\left(t, y - \tilde{y} \right) \right| \left\| \left( P^y_{\xi(t),  2^{ l_1} }   v \right)  \left(t, \tilde{y}, x \right ) \right\|_{L_x^2}  \,\mathrm{d} \tilde{y} \right\|_{L_{t,y}^4 \left(G_\beta^{l_2} \times \mathbb{R}^2 \right)} \\
	\lesssim \sum\limits_{ 0 \le l_1 \le l_2 } &2^{l_1} \left\|P^y_{\xi(t),  2^{ l_1} }  v  \right\|_{L_{t,y}^4 L_x^2 \left(G_\beta^{l_2} \times \mathbb{R}^2 \times \mathbb{R} \right)}
	\lesssim \sum\limits_{0 \le l_1 \le l_2} 2^{l_1} 2^{\frac{l_2 - l_1} 4} \| v  \|_{\tilde{X}_{l_2} \left(G_\beta^{l_2} \right) }
	\lesssim 2^{l_2} \| v   \|_{\tilde{X}_i \left(G_\alpha^i \right)}.
\end{align*}
Thus, it implies
\begin{align*}
	| \eqref{eq6.3v51} | \lesssim 2^{2l_2 - 2i} \left\| w_0 \right\|_{L_{y,x}^2 }^2 \| v \|_{\tilde{X}_i \left(G_\alpha^i \right)} \left\|N_1 \right\|_{L_{t,y}^\frac43 L_x^2 \left(G_\beta^{l_2} \times \mathbb{R}^2 \times \mathbb{R} \right)}
	+ 2^{2l_2 - 2i}  \left\|  w_0 \right\|_{L_{y,x}^2 }^2 .
\end{align*}
Let $v^l = P^y_{\xi(t), \le 2^{ l_2 - 5}  }  v $ and $v^h = P^y_{\xi(t), > 2^{  l_2 - 5} } v $.
We can then decompose $N_1$ as
\begin{align}
	N_1 & = \notag \\
	& \quad  P^y_{\xi(t), \le 2^{ l_2 } } \bigg( \sum\limits_{ \substack{ n_1,n_2,n_3,n\in \mathbb{N},\\
			n_1 - n_2 + n_3 = n } } \Pi_n \left(  v_{n_1}^l  \overline{ { v }_{n_2}^l }   v_{n_3}\right ) \bigg)
	- \sum\limits_{ \substack{ n_1,n_2 ,n_3 ,n \in \mathbb{N},\\
			n_1 - n_2+ n_3  = n} }
	\Pi_n \left( P^y_{\xi(t), \le 2^{  l_2} }   v_{n_1}^l \overline{ P^y_{\xi(t), \le  2^{ l_2} }   { v }_{n_2}^l }  P^y_{\xi(t), \le  2^{ l_2}  } v_{n_3}^l \right)  \label{eq6.9v51}\\
	& + 2  \bigg( P^y_{\xi(t), \le  2^{ l_2} }  \bigg( \sum\limits_{ \substack{ n_1,n_2,n_3,n \in \mathbb{N},\\
			n_1 - n_2 + n_3 = n}  } \Pi_n \left(  v_{n_1}^l \overline{ { v }_{n_2}^l }   v_{n_3}^h \right) \bigg)
	- \sum\limits_{ \substack{ n_1,n_2,n_3 ,n \in \mathbb{N},\\ n_1 - n_2 + n_3 = n} } \Pi_n \left( P^y_{\xi(t), \le  2^{ l_2} }  v_{n_1}^l \overline{
		P^y_{\xi(t), \le 2^{  l_2} }   { v }_{n_2}^l }  P^y_{\xi(t), \le  2^{ l_2} }  v_{n_3}^h  \right) \bigg) \label{eq6.10v51}\\
	& + P^y_{\xi(t), \le l_2} \bigg( \sum\limits_{ \substack{ n_1,n_2,n_3 , n \in \mathbb{N}, \\ n_1 -n_2 + n_3 =n } }  \Pi_n \left(  v_{n_1}^l \overline{ { v }_{n_2}^h }   v_{n_3}^l  \right) \bigg)
	- \sum\limits_{ \substack{ n_1,n_2,n_3,n \in \mathbb{N},\\ n_1 - n_2 +n_3 = n} }  \Pi_n \left( P^y_{\xi(t), \le 2^{  l_2} }  v_{n_1}^l  \overline{
		P^y_{\xi(t), \le  2^{ l_2} }   { v }_{n_2}^h }
	P^y_{\xi(t), \le 2^{  l_2} }  v_{n_3}^l \right) \label{eq6.11v51}\\
	&+ \mathcal{O} \bigg( P^y_{\xi(t), \le  2^{ l_2} }   \bigg( \sum\limits_{ \substack{ n_1,n_2,n_3,n \in \mathbb{N}, \\ n_1 - n_2 + n_3 =n} }  \Pi_n \left(   v_{n_1}^h \overline{ { v }_{n_2}^h }  v_{n_3} \right) \bigg)
	- \sum\limits_{ \substack{ n_1,n_2,n_3,n \in \mathbb{N}, \\ n_1 - n_2 + n_3 = n} }
	\Pi_n \left( P^y_{\xi(t), \le 2^{  l_2} }   v_{n_1}^h \overline{ P^y_{\xi(t), \le  2^{ l_2}  }  { v }_{n_2}^h }
	P^y_{\xi(t), \le 2^{  l_2} }  v_{n_3}  \right) \bigg), \label{eq6.12v51}
\end{align}
where the $``\mathcal{O}"$ in \eqref{eq6.12v51} means there are two high frequency factors in it. Observe that
\begin{align*}
	\eqref{eq6.9v51} = 0.
\end{align*}
We next consider \eqref{eq6.10v51} and \eqref{eq6.11v51}. Since their estimates are very similar, we only prove \eqref{eq6.10v51}. Since $\left( \mathcal{F}_y\,  v_{n_3}^h \right)(t, \sigma, x) $ is supported on $\left\{ \sigma: \left|\sigma  - \xi(t) \right| \le 2^{l_2 + 10} \right\}$, we have
\begin{align}\label{eq6.13v51}
	& P^y_{\xi(t), \le 2^{  l_2} }
 \bigg( \sum\limits_{\substack{ n_1,n_2,n_3,n \in \mathbb{N},\\
			n_1 - n_2 + n_3 = n} }
	\Pi_n \left(  v_{n_1}^l \overline{ { v }_{n_2}^l  }
	v_{n_3}^h  \right) \bigg) - \sum\limits_{ \substack{ n_1,n_2,n_3, n \in \mathbb{N}, \\
			n_1 - n_2 + n_3 = n} }
	\Pi_n \left( P^y_{\xi(t), \le 2^{  l_2} }   v_{n_1}^l \overline{ P^y_{\xi(t), \le 2^{  l_2} }   { v }_{n_2}^l } P^y_{\xi(t), \le  2^{ l_2} }   v_{n_3}^h \right)\\
	= & \notag \\
	\qquad  & \sum\limits_{l_1 \le l_2} \sum\limits_{ \substack{ n_1,n_2,n_3,n \in \mathbb{N}, \\
			n_1 -n_2 + n_3 = n} }
	\iiint e^{-i \tilde{y} \xi(t)} \Pi_n \left( \left( P^y_{\xi(t),  2^{ l_1} }   v_{n_1}^l \right) \left(\tilde{y}, x \right)  \overline{ e^{ -iz \xi(t)}  \left( P^y_{\xi(t), \le  2^{ l_1}  }  { v }_{n_2}^l \right) (z,x) }  v_{n_3}^h( \theta  ,x) \right)\label{eq6.68v71}\\
	\cdot
	\iiint &  e^{i\xi(y - \tilde{y})+ i \eta(\tilde{y} - z) + i \sigma (z -  \theta )}
	\left( \left( \phi \left( \frac{\xi- \xi(t)}{ 2^{l_2}} \right) - \phi \left( \frac{\sigma - \xi(t)}{2^{l_2}} \right)  \right) \phi \left( \frac{\sigma - \xi(t)}{ 2^{l_2 + 10}} \right) \psi_{l_1}(\xi- \eta) \phi \left( \frac{\eta- \sigma}{2^{l_1}} \right)  \right) \,\mathrm{d}\sigma \mathrm{d}\eta \mathrm{d}\xi \mathrm{d} z \mathrm{d} \tilde{y} \mathrm{d} \theta  \notag \\
	+ & \sum\limits_{l_1 \le l_2} \sum\limits_{ \substack{ n_1,n_2,n_3,n \in \mathbb{N}, \\
			n_1 - n_2 + n_3 = n } } \iiint \Pi_n \left(e^{-i \tilde{y} \xi(t)} \left( P^y_{\xi(t), \le 2^{  l_1} }   v_{n_1}^l \right) \left(\tilde{y}, x \right) \overline{  e^{- iz \xi(t)} \left( P^y_{\xi(t), 2^{  l_1 } } v_{n_2}^l \right) (z,x)} v_{n_3}^h( \theta ,  x) \right)\label{eq6.69v71}\\
	\cdot \iiint & e^{i \xi(y - \tilde{y} )+ i \eta ( \tilde{y} - z) + i\sigma ( z -  \theta )}
	\left( \left( \phi \left( \frac{ \xi - \xi(t)}{2^{l_2} } \right)
	- \phi \left( \frac{\sigma- \xi(t)} { 2^{l_2}} \right) \right) \phi  \left( \frac{\sigma - \xi(t)}{ 2^{l_2 + 10}}  \right)  \right) \psi_{l_1}( \eta - \sigma) \phi \left( \frac{\xi - \eta}{2^{l_1}} \right) \,\mathrm{d}\sigma \mathrm{d}\eta \mathrm{d}\xi \mathrm{d}\tilde{y} \mathrm{d}z \mathrm{d} \theta . \notag
\end{align}
We shall only prove estimate \eqref{eq6.68v71}, as the proof of \eqref{eq6.69v71} is similar.
\begin{align*}
	\eqref{eq6.68v71}
	= & \sum\limits_{l_1 \le l_2} \sum\limits_{ \substack{ n_1,n_2,n_3,n \in \mathbb{N},\\
			n_1 - n_2 + n_3 = n} }
	\iiint K(t; \tilde{y}, z,  \theta )
	\Pi_n \Big( e^{-i (y - \tilde{y} ) \xi(t)}
	\left( P^y_{\xi(t), 2^{ l_1} }   v_{n_1}^l \right) \left( y - \tilde{y} , x \right) \\
	& \quad \cdot e^{i (y- \tilde{y} - z) \xi(t)} \overline{ \left( P^y_{\xi(t), \le 2^{  l_1} }   { v }_{n_2}^l \right)  \left(y - \tilde{y} - z, x \right) }
	v_{n_3}^h \left(y - \tilde{y} - z -  \theta , x \right) \Big ) \,\mathrm{d} \tilde{y} \mathrm{d}z \mathrm{d} \theta ,
\end{align*}
where
\begin{align*}
	K\left(t; \tilde{y}, z,  \theta \right)
	= \iiint e^{i \xi \tilde{y} + i \eta z + i \sigma  \theta } \left( \left( \phi \left( \frac{ \xi- \xi(t)}{ 2^{l_2}} \right) - \phi\left( \frac{\sigma - \xi(t)}{ 2^{l_2}} \right)  \right)
	\phi \left( \frac{\sigma - \xi(t)}{2^{l_2}}\right ) \psi_{l_1} ( \xi- \eta) \phi \left( \frac{\eta - \sigma}{2^{l_1}}\right ) \right) \,\mathrm{d} \xi \mathrm{d}\eta \mathrm{d} \sigma.
\end{align*}
By the estimates
$|\xi- \eta| \sim 2^{l_1} , | \eta - \sigma| \lesssim 2^{l_1},$ and the fundamental theorem of calculus, we obtain
\begin{align}\label{eq6.14v51}
	\left|\phi \left( \frac{\xi -\xi(t)}{2^{l_2}} \right) - \phi\left( \frac{\sigma - \xi(t)}{2^{l_2}} \right) \right|
	\lesssim 2^{- l_2} |\xi - \sigma | \lesssim 2^{l_1 - l_2}.
\end{align}
This implies
\begin{align}\label{eq6.15v51}
	\sup\limits_t \int \left|K \left(t;\tilde{y}, z,  \theta \right) \right| \,\mathrm{d} \tilde{y} \mathrm{d}z \mathrm{d} \theta  \lesssim 2^{l_1 - l_2}.
\end{align}
Thus, by Minkowski's inequality, H\"older's inequality, \eqref{eq6.15v51}, Lemma \ref{le5.7v51} and the conservation of mass, we infer
\begin{align*}
	& \bigg\| P^y_{\xi(t), \le 2^{ l_2} }  \Pi_n \bigg( \sum\limits_{ \substack{ n_1,n_2,n_3,n \in \mathbb{N},\\
			n_1-n_2+n_3 =n }  }
	v_{n_1}^l\overline{ { v }_{n_2}^l}
	v_{n_3}^h \bigg)
	- \sum\limits_{ \substack{ n_1,n_2,n_3,n \in \mathbb{N}, \\ n_1 -n_2 + n_3 =n } } \Pi_n \left( P^y_{\xi(t), \le  2^{ l_2}  } v_{n_1}^l
	\overline{ P^y_{\xi(t), \le 2^{  l_2} }   { v }_{n_2}^l }
	P^y_{\xi(t), \le 2^{  l_2} }   v_{n_3}^h \right) \bigg\|_{L_{t,y}^\frac43 L_x^2 \left(G_\beta^{l_2} \times \mathbb{R}^2 \times \mathbb{R} \right)} \\
	\lesssim & \sum\limits_{l_1 \le l_2}  \bigg\| \iiint  \left|K \left(t; \tilde{y}, z, \theta \right) \right|
	\bigg\| \left\| \left(  P^y_{\xi(t), 2^{  l_1} }  v^l \right) \left(y - \tilde{y}, x \right) \right\|_{L_x^2} \left\| \left( P^y_{\xi(t), \le  2^{ l_1} }   v^l \right) \left(y - \tilde{y} - z, x \right) \right\|_{L_x^2} \\
	& \qquad \cdot \left\|  v^h \left(y - \tilde{y} - z -  \theta , x \right) \right\|_{L_x^2} \bigg\|_{L_y^\frac43 } \mathrm{d} \tilde{y} \mathrm{d}z \mathrm{d} \theta \bigg\|_{L_t^\frac43} \\
	\lesssim & \sum\limits_{l_1 \le   l_2} 2^{l_1 - l_2} \left\|P^y_{\xi(t), \le 2^{  l_1} }  v^l \right\|_{L_t^\infty L_{y,x}^2
	} \left\|P^y_{\xi(t),  2^{ l_1} }  u^l \right\|_{L_t^\frac83 L_y^8 L_x^2} \left\| v^h \right\|_{L_t^\frac83 L_y^8 L_x^2}
	\lesssim    \| v \|_{\tilde{X}_i \left(G_\alpha^i \right)}^2.
\end{align*}
We now consider \eqref{eq6.12v51}. Since
\begin{align*}
	& \bigg\| \mathcal{O} \bigg( \sum\limits_{ \substack{ n_1,n_2,n_3,n \in \mathbb{N},\\  n_1 -n_2 +n_3 =n}  }
	\Pi_n \left( P^y_{\xi(t), \le 2^{ l_2} }  \left( v_{n_1}^h \overline{  { v }_{n_2}^h }   v_{n_3} \right)
	- P^y_{\xi(t), \le 2^{  l_2} }  v_{n_1}^h \overline{ P^y_{\xi(t), \le 2^{  l_2}  } { v }_{n_2}^h }  P^y_{\xi(t), \le 2^{  l_2} }   v_{n_3 } \right) \bigg) \bigg\|_{L_{t,y}^\frac43 L_x^2 \left(G_\beta^{l_2} \times \mathbb{R}^2 \times \mathbb{R} \right)} \\
	\lesssim & \| v  \|_{L_t^\infty L_{y,x}^2} \| v^h \|_{L_t^\frac83 L_y^8 L_x^2}^2
	\lesssim \|  v  \|_{\tilde{X}_i \left(G_\alpha^i \right)}^2,
\end{align*}
it follows that
\begin{align}\label{eq6.17v51}
	\|N_1\|_{L_{t,y}^\frac43 L_x^2 \left(G_\beta^{l_2} \times \mathbb{R}^2 \times \mathbb{R} \right)} \lesssim
	\| v  \|_{\tilde{X}_i \left(G_\alpha^i \right)}^2,
\end{align}
and therefore, we have
\begin{align*}
	\eqref{eq6.3v51}
	\lesssim 2^{2l_2 - 2i } \|w_0\|_{L_{y,x}^2 (\mathbb{R}^2 \times \mathbb{R})}^2 \left( 1 + \|  v  \|_{\tilde{X}_i \left(G_\alpha^i \right)}^3 \right).
\end{align*}
\medskip
\textbf{Estimate of \eqref{eq6.4v51}.}
\medskip

Applying integration by parts, we have
\begin{align*}
	\eqref{eq6.4v51}\lesssim  \eqref{eq6.3v51}
	+ 2^{l_2 - 2i} \bigg| \int_{\mathbb{R}} \int_{\mathbb{R}} \int_{G_\beta^{l_2}} \iint | w (t, \tilde{y}, \tilde{x})|^2 \frac1{ |y - \tilde{y}|}
	\Re \left( \overline{ \tilde{ w } } N \right)(t,y,x) \,\mathrm{d}y \mathrm{d}\tilde{y} \mathrm{d}x \mathrm{d}\tilde{x}  \mathrm{d}t \bigg| .
\end{align*}
By the Strichartz estimate, \eqref{eq6.17v51}, \eqref{eq5.1v51}, \eqref{eq5.3v51}, Bernstein's inequality, and the conservation of mass, we have
\begin{align*}
	&  2^{l_2 - 2i} \bigg|\int_{G_\beta^{l_2}} \iint \iint |  w (t, \tilde{y}, \tilde{x})|^2 \frac1{ | y - \tilde{y}|}
	\Re\left( \overline{ \tilde{ w } } N \right)(t,y,x) \,\mathrm{d}y \mathrm{d}\tilde{y} \mathrm{d}x \mathrm{d}\tilde{x} \mathrm{d}t \bigg| \\
	\lesssim & 2^{l_2 - 2i} \| w_0 \|_{L_{y,x}^2 }^2 2^{- \frac{l_2}2} \int_{G_\beta^{l_2}} |\xi'(t)| \left\|P^y_{\xi(t), 2^{  l_2 - 3 }  \le \cdot \le  2^{ l_2 + 3} }
 v \right\|_{L_{y,x}^2 } \left\|(\nabla_y  - i \xi(t))^\frac12 P^y_{\xi(t), \le  2^{ l_2} }   v  \right\|_{L_{y,x}^2} \,\mathrm{d}t\\
	& \ + 2^{l_2 - 2i} \|w \|_{L_{t,y}^4 L_x^2 \left(G_\beta^{l_2} \times \mathbb{R}^2 \times \mathbb{R} \right)}^2 \|N_1\|_{L_{t,y}^\frac43 L_x^2 \left(G_\beta^{l_2} \times \mathbb{R}^2 \times \mathbb{R} \right)} \|  \tilde{ w } \|_{L_{t,y}^\infty L_x^2(G_\beta^{l_2} \times \mathbb{R}^2 \times \mathbb{R})}
	\lesssim   2^{2l_2 - 2i} \|  w_0\|_{L_{y,x}^2}^2 \left( 1 + \| v \|_{\tilde{X}_i \left(G_\alpha^i \right)}^3 \right).
\end{align*}
Thus
\begin{align*}
	\eqref{eq6.4v51} \lesssim 2^{2l_2 - 2i} \|  w_0\|_{L_{y,x}^2}^2 \left( 1 + \| v \|_{\tilde{X}_i \left(G_\alpha^i \right)}^3 \right).
\end{align*}
\medskip

\textbf{Estimate of \eqref{eq6.5v51}.}
\medskip

By Bernstein's inequality and the conservation of mass, we have
\begin{align*}
	& 2^{l_2 - 2i} \bigg | \int_{G_\beta^{l_2}}   \int_{\mathbb{R}^2} \int_{\mathbb{R}^2} \int_{\mathbb{R} }  \int_{\mathbb{R} }  \Im( \bar{ w }( \nabla_y  - i \xi(t) ) w )(t,y, x) \frac{ y - \tilde{y}}{ | y - \tilde{y}|} \Im(   \overline{ \tilde{ w} } N_2)(t, \tilde{y}, \tilde{x}) \,\mathrm{d}y \mathrm{d} \tilde{y} \mathrm{d}x \mathrm{d}\tilde{x} \mathrm{d}t \bigg|\\
	\lesssim & \| w_0\|_{L_{y,x}^2}^2  \left( 2^{-i - l_2} \int_{G_\beta^{l_2}}  |\xi'(t)| \left\|P^y_{\xi(t),  2^{ l_2 - 3 }  \le \cdot \le  2^{ l_2 + 3}  }  v  \right\|_{L_{y,x}^2}  \left\| ( \nabla_y  - i \xi(t)) P^y_{\xi(t), \le 2^{  l_2 } }  v  \right\|_{L_{y,x}^2 } \,\mathrm{d}t \right)\\
	& + \| w_0 \|_{L_{y,x}^2}^2 \left(2^{-i - \frac{l_2}2} \int_{G_\beta^{l_2}} |\xi'(t)| \left\|P^y_{\xi(t), l_2 - 3 \le \cdot \le l_2 +3 }  v  \right\|_{L_{y,x}^2}\left\| ( \nabla_y  - i \xi(t))^\frac12 P^y_{\xi(t), \le  2^{ l_2} }   v  \right\|_{L_{y,x}^2 } \,\mathrm{d}t \right)
	\lesssim \|  w_0 \|_{L_{y,x}^2 }^2.
\end{align*}
We now turn to the estimate of
\begin{align*}
	2^{l_2 - 2i} \bigg | \int_{G_\beta^{l_2}}   \int_{\mathbb{R}^2} \int_{\mathbb{R}^2} \int_{\mathbb{R} }  \int_{\mathbb{R} }  \Im( \bar{ w }( \nabla_y  - i \xi(t) ) w )(t,y, x) \frac{ y - \tilde{y}}{ | y - \tilde{y}|} \Im(   \overline{ \tilde{ w} } N_1)(t, \tilde{y}, \tilde{x}) \,\mathrm{d}y \mathrm{d} \tilde{y} \mathrm{d}x \mathrm{d}\tilde{x} \mathrm{d}t \bigg|.
\end{align*}
Since
\begin{align}\label{eq6.24v51}
	\int_{\mathbb{R}} \Im \bigg( \sum\limits_{ \substack{ n_1,n_2,n_3,n \in \mathbb{N}, \\ n_1 - n_2+ n_3 = n} }
	\overline{ \tilde{w} }  \Pi_n \left(  \tilde{w}_{n_1} \overline{ { \tilde{ w } }_{n_2}}
	\tilde{w}_{n_3} \right) \bigg)(\tilde{x})  \,\mathrm{d}  \tilde{x}
	= 0,
\end{align}
we see
\begin{align*}
	\int_{\mathbb{R}}
	\sum\limits_{n \in \mathbb{N}} \Im \left( \overline{ \tilde{{w}}_{n  }}  N_{1,n} \right)(\tilde{x}) \,\mathrm{d} \tilde{x}
	= &  \int_{\mathbb{R}}
	\Im \bigg( \sum\limits_{n \in \mathbb{N}}  \overline{ { \tilde{ w}}_n }  P^y_{\xi(t), \le 2^{ l_2} } \sum\limits_{ \substack{ n_1,n_2,n_3 \in \mathbb{N}, \\
			n_1 -n_2 + n_3 =n }} \Pi_n \left( v_{n_1} \overline{ { v }_{n_2}  } v_{n_3} \right)  \bigg)(\tilde{x}) \,\mathrm{d} \tilde{x} .
\end{align*}
Using the decomposition $ v  =  v^h +  v^l$, where $ v^l = P^y_{\xi(t), \le  2^{ l_2 - 5} }   v $, together with the above equality, we have
\begin{align}\label{eq6.25v51}
	\int_{\mathbb{R}} \sum\limits_{ n \in \mathbb{N}}  \Im \left(\overline{{ \tilde{w} }_n } N_{1,n} \right)(\tilde{x})
	= &
	\int_{\mathbb{R}}
	\sum\limits_{ n \in \mathbb{N}}  \left(F_{0,n}  + F_{1,n}  + F_{2,n}  + F_{3,n}  + F_{4,n}  \right)(\tilde{x}) \,\mathrm{d} \tilde{x}
	,
\end{align}
where $F_{j ,n}$ consists of $j $ $ v_n^h-$terms and $4-j $ $ v_n^l-$terms, for $j  = 0, 1, 2, 3, 4$, in
\begin{align*}
	\Im \bigg( \overline{ { \tilde{ w } }_n  }
	P^y_{\xi(t), \le  2^{l_2}}
 \sum\limits_{ \substack{ n_1,n_2,n_3 \in \mathbb{N},\\
			n_1 -n_2 + n_3 = n}  } \Pi_n \left(  v_{n_1} \overline{ { v }_{n_2} }
	v_{n_3} \right) \bigg).
\end{align*}
We now consider the estimate of the $F_{j}$ terms, $j = 0,1,2,3,4$ as follows.

By \eqref{eq6.24v51}, we have
\begin{align*}
	\int_{\mathbb{R}} \sum\limits_{n \in \mathbb{N}}  F_{0,n} \left(t,\tilde{y}, \tilde{x} \right) \,\mathrm{d} \tilde{x}
	= 0.
\end{align*}
By Bernstein's inequality, \eqref{eq3.8v41} and Lemma \ref{le5.7v51}, we have
\begin{align*}
	& 2^{l_2 - 2i} \bigg|\int_{G_\beta^{l_2}} \int_{\mathbb{R}} \int_{\mathbb{R}}  \int_{\mathbb{R}^2} \int_{\mathbb{R}^2} \Im \left( \bar{ w } ( \nabla_y  - i \xi(t))  w  \right)(t,y,x) \frac{ y - \tilde{y}}{ |y - \tilde{y}|} (F_3 + F_4)\left(t,\tilde{y}, \tilde{x} \right) \,\mathrm{d}y \mathrm{d}\tilde{y} \mathrm{d}x \mathrm{d}\tilde{x} \mathrm{d}t \bigg| \\
	\lesssim & 2^{l_2 - i} \| w_0\|_{L_{y,x}^2 }^2  \left\| v^h \right\|_{L_t^3 L_y^6 L_x^2 \left(G_\beta^{l_2} \times \mathbb{R}^2 \times \mathbb{R} \right)}^3 \| v \|_{L_t^\infty L_{y,x}^2 \left(G_\beta^{l_2} \times \mathbb{R}^2 \times \mathbb{R} \right)}
	\lesssim   2^{l_2 - i} \| w_0 \|_{L_{y,x}^2 }^2 \| v \|_{\tilde{X}_i \left(G_\alpha^i \right)}^3.
\end{align*}
By a direct calculation, we have
\begin{align} \label{eq6.27v51}
	\sum\limits_{ n \in \mathbb{N}}
	F_{1,n}= &  \sum\limits_{ n \in \mathbb{N}}  \Im \bigg(  \overline{ P^y_{\xi(t),\le 2^{ l_2} }  P^y_{\xi(t), \ge  2^{ l_2 - 2} }   { v }_n^h }  P^y_{\xi(t), \le  2^{ l_2} }  \sum\limits_{ \substack{n_1,n_2,n_3 \in \mathbb{N}, \\
			n_1 -n_2 + n_3 =n } } \Pi_n \left(  v_{n_1}^l \overline{ { v }_{n_2}^l }   v_{n_3}^l \right)   \\
	& + \overline{ { v }_n^l } P^y_{\xi(t), \le 2^{ l_2}  } \bigg( \sum\limits_{ \substack{ n_1,n_2,n_3 \in \mathbb{N}, \\
			n_1 - n_2 + n_3 =n} }  \Pi_n \left(  v_{n_1}^l\overline{  P^y_{\xi(t), \ge  2^{ l_2 - 2} }  { v }_{n_2}^h }  v_{n_3}^l \right)
	+ 2 \sum\limits_{ \substack{ n_1,n_2,n_3 \in \mathbb{N},\\  n_1 - n_2 +n_3 =n} }
	\Pi_n \left(  v_{n_1}^l \overline{ { v }_{n_2}^l }
	P^y_{\xi(t), \ge 2^{  l_2 - 2} }  v_{n_3}^h \right) \bigg) \bigg).  \notag
\end{align}
Since the support of the partial Fourier transform with respect to $\tilde{y}$ of $\sum\limits_{ n \in \mathbb{N}}  F_{1,n} \left(t, \tilde{y}, \tilde{x} \right)$ is contained in $\left\{ \xi:  |\xi| \ge 2^{l_2 - 4}  \right\}$, we can apply the integration by parts with respect to $\tilde{y}$, the Hardy-Littlewood-Sobolev inequality, Bernstein's inequality, the Strichartz estimate, \eqref{eq5.11v51}, and \eqref{eq5.10v51} to give the following estimate
\begin{align}
	& 2^{l_2 - 2i} \bigg| \sum\limits_{n \in \mathbb{N}} \int_{G_\beta^{l_2}} \int_{\mathbb{R}} \int_{\mathbb{R}} \int_{\mathbb{R}^2} \int_{\mathbb{R}^2} \Im \left( \overline{{ w }_n }
	( \nabla_y  - i \xi(t)) w_n \right)  (t,y,x) \frac{ y - \tilde{y}}{ | y - \tilde{y} |}
	\Big( \sum\limits_{n' \in \mathbb{N} } F_{1,n'}  \Big)(t, \tilde{y}, \tilde{x}) \,\mathrm{d} y \mathrm{d} \tilde{y}  \mathrm{d}x \mathrm{d} \tilde{x} \mathrm{d}t\bigg| \notag \\
	\lesssim & 2^{l_2 - 2i} \int_{G_\beta^{l_2}} \int_{\mathbb{R}^2} \int_{\mathbb{R}^2} \int_{\mathbb{R}}  \left| \sum\limits_{n \in \mathbb{N}}  (\overline{{w}_n }( \nabla_y  - i\xi(t)) w_n)(t,y,x) \right|
	\frac1{ |y - \tilde{y}|} \left|{\partial_{\tilde{y}}} \, {\left( -  \Delta_{\tilde{y}} \right)^{- 1} } \int_{\mathbb{R}}
	\left( \sum\limits_{n' \in \mathbb{N} } F_{1,n'} \right) \left(t, \tilde{y}, \tilde{x} \right)\mathrm{d} \tilde{x}
	\right| \,\mathrm{d} y \mathrm{d}\tilde{y} \mathrm{d}x  \mathrm{d}t \notag \\
	\lesssim & 2^{l_2 - 2i } \| w  \|_{L_t^3 L_y^6 L_x^2 \left(G_\beta^{l_2} \times \mathbb{R}^2 \times \mathbb{R} \right)}
	\left\| \left( \nabla_y  - i \xi(t) \right) w  \right\|_{L_t^\infty L_{y,x}^2
		\left(G_\beta^{l_2} \times \mathbb{R}^2 \times \mathbb{R} \right)}
	\left\|{ \partial_{\tilde{y}}}\, { \left( -  \Delta_{\tilde{y}} \right)^{- 1} } \int_{\mathbb{R}} \left( \sum\limits_{n' \in \mathbb{N} } F_{1,n'} \right) \,\mathrm{d}\tilde{x} \right\|_{L_t^\frac32 L_{\tilde{y}}^\frac65 \left(G_\beta^{l_2} \times \mathbb{R}^2 \right)} \notag \\
	\lesssim & 2^{-i} \| w_0 \|_{L_{y,x}^2 }^2 \| v^l \|_{L_t^9 L_y^\frac92 L_x^2 \left(G_\beta^{l_2} \times \mathbb{R}^2 \times \mathbb{R} \right)}^3 \|  v^h \|_{L_t^3 L_y^6 L_x^2 \left(G_\beta^{l_2} \times \mathbb{R}^2 \times \mathbb{R} \right)}
	\lesssim 2^{l_2 - i} \| w_0 \|_{L_{y,x}^2 }^2 \| v \|_{\tilde{X}_i \left(G_\alpha^i \right)}^4.
	\notag
\end{align}
We are now left to show
\begin{align}\label{eq6.31v51}
	& 2^{l_2 - 2i} \bigg| \sum\limits_{n,n'\in \mathbb{N}} \int_{G_\beta^{l_2}}  \int_{\mathbb{R}} \int_{\mathbb{R}} \int_{\mathbb{R}^2} \int_{\mathbb{R}^2} \Im \left( \overline{ {w}_n  }
	( \nabla_y  - i \xi(t)) w_n \right)(t,y,x) \frac{ y - \tilde{y}} {|y - \tilde{y}|} F_{2,n'} \left(t,\tilde{y}, \tilde{x} \right) \,\mathrm{d}y \mathrm{d} \tilde{y} \mathrm{d} x \mathrm{d} \tilde{x} \mathrm{d}t \bigg| \\
	\lesssim & \| w_0 \|_{L_{y,x}^2 }^2 \left( 1 +  \|  v \|_{\tilde{X}_i \left(G_\alpha^i \right)}^4 \right). \notag
\end{align}
Similar to the estimate on the term involved $F_{1 }$ above, from integration by parts, Bernstein's inequality and \eqref{eq5.11v51}, we conclude
\begin{align*}
	& 2^{l_2- 2i} \bigg| \sum\limits_{n \in \mathbb{N}} \int_{G_\beta^{l_2}}  \int_{\mathbb{R}^2} \int_{\mathbb{R}^2} \int_{\mathbb{R}}
	\Im ( \overline{{ w}_n }
	( \nabla_y  - i \xi(t)) w_n)(t,y,x) \frac{ y - \tilde{y}}{ | y - \tilde{y}|}
	\int_{\mathbb{R}} \left( \sum\limits_{n' \in \mathbb{N} } P^y_{\ge l_2 - 10}   F_{2,n'} \right) \left(t, \tilde{y}, \tilde{x} \right) \, \mathrm{d}x \mathrm{d} \tilde{x}
	\mathrm{d}y  \mathrm{d} \tilde{y} \mathrm{d}t \bigg| \\
	\lesssim & 2^{l_2 - 2i} \| w \|_{L_t^3L_y^6 L_x^2 \left(G_\beta^{l_2} \times \mathbb{R}^2 \times \mathbb{R} \right)} \left\|( \nabla_y  - i \xi(t) )  w \right\|_{L_t^\infty L_{y,x}^2  \left(G_\beta^{l_2} \times \mathbb{R}^2 \times \mathbb{R} \right)}
	\left\| {\partial_{\tilde{y}}}\, {  \left( - \Delta_{\tilde{y}} \right)^{- 1} } \left( \int_{\mathbb{R}}  \sum\limits_{n' \in\mathbb{N} } P^y_{\ge l_2 - 10} F_{2,n'} \,\mathrm{d} \tilde{x} \right) \right\|_{L_t^\frac32 L_y^\frac65  \left(G_\beta^{l_2} \times \mathbb{R}^2 \right)}\\
	\lesssim & 2^{-i} \| w_0 \|_{L_{y,x}^2 }^2 \left\| v^l \right\|_{L_t^\infty L_y^4 L_x^2 \left(G_\beta^{l_2} \times \mathbb{R}^2 \times \mathbb{R} \right)}^2 \left\| v^h \right\|_{L_t^3 L_y^6 L_x^2 \left(G_\beta^{l_2} \times \mathbb{R}^2 \times \mathbb{R} \right)}^2
	\lesssim   2^{l_2 - i} \| w_0 \|_{L_{y,x}^2 }^2 \|  v  \|_{\tilde{X}_i \left(G_\alpha^i \right)}^2.
\end{align*}
We now turn to the estimate of the low frequency part of $F_2$. First of all, we can decompose $F_{2,n'}$ as
\begin{align}\label{eq6.33v51}
	& F_{2,n'} \left(t, \tilde{y}, \tilde{x} \right) \\
	= &  \notag \\
	&   \Im \bigg( 2 \overline{ { v }_{n'}^l }  P^y_{\xi(t), \le 2^{  l_2} } \bigg( \sum\limits_{ \substack{ n'_1,n'_2,n'_3 \in \mathbb{N}, \\
			n'_1 -n'_2 + n'_3 =n'} }  \Pi_{n'} \left( v_{n_1'}^h \overline{ { v }_{n_2'}^h }
	v_{n_3'}^l \right) \bigg)
	+\overline{ { v }_{n'}^l } P^y_{\xi(t), \le 2^{ l_2}}  \bigg( \sum\limits_{ \substack{ n_1',n_2',n_3' \in \mathbb{N}, \\
			n'_1 - n'_2 + n'_3 =n' } }
	\Pi_{n'}
	\left(  v_{n_1'}^h \overline{ { v }_{n_2'}^l }   v_{n'_3}^h \right) \bigg) \notag \\
	&    + 2 \overline{ P^y_{\xi(t), \le 2^{  l_2}}    { v }_{n'}^h }
	P^y_{\xi(t), \le 2^{  l_2}}  \sum\limits_{ \substack{ n'_1,n'_2,n'_3 \in \mathbb{N}, \\
			n'_1 -n'_2 + n'_3 = n'}  }\Pi_{n'}  \left(  v_{n_1'}^l \overline{ { v }_{n'_2}^l }
	v_{n'_3}^h \right)
	+ \overline{ P^y_{\xi(t), \le 2^{  l_2} }  { v }_{n'}^h }
	P^y_{\xi(t), \le  2^{ l_2} }  \sum\limits_{ \substack{ n'_1,n'_2,n'_3 \in \mathbb{N},\\
			n'_1 - n'_2 +n'_3 = n' } }
	\Pi_{n'}
	\left(  v_{n'_1}^l \overline{ { v }_{n'_2}^h }   v_{n'_3}^l \right) \bigg). \notag
\end{align}
Since
\begin{align*}
	\Im \bigg( \sum\limits_{n' \in \mathbb{N} } \sum\limits_{ \substack{ n_1',n_2',n_3' \in \mathbb{N}, \\
			n_1' - n_2' + n_3' = n'} }
	\overline{  P^y_{\xi(t), \le  2^{ l_2} }   { v }_{n'}^h  } \Pi_{n'} \left(  v_{n_1'}^l  \overline{ v_{n_2'}^l }
	P^y_{\xi(t), \le  2^{ l_2} }   v_{n_3'}^h \right) \bigg)= 0,
\end{align*}	
and
\begin{align*}	
	\Im \bigg(  \sum\limits_{n' \in \mathbb{N} } \sum\limits_{ \substack{ n_1',n_2',n_3' \in \mathbb{N},\\
			n_1' -n_2' + n_3' = n'} }  \overline{ { v }_{n'}^l }
	\Pi_{ n'}
	\left(   P^y_{\xi(t), \le 2^{  l_2} }  \left(  v_{n_1'}^h \overline{ { v }_{n_2'}^h  } \right)v_{n_3'}^l \right)  \bigg) = 0,
\end{align*}
we obtain
\begin{align}
	& 2 \Im \bigg( \sum\limits_{n' \in \mathbb{N} } \overline{ { v }_{n'}^l }
	P^y_{\xi(t), \le 2^{  l_2 } } \bigg( \sum\limits_{ \substack{ n_1',n_2',n_3' \in \mathbb{N}, \\
			n_1' - n_2' + n_3' = n'} } \Pi_{n'}
	\left(  v_{n_1'}^h \overline{ { v }_{n_2'}^h }  v_{n_3'}^l \right) \bigg)
	+ \sum\limits_{n' \in \mathbb{N} } \overline{ P^y_{\xi(t), \le 2^{  l_2} }   { v }_{n'}^h }
	P^y_{\xi(t), \le  2^{ l_2} }
 \sum\limits_{ \substack{ n_1',n_2',n_3' \in \mathbb{N} , \\
			n_1' -n_2'+n_3' = n'} }  \Pi_{n'} \left(  v_{n_1'}^l \overline{ {v}_{n_2'}^l }   v_{n_3'}^h \right) \bigg)  \notag \\
	= & \notag \\
	\quad &  2\Im \bigg( \sum\limits_{n' \in \mathbb{N} } \sum\limits_{ \substack{ n_1',n_2',n_3' \in \mathbb{N},\\
			n_1' -n_2' + n_3' =n' }
	} \left(  \overline{ P^y_{\xi(t), \le 2^{  l_2} }   { v }_{n'}^h }    P^y_{\xi(t), \le  2^{ l_2} }  \Pi_{n'} \left(  v_{n_1'}^l \overline{ { v }_{n_2'}^l }   v_{n_3'}^h   \right)
	- \overline{ P^y_{\xi(t), \le  2^{ l_2} }   { v }_{n'}^h } \Pi_{n'} \left(  v_{n_1'}^l \overline{ { v }_{n_2'}^l }
	P^y_{\xi(t), \le  2^{ l_2} }   v_{n_3'}^h  \right)   \right)  \bigg)  \label{eq6.81v71}\\
	& +   2 \Im \bigg( \sum\limits_{n' \in \mathbb{N} } \sum\limits_{ \substack{ n_1',n_2',n_3' \in \mathbb{N}, \\
			n_1' -n_2' +n_3' =n' } } \left( \overline{ { v }_{n'}^l } P^y_{\xi(t), \le 2^{  l_2} }  \Pi_{n'} \left(  v_{n_1'}^h \overline{ { v }_{n_2'}^h }  v_{n_3'}^l \right)
	- \overline{ { v }_{n'}^l } \Pi_{n'} \left(   P^y_{\xi(t), \le 2^{  l_2} }  \left(  v_{n_1'}^h \overline{ { v }_{n_2'}^h } \right) v_{n_3'}^l\right)  \right) \bigg)
	\label{eq6.80v71}.
\end{align}
For \eqref{eq6.81v71}, by \eqref{eq6.13v51}, \eqref{eq6.15v51}, Lemma \ref{le5.7v51} and the conservation of mass, we have
\begin{align}\label{eq6.35v51}
	& 2^{l_2 - i} \left\| \int_{\mathbb{R}} \eqref{eq6.81v71} (x) \,\mathrm{d}x \right\|_{L_{t,y}^1 \left(G_\beta^{l_2} \times \mathbb{R}^2 \right)}\\
	\lesssim  & 2^{l_2 - i}  \left\|  v^h \right\|_{L_{t,y}^4 L_x^2 \left(G_\beta^{l_2} \times \mathbb{R}^2 \times \mathbb{R} \right)} \notag \\
& \qquad  \cdot \bigg\| \sum\limits_{ \substack{ n_1,n_2,n_3,n \in \mathbb{N},\\
			n_1 -n_2 +n_3 =n } } \Pi_n P^y_{\xi(t), \le  2^{ l_2} } \left (  v_{n_1}^l \overline{{ v }_{n_2}^l}   v_{n_3}^h  \right)
	- \sum\limits_{ \substack{ n_1,n_2,n_3 , n \in \mathbb{N}, \\
			n_1 -n_2 +n_3 =n}  }\Pi_n \left(  v_{n_1}^l\overline{{ v }_{n_2}^l }
	P^y_{\xi(t), \le 2^{  l_2 } }  v_{n_3}^h \right)  \bigg\|_{L_{t,y}^\frac43 L_x^2 \left(G_\beta^{l_2} \times \mathbb{R}^2 \times \mathbb{R} \right)} \notag  \\
	\lesssim & 2^{l_2 - i} \| v \|_{\tilde{X}_i(G_\alpha^i)} \sum\limits_{l_1 \le l_2} 2^{l_1 - l_2} \left\|P^y_{\xi(t), \le 2^{  l_1}  } v \right\|_{L_t^\infty L_{y,x}^2}
	\left\|P^y_{\xi(t),  2^{ l_1} }   v^l \right\|_{L_t^\frac83 L_y^8 L_x^2} \left\| v^h \right\|_{L_t^\frac83 L_y^8 L_x^2}
	\lesssim 2^{l_2 - i} \| v \|_{\tilde{X}_i \left(G_\alpha^i \right)}^3. \notag
\end{align}
To estimate \eqref{eq6.80v71}. We note similar to \eqref{eq6.14v51}, we have
\begin{align*}
	\left| \phi \left( \frac{ \xi_1 + \xi_2 - \xi(t)}{2^{l_2}} \right) - \phi\left( \frac{\xi_1}{2^{l_2}} \right) \right| \lesssim 2^{-l_2} |\xi_2 - \xi(t)|.
\end{align*}
Then by Lemma \ref{le5.7v51} and the conservation of mass, we have
\begin{align*}
	& 2^{l_2 - i} \left\| \int_{\mathbb{R}} \eqref{eq6.80v71} \,\mathrm{d}x \right\|_{L_{t,y}^1 \left(G_\beta^{l_2} \times \mathbb{R}^2 \right)} \\
	\lesssim & 2^{l_2 - i} \left\| v^l \right\|_{L_t^\infty L_{y,x}^2  \left(G_\beta^{l_2} \times \mathbb{R}^2 \times \mathbb{R}  \right)} \\
	& \qquad \bigg\| \sum\limits_{ \substack{ n_1,n_2,n_3,n \in \mathbb{N},\\
			n_1 -n_2 + n_3 =n } } P^y_{\xi(t), \le 2^{ l_2} } \Pi_n \left(  v_{n_1}^h \overline{{ v }_{n_2}^h }
	v_{n_3}^l \right)
	- \sum\limits_{ \substack{ n_1,n_2,n_3,n\in \mathbb{N},\\
			n_1 -n_2 + n_3 =n} }  \Pi_n \left( P^y_{\xi(t),
		\le 2^{  l_2} }  \left(  v_{n_1}^h\overline{ {  v }_{n_2}^h  }
	\right)  v_{n_3}^l \right) \bigg\|_{L_t^1 L_{y,x}^2 \left(G_\beta^{l_2} \times \mathbb{R}^2 \times \mathbb{R} \right)} \\
	\lesssim & 2^{-i} \left\|  v^h \right\|_{L_t^3 L_y^6 L_x^2 \left(G_\beta^{l_2} \times \mathbb{R}^2 \times \mathbb{R} \right)}^2 \sum\limits_{l_1 \le l_2} 2^{l_1} \left\|P^y_{\xi(t),  2^{ l_1} }  v^l \right\|_{L_t^3 L_y^6 L_x^2 \left(G_\beta^{l_2} \times \mathbb{R}^2 \times \mathbb{R} \right)}
	\lesssim 2^{l_2 -i} \|  v \|_{\tilde{X}_i \left(G_\alpha^i \right)}^3.
\end{align*}
Now we turn to the remaining terms in \eqref{eq6.33v51}. Observe that
\begin{align}
	&  \Im \bigg( \sum\limits_{n' \in \mathbb{N} } \sum\limits_{ \substack{ n_1',n_2',n_3' \in \mathbb{N}, \\
			n_1' - n_2' + n_3' =n'} } \overline{{ v }_{n'}^l }
	P^y_{\xi(t), \le  2^{ l_2} } \Pi_{n'} \left(  v_{n_1'}^h \overline{ { v }_{n_2'}^l }
	v_{n_3'}^h  \right)
	+ \sum\limits_{n' \in \mathbb{N} } \sum\limits_{ \substack{ n_1' ,n_2',n_3' \in \mathbb{N},\\
			n_1' - n_2' + n_3' =n'}  }
	\overline{ P^y_{\xi(t), \le  2^{ l_2} }  { v }_{n'}^h }
	P^y_{\xi(t), \le 2^{ l_2} }  \Pi_{n'} \left(  v_{n_1'}^l \overline{ { v }_{n_2'}^h }
	v_{n_3'}^l \right)  \bigg) \notag \\
	= & \notag\\
	\quad & \sum\limits_{n' \in \mathbb{N} } \sum\limits_{ \substack{ n_1',n_2' ,n_3' \in \mathbb{N},\\
			n_1' -n_2' + n_3' =n'} }
	\Im \left( \overline{ P^y_{\xi(t), \le 2^{  l_2} }   { v }_{n'}^h }
	P^y_{\xi(t), \le  2^{ l_2} } \Pi_{n'} \left( v_{n_1'}^l \overline{ { v }_{n_2'}^h }
	v_{n_3'}^l \right)
	- \overline{ P^y_{\xi(t), \le 2^{  l_2}  } { v }_{n'}^h  }
	\Pi_{n'} \left(  v_{n_1'}^l \overline{  P^y_{\xi(t), \le  2^{ l_2} }   { v }_{n_2'}^h }   v_{n_3'}^l \right) \right) \label{eq6.85v71}\\
	& + \sum\limits_{n' \in \mathbb{N} } \sum\limits_{ \substack{ n_1',n_2',n_3' \in \mathbb{N}, \\
			n_1' -n_2' +n_3' =n'}  }
	\Im \left( \overline{ { v }_{n'}^l  } P^y_{\xi(t), \le  2^{ l_2} }  \Pi_{n'} \left(  v_{n_1'}^h \overline{{ v }_{n_2'}^l }  v_{n_3'}^l \right)
	+ \overline{ P^y_{\xi(t), \le 2^{  l_2}  }  { v }_{n'}^h }
	\Pi_{n'} \left(  v_{n_1'}^l P^y_{\xi(t), \le  2^{ l_2} }  \overline{ { v }_{n_2'}^h }
	v_{n_3'}^l \right) \right) .  \label{eq6.84v71}
\end{align}
Similar to the arguments for \eqref{eq6.35v51}, we have
\begin{align*}
	2^{l_2 - i} \left\| \int_{\mathbb{R}} \eqref{eq6.85v71} \,\mathrm{d}x \right\|_{L_{t,y}^1 \left(G_\beta^{l_2} \times \mathbb{R}^2 \right)}
	\lesssim 2^{l_2 - i} \|v\|_{\tilde{X}_i \left(G_\alpha^i \right)}^3.
\end{align*}
Thus, to show \eqref{eq6.31v51}, we just need to consider the term that contains \eqref{eq6.84v71}. By direct calculation, we get
\begin{align}
	& 2^{l_2 - 2i} \int_{G_\beta^{l_2}} \iiint \sum\limits_{n \in \mathbb{N}}
	\Im ( \overline{ { w}_n  }
	( \nabla_y  - i \xi(t)) w_n)(t,y,x) \frac{ y - \tilde{y}}{ | y - \tilde{y}|}
	\sum\limits_{n' \in \mathbb{N} } \sum\limits_{ \substack{ n_1',n_2',n_3' \in  \mathbb{N},\\
			n_1' -n_2' + n_3' =n'}  }
	P^y_{\xi(t), \le 2^{  l_2 - 10} } \Im \Big( \overline{ { v }_{n'}^l }
	P^y_{\xi(t), \le  2^{ l_2} }  \Pi_{n'} \left(  v_{n_1'}^h \overline{{ v }_{n_2'}^l }
	v_{n_3'}^h \right) \notag \\
	& + \overline{ P^y_{\xi(t), \le  2^{ l_2 } }  { v }_{n'}^h }
	\Pi_{n'} \left(  v_{n_1'}^l \overline{ P^y_{\xi(t), \le  2^{ l_2} }   { v }_{n_2'}^h }
	v_{n_3'}^l \right) \Big)(t,\tilde{y}, \tilde{x}) \,\mathrm{d}y \mathrm{d} \tilde{y} \mathrm{d}x \mathrm{d} \tilde{x} \mathrm{d}t  \notag \\
	= &  2^{l_2 - 2i} \int_{G_\beta^{l_2}} \iint \sum\limits_{ n \in \mathbb{N}}
	\Im \left( \overline{{ w }_n }
	( \nabla_y - i \xi(t)) w_n \right) \left( t, y + 2 \xi \left(G_\beta^{l_2} \right)t , x \right)
	\frac{y - \tilde{y}}{ | y - \tilde{y}|} \sum\limits_{n' \in \mathbb{N} } \sum\limits_{ \substack{ n_1',n_2',n_3' \in \mathbb{N} ,\\
			n_1' -n_2' + n_3' = n' } } P^y_{\xi(t),   \le  2^{ l_2 - 10} }  \notag  \\
	& \quad \Im \left( \overline{ { v }_{n_1'}^l } \overline{   { v }_{n_3'}^l }
	  \left(  v_{n'}^h  v_{n_2'}^h- P^y_{\xi(t), \le  2^{ l_2} }   v_{n'}^h P^y_{\xi(t), \le 2^{  l_2}  }  v_{n_2'}^h \right) \right) \left(t, \tilde{y} + 2 \xi(G_\beta^{l_2})t,\tilde{x} \right) \,\mathrm{d} {x} \mathrm{d} \tilde{ x}  \mathrm{d}y \mathrm{d} \tilde{y} \mathrm{d}t .  \notag
\end{align}
We may take $\xi(G_\beta^{l_2}) = 0$ in the right hand side of the above equality by the invariance of the Galilean transformation.
 By the inverse Fourier transform, we have
\begin{align*}
	& \sum\limits_{n,n' \in \mathbb{N} } \sum\limits_{ \substack{ n_1',n_2',n_3' \in \mathbb{N},\\
			n_1' -n_2' + n_3' =n'} }
	\int_{G_\beta^{l_2}} \iiiint \Im \left( \overline{ {w }_n } ( \nabla_y  - i \xi(t)) w_n \right)(t,y,x) \frac{ y - \tilde{y}}{ |y - \tilde{y}|}\\
	& \quad P^y_{\le 2^{  l_2 - 10} }  \Im \left( \left( \overline{ {v}_{n_1'} }
	v_{n_3'}^l \right)\left(  v_{n'}^h  v_{n_2'}^h - P^y_{\xi(t), \le  2^{ l_2}  }  v_{n'}^h P^y_{\xi(t), \le 2^{  l_2}  }  v_{n_2'}^h \right ) \right) \left(t, \tilde{y}, \tilde{x} \right) \,\mathrm{d}y \mathrm{d} \tilde{y} \mathrm{d}x \mathrm{d} \tilde{x} \mathrm{d}t \\
	= & \sum\limits_{n,n' \in \mathbb{N} } \sum\limits_{ \substack{ n_1',n_2', n_3' \in \mathbb{N}, \\
			n_1' -n_2' +n_3' =n'} }
	\int_{G_\beta^{l_2}} \iiiint \Im \left( \overline{ { w }_n }
	( \nabla_y  - i \xi(t))  w_n \right)(t,y,x)
	\frac{ y - \tilde{y}}{ | y - \tilde{y}|} \bigg(\iiint \phi \left( \frac{ \eta_1 + \eta_2 + \eta_3 + \eta_4}{ 2^{l_2 - 10}} \right)
	e^{i \tilde{y} \left( \eta_1 + \eta_2 + \eta_3 + \eta_4 \right)}\\
	& \qquad \Pi_{n'}
	\Big( \left( \mathcal{F}_{\tilde{y}}  v_{n_1'}^l \right) \left(t, \eta_3, \tilde{x} \right)
	\left( \mathcal{F}_{\tilde{y}}
	v_{n_3'}^l \right)(t, \eta_4, \tilde{x}) \left( \overline{\mathcal{F}_{\tilde{y}}  v_{n'}^h}  \right)(t,\eta_1, \tilde{x}) \left( \overline{\mathcal{F}_{\tilde{y}}  v_{n_2'}^h}  \right)(t , \eta_2 , \tilde{x}) \Big)\\
	& \qquad \left( 1 - \phi\left( \frac{\eta_1 - \xi(t)}{ 2^{l_2}} \right) \phi\left( \frac{\eta_2 - \xi(t)}{ 2^{l_2}} \right) \right) \,\mathrm{d}\eta_1 \mathrm{d} \eta_2 \mathrm{d} \eta_3 \mathrm{d} \eta_4 \bigg) \,\mathrm{d}y \mathrm{d}\tilde{y} \mathrm{d} x \mathrm{d}\tilde{x} \mathrm{d}t .
\end{align*}
Let
\begin{align*}
	q(\eta) = |\eta_1|^2 + |\eta_2 |^2 - |\eta_3|^2 - |\eta_4|^2,
\end{align*}
as in \cite{D1}, we have $\frac1{q(\eta)}$ is a convergent sum of terms with  operator norm being dominated by $  \frac1{ |\eta_1|^2 + |\eta_2|^2} \sim \frac1{|\eta_1| |\eta_2|}$ on the support of $\left(1 - \phi \left( \frac{\eta_1 - \xi(t)}{2^{l_2}}\right) \phi\left( \frac{\eta_2 - \xi(t)}{2^{l_2} } \right) \right) \phi\left( \frac{\eta_3}{ 2^{l_2 - 4}} \right) \phi\left( \frac{\eta_4}{ 2^{l_2 - 4}} \right)$.

Let $G_\beta^{l_2} = \left[ t_0 ,t_1 \right]$. Applying integration by parts (with respect to time), we have
\begin{align}
	& \int_{G_\beta^{l_2}} \iiint \iiint \iint \frac1{i q(\eta)} \left( \frac{d}{dt} e^{it q(\eta)} \right) \Im \left(\overline{ { w }_n }
	( \nabla_y  - i \xi(t))  w_n \right)(t,y,x)  \frac{y - \tilde{y}}{ \left| y - \tilde{y} \right|} \phi\left( \frac{ \eta_0+ \eta_1 + \eta_2 + \eta_3   }{2^{l_2 - 10}} \right) \notag \\
	& \quad e^{i \tilde{y} ( \eta_1 + \eta_2 + \eta_3 + \eta_0)}
	\left(1 - \phi\left(\frac{\eta_0 - \xi(t)}{2^{l_2}} \right)
	\phi\left( \frac{ \eta_2 - \xi(t)}{2^{l_2}} \right) \right)
	\bigg( e^{-it |\eta_0|^2 } \left( \overline{ \mathcal{F}_{\tilde{y}}  v_{n'}^h} \right)\left(t, \eta_0, \tilde{x} \right) e^{-it |\eta_2|^2 } \left( \overline{ \mathcal{F}_{\tilde{y}}  v_{n_2'}^h} \right)\left(t, \eta_2, \tilde{x} \right) \notag\\
	& \qquad e^{it |\eta_1|^2} \left( \mathcal{F}_{\tilde{y}}  v_{n_1'}^l \right)\left(t, \eta_1, \tilde{x} \right) e^{i t |\eta_3|^2} \left( \mathcal{F}_{\tilde{y}}  v_{n_3'}^l \right)\left(t, \eta_3, \tilde{x}\right) \,\mathrm{d}\eta_1 \mathrm{d}\eta_2 \mathrm{d}\eta_3 \mathrm{d} \eta_0 \bigg) \,\mathrm{d}{y} \mathrm{d}\tilde{y} \mathrm{d}x \mathrm{d}\tilde{x} \mathrm{d}t  \notag\\
	& := B_1 + B_2 + B_3 + B_4, & \notag
\end{align}
where
\begin{align}
B_1  = &  \iiint \iiint \iint \frac1{i q(\eta)} e^{it q(\eta)} \Im \left( \overline{ { w}_n }
	( \nabla_y  - i \xi(t)) w_n \right)(t,y,x) \frac{ y - \tilde{y}}{ \left| y - \tilde{y} \right|}  \phi \left( \frac{ \eta_1 + \eta_2 + \eta_3 + \eta_0}{2^{l_2- 10}} \right) e^{i \tilde{y} ( \eta_1 + \eta_2 + \eta_3 + \eta_0  )}  \label{eq6.47v51}\\
& \left( 1 - \phi\left( \frac{\eta_0 - \xi(t)}{2^{l_2}} \right) \phi\left( \frac{\eta_2 - \xi(t)}{2^{l_2}} \right) \right)
	\bigg( e^{-it |\eta_0|^2} \left( \overline{ \mathcal{F}_{\tilde{y}}  v_{n'}^h} \right)\left(t, \eta_0, \tilde{x} \right) e^{-it |\eta_2|^2} \left( \overline{ \mathcal{F}_{\tilde{y}}  v_{n_2'}^h }\right ) \left(t, \eta_2, \tilde{x} \right)\notag \\
	& \qquad e^{it |\eta_1|^2 } \left( \mathcal{F}_{\tilde{y}}  v_{n_1'}^l \right)\left(t, \eta_1, \tilde{x} \right) e^{it |\eta_3|^2 } \left( \mathcal{F}_{\tilde{y}}  v_{n_3'}^l \right)\left(t, \eta_3, \tilde{x} \right) \, \mathrm{d}\eta_0\mathrm{d}\eta_1 \mathrm{d} \eta_2 \mathrm{d} \eta_3
	\bigg) \,\mathrm{d}y \mathrm{d} \tilde{y} \mathrm{d} {x} \mathrm{d} \tilde{x} \bigg|_{t_0}^{t_1}, \notag
\end{align}
\begin{align}
B_2=	& - \int_{t_0}^{t_1}
	\iiint \iiint \iint \frac1{i q(\eta)} e^{it q(\eta)} \frac{\partial}{\partial t} \Im \left( \overline{ { w}_n  }
	( \nabla_y  - i \xi(t)) w_n \right)(t,y,x)
	\frac{ y - \tilde{y}}{  \left| y - \tilde{y} \right|} \phi\left( \frac{ \eta_1 + \eta_2 + \eta_3 + \eta_0 } { 2^{l_2- 10}} \right) \label{eq6.48v51} \\
&  e^{i \tilde{y} ( \eta_1 + \eta_2 + \eta_3 + \eta_0 )}
\left( 1 - \phi\left( \frac{\eta_0 - \xi(t)}{2^{l_2}} \right) \phi\left( \frac{\eta_2 - \xi(t)}{2^{l_2}} \right) \right)
	\bigg( e^{-it |\eta_0|^2} \left( \overline{ \mathcal{F}_{\tilde{y}}  v_{n'}^h} \right)\left(t, \eta_0, \tilde{x} \right) e^{-i t |\eta_2|^2 } \left( \overline{ \mathcal{F}_{\tilde{y}}  v_{n_2'}^h} \right) \left(t, \eta_2, \tilde{x} \right) \notag \\
	& \qquad e^{it |\eta_1|^2} \left(\mathcal{F}_{\tilde{y}}  v_{n_1'}^l \right) \left(t, \eta_1 , \tilde{x} \right) e^{it |\eta_3|^2 } \left(   \mathcal{F}_{\tilde{y}} v_{n_3'}^l \right) \left(t, \eta_3, \tilde{x} \right)
	\,\mathrm{d} \eta_1 \mathrm{d} \eta_2 \mathrm{d}\eta_3 \mathrm{d} \eta_0 \bigg) \,\mathrm{d} {y} \mathrm{d}\tilde{ y}
	\mathrm{d} x \mathrm{d} \tilde{x} \mathrm{d}t, \notag
\end{align}
\begin{align}
B_3 = 	& - \int_{t_0}^{t_1}
	\iiint \iiint \iint \frac1{i q(\eta)} e^{it q(\eta)} \Im \left( \overline{ { w}_n }
	( \nabla_y  - i \xi(t)) w_n \right)(t,y,x) \frac{ y - \tilde{y}} { | y - \tilde{y} |}
	\label{eq6.49v51}  \\
	& \frac{ \partial}{ \partial t}  \left( \phi \left( \frac{ \eta_1 + \eta_2 + \eta_3 + \eta_0} {2^{l_2 - 10}}  \right) e^{ i \tilde{y} ( \eta_1 + \eta_2 + \eta_3 + \eta_0 )}
	\left(1 - \phi \left( \frac{\eta_0 - \xi(t)} {2^{l_2}} \right) \phi \left( \frac{\eta_2 - \xi(t)}{2^{l_2}} \right)  \right)  \right)
	\bigg(e^{-it |\eta_0|^2}  \left( \overline{\mathcal{F}_{\tilde{y}}  v_{n'}^h} \right)\left(t, \eta_0, \tilde{x} \right) \notag\\
	& \qquad e^{-it |\eta_2|^2 } \left( \overline{ \mathcal{F}_{\tilde{y} }  v_{n_2'}^h} \right) \left(t, \eta_2, \tilde{x} \right) e^{it |\eta_1|^2} \left( \mathcal{F}_{\tilde{y}}  v_{n_1'}^l \right) \left(t, \eta_1, \tilde{x} \right)
	e^{it |\eta_3 |^2} \left(  \mathcal{F}_{\tilde{y}}  v_{n_3'}^l \right)\left(t, \eta_3 , \tilde{x} \right) \,\mathrm{d} \eta_1 \mathrm{d} \eta_2 \mathrm{d}\eta_3 \mathrm{d}\eta_0   \bigg) \,\mathrm{d} y \mathrm{d} \tilde{y}\mathrm{d} x \mathrm{d} \tilde{x}  \mathrm{d}t, \notag
\end{align}
\begin{align}
B_4 = 	& - \int_{t_0}^{t_1}
	\iiint \iiint \iint \frac1{i q(\eta)} e^{it q(\eta)} \Im \left( \overline{ { w}_n }  ( \nabla_y  - i \xi(t)) w_n \right)(t,y,x) \frac{ y - \tilde{y}}{ | y - \tilde{y}|}
	\phi\left( \frac{\eta_1 + \eta_2 + \eta_3 + \eta_0 }{2^{l_2 - 10}} \right) \label{eq6.50v51}  \\
	& e^{i\tilde{y} ( \eta_1 + \eta_2 + \eta_3 + \eta_0)} \left( 1 - \phi \left( \frac{\eta_0 - \xi(t)}{2^{l_2}} \right) \phi \left( \frac{\eta_2 - \xi(t)} {2^{l_2}} \right) \right)
	\frac{\partial}{\partial t} \bigg( e^{-it |\eta_0|^2} \left( \overline{\mathcal{F}_{\tilde{y}}  v_{n'}^h} \right) \left(t, \eta_0, \tilde{x} \right)
	e^{-it |\eta_2|^2} \left( \overline{ \mathcal{F}_{\tilde{y}}  v_{n_2'}^h} \right) \left(t, \eta_2, \tilde{x} \right)\notag \\
	& \qquad e^{it | \eta_1 |^2} \left( \mathcal{F}_{\tilde{y}}  v_{n_1'}^l \right) \left(t, \eta_1, \tilde{x} \right)
	e^{it | \eta_3 |^2} \left( \mathcal{F}_{\tilde{y}}  v_{n_3'}^l \right)\left(t, \eta_3 , \tilde{x} \right) \,\mathrm{d}\eta_1 \mathrm{d} \eta_2 \mathrm{d} \eta_3 \mathrm{d} \eta_0  \bigg)
	\,\mathrm{d} y \mathrm{d} \tilde{y} \mathrm{d}x \mathrm{d}\tilde{x} \mathrm{d} t. \notag
\end{align}
For \eqref{eq6.47v51}, set
\begin{align*}
	m(t; \eta_0, \eta_1, \eta_2, \eta_3 )
	= \frac1{q(\eta)} \phi \left( \frac{ \eta_1 + \eta_2 + \eta_3 + \eta_0 }{ 2^{l_2 - 10}}  \right) \left ( 1 - \phi \left( \frac{ \eta_0 - \xi(t)}{2^{l_2}} \right)  \phi\left( \frac{ \eta_2 - \xi(t)} { 2^{l_2}} \right) \right).
\end{align*}
Then we have
\begin{align*}
	\eqref{eq6.47v51} = & - i \iiint \iiint \iint \Im \left( \overline{{ w }_n}
	( \nabla_y  - i \xi(t))  w_n \right)(t,y,x) \frac{y - \tilde{y}}{ | y - \tilde{y}|} K(t; z_0, z_1, z_2,z_3 )  \overline{ { v }_{n'}^h
		\left(t, \tilde{y} - z_0 , \tilde{x} \right) }   \\
	& \quad \overline{ { v }_{n_2'}^h \left(t, \tilde{y} - z_2, \tilde{x} \right) }
	v_{n_1'}^l \left(t, \tilde{y} - z_1 , \tilde{x} \right)  v_{n_3'}^l \left(t, \tilde{y} - z_3 , \tilde{x} \right) \,\mathrm{d}z_1 \mathrm{d} z_2 \mathrm{d}z_3 \mathrm{d}z_0 \mathrm{d}x \mathrm{d}\tilde{x} \mathrm{d} y \mathrm{d} \tilde{y} \bigg|_{t_0}^{t_1}  ,
\end{align*}
where
\begin{align}\label{eq6.51v51}
	K(t; z_0, z_1,z_2,z_3 )
	= \iiiint m(t; \eta_0, \eta_1, \eta_2, \eta_3 ) e^{i z_1 \eta_1} e^{i z_2 \eta_2} e^{i z_3 \eta_3} e^{i z_0  \eta_0 } \,\mathrm{d} \eta_1 \mathrm{d} \eta_2 \mathrm{d} \eta_3 \mathrm{d} \eta_0 ,
\end{align}
which satisfies
\begin{align}\label{eq6.52v51}
	\sup\limits_t \int \left|K(t; z_0, z_1, z_2, z_3 ) \right| \,\mathrm{d} z_1 \mathrm{d}z_2 \mathrm{d}z_3 \mathrm{d} z_0  \lesssim 2^{- 2l_2},
\end{align}
by the Coifman-Meyer theorem \cite{GMS}. Thus, by Bernstein's inequality, \eqref{eq6.52v51} and the conservation of mass, we have
\begin{align*}
	& 2^{l_2 - 2i} \bigg| \sum\limits_{n,n' \in \mathbb{N} } \sum\limits_{ \substack{ n_1',n_2',n_3' \in \mathbb{N},\\
			n_1' -n_2' + n_3' =n'}  }
	\eqref{eq6.47v51} \bigg| \\
	\lesssim & 2^{l_2 - 2i} \| w  \|_{L_t^\infty L_{y,x}^2 } \left\| ( \nabla_y  - i\xi(t))  w  \right\|_{L_t^\infty L_{y,x}^2 } \iiiint \left|K(t; z_0,  z_1,z_2,z_3 ) \right|
	\left\|  v^h \left(t, \tilde{y} - z_0 , \tilde{x} \right) \right\|_{L_{\tilde{x}}^2} \left\|  v^h \left( t, \tilde{y} - z_2, \tilde{x} \right) \right\|_{L_{\tilde{x}}^2}\\
	& \quad \left\|  v^l \left(t, \tilde{y} - z_1  , \tilde{x} \right) \right\|_{L_{\tilde{x}}^2} \left\|  v^l\left(t, \tilde{y} - z_3 , \tilde{x} \right) \right\|_{L_{\tilde{x}}^2} \,\mathrm{d}z_1 \mathrm{d} z_2 \mathrm{d} z_3 \mathrm{d} z_0  \mathrm{d} \tilde{y} \\
	\lesssim & 2^{l_2 - i } \| w_0 \|_{L_{y,x}^2 }^2 2^{- 2l_2} \left\|  v^h \right\|_{L_t^\infty L_{y,x}^2 }^2  \left\|  v^l \right\|_{L_{t,y}^\infty  L_x^2 }^2
	\lesssim \| w_0 \|_{L_{y,x}^2}^2.
\end{align*}
Next, we turn to the estimate of \eqref{eq6.48v51}. By a direct computation, we have
\begin{align*}
	& \frac{\partial} { \partial t} \Im  \int_{\mathbb{R}} \left( \overline{ { w}_n }
	( \partial_{y_k}  - i \xi_k(t)) w_n \right)(t,y,x) \,\mathrm{d}x  \\
	= \ & \xi_k'(t) \| w_n(t,y,x) \|_{L_x^2}^2
	+ \sum\limits_{k'=1}^2 \partial_{y_{k'} }\Re \int_{\mathbb{R}} \left( \overline{ { w }_n }
	( \partial_{y_k}  - i \xi_k(t)) \partial_{y_{k'} }  w_n \right)(t,y,x) \,\mathrm{d}x\\
	& - \sum\limits_{k'= 1}^2 \partial_{y_{k'} } \Re \int_{\mathbb{R}} \left( \overline{ \partial_{y_{k'} }  { w }_n  }
	( \partial_{y_k}  - i \xi_k(t) )  w_n \right)(t,y,x) \,\mathrm{d}x.
\end{align*}
Thus, we get
\begin{align}
	\eqref{eq6.48v51}
	& = \notag \\
	&  \int_{G_\beta^{l_2}} \iiint \iiint \iint \frac{ y - \tilde{y}}{ | y - \tilde{y} |}  \xi'(t) \| w_n(t,y,x) \|_{L_x^2}^2 K(t; z_0, z_1, z_2, z_3 ) \overline{ {v }_{n'}^h
		\left(t, \tilde{y} - z_0 , \tilde{x} \right) } \overline{   { v }_{n_2'}^h \left(t, \tilde{y} - z_2, \tilde{x} \right) }
	\label{eq6.54v51}  \\
	& \qquad  v_{n_1'}^l \left(t, \tilde{y } - z_1 , \tilde{x} \right)
	v_{n_3'}^l \left(t, \tilde{y} - z_3  , \tilde{x} \right) \,\mathrm{d}z_1 \mathrm{d} z_2 \mathrm{d} z_3 \mathrm{d}z_0  \mathrm{d} {y} \mathrm{d}\tilde{y}  \mathrm{d}x \mathrm{d}\tilde{x}\mathrm{d}t
	\notag \\
	& + \sum\limits_{k,k'= 1}^2 \int_{G_\beta^{l_2}} \iiint \iiint \iint \frac{ ( y - \tilde{y} )_k} { | y - \tilde{y}|}
	\partial_{y_{k'} } \Re \left(\overline{ { w }_n  }
	( \partial_{y_k}  - i \xi_k(t))  \partial_{y_{k'}}   w_n \right)(t,y,x) K(t; z_0, z_1 ,z_2, z_3 )  \label{eq6.55v51} \\
	& \qquad \overline{ {  v }_{n'}^h \left(t, \tilde{y} - z_0 , \tilde{x} \right) }  \overline{ { v }_{n_2'}^h \left(t, \tilde{y} - z_2, \tilde{x} \right) }
	v_{n_1'}^l \left(t, \tilde{y} - z_1 , \tilde{x} \right)  v_{n_3'}^l \left(t, \tilde{y} - z_3 , \tilde{x} \right) \,\mathrm{d}z_1 \mathrm{d}z_2 \mathrm{d} z_3 \mathrm{d} z_0
	\,\mathrm{d}y \mathrm{d} \tilde{y} \mathrm{d}x \mathrm{d} \tilde{x}  \mathrm{d}t
	\notag \\
	& - \sum\limits_{k,k'= 1}^2 \int_{G_\beta^{l_2}} \iiint \iiint \iint \frac{ ( y - \tilde{y} )_k }{ | y - \tilde{y}|} \partial_{y_{k' } } \Re \left( \partial_{y_{k' } } \bar{ w }_n  ( \partial_{y_k}  - i \xi_k(t))  w_n \right)(t,y,x) K(t; z_0, z_1,z_2,z_3 ) \label{eq6.56v51}  \\
	& \qquad \overline{ { v }_{n'}^h \left(t, \tilde{y} - z_0  , \tilde{x} \right) }
	\overline{ { v }_{n_2'}^h \left(t, \tilde{y} - z_2,\tilde{x} \right) }
	v_{n_1'}^l \left(t, \tilde{y} - z_1 , \tilde{x} \right) v_{n_3'}^l \left(t, \tilde{y} - z_3 , \tilde{x} \right) \,\mathrm{d} z_1\mathrm{d} z_2 \mathrm{d}z_3 \mathrm{d}z_0  \mathrm{d} y \mathrm{d}\tilde{y} \mathrm{d}x \mathrm{d} \tilde{x} \mathrm{d}t, \notag
\end{align}
where $K(t; z_0, z_1,z_2,z_3 )$ is given in \eqref{eq6.51v51}.

By \eqref{eq6.52v51}, \eqref{eq5.1v51}, \eqref{eq5.4v51}, Bernstein's inequality and the conservation of mass, we have
\begin{align*}
	& 2^{l_2 - 2i} \bigg|
	\sum\limits_{n,n' \in \mathbb{N} } \sum\limits_{ \substack{ n_1',n_2',n_3' \in \mathbb{N},\\
			n_1' -n_2' +n_3' =n'}  }  \eqref{eq6.54v51} \bigg|
	\lesssim  2^{l_2 - 2i} 2^{-2l_2} \| w  \|_{L_t^\infty L_{y,x}^2 }^2 \left\| v^h \right\|_{L_t^\infty L_{y,x}^2 }^2  \left\| v^l \right\|_{L_{t,y}^\infty  L_x^2}^2
	\left( \int_{G_\beta^{l_2}} |\xi'(t)| \,\mathrm{d}t \right)
	\lesssim \|  w_0 \|_{L_{y,x}^2 }^2.
\end{align*}
Integrating \eqref{eq6.55v51} by parts in space, we derive
\begin{align*}
	& \ \eqref{eq6.55v51}\\
	= & - \sum\limits_{k,k'= 1}^2 \int_{G_\beta^{l_2}} \iiint \iiint \iint \left( \frac{ \delta_{kk'} }{ | y - \tilde{y} |} + \frac{ ( y - \tilde{y} )_{k'} ( y - \tilde{y} )_k} { | y - \tilde{y}|^3} \right)
	\Re \left( \overline{{ w }_n }
	( \partial_{y_k}  - i \xi_k(t)) \partial_{y_{ k' } }  w_n \right)(t,y,x) \\
	& K(t; z_0, z_1,z_2,z_3 ) \overline{ { v }_{n'}^h \left(t,  \tilde{y} - z_0 , \tilde{x} \right) }
	\overline{ { v }_{n_2'}^h \left(t, \tilde{y} - z_2, \tilde{x} \right) }
	v_{n_1'}^l \left(t, \tilde{y} - z_1 , \tilde{x} \right)    v_{n_3'}^l \left(t, \tilde{y} - z_4 , \tilde{x} \right) \,\mathrm{d} z_1 \mathrm{d} z_2 \mathrm{d} z_3 \mathrm{d} z_0  \mathrm{d}x \mathrm{d} \tilde{x} \mathrm{d} y \mathrm{d} \tilde{y} \mathrm{d}t.
\end{align*}
Therefore, by the Hardy-Littlewood-Sobolev inequality, \eqref{eq6.52v51}, Lemma \ref{le5.7v51}, the Sobolev embedding theorem, the fact $G_\beta^{l_2} \subseteq G_\alpha^i$, $ | \xi(t)|  << 2^{l_2}$ and $l_2 \le i$, we have
\begin{align*}
	2^{l_2 - 2i} \bigg| \sum\limits_{n,n' \in \mathbb{N} } \sum\limits_{ \substack{ n_1',n_2',n_3' \in \mathbb{N},\\
			n_1' - n_2' +n_3' =n' } } \eqref{eq6.55v51} \bigg|
	\lesssim 2^{l_2 - 2i} 2^{2i} 2^{-2l_2} \| w \|_{L_t^6 L_y^3 L_x^2}^2 \left\|  v^h \right\|_{L_t^3 L_y^6 L_x^2}^2 \left\|  v^l \right\|_{L_t^\infty L_y^4 L_x^2}^2
	\lesssim \| w_0 \|_{L_{y,x}^2 }^2 \| v \|_{\tilde{X}_i \left(G_\alpha^i \right)}^2.
\end{align*}
By a similar argument, we infer
\begin{align*}
	2^{l_2 - 2i } \bigg| \sum\limits_{ n,n' \in \mathbb{N}} \sum\limits_{ \substack{ n_1',n_2',n_3' \in \mathbb{N},\\
			n_1' -n_2' + n_3' =n'}  } \eqref{eq6.56v51} \bigg|
	\lesssim \| w_0 \|_{L_{y,x}^2 }^2 \|  v  \|_{\tilde{X}_i(G_\alpha^i)}^2.
\end{align*}
Now we turn to \eqref{eq6.49v51}. As \eqref{eq6.47v51}, we have the corresponding integral kernel
\begin{align*}
	\tilde{K}(t; z_0, z_1,z_2,z_3 ) = \iiiint \tilde{m}(t;  \eta_0, \eta_1, \eta_2, \eta_3 ) e^{i z_1 \eta_1} e^{i z_2 \eta_2} e^{i \eta_3 z_3} e^{iz_0 \eta_0 } \,\mathrm{d} \eta_1 \mathrm{d} \eta_2 \mathrm{d}\eta_3 \mathrm{d} \eta_0 ,
\end{align*}
where
\begin{align*}
	& \tilde{m}(t; \eta_0, \eta_1, \eta_2, \eta_3 ) \\
	= & - \frac{2^{-l_2}}{ q(\eta)} \phi\left( \frac{ \eta_1 + \eta_2 + \eta_3 + \eta_0 }{ 2^{l_2 - 10}} \right)
	\left( ( \nabla \phi) \left( \frac{\eta_0  - \xi(t)}{2^{l_2}} \right) \phi \left( \frac{\eta_2 - \xi(t)} { 2^{l_2}} \right) + \phi \left( \frac{\eta_0 - \xi(t)}{2^{l_2}} \right) ( \nabla \phi) \left(\frac{\eta_2 - \xi(t)}{2^{l_2}} \right) \right).
\end{align*}
The kernel function $\tilde{K}(t; z_0, z_1, z_2, z_3 )$ satisfies
\begin{align*}
	\sup\limits_t \int \left|\tilde{K}(t; z_0, z_1, z_2,z_3 ) \right| \,\mathrm{d}z_1 \mathrm{d} z_2 \mathrm{d}z_3 \mathrm{d} z_0  \lesssim 2^{- 3 l_2}.
\end{align*}
Thus
\begin{align*}
	2^{l_2 - 2i} \bigg| \sum\limits_{n,n' \in \mathbb{N}} \sum\limits_{ \substack{ n_1',n_2' ,n_3' \in \mathbb{N}, \\
			n_1' - n_2' + n_3' =n'} } \eqref{eq6.49v51} \bigg|
	& \lesssim 2^{-2 l_2 - i } \| w \|_{L_t^\infty L_{y,x}^2 }^2 \left\|  v^h \right\|_{L_t^\infty L_{y,x}^2 }^2 \left\|  v^l \right\|_{L_{t,y}^\infty  L_x^2}^2 \left( \int_{G_\beta^{l_2}} |\xi'(t)| \,\mathrm{d}t \right)
	\lesssim \| w_0 \|_{L_{y,x}^2 }^2.
\end{align*}
Finally, we consider the term \eqref{eq6.50v51}. Following the argument for the estimates \eqref{eq6.47v51} and \eqref{eq6.49v51}, by the Bernstein inequality, the conservation of mass and Lemma \ref{le5.7v51}, we deduce
\begin{align*}
	& 2^{l_2 - 2i} \bigg| \sum\limits_{n,n' \in \mathbb{N} } \sum\limits_{ \substack{ n_1',n_2',n_3' \in \mathbb{N}, \\
			n_1' -n_2' +n_3' =n'} } \eqref{eq6.50v51} \bigg|
	\lesssim   2^{l_2 -i} \| w_0 \|_{L_{y,x}^2}^2 \left( 1 + \|  v  \|_{\tilde{X}_i \left(G_\alpha^i \right)}^4 \right).
\end{align*}
Therefore, we eventually arrive at
\begin{align*}
	\eqref{eq6.5v51}
	\lesssim \|w_0\|_{L_{y,x}^2 }^2 \left( 1 + \|  v \|_{\tilde{X}_i \left(G_\alpha^i \right)}^4 \right).
\end{align*}
The proof of Theorem \ref{th5.16v51} is complete.
\end{proof}

\subsection{Proof of Theorem \ref{th7.1v51}}\label{subsec:theorem 2}

\begin{proof}[Proof of Theorem \ref{th7.1v51}]
	By Theorem \ref{th5.8v51}, we have
	\begin{align}\label{eq7.2v51}
		\|v_\lambda \|_{\tilde{X}_{k_0} \left( \left[0, {\lambda^{-2} } T  \right] \right)}  \lesssim 1,
	\end{align}
	where
	\begin{align}\label{eq6.98v77}
		v_\lambda(t,y,x)  = \lambda  v  \left(\lambda^2 t , \lambda y, x \right) \text{ with } \lambda = \frac{\epsilon_3 2^{k_0}}K.
	\end{align}
	Let $ \tilde{ w}  = P^y_{\le 2^{ k_0} }
	v_\lambda $. Then $\tilde{w}$ satisfies
	\begin{align*}
		i \partial_t  \tilde{ w }  + \Delta_y  \tilde{ w }  = F( \tilde{ w } ) + N,
	\end{align*}
	where
	\begin{align*}
		N= P^y_{\le 2^{ k_0} }  F( v_\lambda  ).
	\end{align*}
Let
\begin{align*}
		M(t) = \int_{\mathbb{R}} \int_{\mathbb{R}} \int_{\mathbb{R}^2} \int_{\mathbb{R}^2} \left| \tilde{ w }  \left(t, \tilde{y}, \tilde{x} \right) \right|^2 \frac{ y - \tilde{y}}{ | y - \tilde{y}|} \Im \left( \overline{ \tilde{ w } } \nabla_y   \tilde{ w } \right)(t,y,x) \,\mathrm{d} y \mathrm{d} \tilde{y}  \mathrm{d} x \mathrm{d} \tilde{x}.
\end{align*}
	Then a direct calculation similar to \cite{D1,D2,D3,PV} gives
	\begin{align*}
		\left\| \int_{\mathbb{R}} |\nabla_y|^\frac12 \left( \left|  \tilde{ w } (t,y,x) \right|^2 \right) \,\mathrm{d}x \right\|_{L_{t,y}^2([ 0,  {\lambda^{-2} }T ] \times \mathbb{R}^2  )}^2
		\lesssim \sup\limits_{t \in [ 0, {\lambda^{- 2} } T ]} |M(t)|  + \mathscr{E},
	\end{align*}
	where
	\begin{align}
		\mathscr{E}
		& = \notag \\
		&  2 \bigg| \int_0^{  {\lambda^{- 2}} T} \int_{\mathbb{R}} \int_{\mathbb{R}} \int_{\mathbb{R}^2} \int_{\mathbb{R}^2}
		\Im \left( \overline{ \tilde{ w } } ( \nabla_y  - i \xi(t)) \tilde{ w} \right)(t,y,x) \frac{ y - \tilde{y}}{ | y - \tilde{y}|} \Im \left( \overline{ \tilde{ w } } N \right)(t, \tilde{y}, \tilde{x})\,\mathrm{d} y \mathrm{d} \tilde{y} \mathrm{d}x \mathrm{d} \tilde{x} \mathrm{d} t \bigg|  \label{eq7.5v51}\\
		+  & \bigg| \int_0^{  {\lambda^{- 2}} T} \int_{\mathbb{R}} \int_{\mathbb{R}}  \int_{\mathbb{R}^2} \int_{\mathbb{R}^2} |  \tilde{ w } \left(t, \tilde{y}, \tilde{x} \right)|^2 \frac{ y - \tilde{y}}{ | y  - \tilde{y}|} \Im \left( \bar{N}( \nabla_y  - i \xi(t)) \tilde{ w }  \right)(t,y, x) \,\mathrm{d} y \mathrm{d} \tilde{y} \mathrm{d} x \mathrm{d}\tilde{x} \mathrm{d}t \bigg| \label{eq7.6v51} \\
		+ & \bigg|\int_0^{  {\lambda^{- 2}} T}  \int_{\mathbb{R}} \int_{\mathbb{R}}   \int_{\mathbb{R}^2} \int_{\mathbb{R}^2} |  \tilde{ w } \left(t, \tilde{y}, \tilde{x} \right)|^2 \frac{ y - \tilde{y}} { | y - \tilde{y} |} \Im  \left( \overline{ \tilde{ w} }( \nabla_y  - i \xi(t)) N \right)(t,y,x) \,\mathrm{d} y \mathrm{d} \tilde{y} \mathrm{d}x \mathrm{d}\tilde{x} \mathrm{d}t \bigg|.\label{eq7.7v51}
	\end{align}
	Since $N(t) \le 1$, we have $N_\lambda(t) \le \frac{\epsilon_3 2^{k_0}}K$. By Theorem \ref{th3.3v51} and the Bernstein inequality, for any $\eta> 0$, if $K  \ge C(\eta)$, we have
	\begin{align}\label{eq7.8v51}
		\left\| \left( \nabla_y  - i \xi(t) \right) \tilde{  w } \right\|_{L_t^\infty L_{y,x}^2 ([ 0,  {\lambda^{-2} }T ] \times \mathbb{R}^2  \times \mathbb{R}  )} \lesssim \eta 2^{k_0}.
	\end{align}
	Therefore, by the Galilean transformation and the conservation of mass, we get
	\begin{align*}
		\sup\limits_{t \in [ 0, {\lambda^{- 2} }T]} |M(t)| \lesssim \eta 2^{k_0}.
	\end{align*}
	We now consider \eqref{eq7.5v51}. As in \eqref{eq6.25v51}, let $v_\lambda^l = P^y_{\le 2^{ k_0 - 3 } } v_\lambda $ and $v_\lambda^h = P_{>  2^{  k_0 - 3 } }^y v_\lambda$, then we have the decomposition
	\begin{align*}
		\int_{\mathbb{R}} \Im \left( \overline{ \tilde{ w } } N \right)(t, \tilde{y}, \tilde{x})  \mathrm{d} \tilde{x}  = \int_{\mathbb{R}} F_0 \left(t, \tilde{y} , \tilde{x} \right) + F_1\left(t, \tilde{y}, \tilde{x} \right) + F_2\left(t, \tilde{y}, \tilde{x} \right) + F_3 \left(t, \tilde{y}, \tilde{x} \right) + F_4 \left(t, \tilde{y}, \tilde{x} \right) \,\mathrm{d} \tilde{x}.
	\end{align*}
We can see
	\begin{align*}
		\int_{\mathbb{R}} F_0 \left(t, \tilde{y}, \tilde{x} \right) \,\mathrm{d} \tilde{x}  = 0 .
	\end{align*}
	Following the same argument as the proof of \eqref{eq6.5v51}, we may obtain
	\begin{align*}
		\int_{\mathbb{R}} \left\|F_2 \left(t, \tilde{y}, \tilde{x} \right) + F_3 \left(t, \tilde{y}, \tilde{x} \right) + F_4 \left(t, \tilde{y}, \tilde{x} \right)  \right\|_{L_{t,\tilde{y}}^1([ 0,  {\lambda^{-2} }T ] \times \mathbb{R}^2  ) } \,\mathrm{d} \tilde{x} \lesssim 1.
	\end{align*}
	Then by \eqref{eq7.8v51} and the conservation of mass, we have
	\begin{align*}
		& \ \left| \int_0^{  {\lambda^{- 2}} T} \int_{\mathbb{R}} \int_{\mathbb{R}}   \int_{\mathbb{R}^2} \int_{\mathbb{R}^2}
		\Im \left( \overline{ \tilde{ w } } ( \nabla_y  - i \xi(t)) \tilde{  w } \right)(t,y,x)
		\frac{y - \tilde{y}}{ | y - \tilde{y}|} \left(F_2  + F_3  + F_4 \right) \left(t, \tilde{y}, \tilde{x}  \right) \,\mathrm{d} y \mathrm{d} \tilde{y} \mathrm{d}x \mathrm{d} \tilde{x} \mathrm{d}t \right|
		\lesssim   \eta 2^{k_0}.
	\end{align*}
To estimate the contribution of the term with $F_1$ in \eqref{eq7.5v51}, we see the support of the spatial Fourier transform of $\int_{\mathbb{R}} F_1(t, \tilde{y} , \tilde{x} ) \,\mathrm{d} \tilde{x}$ is in $\left\{ \xi: |\xi| \ge 2^{k_0 - 4} \right\} $ as in \eqref{eq6.27v51}. Therefore, by integration by parts, the Hardy-Littlewood-Sobolev inequality, the Bernstein inequality, Lemma \ref{le5.7v51} and \eqref{eq7.2v51}, we have
	\begin{align*}
		&  \bigg| \int_0^{  {\lambda^{- 2}} T}  \int_{\mathbb{R}}   \int_{\mathbb{R}^2} \int_{\mathbb{R}^2} \Im \left( \overline{ \tilde{ w } } ( \nabla_y  - i \xi(t) )  \tilde{ w }
		\right)(t, y, x) \frac{y - \tilde{y}}{ |y - \tilde{y}|} \left(
		\int_{\mathbb{R}} F_1 \left(t, \tilde{y}, \tilde{x} \right) \,\mathrm{d}\tilde{x} \right)
		\,\mathrm{d} y \mathrm{d} \tilde{y} \mathrm{d}x \mathrm{d}t \bigg| \\
		\lesssim  & \int_0^{ {\lambda^{-2} } T } \int_{\mathbb{R}^2} \int_{\mathbb{R}^2} \left| \int_{\mathbb{R}} \left( \overline{ \tilde{ w } } ( \nabla_y  - i \xi(t)) \tilde{ w} \right)(t,y, x) \,\mathrm{d}x \right|\ \cdot
		\frac1{ | y - \tilde{y} |} \ \cdot  \left| {\partial_{\tilde{y}}}\  {\left( -  \Delta_{\tilde{y}} \right)^{- 1} }
		\int_{\mathbb{R}} F_1(t, \tilde{y}, \tilde{x}) \,\mathrm{d} \tilde{x} \right| \,\mathrm{d} y \mathrm{d} \tilde{y} \mathrm{d}t \\
		\lesssim & \|  \tilde{ w } \|_{L_t^\infty L_{y,x}^2 \left([0, \lambda^{-2} {T} ] \times \mathbb{R}^2 \times \mathbb{R} \right)} \left\| ( \nabla_y - i \xi(t)) \tilde{ w } \right\|_{L_{t,y}^4 L_x^2\left([0, \lambda^{-2} {T} ] \times \mathbb{R}^2 \times \mathbb{R}  \right)}
		\left\|{ \partial_{\tilde{y}}}\ {\left(  -\Delta_{\tilde{y}}\right)^{- 1} }
		\int_{\mathbb{R}} F_1 \left(t, \tilde{y}, \tilde{x} \right) \,\mathrm{d} \tilde{x} \right \|_{L_{t,y}^\frac43 \left([0, \lambda^{-2} {T} ] \times \mathbb{R}^2 \right)} \\
		\lesssim & 2^{-k_0} \left\| ( \nabla_y  - i \xi(t)) \tilde{ w}  \right\|_{L_{t,y}^4 L_x^2\left([0, \lambda^{-2} {T} ] \times \mathbb{R}^2 \times \mathbb{R}  \right)} \left\| v_\lambda^l \right\|_{L_{t,y}^6 L_x^2\left([0, \lambda^{-2} {T} ] \times \mathbb{R}^2 \times \mathbb{R} \right)}^3.
	\end{align*}
	By the Bernstein inequality, Lemma \ref{le5.7v51}, and \eqref{eq7.2v51}, we have
	\begin{align*}
		\left\| v_\lambda^l \right\|_{L_{t,y}^6 L_x^2\left([0, \lambda^{-2} {T} ] \times \mathbb{R}^2 \times \mathbb{R} \right) }
		\lesssim \sum\limits_{0 \le l \le k_0} 2^{\frac{l}3} 2^{\frac{k_0 - l}6} \lesssim 2^{\frac{k_0}3}.
	\end{align*}
	Note that
	\begin{align}\label{eq7.12v51}
		\left\|\left(\nabla_y  - i \xi(t) \right) \tilde{ w} \right\|_{L_t^\frac52 L_y^{10} L_x^2\left([0, \lambda^{-2} {T} ] \times \mathbb{R}^2 \times \mathbb{R} \right)}
		\lesssim \sum\limits_{0 \le l \le k_0} 2^l 2^{\frac25(k_0 - l)} \lesssim 2^{k_0}.
	\end{align}
	Interpolating \eqref{eq7.12v51} and \eqref{eq7.8v51}, we obtain
	\begin{align}\label{eq7.13v51}
		\left \| \left( \nabla_y  - i \xi(t) \right) \tilde{ w}  \right \|_{L_{t,y}^4 L_x^2\left([0, \lambda^{-2} {T} ] \times \mathbb{R}^2 \times \mathbb{R} \right)} \lesssim \eta^\frac38 2^{k_0}.
	\end{align}
	Thus, by the above estimates, we have
	\begin{align*}
 \eqref{eq7.5v51} \lesssim \eta^\frac38 2^{k_0}.
	\end{align*}
	Now, we turn to \eqref{eq7.6v51}. By \eqref{eq6.17v51} and \eqref{eq7.2v51}, we have
	\begin{align*}
 \eqref{eq7.6v51}
 \lesssim  & \left\| \tilde{ w} \right\|_{L_t^\infty L_{y,x}^2 \left([0, \lambda^{-2} {T} ] \times \mathbb{R}^2 \times \mathbb{R}  \right)}^2 \|N\|_{L_{t,y}^\frac43 L_x^2\left([0, \lambda^{-2} {T} ] \times \mathbb{R}^2 \times \mathbb{R} \right)} \left\| ( \nabla_y  - i \xi(t)) \tilde{ w}  \right\|_{L_{t,y}^4 L_x^2\left([0, \lambda^{-2} {T} ] \times \mathbb{R}^2 \times \mathbb{R} \right)}  \lesssim \eta^\frac38 2^{k_0}.
	\end{align*}
	Finally, we consider \eqref{eq7.7v51}. Applying integration by parts, we have
	\begin{align}
 \eqref{eq7.7v51}
		\le &  \eqref{eq7.6v51}  \notag \\
		& + \left| \int_0^{  {\lambda^{- 2}} T}  \int_{\mathbb{R}} \int_{\mathbb{R}}  \int_{\mathbb{R}^2} \int_{\mathbb{R}^2} \left|\tilde{ w} \left(t, \tilde{y}, \tilde{x} \right)  \right|^2 \frac1{ | y - \tilde{y}|} \Re \left( \overline{\tilde{w} }
		N \right)(t,y, x) \,\mathrm{d}y \mathrm{d}\tilde{y} \mathrm{d}x \mathrm{d} \tilde{x} \mathrm{d}t \right|.\label{eq6.105v77}
	\end{align}
	By \eqref{eq3.8v41} and \eqref{eq6.24v51}, we see
	\begin{align}
		\eqref{eq6.105v77}\lesssim & \notag\\
		&   \int_0^{{\lambda^{- 2} } T} \int_{\mathbb{R}^2} \int_{\mathbb{R}^2} \left\| \tilde{ w } \left(t, \tilde{y}, \tilde{x} \right) \right\|_{L_{\tilde{x}}^2}^2 \frac1{ | y - \tilde{y} |}
		\| \tilde{ w} (t, y,x) \|_{L_x^2} \left\| v_\lambda^h(t,y,x) \right\|_{L_x^2}^3\,\mathrm{d} y \mathrm{d}\tilde{y} \mathrm{d}t \label{eq7.15v51} \\
		& + \int_0^{ {\lambda^{-2}} T } \int_{\mathbb{R}^2} \int_{\mathbb{R}^2} \| \tilde{ w} \left(t, \tilde{y}, \tilde{x} \right)\|_{L_{\tilde{x}}^2}^2 \frac1{ | y - \tilde{y}|}
		\| \tilde{ w} (t,y,x) \|_{L_x^2} \left\| v_\lambda^l(t,y,x) \right\|_{L_x^2}^2  \left\|  v_\lambda^h(t,y,x) \right\|_{L_x^2}\,\mathrm{d} y \mathrm{d} \tilde{y} \mathrm{d}t. \label{eq7.16v51}
	\end{align}
	By the Hardy-Littlewood-Sobolev inequality, \eqref{eq7.8v51}, \eqref{eq7.13v51}, Lemma \ref{le5.7v51}, the Sobolev embedding theorem, the conservation of mass, and interpolation, we have
	\begin{align*}
		\eqref{eq7.15v51} & \lesssim \left\|  v_\lambda^h \right\|_{L_{t,y}^4 L_x^2\left([0, \lambda^{-2} {T} ] \times \mathbb{R}^2 \times \mathbb{R} \right)}^3 \left\| \tilde{ w } \right\|_{L_{t,y}^4 L_x^2\left([0, \lambda^{-2} {T} ] \times \mathbb{R}^2 \times \mathbb{R} \right)}^3
		\lesssim \eta^\frac38 2^{k_0},
		\intertext{ and }
		\eqref{eq7.16v51} & \lesssim \left\|  v_\lambda^h \right\|_{L_t^3 L_y^6 L_x^2\left([0, \lambda^{-2} {T} ] \times \mathbb{R}^2 \times \mathbb{R} \right)} \left\| \tilde{ w}  \right\|_{L_t^9 L_y^\frac{90}{29} L_x^2\left([0, \lambda^{-2} {T} ] \times \mathbb{R}^2 \times \mathbb{R} \right)}^3
		\left\| v_\lambda^l \right\|_{L_t^6 L_y^\frac{60}{11} L_x^2\left([0, \lambda^{-2} {T} ] \times \mathbb{R}^2 \times \mathbb{R} \right)}^2
		\lesssim \eta^\frac16 2^{k_0}.
	\end{align*}
	Thus, by the above estimates, we have
	\begin{align*}
		\left\| \int_{\mathbb{R}} |\nabla_y |^\frac12 \left( \left| \tilde{ w} (t,y,x) \right|^2 \right) \,\mathrm{d} x \right\|_{L_{t,y}^2([0, \lambda^{-2} {T} ] \times \mathbb{R}^2)}^2 \lesssim \eta^\frac16 2^{k_0}.
	\end{align*}
	Undoing the scaling in \eqref{eq6.98v77}, we finally reach the desired estimate \eqref{eq6.97v77}.
\end{proof}

\section{Appendix: well-posedness theory for \eqref{eq1.1} }\label{sec:appendix}

In this appendix, we present the proofs of the recorded results in Section \ref{se2} for self-contained.
Let
\begin{align*}\label{eq8.1v127}
X_1(t) & = x \sin(t) - i \cos(t) \partial_x \text{ and }
X_2(t)   = x \cos(t)+ i\sin(t) \partial_x.
\end{align*}
We have the point-wise identity: for any $f\in \mathcal{S}(\mathbb{R}^3)$,
\begin{align}
\left|X_1(t) f(y,x ) \right|^2 + \left|X_2(t) f(y, x ) \right|^2 = \left|x f(y, x ) \right|^2 + \left|\partial_x f(y,x ) \right|^2, \ \forall t\in \mathbb{R}.
\end{align}
The next result follows by direct computation. We refer to \cite{C} for more explanation.
\begin{lemma}\label{le8.2v126}
The operators $X_1(t)$ and $X_2(t)$ satisfy the following properties:

(1) They correspond to the conjugation of gradient and momentum by the free flow,
\begin{align*}
X_1(t) &  = e^{it \left(\Delta_{\mathbb{R}^3} - x^2 \right)} \left(-i\partial_x \right) e^{-it \left(\Delta_{\mathbb{R}^3}- x^2 \right)},\\
X_2(t) &  = e^{it \left(\Delta_{\mathbb{R}^3}-x^2 \right)} x e^{-it \left(\Delta_{\mathbb{R}^3}-x^2 \right)}.
\end{align*}
(2) They act on the nonlinearity like derivatives, that is for $j  = 1, 2$, we have
\begin{align*}
\left|X_j(t) \left(|u|^2 u \right) \right| \lesssim |u|^2 \left|X_j(t) u \right|.
\end{align*}
\end{lemma}
As a consequence, we have
\begin{align*}
& \left\|e^{-it \left(\Delta_{\mathbb{R}^3}-x^2 \right)} u(t) - u_{\pm} \right\|_{\Sigma} \\
= & \left\|u(t)-e^{it \left(\Delta_{\mathbb{R}^3}- x^2 \right)} u_{\pm} \right\|_{L_x^2 H_y^1}
+ \left\|X_1(t) \left( u(t)-e^{it \left(\Delta_{\mathbb{R}^3}- x^2 \right)} u_{\pm} \right) \right\|_{L_{y,x }^2 }
+ \left\|X_2(t) \left(u(t)-e^{it \left(\Delta_{\mathbb{R}^3}- x^2 \right)} u_{\pm} \right) \right\|_{L_{y, x }^2}.
\end{align*}
We now show the local well-posedness part of Theorem \ref{th2.3} in the following formulation. This is essentially following the argument in \cite{Cazenave,T2}.
\begin{theorem}[Local well-posedness]
For any $E>0$ and $u_0$ with $ \left\|u_0 \right\|_{L_y^2 \mathcal{H}_{x}^{1 }(\mathbb{R}^2  \times \mathbb{R})} \le E$, there exists $\delta_0 = \delta_0(E)>0$ such that if
\begin{equation*}
\left\|e^{it \left(\Delta_{\mathbb{R}^3} -x^2 \right)  } u_0 \right\|_{L_{t,y}^4   L_x^{2}  (I\times \mathbb{R}^2 \times \mathbb{R})}
+ \left\|X_1(t) e^{it \left(\Delta_{\mathbb{R}^3} -x^2 \right)  } u_0 \right\|_{L_{t,y}^4  L_x^{2}  (I\times \mathbb{R}^2 \times \mathbb{R})}
+ \left\|X_2(t) e^{it \left(\Delta_{\mathbb{R}^3} -x^2 \right)  } u_0 \right\|_{L_{t,y}^4   L_x^{2}  (I\times \mathbb{R}^2 \times \mathbb{R})}
 \le   \delta_0,
\end{equation*}
where $I$ is the time interval, there exits a unique solution $u\in C_t^0  L_y^2 \mathcal{H}_{ x}^{ 1 } \left(I\times \mathbb{R}^2 \times \mathbb{R} \right)$ of \eqref{eq1.1}
satisfying
\begin{align*}
\|u\|_{L_{t,y}^4    \mathcal{H}_x^{1}\left(I\times \mathbb{R}^2 \times \mathbb{R} \right) } & \le 2 \left\|e^{it \left(\Delta_{ \mathbb{R}^3}-x^2 \right) } u_0 \right\|_{L_{t,y}^4   \mathcal{H}_x^{1  } \left(I\times \mathbb{R}^2 \times \mathbb{R} \right)}, \text{ and }
\|u\|_{L_t^\infty L_y^2 \mathcal{H}_{ x}^{ 1 }\left(I\times \mathbb{R}^2 \times \mathbb{R} \right) }    \le C \left\|u_0 \right\|_{ L_y^2 \mathcal{H}_{x}^{ 1 }  }.
\end{align*}
\end{theorem}
\begin{proof}
Let
\begin{align*}
\Phi(u) = e^{it \left(\Delta_{\mathbb{R}^3}- x^2 \right)} u_0 - i\int_0^t e^{i(t-s) \left(\Delta_{\mathbb{R}^3} - x^2 \right)}  \left(|u|^2 u \right)(s) \,\mathrm{d}s,
\end{align*}
and set the space $X$ to be
\begin{align*}
X = \left\{ u\in C_t^0 L_y^2 \mathcal{H}_x^1: \|u\|_{L_t^\infty L_y^2 \mathcal{H}_x^1  } \le 2 E,
\|u\|_{L_{t,y}^4 \mathcal{H}_x^1} \le 2 C\delta_0\right\},
\end{align*}
or
\begin{align*}
X = \left\{ u \in C_t^0 L_y^2 \mathcal{H}_x^1 : \,   \|  u\|_{L_t^\infty L_{y,x}^2} \le 2E, \ \|  u\|_{L_{t,y}^4 L_x^2} \le 2C\delta_0,
\left\|X_j(t) u \right\|_{L_t^\infty L_{y,x}^2} \le 2E, \  \left\|X_j(t) u \right\|_{L_{t,y}^4 L_x^2} \le 2C \delta_0,\ j= 1,2\right\}.
\end{align*}
For any $u\in X$, by Proposition \ref{pr2.2v31}, H\"older's inequality, Sobolev's inequality, Lemma \ref{le8.2v126}, and \eqref{eq8.1v127}, we have
\begin{align*}
 & \left\|\Phi(u) \right\|_{L_t^\infty L_{y,x}^2}
\lesssim \left\|u_0 \right\|_{L_{y,x}^2} + \|u\|_{L_{t,y}^4 L_x^2} \|u\|_{L_{t,y}^4 H_x^1}^2,
\intertext{ and }
& \ \left\|X_1(t) \Phi(u) \right\|_{L_t^\infty L_{y,x}^2}+ \left \|X_2(t) \Phi(u) \right\|_{L_t^\infty L_{y,x}^2}\\
&
\lesssim \left\|\nabla_x u_0 \right\|_{L_{y,x}^2} + \left\|x u_0 \right\|_{L_{y,x}^2}
+ \|u\|_{L_{t,y}^4 H_x^1}^2 \left(  \left\|X_1(t) u\right\|_{L_{t,y}^4 L_x^2}
+ \left\|X_2(t) u \right\|_{L_{t,y}^4 L_x^2} \right) .
\end{align*}
Thus
\begin{align}\label{eq8.2v127}
\left\|\Phi(u) \right\|_{L_t^\infty L_{y,x}^2} + \left\|X_1(t) \Phi(u) \right\|_{L_t^\infty L_{y,x}^2} + \left\|X_2(t) \Phi(u) \right\|_{L_t^\infty L_{y,x}^2} \le E + \left(2C\delta_0 \right)^3 \le 2E.
\end{align}
Similarly, we can obtain
\begin{align}\label{eq8.3v127}
\left\|\Phi(u) \right\|_{L_{t,y}^4 L_x^2} + \left \|X_1(t) \Phi(u) \right\|_{L_{t,y}^4 L_x^2} + \left\|X_2(t) \Phi(u) \right\|_{L_{t,y}^4 L_x^2}
\le \delta_0 +  \left(2C\delta_0 \right)^3 \le 2C \delta_0.
\end{align}
In the same time, for any $u, v\in X$, by the Strichartz estimate, H\"older's inequality, and Sobolev's inequality, we have
\begin{align}\label{eq8.4v127}
\left\|\Phi(u)- \Phi(v)\right\|_{L_{t,y}^4 L_x^2} \lesssim \big\| |u|^2 u - |v|^2 v\big\|_{L_{t,y}^\frac43 L_x^2}
& \lesssim \|u-v\|_{L_{t,y}^4 L_x^2} \left(\|u\|_{L_{t,y}^4 H_x^1}^2 + \|v\|_{L_{t,y}^4 H_x^1}^2\right)\\
& \lesssim  \left(2C\delta_0 \right)^2 \|u-v\|_{L_{t,y}^4 L_x^2}. \notag
\end{align}
%This
Combining \eqref{eq8.2v127}, \eqref{eq8.3v127}, and \eqref{eq8.4v127}, we have for $\delta_0$ small enough, $\Phi: X \to X $ is a contractive map. Therefore, the theorem follows from the fixed point theorem.
\end{proof}

We now turn to the proof of the scattering norm in Theorem \ref{th2.3}.
\begin{proof}[Proof of the scattering norm part of Theorem \ref{th2.3}]
We need to show
\begin{align}\label{eq2.2046}
\|u\|_{L_{t,y}^4 \mathcal{H}_x^{1 } \cap L_t^4 W_y^{1,4} L_x^2(\mathbb{R} \times \mathbb{R}^2  \times \mathbb{R})} \le C(M),
\end{align}
then by the scattering theory of the nonlinear Schr\"odinger equations \cite{ACS,C2,T2}, we have scattering in \eqref{eq3.1v131}.
By well-posedness part of Theorem \ref{th2.3}, it suffices to prove \eqref{eq2.2046} as an a priori bound.

Divide the time interval $\mathbb{R}$ into $N \sim \left(1+ \frac{L}\delta \right)^{4}$ subintervals $I_j = [t_j,t_{j+1}]$ such that
\begin{align}\label{eq2.2146}
\|u\|_{L_{t,y}^4 \mathcal{H}_x^{1- \epsilon_0 } (I_j \times \mathbb{R}^2  \times \mathbb{R})} \le \delta,
\end{align}
where $\delta >0 $ will be chosen later.

On each $I_j$, by \eqref{eq8.1v127}, the Strichartz estimate, the Sobolev embedding and \eqref{eq2.2146}, we have
\begin{align*}
 &\  \|u\|_{L_t^4 W_y^{1,4} L_x^2 \cap L_{t,y}^4 \mathcal{H}_x^1  (I_j \times \mathbb{R}^2  \times \mathbb{R})} \\
%\le & \| u\|_{L_{t,y}^2 L_x^2} + \left\|X_1(t) u \right\|_{L_{t,y}^4 L_x^2} + \left\|X_2(t) u \right\|_{L_{t,y}^4 L_x^2}\\
\le &  C  \left( \left\|u \left(t_j \right) \right\|_{\Sigma}  + \left\||u|^2 u \right\|_{L_{t,y}^\frac43 L_x^2} + \left\|X_1(t) \left(|u|^2 u \right) \right\|_{L_{t,y}^\frac43 L_x^2} + \left\|X_2(t) \left(|u|^2 u \right) \right\|_{L_{t,y}^\frac43 L_x^2} \right)\\
\le & C \left( \left\|u \left(t_j \right) \right\|_{\Sigma} + \|u\|_{L_{t,y}^4 H_x^{ 1 - \epsilon_0} }^2  \left( \|u\|_{L_{t,y}^4 L_x^2} +  \left\|X_1(t) u \right\|_{L_{t,y}^4 L_x^2} + \left \|X_2(t) u \right\|_{L_{t,y}^4 L_x^2} \right)+ \|u\|_{L_t^4 W_y^{1,4} L_x^2} \|u\|_{L_{t,y}^4 H^{1- \epsilon_0}_x}^2 \right)\\
\le & C \left( \left\|u \left(t_j \right) \right\|_{\Sigma} + \|u\|_{L_{t,y}^4 H_x^{ 1 - \epsilon_0} }^2 \left(\|u\|_{L_{t,y}^4 L_x^2} + \|\nabla_x u\|_{L_{t,y}^4 L_x^2} + \||x| u \|_{L_{t,y}^4 L_x^2}+ \|u\|_{L_t^4 W_y^{1,4} L_x^2} \right) \right)\\
\le & C \left( \left\|u \left(t_j \right) \right\|_{\Sigma} + \delta^2 \|u\|_{L_t^4 W_y^{1,4} L_x^2 \cap L_{t,y}^4 \mathcal{H}_x^1 } \right).
\end{align*}
Choosing $\delta \le \left(\frac1{2C} \right)^\frac14$ leads to the estimate
\begin{align*}
\|u\|_{L_t^4 W_y^{1,4} L_x^2 \cap L_{t,y}^4 \mathcal{H}_x^1 (I_j \times \mathbb{R}^2   \times \mathbb{R})} \le 2C \left\|u(t_j) \right\|_{\Sigma_{y,x} }.
\end{align*}
The desired bound \eqref{eq2.2046} now follows by adding up the bounds on each subintervals $I_j$.
\end{proof}
We now provide the proof of Theorem \ref{th2.6}. First, we show the following short-time version.
\begin{lemma}[Short-time stability theorem]\label{le3.6}
Let $I$ be a compact interval and let $\tilde{u}$ be an approximate solution to \eqref{eq1.1} in the sense that
$i\partial_t \tilde{u} + \Delta_{\mathbb{R}^3 } \tilde{u} - x^2 \tilde{u} = |\tilde{u}|^2 \tilde{u} + e$
for some function $e$.
Assume that
\begin{equation}\label{eq2.13new}
\|\tilde{u}\|_{  L_t^\infty L_y^2 \mathcal{H}_x^{1}(I\times \mathbb{R}^2   \times \mathbb{R})} \le M
\end{equation}
for some positive constant $M$.
Let $t_0\in I$ and $u(t_0)$ be such that
\begin{equation}\label{eq3.20}
\left\|u(t_0)- \tilde{u}(t_0) \right\|_{    L_y^2 \mathcal{H}_x^{1}} \le M'
\end{equation}
for some $M' >0$.

Assume also the smallness conditions hold:
\begin{align}
&  \left\|\tilde{u} \right\|_{L_t^4 L_y^4 \mathcal{H}_x^{1} (I\times \mathbb{R}^2 \times \mathbb{R})}   \le \epsilon, \label{eq3.21}\\
& \left\|e^{i(t-t_0) \left(\Delta_{\mathbb{R}^3} - x^2 \right) } \left(u(t_0) - \tilde{u}(t_0) \right) \right\|_{L_t^4 L_y^4 \mathcal{H}_x^{1}  } +
 \|e\|_{L_t^\frac43 L_y^\frac43 \mathcal{H}_x^{1} }    \le \epsilon, \label{eq3.22}
\end{align}
for some $0 < \epsilon \le \epsilon_1$, where $\epsilon_1 = \epsilon_1(M,M') > 0$ is a small constant.
Then, there exists a solution $u$ to \eqref{eq1.1} on $I\times \mathbb{R}^2  \times \mathbb{R}$ with initial data $u(t_0)$ at time $t= t_0$ satisfying
\begin{align}
\left\|u-\tilde{u} \right\|_{L_{t,y}^4   \mathcal{H}_x^{1} }  & \lesssim \epsilon, \label{eq3.24}\\
\left\|u-\tilde{u} \right\|_{  L_t^\infty L_y^2 \mathcal{H}_x^{1}}   &  \lesssim M', \label{eq3.25}\\
\|u\|_{ L_t^\infty L_y^2 \mathcal{H}_x^{1}}  &  \lesssim M + M',\label{eq3.26} \\
\big\||u|^2 u - |\tilde{u}|^2\tilde{u}\big\|_{L_t^\frac43 L_y^\frac43 \mathcal{H}_x^{1} }  &  \lesssim \epsilon. \label{eq3.27}
\end{align}
\end{lemma}
\begin{proof}
By symmetry, we may assume $t_0 = \inf I$. Let $w =  u - \tilde{u}$, then $w$ satisfies
\begin{equation*}
i\partial_t w + \Delta_{\mathbb{R}^3} w - x^2 w = |\tilde{u} + w|^2 (\tilde{u} + w) - |\tilde{u}|^2 \tilde{u} - e.
\end{equation*}
For $t\in I$, we define
\begin{equation*}
D (t) = \big\||\tilde{u} + w|^2(\tilde{u} + w) - |\tilde{u}|^2 \tilde{u}\big\|_{L_{t,y}^\frac43 \mathcal{H}_x^{1}([t_0,t]\times \mathbb{R}^2  \times \mathbb{R})}.
\end{equation*}
By \eqref{eq3.21}, we have
\begin{align}
D (t) \lesssim \|w\|_{L_{t,y}^4  \mathcal{ H}_x^{1  } }\left( \left\|\tilde{u} \right\|_{L_{t,y}^4  \mathcal{H}_x^{1} }^2 + \|w\|_{L_{t,y}^4   \mathcal{H}_x^{1} }^2\right)
  \lesssim \|w\|_{L_{t,y}^4  \mathcal{H}_x^{1}([t_0,t]\times \mathbb{R}^2   \times \mathbb{R})}^3 + \epsilon_1^2 \|w\|_{L_{t,y}^4  \mathcal{ H}_x^{1} ([t_0,t]\times \mathbb{R}^2  \times \mathbb{R})}.  \label{eq3.28}
\end{align}
On the other hand, by the Strichartz estimate and \eqref{eq3.22}, we get
\begin{align}
 \|w\|_{L_{t,y}^4   \mathcal{ H}_x^{1} ([t_0,t]\times \mathbb{R}^2   \times \mathbb{R})}
 \lesssim    \left\|e^{i(t-t_0) \left(\Delta_{\mathbb{R}^3} - x^2 \right)} w(t_0) \right\|_{L_{t,y}^4  \mathcal{H}_x^{1}   ([t_0,t]\times \mathbb{R}^2  \times \mathbb{R})}
+  D  (t) + \|e\|_{L_{t,y}^\frac43 \mathcal{H}_x^{1}   ([t_0,t]\times \mathbb{R}^2  \times \mathbb{R})}
\lesssim     D  (t) + \epsilon. \label{eq3.29}
\end{align}
Combining \eqref{eq3.28} and \eqref{eq3.29}, we obtain
\begin{equation*}
D  (t) \lesssim \left( D  (t) + \epsilon \right)^3 + \epsilon_1^2 \left( D  (t) + \epsilon \right).
\end{equation*}
A standard continuity argument then shows that if $\epsilon_1$ is taken sufficiently small, then
\begin{align*}
 D  (t) \lesssim \epsilon, \ \forall\, t\in I,
\end{align*}
which implies \eqref{eq3.27}.

Using \eqref{eq3.27} and \eqref{eq3.29}, one easily derives \eqref{eq3.24}. Moreover, by the Strichartz estimate, \eqref{eq3.20} and \eqref{eq3.27},
\begin{align*}
 \|w\|_{  L_t^\infty L_y^2 \mathcal{H}_x^{1} (I\times \mathbb{R}^2 \times  \mathbb{R})}
  \lesssim    \left\|w(t_0) \right\|_{   L_y^2 \mathcal{H}_x^{1}  } + \big\||\tilde{u} + w|^2(\tilde{u} + w) - |\tilde{u}|^2 \tilde{u}\big\|_{L_{t,y}^\frac43   \mathcal{ H}_x^{1}  } + \|e\|_{L_{t,y}^\frac43  \mathcal{ H}_x^{1}     }
\lesssim &  M' + \epsilon,
\end{align*}
which establishes \eqref{eq3.25} for $\epsilon_1= \epsilon_1(M')$ sufficiently small.

To prove \eqref{eq3.26}, we use the Strichartz estimate, \eqref{eq2.13new}, \eqref{eq3.20}, \eqref{eq3.27} and \eqref{eq3.21},
\begin{align*}
  \|u\|_{  L_t^\infty L_y^2 \mathcal{H}_x^{1}  (I\times \mathbb{R}^2  \times  \mathbb{R})}
 \lesssim &   \left\|\tilde{u}(t_0) \right\|_{    L_y^2 \mathcal{H}_x^{1} } +  \|u(t_0) - \tilde{u}(t_0) \|_{    L_y^2 \mathcal{H}_x^{1}   }   + \big\||u|^2 u - |\tilde{u}|^2 \tilde{u}\big\|_{L_{t,y}^\frac43  \mathcal{ H}_x^{1}     } + \big\||\tilde{u}|^2 \tilde{u}\big\|_{L_{t,y}^\frac43
\mathcal{ H}_x^{1}   }\\
 \lesssim &  M + M' + \epsilon+   \left\|\tilde{u} \right\|_{L_{t,y}^4\mathcal{H}_x^{1} }^3
 \lesssim M + M'+ \epsilon + \epsilon_1^3.
\end{align*}
The proof is complete by choosing $\epsilon_1 = \epsilon_1(M,M')$ sufficiently small.
\end{proof}
We now show the proof of Theorem \ref{th2.6}.
\begin{proof}[Proof of Theorem \ref{th2.6}]
We divide the interval $I$ into $ N \sim  \left( 1 + \frac{L}{\epsilon_0} \right)^4$ subintervals $I_j = [t_j, t_{j+1}]$, $0\le j\le N -1$ such that
\begin{equation*}
	\left\|\tilde{u} \right\|_{L_{t,y}^4   \mathcal{H}_x^{1}   (I_j\times \mathbb{R}^2  \times\mathbb{R})} \le \epsilon_1,
\end{equation*}
where $\epsilon_1 = \epsilon_1(M, 2M')$ is given by Lemma \ref{le3.6}.

By choosing $\epsilon_1$ sufficiently small depending on $J$, $M$ and $M'$, we can apply Lemma \ref{le3.6} to obtain for each $j$ and all $0 < \epsilon < \epsilon_1$,
\begin{align*}
	& \left\|u-\tilde{u} \right\|_{L_{t,y}^4  \mathcal{H}_x^{1}     (I_j\times \mathbb{R}^2  \times\mathbb{R})}    \le C(j) \epsilon,    \quad
	\left\|u-\tilde{u} \right\|_{ L_t^\infty L_y^2 \mathcal{H}_x^{1}    (I_j\times \mathbb{R}^2  \times\mathbb{R})}   \le C(j) M',\\
	& \|u\|_{  L_t^\infty L_y^2 \mathcal{H}_x^{1}   (I_j\times \mathbb{R}^2 \times\mathbb{R})}     \le C(j)(M + M'),   \quad
	\big\||u|^2 u - |\tilde{u}|^2 \tilde{u}\big\|_{L_{t,y}^\frac43  \mathcal{H}_x^{1}    (I_j\times \mathbb{R}^2  \times\mathbb{R})}     \le C(j) \epsilon,
\end{align*}
provided we can prove that analogues of \eqref{eq3.32} and \eqref{eq3.33} hold with $t_0$ replaced by $t_j$.

In order to verify this, we use an inductive argument. By the Strichartz estimate, \eqref{eq3.32}, and the inductive hypothesis,
\begin{align*}
	\left\|u(t_j)-\tilde{u}(t_j) \right\|_{  L_y^2 \mathcal{H}_x^{1}  }
	\lesssim & \  \left\|u(t_0)- \tilde{u}(t_0) \right\|_{  L_y^2 \mathcal{H}_x^{1}   } + \big\||u|^2 u - |\tilde{u}|^2 \tilde{u}\big\|_{L_{t,y}^\frac43  \mathcal{H}_x^{1} ([t_0,t_j]\times \mathbb{R}^2   \times \mathbb{R})} + \|e\|_{L_{t,y}^\frac43  \mathcal{H}_x^{1}    ([t_0,t_j]\times \mathbb{R}^2   \times \mathbb{R})}\\
	\lesssim  &\  M' + \sum\limits_{k=0}^{j-1} C(k)\epsilon + \epsilon.
\end{align*}
Similarly, by the Strichartz estimate, \eqref{eq3.33}, and the inductive hypothesis,
\begin{align*}
	& \left\|e^{i(t-t_j) \left(\Delta_{\mathbb{R}^3} -x^2 \right) } \left(u(t_j) - \tilde{u}(t_j) \right) \right\|_{L_{t,y}^4  \mathcal{ H}_x^{1}     (I_j\times \mathbb{R}^2 \times \mathbb{R})}\\
	\lesssim &  \left\|e^{i(t-t_0) \left( \Delta_{\mathbb{R}^3} - x^2 \right) } \left(u(t_0) - \tilde{u}(t_0) \right) \right\|_{L_{t,y}^4 \mathcal{H}_x^{1}    (I_j \times \mathbb{R}^2  \times \mathbb{R})}
	+ \|e\|_{L_{t,y}^\frac43 \mathcal{ H}_x^{1}  ([t_0,t_j]\times \mathbb{R}^2   \times \mathbb{R})}
	+ \big\||u|^2 u - |\tilde{u}|^2 \tilde{u}\big\|_{L_{t,y}^\frac43  \mathcal{ H}_x^{1}    ([t_0,t_j] \times \mathbb{R}^2 \times \mathbb{R})}\\
	\lesssim & \epsilon + \sum\limits_{k=0 }^{j-1}  C(k) \epsilon.
\end{align*}
It is clear now we may choose $\epsilon_1$ sufficiently small, depending on $N, M$ and $M'$, such that the hypotheses of Lemma \ref{le3.6} continue to hold as $j$ varies. This completes the proof of Theorem \ref{th2.6}.
\end{proof}

\noindent \textbf{Acknowledgments.}
Xing Cheng is grateful to Professor R\'emi Carles for explaining his work on the nonlinear Schr\"odinger equation with quadratic potentials and for pointing out some typos and misunderstanding in an initial version of this manuscript. Xing Cheng also thanks Binghua Feng for calling his attention to \cite{FO,FL}.


\begin{thebibliography}{99}
\bibitem{ACS} A. P. Antonelli, R. Carles, and J. D. Silva, \emph{Scattering for nonlinear Schr\"odinger equation under partial harmonic confinement},
Commun. Math. Phys. {\bf 334} (2015), no. 1, 367-396.

\bibitem{AC} A. H. Ardila and R. Carles, \emph{Global dynamics below the ground states for NLS under partial harmonic confinement}, arXiv: 2006.06205.

\bibitem{B} A. Barron, \emph{On global-in-time Strichartz estimates for the semiperiodic Schr\"odinger equation}, arXiv: 1901.01663.

\bibitem{BBJV} J. Bellazzini, N. Boussaid, L. Jeanjean, and N. Visciglia, \emph{Existence and stability of standing waves for supercritical NLS with a partial confinement}, Comm. Math. Phys. {\bf 353} (2017), no. 1, 229-251.

\bibitem{BBCE} A. Biasi, P. Bizo\'n, B. Craps, and O. Evnin, \emph{Two infinite families of resonant solutions for the Gross-Pitaevskii equation}, Phys. Rev. E {\bf98} (2018), no. 3, 032222.

\bibitem{BG} H. Bahouri and P. G\'erard, \emph{High frequency approximation of solutions to critical nonlinear wave equations}. Amer. J. Math. 121 (1999), no. 1, 131-175.

\bibitem{BV} P. B\'egout and A. Vargas, \emph{Mass concentration phenomena for the $L^2$-critical nonlinear Schr\"odinger
equation}, Trans. Amer. Math. Soc. {\bf359} (2007), no. 11, 5257-5282.

\bibitem{Bo} J. Bourgain, \emph{Refinements of Strichartz' inequality and applications to 2D-NLS with critical nonlinearity}, Internat. Math. Res. Notices 1998, no. 5, 253-283.

\bibitem{BGHS2}  T. Buckmaster, P. Germain, Z. Hani, and J. Shatah, \emph{Analysis of (CR) in higher dimension}, Int. Math. Res. Not. IMRN 2019, no. 4, 1265-1280.

\bibitem{BTT} N. Burq, L. Thomann, and N. Tzvetkov, \emph{Long time dynamics for the one dimensional nonlinear Schr\"odinger equation}, Ann. Inst. Fourier(Grenoble) {\bf 63} (2013), no. 6, 2137-2198.


\bibitem{Candy} T. Candy, \emph{Multi-scale bilinear restriction estimates for general phases}, Math. Ann. {\bf 375 } (2019), no. 1-2, 777-843.

\bibitem{FL} D. Cao, B. Feng, and T. Luo, \emph{On the standing waves for the X-ray free electron laser Schr\"odinger equation}, arXiv: 2005. 01516.


\bibitem{C4} R. Carles, \emph{Critical nonlinear Schr\"odinger equations with and without harmonic potential}, Math. Models Methods Appl. Sci. { \bf 12 } (2002), no. 10, 1513-1523.


\bibitem{C} R. Carles, \emph{Remarks on nonlinear Schr\"odinger equations with harmonic potential}, Ann. Henri Poincar\'e {\bf 3} (2002), no. 4, 757-772.

\bibitem{C5} R. Carles, \emph{Nonlinear Schr\"odinger equations with repulsive harmonic potential and applications}, SIAM J. Math. Anal. {\bf 35 } (2003), 823-843.


\bibitem{C3} R. Carles, \emph{Semi-classical analysis for nonlinear Schr\"odinger equations}, World Scientific Publishing Co. Pte. Ltd., Hackensack, NJ, 2008. xii+243 pp. ISBN: 978-981-279-312-6; 981-279-312-7.

\bibitem{C2} R. Carles, \emph{Nonlinear Schr\"odinger equation with time dependent potential}, Commun. Math. Sci. {\bf 9} (2011), no. 4, 937-964.

\bibitem{CG} R. Carles and C. Gallo, \emph{Scattering for the nonlinear Schr\"odinger equation with a general one-dimensional confinement}, J. Math. Phys. {\bf56} (2015), no. 10, 101503, 15 pp.

\bibitem{CK} R. Carles and S. Keraani, \emph{On the role of quadratic oscillations in nonlinear Schr\"odinger equations. II. The $L^2$-critical case}, Trans. Amer. Math. Soc. {\bf 359} (2007), no. 1, 33-62.

\bibitem{Cazenave} T. Cazenave, \emph{Semilinear Schr\"odinger equations}, Courant Lecture Notes in Mathematics, \textbf{10}.
New York University, Courant Institute of Mathematical Sciences, New York; American Mathematical
Society, Providence, RI, 2003.


\bibitem{CGYZ} X. Cheng, Z. Guo, K. Yang, and L. Zhao, \emph{On scattering for the cubic defocusing nonlinear Schr\"odinger equation on the waveguide $\mathbb{R}^2\times \mathbb{T}$}, Rev. Mat. Iberoam.  {\bf 36} (2020), no. 4, 985-1011.

\bibitem{CGZ} X. Cheng, Z. Guo, and Z. Zhao, \emph{On scattering for the defocusing quintic nonlinear Schr\"odinger equation on the two-dimensional cylinder},  SIAM J. Math. Anal. {\bf 52 } (5) (2020), pp. 4185-4237.

\bibitem{CGT} J. Colliander, M. Grillakis, and N. Tzirakis, \emph{Tensor products and correlation estimates with applications to nonlinear Schr\"odinger equations}, Comm. Pure Appl. Math. {\bf 62 }(2009), no. 7, 920-968.

\bibitem{CKSTT1} J. Colliander, M. Keel, G. Staffilani, H. Takaoka and T. Tao, \emph{Global existence and scattering for rough solutions of a nonlinear Schr\"odinger equation on $\mathbb{R}^3$}, Comm. Pure Appl. Math. {\bf 57 } (2004), no. 8, 987-1014.

 \bibitem{CKSTT0} J. Colliander, M. Keel, G. Staffilani, H. Takaoka and T. Tao, \emph{Global well-posedness and scattering for the energy-critical nonlinear Schr{\"o}dinger equation in $\mathbb{R}^3$}, Ann. of Math. {\bf167} (2008), 767-865.

\bibitem{CKSTT} J. Colliander, M. Keel, G. Staffilani, H. Takaoka, and T. Tao, \emph{Transfer of energy to high frequencies in the cubic defocusing nonlinear Schr\"odinger equation}, Invent. Math. {\bf 181} (2010), no. 1, 39-113.

\bibitem{CS} P. Constantin and J.-C. Saut, \emph{Local smoothing properties of dispersive equations}, J. Amer. Math. Soc. {\bf1} (1988), 413-439.

\bibitem{DELRV} S. Dartois, O. Evnin, L. Lionni, V. Rivasseau, and G. Valette, \emph{Melonic turbulence}, Comm. Math. Phys. {\bf 374 } (2020), no. 2, 1179-1228.

\bibitem{D3} B. Dodson, \emph{Global well-posedness and scattering for the defocusing, $L^2$-critical nonlinear Schr\"odinger equation when $d\ge 3$}, J. Amer. Math. Soc. {\bf25} (2012), no. 2, 429-463.

\bibitem{D4} B. Dodson, \emph{Global well-posedness and scattering for the mass critical nonlinear Schr\"odinger equation with mass below the mass of the ground state}, Adv. Math. {\bf 285} (2015), 1589-1618.

\bibitem{D1} B. Dodson, \emph{Global well-posedness and scattering for the defocusing, $L^2$-critical, nonlinear Schr\"odinger equation when $d=2$}, Duke Math. J. {\bf165} (2016), no. 18, 3435-3516.

\bibitem{D2} B. Dodson, \emph{Global well-posedness and scattering for the defocusing, $L^2$-critical, nonlinear Schr\"odinger equation when $d = 1$},  Amer. J. Math. {\bf 138} (2016), no. 2, 531-569.


\bibitem{D5} B. Dodson, \emph{Global well-posedness and scattering for the defocusing, mass-critical generalized KdV equation}, Ann. PDE {\bf3} (2017), no. 1, Art. 5, 35 pp.

\bibitem{D6} B. Dodson, \emph{Global well-posedness and scattering for the focusing, cubic Schr\"odinger equation in dimension $d=4$}, Ann. Sci. \'Ec. Norm. Sup\'er. (4) {\bf 52} (2019), no. 1, 139-180.

 \bibitem{DL} B. Dodson and A. Lawrie, \emph{Scattering for the radial 3D cubic wave equation}, Anal. PDE {\bf 8} (2015), no. 2, 467-497.

\bibitem{DLMM} B. Dodson, A. Lawrie, D. Mendelson, and J. Murphy, \emph{Scattering for defocusing energy subcritical nonlinear wave equations}, Anal. PDE {\bf 13 } (2020), no. 7, 1995-2090.

\bibitem{BMMZ} B. Dodson, C. Miao, J. Murphy and J. Zheng, \emph{The defocusing quintic NLS in four space dimensions}, Ann. Inst. H. Poincar\'e Anal. Non Lin\'eaire {\bf34} (2017), no. 3, 759-787.

\bibitem{FGH} E. Faou, P. Germain, and Z. Hani, \emph{The weakly nonlinear large-box limit of the 2D cubic nonlinear Schr\"odinger equation}, J. Amer. Math. Soc. {\bf29} (2016), no. 4, 915-982.

\bibitem{Fe} J. Fennell, \emph{Resonant Hamiltonian systems associated to the one-dimensional nonlinear Schr\"odinger equation with harmonic trapping},
Comm. Partial Differential Equations { \bf 44 } (2019), no. 12, 1299-1344.


\bibitem{FO} R. Fukuizumi and M. Ohta, \emph{Instability of standing waves for nonlinear Schr\"odinger equations with potentials}, Differential Integral Equations {\bf 16 } (2003), no. 6, 691-706.

\bibitem{GGT} P. G\'erard,  P. Germain, and L. Thomann, \emph{On the cubic lowest Landau level equation}, Arch. Ration. Mech. Anal. {\bf231} (2019), no. 2, 1073-1128.

\bibitem{GHT} P. Germain, Z. Hani, and L. Thomann, \emph{On the continuous resonant equation for NLS. I. Deterministic analysis},
J. Math. Pures Appl. (9) {\bf105} (2016), no. 1, 131-163.

\bibitem{GHT1} P. Germain, Z. Hani, and L. Thomann, \emph{On the continuous resonant equation for NLS, II: Statistical study},
Anal. PDE {\bf8} (2015), no. 7, 1733-1756.

\bibitem{GMS} P. Germain, N. Masmoudi, and J. Shatah, \emph{Global solutions for 2D quadratic Schr\"odinger equations}, J. Math. Pures Appl. (9) {\bf 97 } (2012), no. 5, 505-543.

\bibitem{GT}  P. Germain and L. Thomann, \emph{On the high frequency limit of the LLL equation}, Quart. Appl. Math. {\bf74} (2016), no. 4, 633-641.

\bibitem{HHK} M. Hadac, S. Herr, and H. Koch, \emph{Well-posedness and scattering for the KP-II equation in a critical space}, Ann. Inst. H. Poincar\'e Anal. Non Lin\'eaire {\bf 26 } (2009), no. 3, 917-941.

\bibitem{HP} Z. Hani and B. Pausader, \emph{On scattering for the quintic defocusing nonlinear Schr\"odinger equation on $\mathbb{R}\times \mathbb{T}^2$}, Comm. Pure Appl. Math. {\bf67} (2014), no. 9, 1466-1542.

\bibitem{HPTV} Z. Hani, B. Pausader, N. Tzvetkov, and N. Visciglia, \emph{Modified scattering for the cubic Schr\"odinger equation on product spaces and applications}, Forum of Mathematics, PI. (2015), Vol. 3, 1-63.

\bibitem{HT} Z. Hani and L. Thomann, \emph{Asymptotic behavior of the nonlinear Schr\"odinger equation with harmonic trapping}, Comm. Pure Appl. Math. {\bf69} (2016), no. 9, 1727-1776.

\bibitem{HHL1} C. Hao, L. Hsiao, and H. Li, \emph{Global well posedness for the Gross-Pitaevskii equation with an angular momentum rotational term in three dimensions}, J. Math. Phys. {\bf48} (2007), no. 10, 102105, 11 pp.

\bibitem{HHL2} C. Hao, L. Hsiao, and H. Li, \emph{Global well posedness for the Gross-Pitaevskii equation with an angular momentum rotational term},
Math. Methods Appl. Sci. {\bf31} (2008), no. 6, 655-664.

\bibitem{He} B. Helffer, \emph{Semi-classical analysis for the Schr\"odinger operators and applications}, Lecture Notes in Mathematics, 1336. Springer, Berlin, 1988.

\bibitem{IMN} S. Ibrahim, N. Masmoudi, and K. Nakanishi, \emph{Scattering threshold for the focusing nonlinear Klein-Gordon equation}, Anal. PDE {\bf4} (2011), no. 3, 405-460.

\bibitem{IP} A. D. Ionescu and B. Pausader, \emph{Global well-posedness of the energy-critical defocusing NLS on $\mathbb{R}\times \mathbb{T}^3$}, Comm. Math. Phys. {\bf 312} (2012), no. 3, 781-831.

\bibitem{J1} C. Jao, \emph{The energy-critical quantum harmonic oscillator}, Comm. Partial Differential Equations {\bf 41} (2016), no. 1, 79-133.

\bibitem{J2} C. Jao, \emph{Energy-critical NLS with potentials of quadratic growth}, Discrete Contin. Dyn. Syst. {\bf38} (2018), no. 2, 563-587.

\bibitem{J4} C. Jao, \emph{Refined mass-critical Strichartz estimates for Schr\"odinger operators}, Anal. PDE {\bf 13 } (2020), no. 7, 1955-1994.


\bibitem{JKV} C. Jao, R. Killip, and M. Visan, \emph{Mass-critical inverse Strichartz theorems for 1d Schr\"odinger operators}, Rev. Mat. Iberoam. {\bf 35 } (2019), no. 3, 703-730.

\bibitem{JP} C. Josserand and Y. Pomeau, \emph{Nonlinear aspects of the theory of Bose-Einstein condensates}, Nonlinearity {\bf 14}(5), R25-R62(2001).

\bibitem{KT0} M. Keel and T. Tao, \emph{Endpoint Strichartz Estimates}, Amer. J. Math. {\bf 120 }(1998), 955-980.

\bibitem{KM} C. E. Kenig and F. Merle, \emph{Global well-posedness, scattering and blow-up for the energy-critical, focusing, non-linear Schr\"odinger equation in the radial case}, Invent. Math. {\bf 166} (2006), no. 3, 645-675.

\bibitem{KM1} C. E. Kenig and F. Merle, \emph{Global well-posedness, scattering and blow-up for the energy-critical focusing non-linear wave equation}, Acta Math. {\bf201} (2008), no. 2, 147-212.

\bibitem{KTV09} R. Killip, T. Tao, and M. Visan, \emph{The cubic nonlinear Schr\"odinger equation in two dimensions with radial data}, J. Eur. Math. Soc. (JEMS) {\bf 11 } (2009), no. 6, 1203-1258.

\bibitem{Killip-Visan1} R. Killip and M. Visan, \emph{Nonlinear Schr\"odinger equations at critical regularity}. Proceedings for the Clay summer school ``Evolution Equations'', Eidgen\"ossische technische Hochschule, Z\"urich, 2008.

\bibitem{Killip-Visan2}  R. Killip and M. Visan, \emph{Global well-posedness and scattering for the defocusing quintic NLS in three dimensions}, Anal. PDE {\bf 5 } (2012), no. 4, 855-885.



\bibitem{KVZ0} R. Killip, M. Visan, and X. Zhang, \emph{The mass-critical nonlinear Schr\"odinger equation with radial data in dimensions three and higher}, Anal. PDE {\bf 1 } (2008), no. 2, 229-266.

\bibitem{KVZ} R. Kllip, M. Visan, and X. Zhang, \emph{Energy-critical NLS with quadratic potentials}, Comm. Partial Differential Equations {\bf 34} (2009), no. 10-12, 1531-1565.

\bibitem{KT3} H. Koch and D. Tataru, \emph{Dispersive estimates for principally normal pseudodifferential operators}, Comm. Pure Appl. Math. {\bf 58 } (2005), no. 2, 217-284.

\bibitem{KT} H. Koch and D. Tataru, \emph{$L^p$ eigenfunction bounds for the Hermite operator}, Duke Math. J. {\bf 128 } (2005), no. 2, 369-392.


\bibitem{KT1} H. Koch and D. Tataru, \emph{A priori bounds for the 1D cubic NLS in negative Sobolev spaces}, Int. Math. Res. Not. IMRN 2007, no. 16, Art. ID rnm053, 36 pp.

\bibitem{KTV} H. Koch, D. Tataru, and M. Visan, \emph{Dispersive equations and nonlinear waves. Generalized Korteweg-de Vries, nonlinear Schr\"odinger, wave and Schr\"odinger maps}, Oberwolfach Seminars, {\bf 45}. Birkh\"auser/Springer, Basel, 2014. xii+312 pp.

\bibitem{Le} E. H. Lieb and M. Loss, \emph{Analysis}. Second edition. Graduate Studies in Mathematics, {\bf14}. American Mathematical Society, Providence, RI, 2001.

\bibitem{Lin-S} J. E. Lin and W. A. Strauss, \emph{Decay and scattering of solutions of a nonlinear Schr\"odinger equation}, J. Funct. Anal. {\bf 30 } (1978), 245-263.


\bibitem{MV} F. Merle and L. Vega, \emph{Compactness at blow-up time for $L^2$ solutions of the critical nonlinear Schr\"odinger equation in 2D}, Internat. Math. Res. Notices (1998), no. 8, 399-425.

\bibitem{Morawetz} C. S. Morawetz, \emph{Time decay for the nonlinear Klein-Gordon equation}, Proc. R. Soc. Ser. A {\bf 306} (1968), 291-296.


\bibitem{Mu} J. Murphy, \emph{Intercritical NLS: critical $H^s$-bounds imply scattering}, SIAM J. Math. Anal. {\bf46} (2014), no. 1, 939-997.


\bibitem{PS} L. Pitaevskii and S. Stringari, \emph{Bose-Einstein condensation}, International Series of Monographs on Physics, vol. {\bf 116}. The Clarendon Press Oxford University Press, Oxford (2003).

\bibitem{PV} F. Planchon and L. Vega, \emph{Bilinear virial identities and applications}, Ann. Sci. \'Ec. Norm. Sup\'er. (4) {\bf 42} (2009), no. 2, 261-290.

\bibitem{R} M. Rosenzweig, \emph{Global well-poseness and scattering for the elliptic-elliptc Davey-Stewartson system at $L^2$-critical regularity}, arXiv: 1808. 01955.

\bibitem{SW} J. Shen and Y. Wu, \emph{Global well-posedness and scattering of 3D defocusing, cubic Schr\"odinger equation}, arXiv: 2008. 10019.


\bibitem{Sta} G. Staffilani, \emph{The theory of nonlinear Schr\"odinger equations}, Evolution equations, 207-267, Clay Math. Proc., {\bf 17}, Amer. Math. Soc., Providence, RI, 2013.

\bibitem{SS} M. Stanislavova and  A. Stefanov, \emph{Ground states for the nonlinear Schr\"odinger equation under a general trapping potential}, arXiv:2002.03822.


\bibitem{T2} T. Tao, \emph{Nonlinear dispersive equations: local and global analysis}, CBMS Regional Conference Series in Mathematics, {\bf 106}. American Mathematical Society, Providence, R.I., 2006.

\bibitem{TVZ0} T. Tao, M. Visan, and X. Zhang, \emph{Global well-posedness and scattering for the defocusing mass-critical nonlinear Schr\"odinger equation for radial data in high dimensions}, Duke Math. J. {\bf 140 }  (2007), no. 1, 165-202.


\bibitem{TVZ1} T. Tao, M. Visan, and X. Zhang, \emph{The nonlinear Schr\"odinger equation with combined power-type nonlinearities}, Comm. Partial Differential Equations {\bf 32 } (2007), no. 7-9, 1281-1343.

\bibitem{TVZ} T. Tao, M. Visan, and X. Zhang, \emph{Minimal-mass blowup solutions of the mass-critical NLS}, Forum Math. {\bf20} (2008), no. 5, 881-919.


\bibitem{Ta} M. Tarulli, \emph{Well-posedness and scattering for the mass-energy NLS on $\mathbb{R}^N \times \mathbb{M}^K$}, Analysis {\bf 37}, Issue 3, 117-132, 2017.

\bibitem{TV} N. Tzvetkov and N. Visciglia, \emph{Small data scattering for the nonlinear Schr\"odinger equation on product spaces}, Comm. Partial Differential Equations. {\bf 37} (2012), no. 1, 125-135.

\bibitem{Ve} L. Vega, \emph{Schr\"odinger equations: pointwise convergence to the initial data}, Proc. Amer. Math. Soc. {\bf102} (1988), 874-878.


\bibitem{Visan} M. Visan, \emph{Global well-posedness and scattering for the defocusing cubic nonlinear Schr\"odinger equation in four dimensions}, Int. Math. Res. Not. IMRN 2012, no. 5, 1037-1067.


\bibitem{YZ} K. Yajima and G. Zhang, \emph{Local smoothing property and Strichartz inequality for Schr\"odinger equations with potentials superquadratic at infinity}, J. Differential Equations {\bf 202} (2004), no. 1, 81-110.

\bibitem{YZ0} K. Yang and L. Zhao, \emph{Global well-posedness and scattering for mass-critical, defocusing, infinite dimensional vector-valued resonant nonlinear Schr\"odinger system}, SIAM J. Math. Anal. {\bf50} (2018), no. 2, 1593-1655.

\bibitem{Z} J. Zhang, \emph{Sharp threshold for blowup and global existence in nonlinear Schr\"odinger equations under a harmonic potential}, Comm. Partial Differential Equations {\bf 30} (2005), no. 10-12, 1429-1443.

\bibitem{Z3} J. Zhang, \emph{Sharp threshold of global existence for nonlinear Schr\"odinger equation with partial confinement}, Nonlinear Anal. {\bf 196 } (2020), 111832.

\bibitem{Z1} Z. Zhao, \emph{Global well-posedness and scattering for the defocusing cubic Schr\"odinger equation on waveguide $\mathbb{R}^2 \times \mathbb{T}^2$},  J. Hyperbolic Differ. Equ. {\bf16} (2019), no. 1, 73-129.

\end{thebibliography}
\end{document}